\documentstyle{amltd}
\begin{document}
\annalsline{153}{2001}
\received{September 8, 1997}
\revised{November 3, 1999}
\startingpage{27}
\def\nmid{\hbox{$\mid \hskip-6.5pt\raise1pt\hbox{$\scriptscriptstyle /$}$}}
\font\tencyr=wncyr10 scaled 1100
\font\smcyr=wncyr8
\def\scyr#1{\hbox{\smcyr #1}}
\def\cyr#1{\hbox{\tencyr #1}}
\def\Sha{{\cyr{Sh}}}
\def\sSha{{\scyr{Sh}}}
\def\nnp{\nonumproclaim}
\def\spn{\specialnumber}
\def\spq{\speqnu}
\def\nn{\nonumber}
\def\ritem#1{\item[{\rm #1}]}
\def\eqref#1{(\ref{#1})}
\def\joinrel{\mathrel{\mkern-4mu}}
\def\relbar{\mathrel{\smash-}}
\def\lrar#1{{\stackrel{#1}{\relbar\joinrel\relbar\joinrel\relbar\joinrel\relbar\joinrel\relbar\joinrel\rightarrow}}}

\font\tenrm=cmr10
\catcode`\@=11
\font\twelvemsb=msbm10 scaled 1100
\font\tenmsb=msbm10
\font\ninemsb=msbm10 scaled 800
\newfam\msbfam
\textfont\msbfam=\twelvemsb  \scriptfont\msbfam=\ninemsb
  \scriptscriptfont\msbfam=\ninemsb
\def\msb@{\hexnumber@\msbfam}
\def\Bbb{\relax\ifmmode\let\next\Bbb@\else
 \def\next{\errmessage{Use \string\Bbb\space only in math
mode}}\fi\next}
\def\Bbb@#1{{\Bbb@@{#1}}}
\def\Bbb@@#1{\fam\msbfam#1}
\catcode`\@=12

 \catcode`\@=11
\font\twelveeuf=eufm10 scaled 1100
\font\teneuf=eufm10
\font\nineeuf=eufm7 scaled 1100
\newfam\euffam
\textfont\euffam=\twelveeuf  \scriptfont\euffam=\teneuf
  \scriptscriptfont\euffam=\nineeuf
\def\euf@{\hexnumber@\euffam}
\def\frak{\relax\ifmmode\let\next\frak@\else
 \def\next{\errmessage{Use \string\frak\space only in math
mode}}\fi\next}
\def\frak@#1{{\frak@@{#1}}}
\def\frak@@#1{\fam\euffam#1}
\catcode`\@=12
\def\bye{ \end{document}}
\def\mod{\, {\rm mod}\, }


 \newcommand{\lb}{\label}

 \newcommand{\BA}{{{\Bbb A}}}
\newcommand{\BB}{{{\Bbb B}}}
\newcommand{\BC}{{{\Bbb C}}}
\newcommand{\BD}{{{\Bbb D}}}
\newcommand{\BE}{{{\Bbb E}}}
\newcommand{\BF}{{{\Bbb F}}}
\newcommand{\BG}{{{\Bbb G}}}
\newcommand{\BH}{{{\Bbb H}}}
\newcommand{\BI}{{{\Bbb I}}}
\newcommand{\BJ}{{{\Bbb J}}}
\newcommand{\BK}{{{\Bbb K}}}
\newcommand{\BL}{{{\Bbb L}}}
\newcommand{\BM}{{{\Bbb M}}}
\newcommand{\BN}{{{\Bbb N}}}
\newcommand{\BO}{{{\Bbb O}}}
\newcommand{\BP}{{{\Bbb P}}}
\newcommand{\BQ}{{{\Bbb Q}}}
\newcommand{\BR}{{{\Bbb R}}}
\newcommand{\BS}{{{\Bbb S}}}
\newcommand{\BT}{{{\Bbb T}}}
\newcommand{\BU}{{{\Bbb U}}}
\newcommand{\BV}{{{\Bbb V}}}
\newcommand{\BW}{{{\Bbb W}}}
\newcommand{\BX}{{{\Bbb X}}}
\newcommand{\BY}{{{\Bbb Y}}}
\newcommand{\BZ}{{{\Bbb Z}}}

\newcommand{\CA}{{{\cal A}}}
\newcommand{\CB}{{{\cal B}}}
\newcommand{\CC  }{{{\cal C}}}
\newcommand{\CD}{{{\cal D}}}
\newcommand{\CE}{{{\cal E}}}
\newcommand{\CF}{{{\cal F}}}
\newcommand{\CG}{{{\cal G}}}
\newcommand{\CH}{{{\cal H}}}
\newcommand{\CI}{{{\cal I}}}
\newcommand{\CJ}{{{\cal J}}}
\newcommand{\CK}{{{\cal K}}}
\newcommand{\CL}{{{\cal L}}}
\newcommand{\CM}{{{\cal M}}}
\newcommand{\CN}{{{\cal N}}}
\newcommand{\CO}{{{\cal O}}}
\newcommand{\CP}{{{\cal P}}}
\newcommand{\CQ}{{{\cal Q}}}
\newcommand{\CR }{{{\cal R}}}
\newcommand{\CS}{{{\cal S}}}
\newcommand{\CT}{{{\cal T}}}
\newcommand{\CU}{{{\cal U}}}
\newcommand{\CV}{{{\cal V}}}
\newcommand{\CW}{{{\cal W}}}
\newcommand{\CX}{{{\cal X}}}
\newcommand{\CY}{{{\cal Y}}}
\newcommand{\CZ}{{{\cal Z}}}

\newcommand{\RA}{{{\rm A}}}
\newcommand{\RB}{{{\rm B}}}
\newcommand{\RC  }{{{\rm C}}}
\newcommand{\RD}{{{\rm D}}}
\newcommand{\RE}{{{\rm E}}}
\newcommand{\RF}{{{\rm F}}}
\newcommand{\RG}{{{\rm G}}}
\newcommand{\RH}{{{\rm H}}}
\newcommand{\RI}{{{\rm I}}}
\newcommand{\RJ}{{{\rm J}}}
\newcommand{\RK}{{{\rm K}}}
\newcommand{\RL}{{{\rm L}}}
\newcommand{\RM}{{{\rm M}}}
\newcommand{\RN}{{{\rm N}}}
\newcommand{\RO}{{{\rm O}}}
\newcommand{\RP}{{{\rm P}}}
\newcommand{\RQ}{{{\rm Q}}}
\newcommand{\RR }{{{\rm R}}}
\newcommand{\RS}{{{\rm S}}}
\newcommand{\RT}{{{\rm T}}}
\newcommand{\RU}{{{\rm U}}}
\newcommand{\RV}{{{\rm V}}}
\newcommand{\RW}{{{\rm W}}}
\newcommand{\RX}{{{\rm X}}}
\newcommand{\RY}{{{\rm Y}}}
\newcommand{\RZ}{{{\rm Z}}}

\newcommand{\Ad}{{{\rm Ad}}}
\newcommand{\Aut}{{{\rm Aut}}}

\newcommand{\Br}{{{\rm Br}}}

\newcommand{\Ch}{{{\rm Ch}}}
\newcommand{\cod}{{{\rm cod}}}
\newcommand{\cont}{{{\rm cont}}}
\newcommand{\cl}{{{\rm cl}}}

\newcommand{\disc}{{{\rm disc}}}
\newcommand{\Div}{{{\rm Div}}}
\renewcommand{\div}{{{\rm div}}}

\newcommand{\End}{{{\rm End}}}

\newcommand{\Frob}{{{\rm Frob}}}

\newcommand{\Gal}{{{\rm Gal}}}
\newcommand{\GL}{{{\rm GL}}}

\newcommand{\Hom}{{{\rm Hom}}}

\renewcommand{\Im}{{{\rm Im}}}
\newcommand{\Ind}{{{\rm Ind}}}
\newcommand{\inv}{{{\rm inv}}}
\newcommand{\Isom}{{{\rm Isom}}}

\newcommand{\Jac}{{{\rm Jac}}}

\newcommand{\Ker}{{{\rm Ker}}}

\newcommand{\Lie}{{{\rm Lie}}}

\newcommand{\new}{{{\rm new}}}
\newcommand{\NS}{{{\rm NS}}}

\newcommand{\ord}{{{\rm ord}}}

\newcommand{\rank}{{{\rm rank}}}

\newcommand{\PGL}{{{\rm PGL}}}
\newcommand{\Pic}{{\rm Pic}}

\renewcommand{\Re}{{{\rm Re}}}
\newcommand{\reg}{{{\rm reg}}}
\newcommand{\Res}{{{\rm Res}}}

\newcommand{\Sel}{{{\rm Sel}}}
\newcommand{\SL}{{{\rm SL}}}
\newcommand{\Spec}{{{\rm Spec}}}
\newcommand{\Sym}{{{\rm Sym}}}
\newcommand{\sgn}{{{\rm sgn}}}

\newcommand{\tr}{{{\rm tr}}}

\newcommand{\ur}{{{\rm ur}}}

\newcommand{\vol}{{{\rm vol}}}

\newcommand{\wt}{\widetilde}
\newcommand{\wh}{\widehat}
\newcommand{\pp}{\frac{\partial\bar\partial}{\pi i}}
\newcommand{\pair}[1]{\langle {#1} \rangle}
\newcommand{\wpair}[1]{\left\{{#1}\right\}}
\newcommand{\intn}[1]{\left( {#1} \right)}
\newcommand{\norm}[1]{\|{#1}\|}
\newcommand{\sfrac}[2]{\left( \frac {#1}{#2}\right)}

\title{Heights of Heegner points\\
on Shimura curves}  
\shorttitle{Height of Heegner points} 
\acknowledgements{This research was supported by a grant from 
the National Science Foundation, and 
a research fellowship of the Sloan Foundation.}
\author{Shouwu Zhang}
   \institutions{ Columbia University,
New York, NY\\
 {\eightpoint {\it E-mail address\/}: szhang@math.columbia.edu}}
 
 \def\sni#1{\smallbreak\noindent{#1}. }
\def\ssni#1{\vglue-1pt\noindent\hskip18pt {#1}.}

\centerline{\bf Contents}
\begin{small}
\bigbreak
\noindent Introduction
\sni{1} Shimura curves
\ssni{1.1} Modular interpretation
\ssni{1.2} Integral models
\ssni{1.3} Reductions of models
\ssni{1.4} Hecke correspondences
\ssni{1.5} Order $R$ and its level structure
\sni{2} Heegner points
\ssni{2.1} CM-points
\ssni{2.2} Formal groups
\ssni{2.3} Endomorphisms
\ssni{2.4} Liftings of distinguished points
\sni{3} Modular forms and $L$-functions
\ssni{3.1} Modular forms
\ssni{3.2} Newforms on $X$
\ssni{3.3} The supercuspidal case
\ssni{3.4} $L$-functions associated to newforms
\ssni{3.5} Eisenstein series and theta series
\sni{4} Global intersections
\ssni{4.1} Height pairing
\ssni{4.2} Computing $T(m)\eta$
\ssni{4.3} Computing $T(m)\hat{\xi}^{\phantom{\sum}}$
\ssni{4.4} Computing $T(m)Z$
\ssni{4.5} A uniqueness theorem
\sni{5} Local intersections
\ssni{5.1} Archimedean intersections
\ssni{5.2} Nonarchimedean intersections
\ssni{5.3} Clifford algebras
\ssni{5.4} The final formula
\sni{6} Derivatives of $L$-series
\ssni{6.1} The Rankin-Selberg method
\ssni{6.2} Fourier coefficients
\ssni{6.3} Functional equations and derivatives
\ssni{6.4} Holomorphic projections
\sni{7} Proof of the main theorems
\ssni{7.1} Proof of Theorem C
\ssni{7.2} Proof of Theorem A
\smallbreak\noindent Index
\smallbreak\noindent References

\end{small}

 \bigbreak

\intro
The purpose of this paper is to  generalize some results of 
Gross-Zagier \cite{G-Z} and Kolyvagin \cite{K} to totally real fields.
 The main result and the plan of its proof are described as follows.  

\demo{Main results} Let $F$ be a totally real number field and $N$ a nonzero
ideal of $\CO _F$. Let $f$ be a newform on
$\GL _2(\BA _F)$, of (parallel) weight $2$, level $   K   _0( N  )$,
and with  trivial central character, where
  $   K   _0(   N  )$ denotes the subgroup of $\GL _2(\wh F)$:
$$   K   _0(   N  )= \left\{ \left( \begin{array}{cc} a&b\\ c&d \end{array}\right)\in 
\GL _2(\wh \CO _F)\bigg| c\in   \wh N  \right\},$$
where for an abelian group $M$, $\wh M$ denotes $M\otimes \prod _p\BZ   _p$.  
Let $\CO _f$ denote the subalgebra of $\BC  $ over $\BZ   $ generated by
eigenvalues $a(f,   m )$ of $f$ under the Hecke operators.
For each embedding $\sigma : \CO_f\to \BC  $,
let $f^\sigma$ denote the newform with the eigenvalues
$a(f^\sigma,   m )=a(f,   m )^\sigma$.
Assume that
either $[F:\BQ]$ is odd or $\ord _v (N)=1$ for at least one finite 
place $v$ of $F$.
 Then
there is an abelian variety $A$ over $F$ of dimension $[\CO_f:\BZ   ]$ such that
$L(s, A)$ equals  $\prod _{\sigma:\CO_f\to \BC  } L(s, f^\sigma)$ modulo 
the factors at places dividing $   N  $.
 Our  main result is
the following:
\enddemo
 
\nnp{Theorem A}
  Assume the   $L$\/{\rm -}\/function
$L(s,f)$ has order $\le 1$ at $s=1${\rm .}
 Then for any $A$ as above\/{\rm :}\/
\begin{itemize}
\ritem{1.} The Mordell\/{\rm -}\/Weil group $A(F)$ has  rank given by
$$\rank\, A(F)=[\CO_f:\BZ   ]\,\ord _{s=1}L(s, f);$$
\ritem{2.} The  Shafarevich\/{\rm -}\/Tate group $\Sha (A)$ is finite{\rm .}
\end{itemize}

\endproclaim

The theorem would hold  with a weaker condition that either
$[F:\BQ]$ is odd or $\ord _v(N)$ is odd for at least one
finite place, provided   we could overcome one technical difficulty.
 That is, it would be enough to prove  Lemma  
5.2.3 (below) without  the assumption  that $\ord _\wp (N)\le 1$.

\demo{Shimura curves} As in the case $F=\BQ$ treated by Gross-Zagier and\break Kolyvagin, 
we will prove the theorem by studying  Heegner 
points over some imaginary quadratic extension. 
 Let $ E$ 
be a totally imaginary quadratic extension of $F$  which is unramified
over places dividing $   N  $. Assume further 
$\varepsilon (   N  )=(-1)^{g-1}$ where $g=[F:\BQ]$ and 
$$\varepsilon=\otimes_v\varepsilon_v :F^\times\backslash \wh F^\times
\to \{\pm 1\}$$ is the  character on $\BA _F^\times/F^\times$ 
associated to the extension
$ E/F$.  Let $\tau $ be a fixed archimedean place, and 
let $  B   $ be a quaternion algebra over $F$ which is nonsplit exactly at all
archimedean places other than $\tau$, and finite places $v$ such that
$\varepsilon _v ( N  )=-1$.
Fix an embedding $\rho:  E\to   B   $ over $F$.
Let $   R   $ be an order of $  B   $ of type $(   N  ,  E)$, that is an order 
of $  B   $ of discriminant $   N  $ which contains $\rho (\CO _ E)$. 
Fix an isomorphism $  B   _\tau\otimes \BR\simeq \RM_2(\BR)$
such that $\rho ( E)\otimes\BR$ is sent to the subalgebra of $\RM_2(\BR)$
of elements $ \left( \begin{array}{cc} a&b\\ -b&a \end{array}\right)$. Then the group
$  B_+$ of the elements in $  B    ^\times$ with totally positive 
reduced norm acts on the Poincar\'e half-plane $\CH$. 
Thus we obtain a Shimura curve
$$X_\tau(\BC  )=  B    _+ \backslash \CH\times 
\wh   B    ^\times /\wh F^\times \wh   R    ^\times \cup\{{\rm cusps}\}$$
where $\{{\rm cusps}\}$ is not empty only if $F=\BQ$
and $\varepsilon_v ( N  )=1$ for any $v| N  $.
By  Shimura's theory \cite{S}, $X_\tau(\BC  )$ has a canonical model
$X$ defined over~$F$.

The curve $X$ over $F$ is connected but not geometrically connected. Let
$\Jac (X)$ denote the connected component subgroup of 
$\Pic (X/F)$. Then, $$\Jac (X)=\Res _{\wt F/F}\Pic ^0 (X/\wt F),$$
where  $\wt F$ denotes  the abelian Galois extension 
of $F$ corresponding to the subgroup $F_+\cdot (\wh F^\times )^2
\cdot \wh\CO _F^\times$ via  class field theory. \enddemo

\nnp{Theorem B}
 There is a unique abelian subvariety $A$ of $\Jac (X)$ defined
over $F$  of dimension $[\CO_f:\BZ   ]$ such that 
$L(s,A)$ is equal to $\prod _{\sigma:\CO_f\to \BC  } L(s, f^\sigma)$ modulo 
the factors at places  dividing $   N  ${\rm . }
\endproclaim

We will prove this theorem in Section~3, by combining the  Eichler-Shimura theory
and a newform theory for $X$  obtained by using
 Jacquet-Langlands theory \cite{J-L}.
The key to the newform theory on $X$ is Proposition 3.3.1. 
I am indebted to H. Jacquet for showing me the proof in  
the supercuspidal case using results of Waldspurger.
(After the paper was submitted, I learned from Gross and the referees that 
some related results have been obtained by Tunnell \cite{Tu}
and Gross \cite{Gr5}.) 

\demo{Heegner points}  Let $  x  $ denote the image on
$X_\tau(\BC  )$ of $\{\sqrt{-1}\}\times \{1\}\in \CH\times \wh B^\times$.
By  Shimura's theory \cite{S},  $  x  $ is defined over the 
Hilbert class field $H$ of $ E$. We call  
$  x  $ a Heegner point on $X$.

In order to construct  a   point in the Jacobian $\Jac (X)$ from $  x  $, 
we need to define a map from $X$ to $\Jac (X)$.
Write  $X_\tau (\BC  )$ as a union $\mathbold{\cup} X_i$ of connected compact Riemann
surfaces of the form 
$$X_i= \Gamma   _i\backslash \CH \cup \{{\rm cusps}\}$$
with $ \Gamma   _i\subset B_+/F^\times\subset {\rm PSL}_2(\BR)$.
Then one has
$\Jac (X)({\Bbb C})=\prod\Jac (X_i)$.
We define a canonical divisor 
class of degree $1$ in $\Pic (X_i)\otimes \BQ$ by the formula
$$\xi_i:=
\left\{[\Omega ^1_{X_i}]+\sum _{p\in X_i}(1-\frac 1{u_p})[p]+[{\rm cusps}]
\right\}\bigg/\int _{X_i}\frac {dxdy}{2\pi y^2},$$
where for any noncuspidal point  $p\in X_i$,
$u_p$ denotes the  cardinality of the group of
stabilizers of  $\wt p$ in $ \Gamma  $,
where $\wt p$ is a point in $\CH$  projecting to $p$.
 Now we define
a map $\phi: X\to\Jac (X)\otimes \BQ$ which sends a point $p\in X_i$ to the class
of $p-\xi _i$.
It is easy to see that some positive multiple of $\phi$ is actually 
defined over $F$.

Let $ z  $ denote the class
$$u_  x   ^{-1}\sum _{\sigma \in \Gal (H/ E)}\phi(  x  ^\sigma)$$
 in $\Jac (X)( E)\otimes \BQ$.
Let $ z   _f$ be the  component of $ z  $ in 
$A\otimes \BQ$. 
\enddemo

\vglue5pt 

\demo{Gross-Zagier formula} 
Now we assume that
a prime $\wp$ is split in $E$ if either 
$\wp$ divides $2$ or $\ord _\wp(N)>1$.\enddemo

\vglue5pt 

\nnp{Theorem C}
  Let 
$L_ E(s, f)$ denote the product $L(s, f)L(s, \varepsilon, f)${\rm ,} where 
$L(s, \varepsilon, f)$ is the $L$\/{\rm -}\/function of $f$ twisted by $\varepsilon${\rm .}  
Then $L_ E(f, 1)=0$ and
$$L_ E'(f, 1)=\frac {(8\pi ^2)^g}{d_F^2\sqrt{d_ E}}
[K   _0(1):   K  _0( N  )](f, f)
\pair{ z   _f,  z   _f},$$
where 
\begin{itemize}
\ritem{1.} $\pair{ z   _f,  z   _f}$ is the N{\rm \'{\it e}}ron\/{\rm -}\/Tate height of $ z  _f ${\rm ;}
\ritem{2.} $d_F$ is the discriminant of $F${\rm ,} and 
$d_E$ is the norm of the relative discriminant of $E/F${\rm ;}
\ritem{3.} $(f, f)$ is the inner product with respect to the standard measure
on\break $Z(\BA _F)\GL _2(F)\backslash \GL _2(\BA _F)${\rm .}
\end{itemize}

\endproclaim

If $F=\BQ$ and every prime factor of $   N  $ is split in $ E$,
this is due to  Gross-Zagier \cite{G-Z1}. Again, the extra
condition that $\wp$ is split in $E$ when $\ord _\wp (N)>1$
  can be eliminated if we know how to
compute local intersections at $\wp$ when the integral model 
of Shimura curves has
some mild singularities over $\wp$.
 
For the proof of the second part of Theorem A, we assume that
 $L(s, f)$ has order less than or equal to $1$.   
By some results in  \cite{B-F-H} and  \cite{W} (the theorem in \cite {B-F-H}
is stated for $\BQ$, but its proof can be easily generalized to 
any number field), there is an 
$ E$ such that $L_ E(f, s)$ has order equal to $1$ at $s=1$. It follows 
from Theorem C that $ z  _f$ has infinite order. Now the second part of
Theorem A follows from  Kolyvagin's method 
\cite{G3}, \cite{K}, \cite{K-L}, \cite{K-L2}, which applies directly to
our case without any new difficulty. The only thing we need is to
give a correct  system of CM-points which we will do at the 
end of this paper.\enddemo

\demo{Plan of proof} Now we sketch the proof of Theorem C. Let $ \Psi  $ and $ \Phi $ be 
two cusp forms on 
$\GL _2(\BA _F)$ of weight $2$ and level $   K  _0( N  )$
 characterized by the following
properties:
\begin{itemize}
\item  The Fourier coefficients of $ \Psi  $
are given by 
$$a( \Psi  ,   m )=\pair{ z  , \RT(  m ) z  }$$ for all
$  m $.
\item  The form $ \Phi $ satisfies the equality
$$L_ E'(f, 1)=c{(f,  \Phi )}
$$
for any newform $f$ on $\GL _2(\BA _F)$ of   weight $2$,  level
 $   K   _0( N   )$, and with  trivial central character,
where $c$ is some constant, and
$(\cdot, \cdot)$ denotes the Weil-Petersson product.
\end{itemize}
Then the equality in Theorem C is equivalent to  $ \Phi \equiv  
{\rm const}\cdot\Psi  $
modulo  old forms, and the proof of Theorem C is reduced to the computations
of Fourier coefficients of $ \Psi  $ and $ \Phi $ respectively.
We will do this by using Arakelov theory and
the Rankin-Selberg  method respectively. (In a separate paper
\cite {Go-Zh}, we will provide a more simple and direct proof
for the Fourier coefficients of $ \Phi $ when $F=\BQ$.) 

The absence 
of  a cuspidal divisor representative for
$\xi _i$ and the absence of Dedekind's $\eta$-function in the
general  case cause some essential
difficulties in our height computation. 
Fortunately, these difficulties can be overcome by using 
Arakelov theory and the strong multiplicity-one argument. See Section~4 for a detailed
explanation of our method.
Even in the case $X=X_0(N)$, our method simplifies the 
computation of Gross and Zagier. 
\enddemo

\demo{Acknowledgment} Modulo the construction of the map from Shimura\break curves to their 
Jacobians, the formula in Theorem C was first conjectured  by Gross
in \cite{G4}. In some cases of $F=\BQ$, K. Keating \cite{Ke} and D.\ Roberts
\cite{R} have made some computations of local intersection numbers of 
some CM-points based on desingularizations.  See also S. Kudla's paper
\cite{Ku}.

I would like to thank D. Blasius, P. Deligne,
B. Gross, and G. Shimura for their helpful   correspondence, and  
K. Keating and D. Roberts for sending me their papers.
I would also like to thank the referees for their 
  kind   comments and suggestions which made the paper more readable.
Finally, I would like to thank  D. Goldfeld and H. Jacquet for their
constant support. \enddemo

\demo{Notation}
\begin{itemize}
\item $\BN _F$: the multiplicative monoid of nonzero ideals of $\CO_F$.
\lb{BN_F}
\item For any ideal $  m $ of $\CO_F$, we define $\varepsilon (  m )$
\lb{epsilon(m)} such that
$\varepsilon $ is multiplicative on $\BN_F$ and such that $\varepsilon (\wp)
=\varepsilon _\wp (\pi )$ if $\wp$ is unramified in $ E$ and $\pi $ is a  uniformizer
of $\wp$ in $\CO _\wp$; otherwise, $\varepsilon (\wp)=0$.
\item Let $D_F$ \lb{D_F} denote the inverse different ideal of $F$, 
$$D_F^{-1}=\{x\in F: \quad \tr _{F/\BQ}(x\CO _F)\subset \BZ\}$$
and let $D_E$ denote the relative discriminant of $E$ over $F$.
\item   Let $d_F$, $d_E$, $d_N$ denote the absolute norms of $N$,
$D_F$, and $D_E$.\lb{d_F}
\item For a quaternion algebra we let $\det$\lb{det} (resp.\ $\tr$) denote 
the reduced norm map (resp.\ reduced trace map). For an order in a
quaternion algebra, we call {\em the reduced discriminant} simply 
a  {\em discriminant}.
\end{itemize}
\enddemo

\section{Shimura curves}

In this section we introduce  some of the theory of Shimura curves
 which will be used in later sections. 
We start from  the construction of the integral model for general 
Shimura curves through a moduli interpretation in subsections 1.1. and  1.2.
Then we give a description of the set of special fibers in  1.3.
In 1.4, we study Hecke operators and their reductions. After some modular interpretations, we  prove the
Eichler-Shimura congruence relation in a special case. 
Finally in  1.5,  we move to the special Shimura curve $X$ constructed in 
the introduction and   define the order $R$ and its corresponding level structure.

\demo{{\rm 1.1.} Modular interpretation} 
\enddemo

 1.1.1. {\it General properties of Shimura curves}.
Let $F$\lb{F} be a totally real field of degree $g$\lb{g}. 
This means that all Archimedean places  of $F$ are real.
Fix a real place $\tau$\lb{tau} which allows us to consider $F$ as a subfield of 
$\BR$ by the embedding which we still denote by $\tau$.
Let  $B$\lb{B} be a quaternion algebra over 
$F$ which is unramified at $\tau$ but not at the other infinite 
places. Then we can fix an isomorphism
\begin{equation}\lb{E:BtR}
B\otimes \BR\simeq \RM_2(\BR)\oplus \BH ^{g-1}\spq{1.1.1}\end{equation}
where the first factor corresponds to $\tau$, and ${\BH}$\lb{BH} is the quaternion division algebra 
over $\BR$. See \cite{V}
for basic properties of quaternion algebras.
Let ${\CH^\pm}$\lb{CHpm} denote the Poincar\'e double-half plane 
equal to $\BC-\BR$ equipped with the usual  action by $\GL _2(\BR)$. 
Thus the first projection in \eqref{E:BtR} gives an action of $B$ on $\CH^\pm$.

For   each open subgroup ${K}$\lb{K} of $\wh  B^\times$ which is compact
modulo $\wh F^\times$, we have a Shimura curve 
\begin{equation}
{M_K(\BC)}= B^\times\backslash \CH ^\pm \times \wh  B^\times/K, \spq{1.1.2}
\end{equation}
where for any abelian group $M$, $\wh M$ denote the completion $M\otimes \prod_p\BZ _p$.
For any $g\in \wh  B ^\times$, and open subgroups $K_1,    K_2$ such that
$g K_1g^{-1}\subset    K_2$, the right multiplication on $(B\otimes\wh F)^\times$ by
$g^{-1}$ induces a morphism  $g: M_{K_1}(\BC)\to M_{K_2}(\BC)$.
By Shimura's theory (see \cite{D}), the   curve $M_K(\BC) $ has a canonical
model $M_K$\lb{MK} defined over $F$ and the   
 morphism  $g: M_{K_1}\to M_{K_2}$ is  also defined over $F$
with respect to these models. 

By   work of Drinfeld and Carayol \cite{B-C}, \cite {C},
\cite {Dr}, one can even define an integral  model ${\CM _K}$ over 
$\Spec \CO_F$ such that  $\CM_K$ is regular if $K$ is sufficiently small.
This is what we need for the computation of heights in Sections~4 and 5.

If $F=\BQ$, then $M_K(\BC)$ parametrizes elliptic curves or abelian surfaces.
The canonical models and integral models can be obtained by extending
the corresponding modular problems to integers. See \cite{K-M} and 
\cite{B-C} for details.

If  $F\ne \BQ$, $M_K(\BC)$ does not parametrize abelian varieties
in a convenient way. But $M_K(\BC)$ has a finite map to another Shimura 
curve $M_{K'}(\BC)$ which apparently parametrizes abelian varieties.
Thus extending the moduli problem to integers gives the integral models.
 In the following we will describe the curve $\CM_{K'}$ and its
moduli interpretation.

Let us fix a quadratic extension ${F'}=F(\sqrt{\lambda})$  of $F$, where 
$\lambda$\lb{lambda} is a negative integer. Consider $\sqrt\lambda$ as an element in
$\BC$. Then  $\tau$\lb{tau'} can be extended to a complex place for $F'$:
\begin{equation}\lb{E:tauF'} \tau (x+y\sqrt \lambda)=
\tau (x)+\tau (y)\sqrt\lambda.\spq{1.1.3}\end{equation}

Let ${B'}$\lb{B'} denote $B\otimes F'$,
let ${J}$\lb{J} be a compact open subgroup of $\wh F^{'\times}$,
and let ${K'}$\lb{K'} denote the subgroup $K\cdot J$ of $\wh {B'}^\times$.
Then we have a Shimura curve
\begin{eqnarray}
M_ {K'}(\BC)&
=&F^{'\times}B^\times\backslash
\CH ^\pm \times B^\times \wh F^{'\times}/K' \spq{1.1.4}\lb{E:MK'}\\
&=&M_ K  (\BC   )\times_{\wh F^\times} 
\left[(F')^\times\backslash \wh{F'}^\times/ J \right].\nn
\end{eqnarray}
Again by Shimura's theory, this curve has a canonical model $M_{K'}$\lb{MK'} over 
$F'$ and  the  morphism 
\begin{equation}\lb{E:MKK'}M_ K (\BC  )\to M_ {K'}(\BC  )\spq{1.1.5}
\end{equation}
is defined over some extension of $F'$. For example, we have the 
abelian extension corresponding to $ J $ via class field theory.
The image of the morphism in \eqref{E:MKK'} is another Shimura curve
$M_{\wt K}$ where 
$$\wt K=K\cdot\left[\wh F^\times \cap \left( F^{'\times}\cdot J\right)\right].$$
\enddemo

\demo{{\rm 1.1.2.} A moduli problem over $F'$}\lb{mpF'}
In the following we will explain how the curve $M_{K'}(\BC)$ 
parametrizes certain
abelian varieties over $F'$ and, therefore, has a model $M'_   K  $
defined over $F'$.

For this we write  ${V}$\lb{V} for  $B'$ as a left $B'$-module and write 
${V_\BR}$\lb{VBR} for $V\otimes \BR$. Then there is a decomposition:
$$V_\BR=(B\otimes \BR)\otimes _\BR{F'\otimes \BR}=(\RM_2(\BR)\otimes\BC)
\oplus (\BH \otimes \BC)^{g-1},$$
where we use \eqref{E:BtR}) and \eqref{E:tauF'} with $\tau$ replaced by all
places of $F$. Now we define a 
complex structure on $V_\BR$ such that $\sqrt {-1}$ acts 
on $V_\BR$ by  right multiplication of the following
element ${j}\in B\otimes \BC  $\lb{j}:
$$j=\left ( \left( \begin{array}{cc} 0&1\\ -1&0 \end{array}\right), 
1\otimes \sqrt{-1}, \cdots, 
1\otimes \sqrt{-1}\right).$$
Then the space $\CH ^\pm$ can be identified with the 
$(B\otimes \BR)^\times \cdot (F'\otimes \BR)^\times$-conjugacy classes of 
$j$:  each $z=x+yi\in \CH ^\pm$ corresponds to
an element given  by 
\begin{equation}j_z\lb{j_z}=\left (\alpha _z \left( \begin{array}{cc} 0&1\\ -1&0 \end{array}\right)\alpha
_z^{-1},  1\otimes \sqrt{-1}, \cdots, 
1\otimes \sqrt{-1}\right)\spq{1.1.6} \end{equation}
where $\alpha _z$ is an element of $\GL _2(BR)$ such that its action on $\CH ^\pm$ gives
$\alpha _z(\sqrt{-1})\break =z$.

Thus $V_\BR$ is a $\BC$ vector space with an action by
$B'$. The traces of elements $\ell$ of $B'$ acting  on 
the $\BC  $-space $V_\BR$ are 
given by the following formula:
$$\tr (\ell, V_\BR/\BC  )={t(\ell)}\lb{tell}$$
where $t$ is a map $t: B'\to F'$ given by
\begin{equation}\lb{E:tell}
t(\ell )=2\tr _{F/\BQ}(x)+2\left(\tr _{F/\BQ}(y)-y\right)\sqrt\lambda \spq{1.1.7}
\end{equation}
if $\tr_{B'/F'}(\ell)=x+y\sqrt\lambda.$
The function $t$ characterizes $(V_\BR, j)$ uniquely in the sense that a complex 
$B'$-module $W$ is $B'$-linearly isomorphic to $(V_\BR, j)$ if and only if 
$\tr (\ell, W)=t(\ell)$ for every $\ell\in B'$.

Let $v\to \bar v$ denote the product of the involutions on both factors
of $B' =B\otimes _FF'$, and let ${\delta}$\lb{delta} be  a symmetric 
($\bar\delta=\delta$) and invertible element in $B'$. 
Let $\ell\to \ell ^*$ be an anticonvolution on $V$ defined by
$$\ell ^*=\delta^{-1} \bar \ell \delta .$$
Notice that every  anticonvolution of $B'$ which extends the convolution on $E$
can be obtained in this manner.

Let ${\psi_F}$\lb{psiF} be a pairing  on $V$ with values in $F$ given by 
$$\psi _F(u,v)=\tr _{B'/F}\left(\sqrt\lambda u\bar v\delta\right).$$
Then for any $\ell \in B'$,
$$\psi _F(\ell u, v)=\psi _F(u, \ell ^* v).$$
One can show that the  similitudes of $\psi _F$ consist  of 
right multiplication on $V=B'$ of elements of $B^\times \cdot F^\times$. 

Choose a $\delta$ such  that $\psi _F(v, vj)\in F\otimes \BR$ 
is  totally positive for $v\in V_\BR$.  \enddemo

\vglue9pt 

 \nnp{Proposition 1.1.3} 
The curve $M_{K}'$  is the coarse moduli space of the 
following moduli functor
${\CF^0_{K'}}$\lb{CF^0K'} over $\BC  ${\rm :} For an $F'$\/{\rm -}\/scheme $S${\rm ,}
 $\CF^0_{K'}(S)$ is the set  of the isomorphism classes
of objects $[\bar A, \iota, \bar \theta, \bar \kappa ]$ where
\begin{itemize}
\ritem{1.} $\bar A$ is an  abelian scheme   over $S$ up to isogeny  with an
action $\iota :B'\to \End_S (\bar A)$ such that 
for any $\ell \in {B'}$ there is the equality
$$\tr (\iota (\ell ), \Lie \bar A)=t(\ell).$$
\ritem{2.} $\bar\theta$ is an $F^\times$\/{\rm -}\/class  of
 polarizations  $\theta: A\to A^\vee$ for $A\in \bar A$  such that for any $\ell \in B'${\rm ,}
 the associated Rosati involution takes $\iota (\ell)$
to $\iota (\ell^*)${\rm .}
\ritem{3.} $\bar \kappa $ is a   $K'$\/{\rm -}\/class  of  $B'$\/{\rm -}\/linear
isomorphisms $\kappa : \wh V\to \wh V(\bar A)$
which are\break $\wh F$\/{\rm -}\/symplectic similitudes{\rm ,} where 
$\wh V(\bar A)=\wh T(\bar A)\otimes \BQ$
with $\wh T (\bar A)=\prod T_p(\bar A)${\rm .}
This means that  each $\kappa\in \bar \kappa$ is symplectic between the 
form $\psi _A$ induced by a polarization $\theta\in \bar \theta${\rm ,}
and the form $\tr _{F/\BQ}(ua\psi _F)$ for some  $u\in \wh
F^\times${\rm ,} $a\in \det K'${\rm .}
\end{itemize}

\endproclaim

\demo{Proof}
Let $x$ be a point of $M_{K'}(\BC)$; we want to construct an element $[A, \iota, \bar\theta,
\bar\kappa]$ in $\CF _{K'}^0(\BC)$ as follows. Assume that $x$ is represented by
$(z, \gamma)$.
\begin{itemize}
\ritem{1.} $\bar A$ is the  abelian variety up to isogeny,
$$\bar A=\Lambda\backslash (V_\BR, j_z)$$
with  $\Lambda$ any lattice of $V$, 
where $j_A$ is the complex structure constructed as in \eqref{j_z}.
Thus,  $\wh V(A)=\wh V$. 
\ritem{2.} $\iota: B'\to \End (\bar A)$ is induced by  left multiplication of $B'$ on $V$.
\ritem{3.}  $\bar\kappa$ is the $K'$-class of the
map $\wh V\to \wh V(\bar A)$ induced by right multiplication of $\gamma$. 
\end{itemize}
It follows from the definition that the isomorphic class $[\bar A, \rho, \bar\theta, \bar\kappa]$
is an element of $\CF _{K'}^0(\BC)$.

Conversely, we can construct a point $x\in M_{K'}(\BC)$ from an element 
 $[A, \iota, \bar\theta, \bar\kappa]$  of $\CF ^0_{K'}(\BC)$ as follows.
Let $V_A$ denote $H_1(A, \BQ)$ and let $\psi_A$ be an alternative form defined 
by one polarization in $\bar \theta$. Then $V_A$ is a $B'$-algebra which is isomorphic
to $V$ at each place of $F'$ by a map in $\bar\kappa$.   It follows that $V_A$ must 
be isomorphic to $V$. We may identify $V_A$ with $V=B'$ by fixing such an isomorphism.
By the second condition, the alternative form $\psi _A$ has the form
$$\psi _A(v_1, v_2)=\tr _{F/\BQ}\psi _F(v_1b, v_2)$$
for some $b'\in B^\times$. If $\kappa\in \bar \kappa$ then $\kappa$ is induced by the map
$v\to v\gamma$ with $\gamma\in \wh B^{'\times}$. Condition 3 implies
\begin {equation}\lb{E:psi_FA} \tr _{F/\BQ}(u\psi_F (v_1x, v_2x))=\tr _{F/\BQ}(\psi _F(v_1b,
v_2)). \spq{1.1.8}
\end{equation}
This is equivalent to the equation $ux\bar x=b$. By Hasse's principal
(see \cite[\S2.2.3]{Kn}), this equation must have a solution  $x\in B^{'\times}$.
After modifying  the  isomorphism $\phi: V_A\to V$, we may assume that
$b=1$. Then equation \eqref{E:psi_FA} implies that $\gamma\in \wh B^\times \cdot \wh {F'}^\times$.
Such a  $\gamma$ is uniquely determined modulo  right multiplication of $K'$
once $\phi$ is fixed. We may replace $\phi$ by $b\phi$ with $b\in B^\times\cdot F^{'\times}$
which acts on $V$ by left multiplication. Then $\gamma$ is changed to $b\gamma$.

Let $j_A\in \GL _\BR (V_\BR)$ be multiplication of $\sqrt{-1}$ in the complex structure on
$\Lie (A)$. Then $j_A$ commutes with the action of $B'_\BR$ so that  it is given by   right
multiplication of an element which we still denote by $j_A$. Since $j_A$ preserves the alternative
form $\psi _A$ or equivalently the form $\psi _F$, we see that $j_A\in B_\BR ^\times \otimes \BR$.
Now,
$$j_A=(\alpha _1\otimes \beta _1, \cdots, \alpha _g\otimes \beta _g)$$
where for each $i$, either $\alpha _i=1, \beta _i=\sqrt{-1}$ or $\alpha _i^2=-1, \beta _i=1$.
By computing the trace of $B_\BR$ over $V_A$ which must satisfy
 condition 1, we see that 
$j_A$ must have the form
$$j_A=(\alpha \otimes 1, 1\otimes \sqrt{-1}, \cdots, 1\otimes \sqrt{-1})$$
with $\alpha ^2=-1$. Now $\alpha$ must be conjugate to $ \left( \begin{array}{cc}0&1\\ -1&0 \end{array}\right)$
so that $j_A$ must be $j_z$ for some $z\in \CH ^\pm$. Again this $z$ is unique if $\phi: V\to V_A$
is fixed. If we change $\phi$ to $b\phi$ then $z$ is changed to $b(z)$. It follows that
the image $x$ of $(z, \gamma)$ in $M_{K'}(\BC)$ is a well-defined point.
\enddemo

\vglue9pt
 
\demo{{\rm 1.1.4.} Second version}
We may also describe $M_{K'}$ as a coarse moduli space of abelian varieties,
rather than abelian varieties up to isogeny.
For simplicity, we assume that $K$ is compact. 
 Then there is a  maximal order  $\CO _B$\lb{COB} of $B$ such that 
$K$ is included in $\wh \CO _B^\times $.
(Notice that this is not the exact case wanted, as the Shimura curve $X$ 
defined in the  introduction is the
compactification of $M_K$ with $K=\wh R^\times \cdot \wh F^\times$.)
 Let $\CO _{B'}$\lb{COB'} be the  order $\CO _B\otimes \CO _{F'}$
of $B'$. Write  $V_\BZ$\lb{VBZ} for  $\CO _{B'}$ as a 
left $\CO _{B'}$-module. 

Let  $\CU$\lb{CU} be a subset of  $\wh F^\times$ representing
$$F\backslash\wh F^\times/\det  K')$$
such that for each $u\in \CU$, the alternating pairing
$u\psi _F$ is integral on $V_\BZ$.
Let $\nu (K')$\lb{nuK'} denote $\det K'\cap F^\times$ of $\CO _F^\times$. \enddemo

\spn{1.1.5} \proclaim{Proposition}\lb{CF}
The functor $\CF^0 _{K'}$ is isomorphic to $\CF _{K'}$\lb{CFK'} defined as follows\/{\rm :} 
 For an $F'$\/{\rm -}\/scheme $S${\rm ,} $\CF_{K'}(S)$ is the set  of  isomorphism classes
of objects $[A, \iota, \bar\theta, \bar \kappa ]$ where\/{\rm :}\/
\begin{itemize}
\ritem{1.} $A$ is an abelian scheme   over $S$   with an action $\iota :\CO _{B'}
\to \End_S (A)$ such that for any $\ell \in \CO_{B'}$ there is the equality
$$\tr (\iota (\ell ): \Lie A)=t(\ell).$$
\ritem{2.} $\bar\theta$ is a $\nu (K')$\/{\rm -}\/class of  polarizations  $\theta: A\to A^\vee$
 such that for any $\ell\in \CO _{B'}${\rm ,}  the associated Rosati involution takes $\iota (\ell)$
to $\iota (\ell^*)${\rm .}
\ritem{3.} $\bar \kappa $ is a $K'$\/{\rm -}\/class of  $\CO _{B'}$\/{\rm -}\/linear
isomorphisms $\kappa : \wh V_\BZ   \to \wh T(A)$
which is symplectic with respect to $\psi _{u, a}:=\tr _{\wh F/\wh \BQ}(ua\psi _F) $
for some $u\in \CU$ and $a\in \det K'${\rm .}
\end{itemize}

\endproclaim

\demo{Proof} There is an obvious morphism from $\CF _{K'}$ to $\CF _{K'}^0$. Now we want to
define its converse. Let $[\bar A, \rho, \bar\theta, \bar\kappa]$ be an object in 
$\CF _{K'}^0(S)$. Then  the lattice $\kappa (\wh V_\BZ)$ does not depend on the choice
of $\kappa\in \bar\kappa$. Let $A$  be the corresponding abelian variety isogenous to
$\bar A$. 
Then $A$ has the action by $\CO _{B'}$ such that condition 1 is satisfied. As
$\kappa $ varies in $\bar\kappa =K'\kappa $ and $\theta$ varies in $\bar\theta=F^\times \theta$,
$u$ in condition 3 in Proposition 1.1.3 varies in a single double-coset of
$F^\times\backslash \wh F^\times /\det K'$. Thus we may choose a $\theta_0\in \bar \theta$
such that $u\in \CU$, and the set of such $\theta$'s forms a class $\nu (K')\theta_0$.
As $\tr _{F/\BQ}(ua\psi _F)$ is always integral, $\psi _A$'s corresponding to
$\theta_0\in \nu (K')\theta _0$  take integral values on $T(A)$. It follows that
every such $\theta _0$ defines a polarization of $A$. Condition 2 in this proposition 
is obviously satisfied. This defines a
morphism $\CF _{K'}^0\to \CF_{K'}$ which is 
obviously the inverse of the obvious morphism $\CF_{K'}\to \CF _{K'}^0$.
\enddemo 

\vglue9pt

\demo{{\rm 1.1.6.} Remark}
Let  $x\in M_{K'}(\BC)$ be represented by $(z, \gamma)$. 
From the proof of Propositions 1.1.3 and \ref{CF}, 
we see that the  object
$[A, \iota, \bar\theta, \bar\kappa]$ in $\CF _{K'}(\BC)$  
parametrized by $x$ has the following form: 
\begin{itemize}
\ritem{1.}
\hfil $A=V_\BZ \gamma^{-1} \backslash (V_\BR, j_z).$\hfill
\ritem{2.} $\iota$ is induced by left multiplication by $\CO_{B'}$
on  $V_\BZ\gamma ^{-1}$.
\ritem{3.} $\bar\theta$ is the unique class induced by  alternative
forms 
$$\left\{ \tr _{F/\BQ}(t\psi _F): 
\quad t\in F^\times\cap (\coprod_{u\in \CU}  u\det \gamma K')
\right\}.$$
\ritem{4.} $\bar\kappa$ is the $K'$-class of the morphism $\wh V_Z\to \wh V_Z\gamma ^{-1}$
induced by right multiplication by $\gamma ^{-1}$.
\end{itemize}

\enddemo

\spn{1.1.7} \proclaim{Proposition}
When $K'$ is sufficiently small{\rm ,} then $\CF _{K'}$ {\rm (}\/therefore $\CF _{K'}^0${\rm )}
is representable{\rm .}
\endproclaim

\demo{Proof}
For each $u\in \CU$, let $\CF _{K', u}$ denote the subfunctor of 
$\CF _{K'}\otimes F'\wt F$ with given  $u$ in condition 4 of
Proposition \ref{CF}, where $\wt F$ is the extension of $F$ corresponding to $F^\times \det K'$
via class field theory. Wishing  to show that $\CF _{K', u}$ is representable, we need the following:

\spn{1.1.8} \proclaim{Lemma} There is   a positive
integer $n$ such that
$$(1+m\CO _F)^\times :=(1+m\wh \CO _F)^\times \cap F^\times 
\subset [(1+m\CO _F)^\times ]^2$$
with some $n\ge 3${\rm .}
\endproclaim

\demo{Proof}
We fix an $n\ge 3$ and let $S$ be a finite subset of 
$(1+n\CO _F)^\times$ which contains $1$ and represents the quotient
$$(1+n\CO _F)^\times/[(1+n\CO _F)^\times]^2.$$
For each $s\in S-\{1\}$, let $p_s$  be a prime not dividing $2m$  such that 
$s$ is not a square in $(\CO _F/p_s\CO _F)^\times$. Then 
$$m:=n\prod _{m\in S-\{1\}}p_s$$
will satisfy the requirement. Indeed, the definition of $m$ implies that the morphism
$$\frac {(1+n\CO _F)^\times}{[(1+n\CO _F)^\times ]^2}
\to \frac {(\CO _F/m\CO _F)^\times }{[(\CO _F/m\CO _F)^\times ]^2}$$
is injective. Thus $(1+m\CO _F)^\times$ is included in $[(1+n\CO _F)^\times ]^2$.
\enddemo

We return to our proof of Proposition \ref{CF}. Assume that $K'$ is sufficiently small 
so that 
$$\det K'\subset (1+m\wh\CO _B)^\times\cdot (1+m\wh \CO _{F'})^\times.$$
Then $\nu (K')\subset (1+n\CO _{F})^{\times 2}$. 
Let $T$ denote the set of elements in $(1+n\CO _F)^\times$ whose squares are in
$\nu (K')$. Let $\wt K'$ denote $K'\cdot T$. If $t\in T$, then 
 multiplication by $t$
induces an isomorphism 
$$[A, \iota, \bar\theta, \bar\kappa]\to [A, \iota, \bar\theta, t\bar\kappa]$$
of objects in $\CF _{K', u}(S)$. Here if $\kappa$ is symplectic with respect to 
$\psi _{u, a}$, then $t\kappa $ is symplectic with respect to $\psi _{u, at^2}$.
As $T^2=\nu (K')=\nu (\wt K')$, it follows that the canonical morphism
$\CF _{K', u}\to \CF _{\wt K', u}$ is an isomorphism.

Let $\wt \CF _{\wt K', u}$ denote the functor defined in the same way as $\CF _{\wt K', u}$
but with $\nu (K')$-class $\bar\theta$ to replace a single $\theta$.
Then   multiplication by $t$ induces an isomorphism 
$$[A, \iota, \theta, \bar\kappa]\to [A, \iota, t^2\theta, \bar\kappa]$$
of objects in $\wt\CF_{\wt K, u}(S)$. So the canonical morphism
$\wt\CF_{\wt K', u}\to \CF _{\wt K', u}$ is also an isomorphism.

In this way we have shown that $\CF _{K', u}$ is isomorphic to $\wt\CF_{\wt K', u}$. 
Now we want to show the representability of $\wt \CF _{\wt K', u}$. 
Let $d$ denote the degree of $\psi _{u, 1}$. Let $\CA$ denote the moduli functor 
which classifies abelian varieties of dimension $4g$, with a full level $n$ structure
and a polarization of degree $d$. As $n\ge 3$, $\CA$ is representable by a scheme
$M_{d, n}$ (\cite[Prop.\ 7.9]{Mu}). The functor $\wt \CF _{\wt K', u}$ has a finite 
morphism to $\CA$. The conditions in the definition of $\wt\CF _{\wt K', u}$
defines a finite scheme $\wt M_{K', u}$ over $M_{d, n}\otimes F'\wt F$ which represents 
$\wt\CF _{\wt K', u}$, also $\CF _{K', u}$. Now the union $\wt M_{K'}$ of $\wt M_{K', u}$
represents $\CF _{K'}\otimes F\wt F$. Notice that $\wt M_{K'}$ has an action
by $$\Gal (\wt F/F)=F^\times \backslash \wh F^\times /\det K'$$
which induces a model $M_{K'}$ of $\wt M_{K'}$ defined over $F'$.
This model represents $\CF _{K'}$. As $M_{K'}(\BC)$ is a smooth Riemann surface,
$M_{K'}$ is a regular scheme.\hfill\qed
\enddemo

\demo{{\rm 1.2.} Integral models}
\enddemo

 1.2.1. {\it A new version of $\CF _{K'}\otimes F_\wp$}.
Let $\wp$\label{wp} be a prime of $F$ of  characteristic~$p$.
Assume that $\lambda $ is prime to $p$ and 
and that $\sfrac \lambda p$
is $1$; we  fix a square root $\mu _p$\label{mup} in $\BQ _p$.
 Then  $F'$ can be embedded into $F_\wp$ over $F$ by sending 
$\sqrt \lambda$ to $\mu _p$. 
We want to extend $M_{K'}\otimes F_\wp$ to a model $\CM_{K',\wp}$
over $\CO _\wp$, the ring of integers in $F_\wp$.
For this we need a new version of  the moduli
problem $\CF_{K'}$ over $F_\wp$-schemes. We start with some 
notation.

The algebra $\CO _{F', p}=\CO _{F'}\otimes \BZ_p$ is the sum of all 
completions  $\CO _{F',  q }$ at its places $ q $ over $p$. 
We have the following decomposition:
$$\CO _{F', p}=\CO _{F', p}^1\label{^1}+ \CO _{F', p}^2\label{^2}$$
where $\CO _{F', p}^1$  (resp.\ $\CO _{F', p}^2$)
denotes the   sum of all completions 
$\CO _{F',  q }$ such that the map $\CO_{F'}\to \CO _{F',q}$
takes $\sqrt \lambda$ to $\mu _p$ (resp.\ $-\mu _p$).
For any $\CO _{F', p}$ module $M$, we let $M^1$ (resp.\ $M^1_\wp$, $M^2$, $M^2_\wp$) 
denote $\CO _{F', p}^1 M$ (resp.\  $\CO _{F', p}^2 M$). Let $\CO _{F', p}^{\wp}$
\lb{cdot^wp} denote  the sum of components of $\CO^{\phantom{|}}_{F', p}$
not over $\wp$ and let $M^{\wp}$ (resp.
$M_\wp$) 
denote $\CO _{F', p}^{\wp}M$
(resp.\ $\CO _{F', \wp}M$).

Choose $\delta$ and $\CU$ such that for each $u\in \CU$, $u\psi _F$  has  degree prime to $p$,
and  $u$ has component $1$ at places dividing $p$.
Then we have  the following: \enddemo

\spn{1.2.2} \proclaim{Proposition}\label{CF'} The functor $\CF _{K'}\otimes F_\wp$
 is equivalent to the following functor
$\CF _{K', \wp}${\rm :} for any $F_\wp$-scheme $S${\rm ,} $\CF _{K', \wp}(S)$ is the isomorphism
classes of objects $[A, \iota, \bar\theta, \bar\kappa _p, \bar\kappa ^p]$ where 
\begin{itemize}
\ritem{1.} $A$ is an abelian scheme  over  $S$ with an 
action $\iota: \CO _{B'}\to \End (A/S)$ such that the following two
condition are satisfied\/{\rm :}\/
\begin{itemize}
\ritem{(a)} $\Lie (A)_\wp^2$ is a locally free $\CO_S$ module of rank $2$
such that the action of $\iota(F)$ is given by the 
inclusion $F\to F_\wp\to \CO _S${\rm ;}
\ritem{(b)} $\Lie (A)^{2, \wp}=0${\rm .}
\end{itemize}
Here $\Lie (A)$ may be viewed as an $\CO _{B', p}$
via the action $\iota${\rm .}
\ritem{2.} $\bar\theta$ is a $\nu (K')$\/{\rm -}\/class of  polarizations on
$A$ of degrees  prime to $p${\rm ,} such that 
 the Rosati involutions take
$\iota (\ell)$ to  $\iota (\ell^*)${\rm .}
\ritem{3.} $\bar \kappa ^2_p$ is a $K_p$\/{\rm -}\/class of  
$\CO _{B'}^p$-linear isomorphisms\/{\rm :}\/
$$\kappa_p ^2: V_{\BZ, p}   ^2\to  T_p(A)^2.$$
\ritem{4.} $\bar \kappa ^p$ is a $K  ^{'p}$\/{\rm -}\/class of $\CO _{B'} ^p$\/{\rm -}\/linear
isomorphisms
$$\kappa ^p: \wh V_\BZ   ^p\to \wh T(A)^p$$
which is symplectic with respect to some $\psi _{u, a}^p${\rm .} 
Here{\rm ,} for a $\wh\CO_F$\/{\rm -}\/module $M=\prod _q M_q ${\rm ,}  let $M^p$ denote the product of components
$M_q$ for $q\!\not| \  p${\rm .}
\end{itemize}

\endproclaim 

\demo{Proof} Let us first define a morphism from $\CF _{K'}\otimes F_\wp$ to $\CF _{K', \wp}$.
 Let
$[A, \iota, \bar\theta, \bar\kappa]$ be an object in $\CF _{K'}(S)$ where
$S$ is an $F_\wp$-scheme. Then we can decompose any $\bar\kappa$ into parts 
\begin{equation}\label {kappapp}\bar\kappa =\bar\kappa _p\oplus \bar\kappa ^p:
V_{\BZ, p}\oplus \wh V_{\BZ}^p\to T_p(A)\oplus \wh T(A)^p.\spq{1.2.1}
\end{equation}
Furthermore we can decompose $\kappa _p$ into two parts: 
$$\kappa _p^1\oplus \kappa _p^2: V_{\BZ, p}^1\oplus V_{\BZ, p}^2\to T_p(A)^1\oplus
T_p^2.$$
We claim that the object $[A, \iota, \bar\theta, \bar\kappa _p, \bar\kappa ^p]$ 
is an object of $\CF_{K' \wp}(S)$. We need only verify condition 1 in Proposition \ref{CF}.
Indeed, by a result of Carayol, condition 1 in Proposition \ref{CF}, 
which states that 
$$\tr (\iota (\ell), \Lie A)=t(\ell),\qquad \ell \in B',$$
can be replaced by the first condition in Proposition \ref{CF'} together with
one further condition that
\smallbreak
\centerline{\em  $A$ is an abelian variety of dimension $4g$.}
\medbreak

 This  is a slight generalization of Carayol's proposition in 
\cite[p.\  171]{C}. In his case $\wp$ is split in $B$. His proof can be generalized to our case
without any difficulty.  In this way, we obtain a morphism $\CF_{K'}\otimes F_\wp\to \CF _{K',
\wp}$. 

Now we want to  construct the  converse of the morphism of functors constructed as above.
 Since  the fact  that $A$ has dimension $4g$ is  implied by condition 4 in Proposition \ref{CF'}, 
 we need only show that
for a given $K'_p$-class $\bar \kappa_p^2$ as in Proposition \ref{CF'}, 
we can find $\bar\kappa _p$ with a decomposition as in \eqref{kappapp}   such that 
{$\bar\kappa _p$ is a $K'_p$-class of isomorphisms $\kappa _p: V_{\BZ, p}\to T_p(A)$ which is 
$\CO _{B', p}$-linear and symplectic with respect to $\tr a \psi _F$ and some $\psi _A$ induced
by a $\theta$ in $\bar\theta$, where $a$ is some element in $\det K'_p$.}

Notice that condition 2 in  Propositions 1.1.5 and 1.2.2 implies that all these subspaces are 
null spaces under symplectic forms. So each pair of these spaces forms a complete
dual. Now we may take  $\kappa _p^1$ to be the dual of $\kappa _p^2$.
\enddemo

\demo{{\rm 1.2.3.} Definitions} Let $S$ be a scheme over $\CO _\wp$,  $\CG$
 an $\CO _{B, \wp}$-module scheme over $S$.
\begin{itemize}
\ritem{1.} We say that  $\CG$ is a {\em special $\CO _{B, \wp}$-module}\label{specialmodule} if the induced
action of
$\CO _{B, \wp}$ on $\Lie (\CG)$ makes $\Lie (\CG)$   a locally free module
of rank one over $\CO _S\otimes _{\CO _\wp}\CO _{E, \wp}$, where $\CO _{E, \wp}$
is any unramified quadratic extension of $\CO _\wp$ contained in $\CO _{B, \wp}$.
\ritem{2.} Let  $n\in \BN$ and $x\in \CG[n] (S)$. We say $x$ is a {\em Drinfeld base of $\CG$
of level $n$}\label{Drinfeldbase} if, as cycles in $\CG$,  there exists the identity:
$$\left[\CG[n]\right]=\sum _{a\in \CO _{B, \wp}/\wp^n}[nx].$$
\end{itemize}

\enddemo

\spn{1.2.4} \proclaim{Proposition}\label{CF''} \hskip-8pt Assume that $J$ is maximal at all places dividing~$p${\rm .}
Then the functor 
$\CF _{K', \wp}$ can be extended to the following functor over $\CO _\wp$ which
is still denoted by $\CF _{K', \wp}${\rm :}  For any $\CO_\wp$-scheme $S${\rm ,} 
$\CF _{K_0',\wp}(S)$ is the set of isomorphism 
classes of objects  $[A, \iota, \theta, \bar x, \bar \kappa ^{2, \wp}, \bar \kappa ^p]$ where
\begin{itemize}
\ritem{1.} $A$ is an abelian scheme  over the scheme $S$ with an 
action $\iota: \CO _{B'}\to \End (A/S)$ such that the following two
condition are satisfied\/{\rm :}\/
\begin{itemize}
\ritem{(a)} $\CG :=A[\wp^\infty]^2$ is a special formal $\CO _{  B   , \wp}$\/{\rm -}\/module{\rm .}
\ritem{(b)} $A[p^\infty]^{2, \wp}$ is an \'etale $\CO _{B', _p}^{2, \wp}$\/{\rm -}\/module.
\end{itemize}
\ritem{2.} $\theta$ is a  $\nu (K')$-class of polarizations on
$A$ of degrees  prime to $p${\rm ,} such that 
 the Rosati involutions take
$\iota (\ell)$ to  $\iota (\ell^*)${\rm .}
\ritem{3.} $\bar x$ is a $K_\wp$\/{\rm -}\/class of Drinfeld bases of $\CG $ of level $n${\rm ,} 
where  $n$ is a
positive integer such that $K_\wp$ contains $1+\wp^n\CO_{B, \wp}${\rm .}
\ritem{4.} $\bar \kappa ^{2, \wp} $ is a $K_p^\wp$\/{\rm -}\/class of  
$\CO _{B', p} ^{2, \wp} $\/{\rm -}\/linear isomorphisms\/{\rm :}\/
$$\kappa ^{2, \wp}: V_{\BZ, p}   ^{2, \wp}\to T_p(A)^{2, \wp}.$$
\ritem{5.} $\bar \kappa ^p$ is a $K^{'p}$\/{\rm -}\/class of $\CO _{B'} ^p$\/{\rm -}\/linear
isomorphisms
$$\kappa ^p: \wh V_\BZ   ^p\to \wh T(A)^p$$
which is symplectic with respect to some $\tr (ua\psi _F)$ 
for some $u\in \CU$ and $a\in \det K'$\/{\rm ,}\/  and some $\psi _A$ induced by
some element in $\bar\theta$\/{\rm .}\/
\end{itemize}

Moreover{\rm ,} when $K_\wp=(1+\wp ^n\CO _{B, \wp})^\times$
and $K^{'p}$ is sufficiently small{\rm ,} the functor $\CF _{K', \wp}$ is representable 
by a regular scheme $\CM _{K',\wp}$\label{CMK'wp} over $\CO _\wp${\rm .} In general{\rm ,} the functor $\CF _{K', wp}$ 
has a coarse moduli space $\CM _{K',\wp}$ over $\CO _{\wp}${\rm .}
\endproclaim 

\demo{Proof} It is easy to see that the conditions in Proposition \ref{CF''} are equivalent to the 
conditions in Proposition \ref{CF'} when $S$ is a $F_\wp$-scheme. The representability can
be proved in the same way as in Proposition \ref{CF}. Here we need to choose $K^{'p}$ to be
sufficiently small and take care of Drinfeld bases. See \cite[\S5.3  and  \S7.3]{C}.
Also the regularity can be proved using the same argument as in \cite[\S5.4 and \S7.4]{C}.

If $K^{'p}$ is not sufficiently small
or if $K_\wp$ does not have the form\break 
$(1+\wp ^n\CO _{B, \wp})^\times$,
 then $\CF _{K', \wp}$ may not be representable.
But we may choose a sufficiently small normal subgroup $\wt K'$ of $K'$ so that the functor 
$\CF_{\wt K', \wp}$ is representable. The quotient $\CM _{\wt K',\wp}/K'$ does not depend on 
the choice of $\wt K'$ and it is actually the coarse moduli space of $\CF _{K', \wp}^{\phantom{|}}$
\enddemo

\demo{{\rm 1.2.5.} Modules $\CM_K$ and modules $\CG_K$}
Recall that $M_ K$ has a finite morphism to $M _{K'}$.
We, therefore, obtain a model $\CM _{   K  , \wp}$\label{CMKwp}
for the Shimura curve $M _   K  $ by taking the normalization
 of  $\CM _{   K  ', \wp}$ in $M_K$. 
 One can show that this model does not depend on the 
choice of $F'$ and $ J $. 
By gluing these models,   one has a model $\CM _K$\label{CMK} over $\Spec \CO _F$
for $M_K$,  which is regular when $   K  $ is sufficiently small.

Assume that  $K'$ is sufficiently small so that $\CF _{K', \wp}$ is represented
 by a regular scheme 
$\CM _{K', \wp}$. Then over $\CM _{K', \wp}$, we have 
a divisible $\CO _{B,\wp}$-module $\CG_{K'}:=\CG_{\CA}$ \label{CGK'}, where $\CA$ is the 
universal abelian variety on $\CM _{K'}$. 
 Let $ \CG_{K}$\label{CGK}  be the pull-back  of $ \CG   _\CA$
on $\CM _{K}$. Then $ \CG  _{K_0}$ does not depend on the choice of $J$.

Let $K_0$ denote $\CO _{B, \wp}^\times \cdot K^\wp$. Then 
the scheme $\CM _{K}$ over $\CM _{K _0}$ classifies  
the $K _\wp$-class of Drinfeld bases
in  $\CG   _{K_0}[\wp ^n]$.

Let $x$ be a geometric point of the special fiber of $\CM _{K_0}$.
Then over the completion of the strict localization $\wh{\CM _{K_0, x}}$,
$ \CG_{K_0}$ is the universal deformation of $ \CG _{K_0}|_x$.\enddemo

\demo{{\rm 1.3.} Reductions of models} We want to study the set of irreducible components of the special
fibers of $\CM _K$.

\demo{{\rm 1.3.1.} Split case}
Assuming that $B$ is split at $\wp$,  we can fix an isomorphism between $\CO _{B, \wp}$ and $\RM_2(\CO _\wp)$,
thus obtaining the decomposition of $\CO _\wp$ modules over 
$\CM _{   K  _0}$:
$$ \CG   _{   K  _0}= \CG   ^1\oplus  \CG   ^2,
\qquad  \CG   ^1= \left( \begin{array}{cc} 1&0\\ 0&0 \end{array}\right) \CG   _{   K   _0},
\quad  \CG   ^2= \left( \begin{array}{cc}0&0\\ 0&1 \end{array}\right) \CG   _{   K   _0}.$$
The two $\CO _\wp$-modules $ \CG   ^1$, $ \CG   ^2$ are isomorphic by 
the element $ \left( \begin{array}{cc}0&1\\ 1&0 \end{array}\right)$. 
It is easy to see that in this setting, $\CM _{   K   }$ over 
$\CM _{   K   _0}$ classifies the $   K   _\wp$-class of morphisms
$$\phi: \left (\CO _\wp/\wp ^n\right) ^2\to  \CG   ^1[\wp ^n]$$
such that this homomorphism is surjective on cycles.

Let $x$ be a geometric point in the special fiber of $\CM _{   K   _0, \wp}$. Then
the\break $\CO _\wp$-module $ \CG   ^1_x$ has two possibilities:
\begin{itemize}
\item[1.] {\em Ordinary case\/}: \label{ocase} The group $ \CG   ^1_x$ is isomorphic to 
the product of $(F_\wp/\CO _\wp)$ and a formal $\CO _\wp$ -module
$\Sigma _1$ of height $1$.
\item[2.] {\em Supersingular case\/}:\label{sscase} The group $ \CG   ^1_x$ is 
isomorphic to a formal $\CO _\wp$-module $\Sigma _2$ of height $2$.  
\end{itemize}

The set of connected geometric components of the special fiber of 
$\CM _{   K  _0}$ over $\wp$ is 
the  same as  that of the generic fiber. Fix a geometrically irreducible 
 component $D$ of the 
special fiber of $\CM _{   K  _0}$ over~$\wp$. Then we have:

\spn{1.3.2} \proclaim{Proposition} Assume that $\wp$ is split in $B${\rm .} Then the set of the irreducible 
geometric components of $\CM _   K  $ over $\wp$ is indexed by $\BP ^1(\CO _\wp)/   K   _\wp${\rm .}
More precisely{\rm ,} for each line $C\subset \CO _\wp^2${\rm ,} the corresponding
component of $\CM _   K  $ over $\wp$ will classify the morphism 
$$\phi : (\CO _\wp /\wp )^2\to  \CG   ^1[\wp^n]$$
such that $\ker \phi$ contains  $C\hskip-4pt \pmod \wp${\rm .}
\endproclaim

\demo{{\rm 1.3.3.} Nonsplit case}
It remains to study the reduction of $\CM _   K  $ in the case that 
$B$ is not split at $\wp$. In this case, one can show that 
$ \CG  _   K  $ is a formal group. It follows that the map
$$\CM _   K  \to \CM _{   K  _0}$$
is purely inseparable at the fiber over $\wp$.
So the set of irreducible components in the special fiber
of $\CM _ K$ over $\wp$ is the same as that of $\CM _{   K   _0}$.

To study the irreducible components of $\CM _{   K  _0}$ over 
$\wp$ we can use the uniformization theorem of Cerednik -- Drinfeld 
\cite{B-C}, but we need some notation.
Let $\wh \CM_{   K   _0}$ denote the formal completion of
$\CM _{   K   _0}$ along its special fiber over~$\wp$.
Let $B(\wp)$ denote the quaternion algebra over $F$ 
obtained by switching the invariants of $B$ at $\tau$ and $\wp$. 
Fix an isomorphism:
$$\wh {B(\wp )}\simeq \RM_2(F _\wp)\cdot \wh B^{\wp}$$
where the superscript $\wp$ means that the component at the place $\wp$
is removed.
Let $\wh \Omega $ denote Deligne's formal scheme over $\CO _\wp$
obtained by blowing-up
$\BP ^1$ along its rational points in the special fiber over the residue
field $k $ of $\CO _\wp$ successively. The generic fiber $\Omega $
of $\wh \Omega $ is a rigid analytic space over $F_\wp$  whose 
$\bar F _\wp$ points  are given by $\BP ^1(\bar F _\wp)-\BP ^1(F _\wp).$
The group $\GL _2(F_\wp)$ has a natural action on $\wh \Omega $. 
The theorem of Cerednik-Drinfeld gives a natural isomorphism 
$$\wh {\CM }_{   K   _0}\simeq 
B(\wp )^\times \backslash \wh \Omega \wh\otimes \wh \CO ^{nr}_\wp
\times \wh B ^{\times, \wp}/   K   ^\wp$$ 
where $\wh \CO ^{nr}_\wp$ denote the completion of the maximal unramified 
extension of $\CO _\wp$.

To obtain a description of the special fiber of $\wh\CM _{   K   _0}$, 
we notice that the irreducible components  of special fibers of $\wh\Omega $
correspond one-to-one to the classes modulo $F^\times$ of 
$\CO _\wp$ lattices in $F_\wp ^2$. Consequently, we have the following: 

\spn{1.3.4} \proclaim{Proposition} Assume that $\wp$ is not split in $B${\rm .} Then the set of 
irreducible geometric  components of $\wh \CM _{K_0}$ over $\wp$ 
is indexed by the set
 \begin{eqnarray*}&&
 B(\wp )^\times \backslash \GL _2(F_\wp)/F_\wp ^\times 
\GL _2(\CO _\wp)\times \wh B ^{\times, \wp}/   K   ^\wp
 \\ &&\hskip1in 
\simeq B(\wp )^\times \backslash \wh {B (\wp )}^\times 
/F_\wp^\times \GL _2(\CO _\wp) K^\wp.  \end{eqnarray*}
\endproclaim

\demo{{\rm 1.4.} Hecke correspondences}\label{S:Hc} 
\enddemo
1.4.1. {\it  Definition}.
Let $M_K$ be a Shimura curve with a compact 
$K$ contained in $\wh \CO _B^\times$.
Let $m$ be an ideal of $\CO_F$  such that  
at  every prime $\wp$ dividing $m$, $K$ has maximal components and $B$ is split.
Let $G_  m $\lb{G_m} (resp.\ $G_1$) be the set of element $g$ of $\wh\CO_  B    $ 
which has component
$1$ at places not dividing $  m $, and such that $\det (g)$ generates $m$ (resp.\ is invertible)
at each place dividing $  m $.
Then we may consider $G_1$\lb{G_1} as a subgroup of 
$K$. 
 The Hecke operator $\RT (m)$\lb{OT(m)} on $M_K$ is defined 
by the formula
\begin{equation}\label{OTm}
\RT (  m )x=\sum _{\gamma \in 
G_m/G_1}
[(z, g\gamma)], \spq{1.4.1}
\end{equation}
where  $(z, g)$ is a representative of $x$ in $\CH\times \wh B ^\times$,
and $[(z, g\gamma )]$ is the projection of $(z, g\gamma )$ on $X$. It is easy to see that
the correspondence 
has the degree
$$\deg \RT(  m )=\sigma _1(  m )=\sum  _{  a   |  m }\RN (    a   )\lb{sigma_1(m)}.$$

To see that $\RT (m)$ is a correspondence 
given by  algebraic cycles,  
decompose $G_  m $ into a union of  double cosets:
$$G_  m =\coprod G_1  g_iG_1.$$
For each $i$,
let $   K  _i$ denote the group 
$g_i   K   g_i^{-1}\cap    K  $.
Then we obtain two morphisms $p_1, p_2$ from $M_{   K  _i}$ to
$M_   K  $, induced by right multiplication on $\wh  B    ^\times$ by $1$ 
and $g_i$ respectively. The  image of 
$M_{   K  _i}$
in $M_   K  \times M_   K  $ by $(p_1, p_2)$, as an algebraic cycle,
 defines a correspondence $T_i$.
Then $\RT (  m )$ is defined to be the sum of $T_i$.

If $\CM  _   K  $ is the integral model constructed as before then
the Hecke correspondences  $\RT (  m )$ can be extended to $\CM _   K  $ by taking Zariski
closure of cycles in $\CM _K\times 
\CM _K$. See moduli interpretation in
the next section. \enddemo

\demo{{\rm 1.4.2.} Moduli interpretation}
Let $F'=F(\sqrt{\lambda})$ be a quadratic extension as in 1.1.1. 
Let $J$ be a compact subgroup of $\wh F^{'\times}$ which has maximal components 
for places dividing $m$. 
Let $K'=K\cdot J$. Then we can use
the same formula \eqref{OTm} to define
Hecke correspondence $\RT(m)$ on 
$M_{K'}$. In the following we want to 
describe a moduli interpretation for 
$\RT(m)$. 
\enddemo
 
\demo{{\rm 1.4.3.} Definition} 
Let $[A, \rho, \bar\theta, \bar\kappa]$ be
an object in $\CF _{K'}(S)$ as in 
Proposition \ref{CF},
let $m$ be an ideal of $\CO_F$,
and let $D$ be an $\CO _{B'}$-submodule of $A[m]$. We say
that $D$ is an {\em admissible submodule of level 
$m$} if the 
following conditions are satisfied:
\begin{itemize}
\ritem{1.} $D$ is its own annihilator under a
Weil pairing  \end{itemize}
\begin{equation}\label{E:Weilpair}
(\cdot, \cdot): A[m]\times A[m]\to   \oplus _{\ell|m}\CO_\ell/m\CO _\ell\spq{1.4.2}
\end{equation} \begin{itemize}
\item[]
induced by a polarization in $\bar\theta$.
\ritem{2.} $D^1$ and $D^2$ have the same order.
\end{itemize}
\enddemo

\spn{1.4.4} \proclaim{Proposition}\label{P:hecke} 
Assume that each prime factor $\ell$ of $m$
is split in both $B$ and $F'${\rm .}
Let $[A, \rho, \bar\theta,
\bar\kappa]$ be an object in $\CF _{K'}(S)${\rm .}
\begin{itemize}
\ritem{1.}
 Let $D$ be an admissible submodule 
of $A$ of level $m${\rm , }
let $A_D$ denote
the abelian variety $A/D${\rm ,}
and let $\rho _D$
denote the action of $\CO _{B'}$ on $A_D$ 
induced from that on $A${\rm .} Then there  are a unique
$\nu(K)$\/{\rm -}\/class $\bar\theta _D$ of polarizations on $A_D$ inside of
$F^\times\bar\theta${\rm ,} and a unique $K'$\/{\rm -}\/class
$\bar\kappa _D$ of level structure 
which have the same components as 
$\bar\kappa$ out side of $\ell$ such that
$[A_D, \rho _D, \bar\theta _D, \bar\kappa _D]$
defines an element in $\CF _{K'}(S)${\rm .}
\ritem{2.} The Hecke operator 
as a correspondence  acting on $M _{K'}$  is given by  the 
following formula\/{\rm :}\/
$$\RT(m)[A, \rho, \bar\theta, \bar\kappa]
=\sum _D[A_D, \rho _D, \bar\theta _D,
\bar\kappa _D]$$
where $D$ runs over all admissible submodule
of $A$ of level $m${\rm .}
\end{itemize}
\endproclaim

\demo{{P}roof} 
Choose a root $\mu _\ell$ of $\lambda$ in 
$F _\ell$ for each $\ell$ dividing $m$.
 Then we have an isomorphism
$$\CO _{F',   \ell }\to \CO _{  \ell }\oplus \CO _{  \ell }, \qquad
\sqrt\lambda\to (\mu _\ell, 
-\mu _\ell).$$
 Any $\oplus _{\ell|m}\CO _{F',   \ell }$-module $M$  has a corresponding decomposition 
$M=M^1+M^2$.\lb{cdot^1ell}

Since $\RT(m)$ is multiplicative 
for coprime $m$'s, we may assume that $m$
is a power of a prime  ideal $\ell$. \enddemo

\demo{{\rm 1.4.5.} Models for $\CO_{B', \ell}$
and $V_{\BZ, \ell}$} In the decomposition 
$$\CO_{B', \ell}=\CO_{B', \ell}^1\oplus 
\CO _{B', \ell}^2,$$
the Rosatti convolution switches two factors.
Now, we can fix  an isomorphism 
\begin{equation}\label{E:COB'-dec}
\CO_{B', \ell}=\CO_{B, \ell}\oplus 
\CO _{B, \ell}, \spq{1.4.3}\end{equation}
such that the following conditions are satisfied:
\begin{itemize}
\item  The second projection 
is the projection onto $\CO _{B', \ell}^2$
composing with the canonical isomorphism
$\CO _{B', \ell}^2\simeq \CO _{B, \ell}$.
\item The Rosatti operator is given
by 
$$(a, b)^*=(\bar b, \bar a).$$
\end{itemize}

Similarly, we fix a model for $V_{\BZ, \ell}$ as follows. First of all, since 
\begin{equation}\label{E:aa*}
\psi _F(ax, y)=\psi _F(x, a^*y), \spq{1.4.4}
\end{equation}
 it follows that in the decomposition 
$$V_{\BZ, \ell}=V_{\BZ,\ell}^1
\oplus V_{\BZ, \ell}^2,$$
$\psi _F$ has the form
$$\psi _F(x^1+x^2, y^1+y^2)
=\psi _F(x^1, y^2)-\psi _F(y^1, x^2)$$
for $x^i, y^i\in V_{\BZ, \ell}^i.$
It follows that $\psi _F$ gives a
perfect pairing between $V_{\BZ, \ell}^i$'s.
So we have an isomorphism 
\begin{equation}\label{E:VBZ-dec}
V_{\BZ, \ell}\simeq \CO_{B, \ell}
\oplus \CO _{B, \ell}\spq{1.4.5}\end{equation}
such that:
\begin{itemize}
\item The second projection is the projection
onto $V_{\BZ, \ell}^2$ composing with the 
canonical isomorphism 
$$V_{\BZ, \ell}^2\simeq \CO _{B, \ell}.$$
(Recall that $V_\BZ=\CO_{B'}$ in its definition.)
\item The pairing $\psi _F$ is given by \end{itemize}
\begin{equation}\label {E:psi-dec}
\psi _F((x^1, x^2), (y^1, y^2))
=\tr _{B/F}(\bar x^1y^2)-\tr_{B/F}(\bar y^1x^2) \spq{1.4.6}
\end{equation} \begin{itemize}\item[]
for $x^i, y^i\in \CO_{B, \ell}.$
\end{itemize}

With respect to the decompositions 
\eqref{E:COB'-dec} and \eqref{E:VBZ-dec},
the action of the second factor of $\CO _{B', \ell}$ 
 on $V_{\BZ, \ell}$ is given by left
multiplication on the second factor 
of $V_{\BZ, \ell}$. It follows   
from \eqref{E:aa*}, that the same is true 
for the first factor. 

Now we want to find a formula for another  action $\CO _{B, \ell}$
on $V_{B, \ell}$ which is originally given
by right multiplication in its definition.
Let us  denote this action by $r$.
Let $a\in \CO _{B, \ell}$. Since $r(a)$ 
is $\CO_{B', \ell}$-linear, $r(a)$  must
be given by right multiplication of some 
element $(a_1, a_2)$ of $\CO _{B', \ell}$
with respect to the decomposition
\eqref{E:COB'-dec}. From the
definitions of the decomposition, $a_2=a$. Recall that $r(a)$ is
a similitude of $\psi _F$:
$$\psi _F(r(a)x, r(a)y)=\det (a)\psi_F(x, y).$$
Combining this with \eqref{E:psi-dec}, we must
have $a_1=a$. So the action $r$ is still given
by right multiplication.
\enddemo

\demo{{\rm 1.4.6.} First statement}
 Let $\kappa$ be  one element
 in $\bar\kappa$. Then tensoring with $\BQ$
we obtain an isomorphism 
$$\kappa :\wh V\to T(A)\otimes \BQ.$$
Notice that
the natural map $A\to A_D$ induces inclusions
$$T(A)\subset T(A_D)
\subset T(A)\otimes \BQ.$$
We want to find  $\gamma\in G_m$
such that $\wh V_\BZ \gamma ^{-1}
=\kappa ^{-1}(T(A_D))$. Notice that 
 such a  $\gamma$ is unique modulo $G_1$
if it exists. We need only work at the place~$\ell$.

Let $W$ denote $\kappa ^{-1}(T_\ell(A_D))$.
Then $D$ is isomorphic to $W/V_{\BZ, \ell}$.
Notice that the Weil pairing on $A[m]$ is 
induced up to an invertible factor by the pairing
\begin{equation}
\alpha^2\psi _F: m^{-1}V_{\BZ, \ell}
\times m^{-1}V_{\BZ, \ell}
\to \CO _\ell,
\spq{1.4.7}
\end{equation}
where $\alpha\in \wh F^\times$ is a generator
of  $m$.
Since $D$ is its own annihilator, 
it follows that the pairing
$$\alpha \psi _F: W\times W\to \CO _\ell$$
is perfect. 

With respect to the decomposition
\eqref{E:VBZ-dec}, $W$ must have the form
$$W=\CO _{B, \ell}\gamma_1^{-1}
\oplus \CO _{B, \ell }\gamma _2^{-1}$$
with $\gamma _i\in \CO _{B, \ell}$. 
Notice that  $D_i$ is  isomorphic to $\CO_{B, \ell}\gamma_i^{-1}/\CO _{B, \ell}$
respectively. So $D_i$ has order
$(\det (\gamma _i))$. Since
$D^1$ and $D^2$ have the same order,
it follows that both $\det \gamma _1$
and $\det \gamma _2$ generate $m$.
Now as $\alpha \psi _F$ is perfect on $W$, 
$\gamma _2$ must be equal to $\gamma _1$ times
a unit. So $W=V_{\BZ, \ell}\gamma _1^{-1}$.

Let $\gamma $ be an element of
$G_m$ which has the component $\gamma _1$ at 
the place $\ell$; then we have 
$\kappa ^{-1}(T(A_D))=\wh V_{\BZ}\gamma ^{-1}$.
Now we can define $\bar\kappa _D$ as the class
of the composition
$$\kappa \circ \gamma^{-1} :
\wh V_{\BZ}\to \wh V_{\BZ}\gamma ^{-1}=
\kappa ^{-1}(T(A_D))
\to T(A_D).$$
As in the proof of Proposition \ref{CF}, there will be a unique 
class $\theta _D$ inside $F^\times \theta$
such that $[A_D, \rho_D, \bar\theta _D,
\bar\kappa _D]$ is an object in 
$\CF _{K'}(S)$.
\enddemo

\demo{{\rm 1.4.7.} Second statement}
Let $\gamma $ be an  element in $G_m$.
Recall that in the proof of Proposition 1.1.3, 
 if $[(z, g)]$ represents an object
$[\bar A, \rho, \bar\theta_{\bar A}, \bar \kappa]$
in $\CF_{K'}^0(\BC)$ then $[z, g\gamma]$ represents
the object $[\bar A, \rho, \bar \theta_{\bar A},
\bar\kappa \gamma^{-1}]$. Also recall that in the 
proof of Proposition \ref{CF} of the 
equivalence $\CF ^0_{K'}$ and $\CF _{K'}$,
these two objects are equivalent to 
$$[A, \rho, \bar\theta, \bar\kappa]
\quad{\rm and}
\quad
[A', \rho, \bar\theta', \bar\kappa\gamma^{-1} ]$$
where $A$ (resp.\ $A'$) is the abelian 
variety isogenous to $\bar A$ such that
$$\kappa(\wh V_\BZ)=T(A)
\quad \left({\rm resp.}\quad
\kappa \circ \gamma^{-1} (\wh V_\BZ))
=T(A')\right).$$
Here $\bar \theta$ (resp.\ $\bar\theta'$ )
is the unique $\nu (K')$-class inside
$\bar\theta _{\bar A}$ to make this  an object in $\CF _{K'}(\BC)$.

The inclusion $\wh V_\BZ\subset \wh V_\BZ\gamma ^{-1}$ induces an isogeny $A\to A'$ with 
kernel $D$ isomorphic to 
$$\wh V_\BZ\gamma ^{-1}/\wh V_\BZ.$$
We want to show that $D$ is admissible of
level $m$.
As the Weil pairing on $A[m]$ up to an invertible scale is induced by
a pairing as in
(1.4.7), it follows  easily that 
$D$ is its own annihilator.  Also $D_1$ and 
$D_2$ have the same cardinality as both of them
are isomorphic to $\CO _{B, \ell}\gamma ^{-1}
/\CO _{B, \ell}.$  Thus $D$ is admissible.

\smallbreak
By the first statement, $[A', \rho, \bar\theta',
\bar\kappa\gamma ^{-1}]$ is equal to
$[A_D, \rho, \bar\theta _D, \bar\kappa _D]$.
From the above arguments, one sees that 
the correspondence between admissible
submodules of level $m$ and $\gamma$'s in
$G_m/G_1$ is bijective. The second statement of
the proposition thus follows.\hfill\qed
\enddemo

\demo{{\rm 1.4.8.} Remarks} 
First of all,  we may 
 extend Definition 1.4.3
and Proposition \ref{P:hecke} to $\CM _{K', \wp}$ 
where $\wp$ is a prime of $F$.
 Indeed,   everything is exactly the same as above except when
 $\ell=\wp$. In this case,  we need  the following assumptions:
\begin{itemize}
\ritem{1.} Assume further that $\lambda$ is split
in $F_\wp$ and choose a square root $\mu _\wp$
in $F_\wp$. 
\ritem{2.} Assume the Weil pairing in \eqref{E:Weilpair} has the values in 
$$\Sigma _1 [m]\oplus\oplus _{\ell\not\mid m}m^{-1}\CO _\ell/m\CO_\ell$$ where $\Sigma _1$ is the formal $\CO
_\ell$-module  of height $1$.
\end{itemize}

Secondly, $D$ is uniquely determined by $D^2$
as $D^1$ is the annihilator of $D^2$ in 
$A[m]^1$.
Actually, the correspondence $D\longmapsto D^2$
gives a bijection between admissible submodules
of $A[m]$ and submodules of $A[m]^2$ of order
$m^2$.

Moreover, if each $\ell |m$ is split in $B$ then we may give a
further decomposition for $M^2$. For this  we fix an 
isomorphism $\CO _{B, \ell}
\simeq \RM_2(\CO _\ell)$. Then
$M^2$ has a decomposition:
$$M^2=M^{2,1} +M^{2,2}$$
where \lb{cdot^{2,1}} \lb{cdot^{2,2}}
$$M^{2,1}= \left( \begin{array}{cc}1&0\\ 0&0 \end{array}\right)M^2,\qquad
M^{2,2}= \left( \begin{array}{cc}0&0\\ 0&1 \end{array}\right)M^2.$$
The element $w= \left( \begin{array}{cc} 0&1\\ -1&0 \end{array}\right)$
switches $M^{2,1}$ and $M^{2,2}$.

If $D$ is an admissible submodule
of $A$ of level $m$, then $D^{2,1}$ is an\break
$\CO _\wp$-submodule of  
$A[m]^{2,1}$ of order $  \RN(m) $. 
The map $M\longmapsto M^{2,1}$ is bijective 
between the set of admissible submodules
of $A$ of level $m$, and the set of submodules
of $A[m]^{2,1}$ of order $\RN(m)$.
Indeed,
 for a given $\CO _F$-submodule $D_1$ of $A[m ]^{2,1}$
of order $  m $, 
we can obtain a module $D_2=D_1+wD_1$ as an $\CO _{B,   m }$
module of $A[  m ]^2$. Let $D_3$ be the annihilator of 
$D_2$ in $A[  m ]$; then $D=D_2+D_3$ is 
an admissible submodule of level~$m$.
\enddemo

\demo{{\rm 1.4.9.} The Eichler\/{\rm -}\/Shimura congruence
relation}
Let $\wp$ be a prime in $\CO _F$ over which 
$K$ has the maximal component and $B$ is split.  Let $\Frob (\wp)$\lb{Frob(wp)}
 be the Frobenius correspondence 
on $\CM_{K,k}$ where $k$ is the residue field of $\CO _{F, \wp}$.
Then we have the following Eichler-Shimura congruence relation: \enddemo 

\nnp{Proposition 1.4.10}
Let $\Frob (\wp)^*$  denote the dual correspondence 
of $\Frob (\wp)${\rm .} Then
$$\RT(\wp)=\Frob(\wp)+\Frob (\wp) ^*.$$
\endproclaim
 
\demo{Proof} Let $F'=F(\sqrt\lambda)$ be as before,
such that $\wp$ is split in $F'$.
We will only give a proof for the special case where $M_K$ can be embedded
 into $M_{K'}$ for some $K'$ which is sufficient to apply to the curve
$X$ defined in the introduction. (The proof of the general case can be
 found in Carayol's paper \cite[\S10.3]{C} where he uses
 a slightly different definition of $M_{K'}$ so that every
$M_K$ can be embedded into his $M_{K'}$.)

It is obvious that we need only prove the 
same identity for $M_{K'}$. Furthermore, it
is true if the identity is true for one $K'$ then
it is true for a  smaller one, as the identity
in the proposition is stable under pushforward
of cycles. So we may assume that $K'$ is compact. Now we need only verify the identity
for  points in $\CM _{K', \wp}(\bar k)$. 
We may only restrict ourselves to the dense subset of smooth and
ordinary points. These points   are reductions of points in $\CM _{K', \wp}(W)$
where $W$ is the completion of the maximal
unramified extension of $F_\wp$.

Let $[A, \rho, \bar\theta, \bar\kappa]$ be 
one object in $\CF _{K'}(W)$. Then 
$${\rm T}(\wp)[A, \rho, \bar\theta, \bar\kappa]
=\sum _{D}[A_D, \rho_D, \bar\theta_D, \bar\kappa_D]$$
where $D$ runs through the set of admissible 
submodules of level $m$. We want to study 
the reduction of this identity module
$\wp$. 

As   explained in 1.4.8,
$D$ is completely determined by a submodule $D^{2,1}$ of $A[\wp]^{2,1}$ of order $\wp$.
Since our object is ordinary, $A[\wp ^\infty]^{1,2}$ is isomorphic to 
$$\Sigma:=\Sigma _1\oplus F_\wp/\CO _\wp$$
where $\Sigma _1$ is a formal $\CO _\wp$-module
on $W$ of height $1$. The generic fiber
of $\Sigma$ isomorphic to $F_\wp/\CO _\wp
\oplus F_\wp/\CO _\wp$. For any $t\in \CO _\wp/\wp$, let $\Sigma ^t$ denote the 
submodule of $\Sigma $ whose generic fiber is a
group of points $(tx, x)$. The submodules of $\Sigma$ order $\wp$ are exactly those of $\Sigma ^t$ and $\Sigma _1$. 

As the universal deformation space of 
$[A, \rho, \bar\theta, \bar \kappa]$ is isomorphic to that of $A[\wp^\infty]^{2,1}$,
 it is easy to see
that the isogeny $A\to A_D$ is purely inseparable
if $D^{2,1}$ corresponds to $\Sigma_1$, and
is \'etale  if $D^{2,1}$ does not correspond to
$\Sigma _1$. Thus in the first case, 
$$
[A_D, \rho _D, \bar\theta _D, \bar\kappa _D]
=\Frob(\wp)[A, \rho , \bar\theta , \bar\kappa]\pmod{\wp}$$
and in the second case
$$[A, \rho , \bar\theta , \bar\kappa]
=\Frob(\wp)[A_D, \rho _D, \bar\theta _D, \bar\kappa _D]\pmod{\wp}.$$
Now the congruence relation in the proposition follows.\enddemo

\demo{{\rm 1.5.} Order $R$ and its level structure}\enddemo 

 1.5.1. {\it Construction of $R$ and
$X$}.
Let $N$\label{N} be a nonzero  ideal of $\CO_F$ and let $E$\label{E} be a totally imaginary quadratic extension of $F$
whose relative discriminant is prime to $N$. Assume that $\varepsilon (N)=(-1)^{g-1}$, where 
$$\varepsilon\label{epsilon} : F^\times \backslash \wh F^\times\to \{\pm 1\}$$
is the character associated to the extension $E/F$. Then up to isomorphisms, there is a unique quaternion algebra $B$ such that $B$ 
is ramified exactly at the
place $\tau$ and  finite places $\wp$ where $\varepsilon _\wp(N)=-1$, as this ramification set has 
even cardinality by our assumption. Also by construction,
every ramification place of $B$ is not split in $E$. So we may
fix an embedding $\rho :E\to B$\label{rho} over $F$. This allows us
to consider $E$ as a subalgebra of $B$. 

In the following we want to construct an order $R$\label{R}
of $B$ of type $(N, E)$\label{typeNE}; this means that $R$ 
contains $\CO_E$ and has discriminant $N$. 
For each prime $\wp$ dividing $N$, let $\wp_K$
\label{wpK} be a prime of $\CO _E$ dividing $\wp$. Let $N_E$\label{NE} be an ideal of $\CO _E$ which is a product of powers of $\wp _E$\label{wpE} and which has relative norm 
$N/N_B$. The existence of such $N_E$ follows easily from our assumptions. 
Indeed, if we write
$$\wt N\label{wtN}=\prod _{\varepsilon (\wp)=1}\wp _E^{\ord _\wp (N)}
\cdot \prod _{\varepsilon (\wp)=-1}\wp _E ^{[\ord _\wp (N)/2]}$$
then
$$N_E=\prod _{\wp}\wp _E^{\ord _\wp(\wt N)}.$$
Let $\CO _B$ be a 
maximal order of $B$ containing $\CO _E$. Then we obtain an order of $B$   by the following formula:
$$R=\CO _E+N_E\CO _B.$$

Conversely, any order of type $(N, E)$
of $B$ has the above form with some choice
of the maximal order $\CO_B$.

As in the introduction, our primary  curve of study
 is the compactification $X$\label{X} of the  Shimura curve associated to the noncompact group $\wh F^\times \cdot \wh
R^\times$.
\enddemo

\demo{{\rm 1.5.2.}  Cyclic submodule
structures}
 Let $K$ be an open subgroup of $\wh R$ 
which has the same components as $\wh R$
over places dividing $N$. Let $J$ be some
compact open subgroup of $\wh F^{'\times}$
which has maximal components at places dividing
$N$. Let $K_0$\label{K0} denote the subgroup of 
$\wh \CO _B^\times$ which is obtained by
replacing components of $K$ over places
diving $N$ with maximal ones. Let
$K'$ denote $K\cdot J$ and $K'_0$\label{K0'} denote 
$K_0\cdot J$. Then we have a morphism of functors
$$\CF _{K'}\to \CF _{K'_0}.$$
In the following we want to show that the 
fiber of this morphism is given by 
so-called cyclic submodule structures.
 
For every prime $p$ which is divided by at 
least one prime factor $\wp$ in $ N  $, 
we assume that $\sfrac \lambda p=1$, and 
fix a square root $\mu _p$ of $\sqrt \lambda$ in $\BQ_p$. In this way any $\CO _{F'}/N$
module $M$ has decomposition $M=M^1\oplus 
M^2$ in the same fashion as before.
\enddemo

\demo{{\rm 1.5.3.} Definition}
\lb{cyc.sub.mod.str}
Let $A$ be an object of $\CF _{K_0}(S)$.
By a cyclic submodule structure on $A$
of level $N_E$,  we mean 
an $\CO _{E}/N_E$-submodule $C$ of $A[N_E]^2$
such that locally there is an element
$x\in A[N_{E}]$ with the following properties:
\begin{itemize}
\ritem{1.} The element $x$ is a Drinfeld base for  $C$. This
means  that as cycles one has:
$$[C]=\sum _{a\in \CO _{E}/N_{E}}[ax].$$
\ritem{2.}  If $\wp$ is a prime of $F$ over which
$B$ is not split, then $x$ is also a
Drinfeld base for $\CO_B/\wt N$-module 
$A[\wt N]^2$.
\end{itemize}

Notice that the second condition here is equivalent to the fact
that $x$ is not divisible 
by uniformizers of $B$ in $A[\wt N]$.
\enddemo

\spn{1.5.4} \proclaim{Proposition}\label{P:cyclicsubmodules}
The functor $\CF _{K'}$ is equivalent
to the functor  which sends an $F'$\/{\rm -}\/scheme $S$ to the set  of objects $[A, C]${\rm ,}  where 
 $A$ is an object in  $\CF _{K'_0}(S)${\rm ,}   and where $C$ is a cyclic submodule structure of level
 $N_{E'}$ on $A${\rm .}
\endproclaim

Before  the proof of this proposition, we need the following crucial lemma.
Let $ E'= E\otimes F'$.\label{E'} Then every prime factor
 $\wp_E $ of $\CO _ E$ can be lifted to a 
prime $\wp _{ E'}$\label{wpE'}  which is the 
preimage of $\wp_E$ via   the map 
$$\CO _{E'}\to \CO _{E, \wp_E}, \qquad \sqrt \lambda \to -\mu _p.$$
Let $N _{E'}$\label{NE'}   be the lifting  of $N_E$  to the ideal in $\CO _{E'}$.
  Now we have the formulas
$$N_{F'}=\prod _\wp
\wp_{F'}^{\ord _\wp (\wt N)}.$$

\spn{1.5.5} \proclaim{Lemma}\label{L:identities} The following identities hold in $\wh B$\/{\rm :}
\begin{eqnarray*}
\wh\CO _  B    ^\times \cdot\wh\CO _{F'}^\times
&=&\left\{g\in\wh   B    ^\times \cdot\wh F^{'\times}:   \wh \CO _{  B   '}g=\wh\CO _{  B   '}\right\},\\
\wh  R  ^\times
&=&\left\{g\in \wh\CO _B^\times: \wh  N_ E ^{-1}g =\wh N_ E ^{-1}\pmod {\wh\CO _B}\right\},\\
\wh   R    ^\times
\cdot \wh\CO _{F'}^\times& =&\left\{g\in \wh\CO _  B    ^\times \cdot \wh\CO _{F'}^\times: 
  \wh  N_{ E'}^{-1}g= \wh N_{ E'} ^{-1} \pmod{\wh\CO _{B'}}\right\}.
\end{eqnarray*}
\endproclaim

{\it Proof}.
\quad 1.5.6. {\it First identity}.
We need only prove the inclusion ``$\supset$" 
for each place $\wp$. Let $a\in B_\wp^\times$
 and $b\in F^{'\times }_\wp$ such that 
$c=ab\in \CO _{B', \wp}^\times$.

If $F'_\wp$ is a field unramified over $F_\wp$, then $b=db'$ with $d\in F_\wp^\times$ and $b'\in \CO _{F', \wp}^\times$.
So we may write $c=a'b'$ with
$$a'=ad=cb^{'-1}\in B_\wp\cap \CO _{B', \wp}^\times=\CO _{B, \wp}^\times.$$

If $F_\wp'$ is split, then we have a decomposition
$$F_\wp'=F_\wp\oplus F_\wp, 
\quad B_\wp'=B_\wp\oplus B_\wp.$$
Write $b=(b_1, b_2)$ with respect to these decompositions. Then $ab_1$ and $ab_2$ are both in $\CO _{B, \wp}^\times$. It follows that
$b_1b_2^{-1}\in \CO _{F, \wp}^\times$. So we may write $c=a'b'$ with 
$$b'=(b_1b_2^{-1}, 1)\in \CO _{F', \wp}^\times
\quad {\rm and} \quad a'=ab_1\in \CO _{B, \wp}^\times.$$

Finally let us assume that $F'_\wp$ is a ramified quadratic extension of $F_\wp$.
Let $\pi'$ be a uniformizer for $F_\wp'$. By 
replacing $b$ with a multiple of elements in 
$F_\wp^\times\cdot \CO _{F'}^\times$, we may assume that $b$ is either $1$ or $\pi'$.  In 
the first case, $a$ must be a unit and we are done. Now we assume that $b=\pi'$. Since 
$\wp$ is ramified in $F'$, $\wp$ must be split
in $B$. By writing $a$ as a 2 by 2 matrix over 
$F_\wp$, we see that the integrality of $a\pi'$ implies that of $a$. But this implies that $c=a\pi'$ cannot be a unit. 
\enddemo

 {\rm 1.5.7.} {\it Second identity}.
This identity follows from the definition because
\begin{eqnarray*}
\wh N_E^{-1} g=\wh N_E^{-1}\pmod {\wh \CO _B}
&\Longleftrightarrow&
\wh N_E^{-1}g=\wh N_E^{-1}+\wh \CO _B\\
&\Longleftrightarrow&
g\in \wh \CO _E+N_E\wh\CO _B=\wh R.
\end{eqnarray*} 
 
\vglue12pt

\demo{{\rm 1.5.8.} Third identity}
For this, we need only show the following:
\begin{equation}\label{E:third}
\wh \CO _B^\times \cap \left(\wh \CO _{E'}
+\wh N_{E'}\wh \CO _{B'}\right)\subset \wh R^\times, \spq{1.5.1}
\end{equation}
since by similar reasoning  to that above,
$$\wh N_{E'}^{-1}g=\wh N_{E'}^{-1}
\pmod{\wh \CO _{B'}}
\Longleftrightarrow
g\in \wh\CO _{E'}+\wh N_{E'}\wh\CO _{B'}.$$
We need only check \eqref{E:third} for each place $\wp$ of $F$. This is clear if $\wp$
does not divide $N$. But if $\wp$ divides $N$, then it is split in $F'$ and we have decompositions:
$$B_\wp'=B_\wp\oplus B_\wp, \quad
\sqrt{\lambda}\to (\mu _p, -\mu _p),$$
$$N_{E', \wp}=\CO _{E, \wp}\oplus 
N_{E, \wp}.$$
It follows that
$$\CO _{E', \wp}+N_{E'}\CO _{B', \wp}
=\CO _{B, \wp}\oplus R_\wp.$$
Thus \eqref{E:third} is proved.\hfill\qed
\enddemo

\demo{{\rm 1.5.9.} Proof of Proposition
{\rm \ref{P:cyclicsubmodules}}}
Let $S$ be an $F'$-scheme, and $[A, \bar\kappa ]$ an element of 
$\CF_{K'} (S)$ with $A\in \CF _{K_0'}(S)$ and $\bar\kappa $ a class
modulo $K'$ isomorphisms $\kappa : \wh V_\BZ   \to \wh T (A)$.
This  $\kappa $ will induce an isomorphism $\kappa : \wh V\to \wh V(A)$ 
and a map
$$\wt\kappa  :\wh V\to \wh V(A)/\wh T (A)=A_{{\rm tor}}.$$
Thus, we  have an $\CO _{ E'}$-submodule
 $C_{\kappa }:=\wt\kappa  (N _{ E'}^{-1}/\CO _{ E'})$
of $A[N_{ E'}]$.
By the above lemma,  $C_\kappa $ does not depend on the choice of $\kappa $
in the class $\bar\kappa $.\break 
 Since $N_{E'}^{-1}/\CO _{E'}$ is a free module 
of rank $1$ over $\CO _{F'}/N_{F'}$ and generates\break $\CO _{B'}/N_{F'}$-module $N_{F'}^{-1}\CO _{B'}/\CO _{B'}$,  it
follows that $C_\kappa $ is generated by a Drinfeld base $x$ of the
  order $N_{F'}$.

Conversely, for any $\CO _{ E'}$-submodule $C$ of $A[  N_{ E'}]$ which 
is generated by 
a Drinfeld base of order $N_{F'}$, and any level structure $\kappa  _0$
for the compact subgroup $  K_0'$, we have a unique level structure
$\kappa $ so that $\kappa _0=\kappa \pmod{K_0'}$
and $C=C_\kappa $.

In a similar manner, we have the following:\enddemo

\spn{1.5.10} \proclaim{Proposition}
Let $\wp$ be a prime of $F$ of characteristic   $p${\rm .}
 Assume that $ J $ is maximal at places over $p${\rm .}
Then the functor $\CF _{K',\wp}$ is equivalent to
the functor  which sends the $W$\/{\rm -}\/scheme $S$ to the 
set of isomorphism classes of objects $[A, C]$
where $A$ is an object in   $\CF _{K  '_0, \wp}(S)$
and $C$ is a cyclic submodule structure on $A$ of order $N_{F'}${\rm .}
\endproclaim

\section{Heegner points}

In this section we study Heegner points.
We start in Section 2.1 with the general definition of CM-points and Heegner points
as complex points,
and their modular interpretations. Then we move to the 
study of their reductions which are so-called distinguished points,
 first  the structure of formal group
in  Section~2.2 and then the structure of endomorphism rings in  Section~2.3
using Honda-Tate theory.
Finally in  Section~2.4,  we study the lifting of distinguished points by 
Serre-Tate theory and Gross's theory.
In this section we assume that every prime factor of $2$
is split in $E$.

\demo{{\rm 2.1. CM-}\/points} 
\enddemo

2.1.1. {\it Definitions and general properties}.
Our primary object of study in this paper is the class
of Heegner
points on the curve $X$ defined in 1.5.1 by the 
noncompact group $\wh F^\times \wh R^\times$. From 
the modular  point of view, it is more natural
to study Heegner points on the  
Shimura curve $Y$\lb{Y} defined by the 
compact group $\wh R^\times$:
$$Y=B^\times \backslash
\CH ^\pm\times \wh B^\times /\wh R^\times.$$
The curve $X$ is then a quotient of $Y$ by the 
action of $\wh F^\times$.
As in the introduction, we fix a splitting
$$B\otimes _\tau \BR =
\RM_2(\BR)$$
such that $\rho (E)\otimes \BR$ is sent to
the subalgebra of $\RM_2(\BR)$ of elements 
$ \left( \begin{array}{cc}a&b\\ -b&a \end{array}\right)$. 
We then extend $\tau: F\to \BR$ to
$\tau: E\to \BC$\lb{tauE} such that
$$\tau (x)=a+bi
\Longleftrightarrow
\rho(x)=
 \left( \begin{array}{cc}a&b\\ -b&a \end{array}\right).$$
 
We say a point  $z$ in $Y$ is a 
CM-\/{\em point {\rm (}by} $ E$)\lb{CMpoint}, if $z$ is represented
by an element of $\CH ^\pm\times \wh  B   ^\times$ of 
the form $(\sqrt{-1}, g)$.

For a CM-point $z$, let $\phi _z$ denote the morphism
$$g^{-1}\rho g:  E\to \wh  B.$$
Then up to conjugation by $\wh    R   ^\times$, $\phi _z$ does not depend on the
choice of $g$. 
The order $\End (z):=\phi _z^{-1}(\wh   R    )$ 
in $ E$, which  does not depend on the 
choice of $g $, is called 
{\em the endomorphism ring  of $z$}.\lb{end.rin}
The ideal   $ c $  of $\CO _F$, such that
$$\End (z)=\CO _ c :=\CO _F+ c \CO _ E,$$
 is called {\em the conductor of $z$}.
\lb{conductor}

For a  place $\wp$ prime to $ c $, the homomorphism  $\phi _z$  defines 
an {\em orientation}\lb{orientation} in 
$$U_\wp=\Hom (\CO _{ E, \wp},   R    _\wp)/   R    _\wp^\times .$$
This  set has only one element if $\wp$ does not divide $   N  $; otherwise
it has two elements: the image of  $\rho $ which we called {\em the positive 
orientation},\lb{pos.ori}  and the image of $\bar\rho $ which we called {\em the negative orientation}.
We say two CM-points 
have {\em the same orientation}\lb{sam.ori}, if they define the same elements in $U_\wp$
for $\wp|   N  $. 
If we write 
\vglue-9pt
$$\CO _{ E, \wp}=\CO _{F, \wp}+\CO _{F, \wp}e$$
\vglue3pt\noindent 
with $e^2\in F$, then two embeddings 
\vglue-9pt
$$\phi_1, \phi _2:\quad \CO _{ E, \wp}\to    R    _\wp$$
\vglue3pt\noindent 
define the same element in $U_\wp$ if and only if
\vglue-9pt 
\begin{equation}\lb{E:equ.con.ori}\left(\phi _1(e)-\phi _2(e)\right)^2\equiv 0\pmod    N  . \spn{2.1.1}
\end{equation}
\vglue3pt\noindent 
Indeed, write $R_\wp =\CO_{E, \wp}+t\CO _{E, \wp}$ with $t\in R_\wp$
such that $\det(t)$ generates $N_\wp$. Then if $e_1$ and $e_2$ have the 
 same orientation, it follows that
 $e_1-e_2\in tR_\wp$. This implies that
\vglue-9pt  $$(e_1-e_2)^2=-\det (e_1-e_2)=0\pmod {N_\wp}.$$ 
If $e_1$ and $e_2$ do not have the same orientation, then 
$e_1-e_2=2e_1\mod {tR_\wp}$. Thus 
\vglue-9pt  $$(e_1-e_2)^2=4e_1^2 \ne 0\pmod {N_\wp}.$$
\vglue3pt 
 
The curve $Y$ admits an action by the group\vglue-6pt  
$$\CW\lb{CW} =\{b\in \wh  B    ^\times: b^{-1}\wh   R    ^\times b=\wh   R    ^\times \}
/\wh   R    ^\times . $$ 
\vglue3pt\noindent 
 This group  has $2^s$
elements, where $s$\lb{s} is the number of prime factors of $   N  $.
The action of  $\CW$  on CM-points does not change
the conductors, and the induced action on  $\prod_{\wp|   N  }U_\wp$ is free 
and transitive. 

Let $Y_ c $\lb{Y_c} denote the subscheme 
 of the positively oriented CM-points of the  conductor
$ c $. Then $Y_ c $    is defined over $ E$
and every point in $Y_ c (\bar  E)=Y_ c (\BC  )$ is defined over the ring class field
 $H_ c $\lb{H_c} of $\CO _ c $. 
Indeed, if $(\sqrt{-1}, g)$ is a CM-point of $Y$ with positive orientation and conductor $c$, then $Y_ c $ is identified
with  the set of points represented by $(\sqrt{-1}, \wh  E^\times g)$. 
The correspondence  which sends $x$ to the class of $(\sqrt{-1}, xg)$, 
therefore, defines a bijection \vglue-9pt  
$$ E^\times\backslash \wh  E^\times / \wh\CO _ c ^\times\simeq Y_c.$$ 
\vglue3pt\noindent The Galois action of $\Gal (H_ c / E)$ on $Y_c$ is given by the inverse of the map,\vglue-9pt  
$$\Gal (H_ c / E)\simeq  E^\times\backslash \wh  E^\times / \wh\CO _ c ^\times,$$
\vglue3pt\noindent via  class field theory.

A CM-point $z$ by $ E$ is called a {\em Heegner point}\lb{Hee.poi} if its conductor
is the trivial ideal $\CO_F$. 
Obviously, the point $(\sqrt{-1}, 1)$ is a Heegner point.
In this paper we only consider CM-points with  conductors prime
to $ N{D_E}$ and with  positive orientation,
where $D_E$\lb{D_E} is the relative discriminant ideal
in $\CO _F$ for the extension $E/F$.

Notice that the property of a point to be a CM-point of conductor
 $c$ is invariant under
the action by $\wh F^\times$. Thus, all the
above discussion is valid  for $X$ or
any Shimura curves between $X$ and $Y$.\enddemo

\demo{{\rm 2.1.2.} Modular interpretation} 
We fix $F'$ as in Section~1.5.
In the following we want to give a modular interpretation
of Heegner points over  $ E'=F'\cdot  E$.
We let $Y'$\lb{Y'} denote the Shimura curve $M_{K'}$
with 
$$K'=\wh R^\times \cdot \wh\CO _{F'}^\times.$$
Then  $Y$ has a finite morphism to $Y'$.

Let $\CF$\lb{CFR} denote the functor $\CF _{K'}$
and let $\CF _0$\lb{CFR0} denote $\CF _{K'_0}$ where 
$$K'_0=\wh \CO _{B}^\times \cdot \wh \CO _{F'}
^\times.$$ Then 
every point $x$ in $Y(\BC)$ represents an object
$[A, C]$, where $A$ stands for an object 
$[A, \rho_A, \bar\theta_A, \bar\kappa _A]$
of $\CF _0(\BC)$ and $C$ is a cyclic 
$\CO _{E}$-submodule
structure of $A$ of level $N_{E}$.
 We need some notation to state our result:
\begin{itemize}
\item For $[A, C]$ in $\CF (S)$, 
\begin{itemize}
\item let $\End _{\CF _0}(A)$\lb{End_0} denote the 
$\CO _{F'}$-subalgebra of $\End _{\CO _{B'}}
(A)$ generated by elements $\phi:A\to A$
such that $\phi\phi ^*\in F^\times$,
where $\phi\to \phi^*$ is a Rosati involution
induced by a polarization in $\bar\theta _A$.
\item let $\End _\CF (A, C)$\lb{End} denote the 
subalgebra of $\End _{\CF _0}(A)$ of elements 
$\phi$ such that $\phi (C)\subset C$.
\end{itemize}
\item Let $t': E'\to E'$\lb{t'} be a map defined by
$$t'(a+b\sqrt\lambda)=\tr _{E/\BQ}(a)
+a-\bar a+\left(\tr _{E/\BQ}(b)-b-\bar b\right)
\sqrt\lambda$$
for any $a, b\in E$.
\end{itemize}

\spn{2.1.3} \proclaim{Proposition} Let $x$ be a point on $Y(\BC)${\rm ,}
and let $[A, , C]$ be an object represented
by $x${\rm .}
Then the point $x$ is a {\rm CM-}\/point by $ E$  if and only if 
$\End _{\CF_0}(A)\otimes\BQ\simeq  E'${\rm .}
Moreover{\rm ,} if $x$ is a {\rm CM-}\/point by $ E${\rm ,} then\/{\rm :}
\begin{itemize} 
\ritem{1.}  There is a unique isomorphism
$$\alpha : E' \simeq \End _{\CF_0}(A)\otimes\BQ$$\lb{alpha}
over $F'$ such that  for any $a\in  E'${\rm ,}
$$\tr (\alpha (a): \Lie A)=2s(a).$$
\ritem{2.} With $\alpha $ as above{\rm ,}
$$\End (x)=\{ a\in  E: \quad\alpha (a)\in \End _{\CF}(A, C)\}.$$
\end{itemize}
\endproclaim

\demo{{P}roof} Let $x$ be represented by $(z, \gamma )$.  With notation as in the proof of Proposition
\ref{CF}, the endomorphism ring $\End _{\CF _0}
(A)$ can be identified with the subring of 
$B$ generated by elements $b\in B^\times 
\cdot F^\times $ such that
\begin{itemize}
\ritem{1.} $V_\gamma b\subset V_\gamma$, or 
equivalently, $b\in \gamma \wh \CO _{B'}
\gamma^{-1};$
\ritem{2.} $bj_z=j_zb$, or equivalently,
$\tau (b)\in a\rho (E)a ^{-1}
\otimes _\tau \BR$, where $a\in \GL _2(\BR)$
such that $a(\sqrt{-1})=z$.
\end{itemize}
It follows that $\End _{\CF _0}(A_x)
\otimes \BQ$ is an $F'$-subalgebra of
$B'$ generated by elements $b\in B^\times $
satisfying the second condition. \enddemo

\demo{{\rm 2.1.4.} Equivalence}
If $\End _{\CF_0}(A)\otimes \BQ\simeq E'$, 
then there is an embedding $\beta :E\to B$ over
$F$ such that in $B'$,
$$\End _{\CF _0}(A)\otimes \BQ
=\beta (E)\otimes F'.$$
As all embeddings of $E$ into $B$ are conjugate,
it follows that $\beta =b\rho b^{-1}$ where $b\in B^\times$
 is uniquely determined by $\beta$ modulo
$\rho (E)^\times$. Now condition 2 implies
that in $B\otimes _\tau \BR$,
$$b\rho  (E)b^{-1}\otimes_\tau \BR
=a\rho (E)a^{-1}\otimes _\tau \BR.$$
It follows that $b=a k$ with some $k\in \rho (E)\otimes _\tau \BR$. As 
$a (\sqrt {-1})=z$,  $k(\sqrt {-1})
=\sqrt{-1}$, one must have $z=\beta (\sqrt {-1}).$
Thus $x$ can be represented by 
$$\beta ^{-1}(z, \gamma )=
(\sqrt{-1}, \beta ^{-1}\gamma ).$$
So $x$ is a CM-point.

Conversely, if $x$ is a CM-point and is represented by $(\sqrt{-1}, g)$,
 then 
in the above description of $\End _{\CF _0}
(A)$,  we may take $a=1$ in condition~2. Now, 
$$\End _{\CF _0}(A)\otimes \BQ
=\rho (E)\otimes F'.$$
This is isomorphic to $E'$ by the following 
map:
$$\alpha : E'=E\otimes F'\to 
\End _{\CF _0}(A)\otimes \BQ,$$
$$\alpha (x\otimes y)=\rho (x)\otimes y.$$ \enddemo

\demo{{\rm 2.1.5.} First property}
It remains to show that $\alpha$ satisfies both
properties in the proposition. If $a\in E$ then
$a$ acts on $V_\BR$ via right multiplication
by $\rho (a)$. Write 
$\rho (a)=(a_1, \cdots, a_g)$
with respect to the decomposition
$$B\otimes \BR=\RM_2(\BR)\oplus 
(\BH)^{g-1}.$$
Then by the  definition of complex structure 
in Section~1.1 on
$$V_\BR=\RM_2(\BR)\otimes \BC
\oplus (\BH\otimes \BC)^{g-1},$$
one has
$$\tr (a+b\sqrt{\lambda})
=4a+2\sum _{i\ge 2}
\left(\tr _{\BH/BR}(a_i)
+\tr _{\BH/BR}(b_i)\sqrt\lambda\right)
=2t'(a).$$
The only other isomorphism between 
$\End _{\CF _0}(A)\otimes \BQ$ and $E'$ is
$\bar\alpha$ defined by 
$$\bar\alpha (x\otimes y)=\alpha (\bar x
\otimes y)$$
which does not satisfies property 1 as
$t'(\bar a)\ne t'(a)$ for $a\in E-F$.\enddemo

\demo{{\rm 2.1.6.} Second property}
Finally we want to prove the second property
in the proposition.
By the proof of Proposition \ref{CF}
and 1.5.4, $C$ is isomorphic to 
$N_{E'}^{-1}\gamma ^{-1}$ modulo
$V_\gamma=\wh \CO _{B'}\gamma ^{-1}$. It follows that
\begin{eqnarray*}
&&\hskip-.25in\left\{a\in E, \quad \alpha (a)\in \End
_{\CF}(A, C)\right\}\\
&&\qquad\quad =\left\{a\in E, \quad  
\begin{array}{c}
\rho (a)\in \gamma \wh \CO _{B'}\gamma ^{-!},\\
\noalign{\vskip4pt}
N_{E'}^{-1}\gamma ^{-1}\rho (a)
\subset N_{E'}^{-1}\gamma ^{-1}
\pmod{\wh \CO _{B'}\gamma ^{-1}}
\end{array}\right\}\\ \noalign{\vskip4pt}
&&\qquad\quad =\left\{a\in E, \quad\gamma ^{-1}\rho (a)
\gamma \in \wh \CO _{E'}+N_{E'}\wh\CO _{B'}
\right\}.\end{eqnarray*}
Similarly, as in 1.5.9, it is easy to see that
$$\wh B\cap \left( \wh \CO _{E'}
+N _{E'}\wh \CO _{B'}\right)=\wh R.$$
Thus we have
\begin{eqnarray*}
\left\{a\in E, \quad \rho (a)\in \End _{\CF}
(A, C)\right\}
&=&\left\{a\in E, \quad \rho (a)\in \gamma 
\wh R\gamma ^{-1}\right\}\\
&=&\End (x).\\
\noalign{\vskip-36pt}
\end{eqnarray*}
\hfill\qed\enddemo

 \phantom{I want to}

\nnp{Proposition 2.1.7} 
Let $x$ and $y$ be  two {\rm CM-}\/points with  conductors prime to $   N  ${\rm , }
and representing the objects $[A, C]$ and 
$[A', C']${\rm .}
Then $x$ and $y$ have the same orientation  if and only if 
there is an $\left(\iota(\CO _{B'})\otimes 
\alpha(\CO_ E)\right)_   N  $-linear symplectic
similitude from $T(A)_   N  $ to $T(A')_   N  $ which takes
$C$ to $C'${\rm .} Here for an $\CO _F$\/{\rm -}\/module $M$, $M_   N  =M\otimes \mathbold{\oplus}_{\ell |N}
\BZ_\ell${\rm .}
\endproclaim
 
\demo{Proof}
We may assume that $x$ is represented by 
$(\sqrt{-1},1)$ and prove only the local statement for each $\wp$ dividing 
$N$. Let $y$ be represented by $(\sqrt {-1},
\gamma)$. Then we have isomorphisms of
$\iota(\CO _{B'})\otimes \alpha(\CO _{E, \wp})$-modules
$$T_\wp (A)\simeq \CO _{B', \wp}
\quad\hbox{and}\quad
T_\wp (A')\simeq \CO _{B', \wp}\gamma _\wp ^{-1}$$
where $\iota(\CO _{B'})$ acts by left multiplication
and $\alpha(\CO _E)$ acts by right multiplication
of $\rho (\CO _E)$,
and isomorphisms of 
$\alpha (\CO_{E'})$-submodules
$$C_\wp\simeq N_{E', \wp}^{-1}
\pmod{\CO _{B', \wp}}
\enspace \hbox{ and }\enspace
C_\wp'\simeq N_{E'}^{-1}\gamma ^{-1}
\pmod {\CO _{B', \wp}\gamma^{-1}}.$$
 
As any $B_\wp'$-linear endomorphism of $B_\wp'$
is given by right multiplication by an element
of $B_\wp'$, the ``if" part of the proposition is,
 therefore, equivalent to the existence 
of $a\in B_\wp'$ such that the following conditions are verified:
\begin{itemize}  
\item[1.] $\gamma _\wp^{-1}a\in \CO _{B', \wp}^\times$;
\item[2.] $a$ commutes with $\rho (E)$;
\item[3.] $a\in B_\wp^\times \cdot F_\wp^{'\times}$;
\item[4.] $N_{E'}\gamma _\wp ^{-1}a\subset N_{E'}
\pmod {\CO _{B', \wp}}$.
\end{itemize} 
By the first identity of Lemma 1.5.5,
conditions 1 and 3 here are equivalent to the fact
that $a$ has the form $\gamma _\wp^{-1}a
=bc$ where $b\in \CO _{B, \wp}^\times$
and $c\in \CO _{F', \wp}^\times$.
Replacing $a$ by $ac^{-1}$, 
 we may assume that $c=1$ and then 
$a\in B_\wp$.  

Now condition 2
is equivalent to $a\in \rho (E)$, and 
condition 4 is equivalent to 
$\gamma _\wp ^{-1}a\in R_\wp$ by
a similar argument to  1.5.8.  It follows
that the ``if" part of the proposition is equivalent 
to $\gamma _\wp\in \rho (E)^\times
\cdot R_\wp ^\times$, or equivalently to
the fact that the map 
$$\gamma _\wp ^{-1}\rho \gamma _\wp:
E_\wp\to B_\wp$$
has positive orientation.
\enddemo

\vglue9pt

\demo{{\rm 2.2.} Formal groups} Let $q$ be a finite place of $ E$
and let $ E ^\ur_q$\lb{E_q^ur} be the 
completion of the maximal unramified extension of $ E_ q$ with 
 ring  of  integers $\CO _ q^\ur$,\lb{CO_q^ur} and  residue field $k$.\lb{k} 
Let ${y} $\lb{y} be a CM-point of
 $Y$ with
 conductor $ c $ prime to $ N   {D_E}   $ and $ q $. Then ${y}$ is defined
over  $ E^{\ur}_ q $.
 Let $\bar{y}$\lb{bary} denote the Zariski closure
of ${y}$ in $\CY\otimes \CO ^{\ur}_q$, where $\CY$\lb{CY} is the integral model of 
$Y$ over $\CO _F$ constructed in Section~1.2. We want to study the reduction
${y} _k$\lb{y_k} of $\bar{y}$ in $\CY _k:=\CY\otimes k$.\lb{Y_k}

Let $p$ denote the characteristic of $k$ and let $\wp$ denote the 
prime of $\CO _F$\break under~$q $.
As usual, we will choose an auxiliary negative integer $\lambda$ 
as in Section~1.1 and work on $F'=F(\sqrt{\lambda })$. 
We   assume that $\left(\frac \lambda p\right)=1$ and choose
a square root $\mu _p$ of $\lambda$ in $\BQ _p$.
Then there is the  usual decomposition $M=M^1\oplus M^2$ for $F'_p$ modules $M$.
  Let $i$ denote the embedding 
$$i: E'= E(\sqrt\lambda)\to  E_ q^{\ur}$$\lb{i}
 which takes $\sqrt \lambda$ to $\mu _p$.

Let $[\bar A, \bar C]$ be the abelian variety represented 
by $\bar {y}$. Then the action of $\End ({y})\otimes \CO _{F'}$
 on $A=\bar A\otimes  E_q^{\ur}$ extends to an action on $\bar A$.
Let $\CG$\lb{CG} denote the divisible  $\CO _{E'}$-module $\bar A[\wp ^\infty]^2$. 
\enddemo

\vglue9pt 

\nnp{Proposition 2.2.1} The action of $\CO _{ E} $ on the $\CO _q^\ur$\/{\rm -}\/module 
$\Lie (\CG)$ induced by the action $\alpha$ on $\bar A$
is given by the canonical embedding $\CO _{ E}\to \CO _q^{\ur}${\rm .} 
\endproclaim

\demo{Proof}
We want to prove the proposition by computing the trace of the action of $\alpha (\CO _{E'})$.
Recall that the action $\alpha: E'\to \End (A)\otimes \BQ$ induces an action of
$E_p'=E'\otimes \BQ_p$ on $\Lie (A)$, therefore an action of $E_\wp'=E'\otimes F_\wp$ on
$\Lie (A[\wp ^\infty])$. This last module has a projection
$$\Lie (A[\wp ^\infty])\to \Lie (\CG),
\quad \sqrt{\lambda}\to -\mu _p.$$
We denote all of these actions by $\alpha$. By Proposition 2.1.3, the action $\alpha$ of 
$E'$ on $\Lie (A)$ has the trace map $i\circ 2t': E'\to E_q^\ur$. It is easy to see
that the action $\alpha $ of $E_\wp'$ on $\Lie (A[\wp]^\infty)$ will have the trace
 $2t'_\wp$, where $t'_\wp: E_\wp'\to E_\wp'$ has the same formula as $t'$ but with
$\tr _{E/\BQ}$  replaced by $\tr _{E_\wp/\BQ}$. The trace $2t''$ of $E$ on
$\Lie (\CG)$ is given by composing $2t'_\wp$  with the embedding
$$E\to E_\wp', \quad
x\to \frac x2-\frac x2\frac {\sqrt{\lambda}}{\mu _p}$$
and the projection 
$$E_\wp '\to E_q^\ur, \quad a+b\sqrt\lambda \to a+b\mu _p.$$
So for $x\in E$, 
\begin{eqnarray*}
t''(x)&=&\tr _{E/\BQ}\left(\frac x2\right)
+\frac x2-\frac {\bar x}2-
\left(\tr _{E_\wp/\BQ_p}\left(\frac x{2\mu _p}\right)
-\frac x{2\mu _p}-\frac {\bar x}{2\mu _p}\right)\mu _p\\
&=&x.
\end{eqnarray*}
Now,  the action $\alpha$ of $E$ on $\Lie (\CG)$
has the trace $2x$. Thus $\Lie (\CG)$ is a two dimensional space of 
$E_q^\ur$ and the action $\alpha$ of $E$ is given by the  usual scalar
multiplication of $E\subset E_q^\ur.$
\enddemo

\demo{{\rm 2.2.2.} The structure of $\CG$}
Let $\CC$ be the component of $C$ in $\CG$.\lb{CC}
In the following we want to identify the structure of 
$[ \CG, \CC]$ as an $\CO _{B, \wp}-\CO _{E, \wp}$-module. First of all let us construct a
special object $[\CG^0,  \CC^0]$.\lb{CG^0CC^0}

Let $\Sigma$\lb{Sigma} denote the following $\CO_\wp$-module:
$$\Sigma=\left\{ \begin{array}{ll} 
\Sigma _2&\hbox{if $\wp$ is not split in $ E$},\\
\Sigma _1\oplus F_\wp/\CO _\wp&\hbox{if $\wp$ is split in $ E$}.
\end{array}\right. 
$$
Here for  any positive integer $h$, let $\Sigma _h$\lb{Sigma_h} denote 
a formal $\CO _\wp$-module of height $h$ over $\CO _ q ^{\ur}$
which is special in the sense  that the induced action on the tangent space is given 
by scalar multiplication, which exists
uniquely up to isomorphism. Let 
$$\CO _{ E, \wp}\times \Sigma \to  \Sigma, \quad (a, x)\to ax $$
be a faithful $\CO _\wp$-linear action such that the induced
 action of $\CO _{ E,\wp}$ on $\Lie (\Sigma)$ is given by 
the reduction $\CO _{ E, \wp}\to \CO _ q    \to k$.

Let us define  a special $\CO _{B, \wp}-\CO _{ E, \wp}$-module 
 $[\CG  ^0,  \CC ^0]$ such that
\begin{itemize}
\item[1.]
$\CG  ^0=\Sigma \oplus \Sigma $; 
\item[2.] The action $\alpha ^0$ of $\CO _{E, \wp}$ is given by the multiplication: 
$$\alpha ^0 (x) (u, v)=(xu, xv);$$
\item[3.] $\CC ^0=\Sigma [N _ E]\oplus 0$ as  $\CO _{ E, \wp}$-modules;
\item[4.] The action of $\CO _{B, \wp}$ is given as  follows:
\begin{itemize}
\item[(a)] If $\wp$ is ramified in $ E$ we fix an 
isomorphism $\CO _{B, \wp}\simeq \RM_2(\CO _\wp)$. Define 
the action  $\iota^0:\CO _{B, \wp}\to \End _{\CO _\wp}(\CG  ^0)$  by
  matrix multiplications.
\item[(b)] If $\wp$ is not ramified in $ E$, then $\CO _{B, \wp}$ is generated
by $\CO _{ E, \wp}$ and an element $\varpi$ such that 
$\varpi x=\bar x\varpi$ and such 
that $\pi:=\varpi ^2$ is a uniformizer of $\CO _\wp$
if $\wp$ is ramified in $B$, and $1$ if $\wp$ is split in $B$.
Then we define the action of $\CO _{B, \wp}$ on $\Sigma^2$ 
by the following formula:
$$\iota ^0 (x)(u,v)=(xu, \bar x v), \qquad \iota ^0 (\varpi) (u, v)=(\pi v, u).$$
\end{itemize}
\end{itemize}
\enddemo

\nnp{Proposition 2.2.3} 
The object  $[ \CG, \CC]$ is  isomorphic 
to $[\CG^0,  \CC]${\rm .} In other words{\rm ,} there is an isomorphism $\phi: \CG\to \CG^0$
such that
\begin{itemize}
\ritem{1.}  $\phi$ is $\CO_{E, \wp}$\/{\rm -}\/linear with respect to the actions $\alpha$, $\alpha ^0${\rm ,} 
\ritem{2.} $\phi$ is $\CO _{B, \wp}$\/{\rm -}\/linear with respect to the actions $\iota$, $\iota^0${\rm ,}
\ritem{3.} $\phi (\CC)=\CC ^0${\rm .}
\end{itemize}
\endproclaim
 
{\it Proof}.
\demo{{\rm 2.2.4.} Case {\rm 1:} $\wp$ is ramified in $ E$}
 In this case $\CC=\CC^0=0$.
Define
$$ \CG  _1=\iota \left( \left( \begin{array}{cc}1&0\\ 0&0 \end{array}\right)\right) \CG  ,\qquad
 \CG  _2=\iota \left( \left( \begin{array}{cc}0&0\\ 0&1 \end{array}\right)\right) \CG  .$$
Then $ \CG  $ is isomorphic to $ \CG  _1\oplus  \CG  _2$ and 
$\iota \left( \left( \begin{array}{cc}0&1\\ 1&0 \end{array}\right)\right)$ 
switches two factors. Since the $ \CG  _i$'s are stable under the action of $\CO _ q^{\phantom{|}}$,
they are isomorphic to $\Sigma_2$. \enddemo

\demo{{\rm 2.2.5.} Case {\rm 2:} $\wp$ is not ramified in $ E$}
Let $ \CG  _1$ (resp.\ $ \CG  _2$) be the maximal $\alpha (\CO _{E, \wp})$-submodule over which 
$\iota(x)=\alpha (x)$ (resp.\ $\iota(x)=\alpha (\bar x)$)
for any $x\in \CO _{ E, \wp}$; then the  $ \CG  _i$'s are $\CO _{ E, \wp}$-modules (via $\alpha $)
of height $1$ and $ \CG  = \CG  _1+ \CG  _2$. The action of $\iota (\varpi)$ gives two
$\CO _{ E, \wp}$-linear morphisms $u:  \CG  _1\to  \CG  _2$ and\break $v:  \CG  _2\to  \CG  _1$ such that
$uv=vu=\pi$. The object $[ \CG  , \iota, \alpha ]$ is completely determined by
$[ \CG  _1,  \CG  _2, u, v]$. As up to isomorphism there is only one special formal 
$\CO _{ E, \wp}$-module of height $1$, so $ \CG  _i$ is isomorphic to $\Sigma $ and
 one of $u$ and $v$ is an isomorphism. In other words,
 up to isomorphism, $[ \CG  _1,  \CG  _2, u, v]$  is isomorphic to 
$[\Sigma , \Sigma , 1, \pi]$ or $[\Sigma , \Sigma , \pi, 1]$. 

By Proposition 2.1.7, the generic fiber of the 
$(\CO _{B, \wp}, \CO _{ E, \wp})$-module
$( \CG , \CC)$ is isomorphic to that corresponding to $(\sqrt{-1}, 1)$; that
is,
$$( \CG , \CC)_{ E_ q ^{\ur}}\simeq (B_\wp/\CO _{B, \wp}, N _ E ^{-1}\CO _{ E, \wp}/\CO _{ E, \wp})$$
 with the action $\iota$ by   multiplication from the
left and the action of $\alpha$ by   multiplication from the right.
It follows that $[ \CG  _1,  \CG  _2, u, v]$  is isomorphic to 
$[\Sigma , \Sigma , 1, \pi]$ and $\CC$ is isomorphic to
$\Sigma [N _ E]\oplus 0$.\hfill\qed
\enddemo

\demo{{\rm 2.3.} Endomorphisms}
Now we assume that ${y}$ is a Heegner point.
We want to study $\End _{\CF} (A_k, C_k)$, where $A_k$ is the reduction 
of $A$ over $\Spec (k)$. 
Let $\BF$ be a finite subfield of $k$ over which 
$[A_k, C_k]$ and $\alpha $ are  defined. In other words, 
$[A_k, C_k]$ is the base change
of some object $[A_\BF, C_\BF]$ with an action of $\CO _E$.  Let 
$\sigma$ be the Frobenius over $\BF$ which acts on $A_\BF$.
By the Honda-Tate theorem and the Tate theorem  (\cite{T} and \cite{W-M}), 
$\End (A_\BF)$ is a semisimple algebra with  center
$\BQ (\sigma)$, and for any prime $\ell$,
$$\End (A_\BF) _\ell\simeq 
\End(A_\BF[{\ell^\infty}])\simeq \End ^\sigma (A_k[\ell ^\infty])$$
where  $\End ^\sigma(\cdot)$ means the commutator of 
$\sigma$ in $\End (\cdot)$.
It follows that\break $\End _{\CF}(A_\BF, C_\BF)$ is also a 
semisimple algebra with  center containing $\CO_{F'(\sigma )}$,  
and  such that for any place $\ell'$ of $F'$, 
\begin{eqnarray*}
\End _{\CF}(A_\BF, C_\BF)_{\ell'}&\simeq &
\End_{\CF}\left(A_\BF[\ell^{'\infty}], C_\BF[\ell ^{'\infty}]\right)\\
&\simeq &
\End_{\CF}^\sigma 
\left(A_k[\ell^{'\infty}], C_k[\ell ^{'\infty}]\right).  
\end{eqnarray*}
Here two $\End _\CF$'s on the right  are defined in the same way as
in \ref{End}.

Fix  $\CO _{B'}$-linear isomorphisms from the level structure of $[A_\BF, C_\BF]$:
$$ \left\{ \begin{array}{ll} 
\kappa  _\wp^2: &[\CG^0, \CC^0]\to [\CG, \CC],\\
 \kappa  _p^{2, \wp}:&  [V_{\BZ, p}^{2, \wp}, \left(N_{ E'} ^{-1}/\CO _{ E'}\right
)_p^{2, \wp}]\to
[T(A)_p^{2, \wp}, C_p^{2, \wp}],\\
\kappa  ^p: &[\wh V_\BZ   ^p, \left(N_{ E'} ^{-1}/\CO _{ E'}\right)^p]
\to [\wh T(A)^p, C^p].
\end{array}\right. \eqno(2.3.1)  $$
Then we obtain isomorphisms:
$$ \End _{\CF}(A_\BF, C_\BF)_{\ell'}
=\left\{ \begin{array}{ll} 
\End  _{\CO _B}^{\wt\sigma}(\CG ^0, \CC^0)
&\hbox{if $\ell'\mid \wp$,}\\
\End _{ \CO_B}^{\wt\sigma} \left(V_{\BZ, p   }^{2, \wp}, 
\left(N_{ E'} ^{-1}/\CO _{ E'}\right
)_p^{2, \wp}\right)_{\ell'}
&\hbox{if $\ell \mid p$ and  $\ell\!\not|\  \wp$}\\
\End _{\CF}^{\wt\sigma} \left(\wh V^p_{\BZ   }, 
\left(N_{ E'} ^{-1}/\CO _{ E'}\right)^p\right)_{\ell'}&\hbox{otherwise,}
\end{array}\right.\eqno(2.3.2)  
$$
where $\wt\sigma$ denotes the  endomorphism  induced from $\sigma$
through the isomorphism  $\kappa $'s. It follows that 
$\End _{\CO _{B'}}(A_\BF, C_\BF)$ is the commutator of $\sigma$ in a
quaternion algebra over $\CO _{F'}$. 
\enddemo  

\nnp{Proposition 2.3.1}
If $\wp$ is split over $ E${\rm ,} then
$$\End_{\CF} (A_k, C_k)=\CO _{ E'}.$$
\endproclaim

\demo{Proof}
From the definition of $\CG ^0$, one sees that 
$$\End _{\CO _B}(\CG ^0)\simeq \End _{\CO_F}(\Sigma)
\simeq \CO _\wp\oplus \CO _\wp.$$
It follows that for $\ell'|\wp$,
$\End_{\CF}(A_k, C_k)_\ell'$ can only be an algebra over $\CO _\wp$ of 
degree at most $2$. Thus $\End _\CF (A_k, C_k)$ is an algebra over $\CO _F$ 
of degree at most $2$.

Obviously, the right side is isomorphic to $\End _{\CF}(\bar A, \bar C)$;
therefore, it is included in the left-hand side.
As $\CO _{ E'}$ is a maximal order, we must have an equality.
\enddemo

\nonumproclaim{Proposition 2.3.2}  
Assume that $\wp$ is not split in $ E${\rm .}
Let $  B    (\wp)$ be the quaternion algebra obtained
by changing invariants at $\tau $ and $\wp${\rm .} Then there is 
 an order  $R(\wp)$   of $B(\wp)$ of type 
$(N(\wp),    E)$ such that 
$$\End_{\CF} (A_k, C_k)\simeq    R(\wp)\otimes \CO _{F'},$$
where $N(\wp)=N\wp ^{1-2\ord _q(N_E)}$.
\endproclaim 

\demo{{P}roof} We need only prove  this identity locally at each place
$\ell'$ of $F'$ using (2.3.2).
In this case, one can show that for $\BF$ sufficiently
large, $\sigma \in F'$. (See [4, \S\S 11.4.4, 11.4.5] for a proof.)

It is easy to check 
that  if $\ell'$ does not divide $\wp$ then $\End_\CF(A, C)_{\ell'}$
is isomorphic to $   R    \otimes \CO _{F', \ell'}.$
It remains to show that 
$\End _{\CO _B}(\CG ^0, \CC^0)$
is isomorphic to $R(\wp)_\wp$.
Notice that $\CC^0$ has only the geometric point $0$, thus does not
play any role in the computation.

Let $D$ denote $\End _{\CO _\wp}(\Sigma)$ which is the maximal order
of a quaternion division algebra over $F_\wp$. The action of $\CO _{E, \wp}$
defines an embedding of $\CO _{E, \wp}$ into $D$.
By a direct computation, we have the following description of
$\End_{\CO _{B, \wp}} (\CG ^0, \CC^0)$:
$$\End _{\CO _\wp}(\CG^0)=\End _{\CO _\wp}(\Sigma\oplus \Sigma)
=\RM_2(D):$$ 
\begin{itemize}
\item[1.] If $\wp$ is ramified in $ E$, then 
$$\End_{\CO _{B, \wp}}(\CG ^0)=D \left( \begin{array}{cc}
1&0\\ 0&1 \end{array}\right).$$
\item[2.] If $\wp$ is ramified in $B $, then
$$\End_{\CO_{B, \wp}} (\CG^0)\simeq
\CO _{ E, \wp}+\CO _{ E, \wp}\wp
 \left( \begin{array}{cc}0& \varpi^{-1}\\ \varpi &0 \end{array}\right),$$
where $\varpi$ is an element in $D$ such that
$\varpi ^2=\wp$ and such that $\varpi x=\bar x\varpi$ for $x\in  E_\wp$.
\item[3.] If $\wp$ is unramified  in both $B$ and  $ E$,
 then 
\medbreak
\hfill ${\displaystyle \End _{\CO _{B, \wp}}(\CG^0, \CC^0)\simeq
\CO _{ E, \wp}+\CO _{ E, \wp}\varpi
 \left( \begin{array}{cc}0&1\\ 1&0 \end{array}\right).}$\hfill\qed
\end{itemize}
\enddemo

\demo{{\rm 2.3.3.} Some remarks and definitions}
Let ${y} _k$ be point in $\CY(k)$. We call $y_k$ a {\em
distinguished point}\lb{dis.poi} in $Y(k)$   if it is the reduction
of Heegner points in $Y(E_q^\ur)$. We can define a
similar concept for any
curve between $Y$ and $X$.

Assume that $y_k$ is a CM-point representing
$(A_k, C_k)$.
If $\wp$ is split in $ E$ then
$$\End_{\CF} (A_k, C_k)\simeq \CO _{ E'}.$$
We write $\End ({y} _k)$  or $\End ^{\rm a}(A_k, C_k)$ for the unique subring of 
$\End_{\CF} (A_k, C_k)$ corresponding to $\CO _ E$ (the superscript
${\rm a}$ stands for the ``admissible endomorphisms").

If $\wp$ is not split in $E$ then 
$$\End_{\CF} (A_k, C_k)\simeq  R(\wp)    \otimes \CO _{F'}$$
with $R(\wp)$ an order of $B(\wp)$ of 
type $( N^ \wp,  E)$.
One may fix the isomorphism   so that
the involution 
$\End_ {\CF_0} (A_k)_\BQ$ induced by the polarization
corresponds to the product of the 
involutions on $  B    (\wp)$ and $F'$ respectively. In this way
the image of $   R(\wp)$  does not depend on  
 the choice of the isomorphism.  Denote this image by $\End (y_k)$ or
$\End^{\rm a} (A_k, C_k)$. 
Notice that two orders in $B(\wp)$ of type
$(N(\wp), E)$ are isomorphic if and only if
they are conjugate under $B(\wp)^\times$.

For a fixed  point $z$ of type $(N\wp,  E)$ with $\End (z)=R(\wp)$, 
the reduction thus defines
a map from the set of CM-points reducing to $z$,
 with conductor $ c $ prime to $N\wp$, to  the set
$$\prod_{v|N  \wp}R(\wp)_v^\times \backslash \Hom (\CO_{ E, v}, R(\wp)_v)$$
 of orientations. This set
 has $2^{s(\wp)}$ elements, where $s(\wp)$ is the number of prime factors
of $   N  \wp$ which do not divide $   {D_E}   $.  Two CM-points $x$ 
and $y$ reducing
to $z$ have the same orientation with respect to $   R   $ if and 
only if they 
have the same orientation with respect to $   R    (\wp)$. We call the {\em orientation}
defined by the reduction of the point $(\sqrt{-1}, 1)$  {\it  the positive orientation}.
\enddemo

\spn{2.3.4} \proclaim{Proposition} Assume that $\wp$ is not split in $ E${\rm .}
Then the map $x\to \End(x)$ gives a bijection between
  the set of distinguished points in $\CX(k)$ and the set of 
conjugacy classes of orders of $ B(\wp)$ of type 
$(N\wp, E)${\rm .} 
\endproclaim

\demo{Proof}
The set of distinguished points on $\CX(k)$ is the set of $\wh F^\times$-orbits
of distinguished points on $\CY (k)$ or any curves between $X$ and $Y$.

Notice that the set of Heegner points in $Y(\BC)$ is represented 
by $(\sqrt{-1}, \wh E)$. Thus the corresponding objects are
$$[A_\gamma, C_\gamma]=\left[\wh V_\BZ\gamma^{-1}\backslash (V_\BR, j),
\quad \left(N_{E'}^{-1}/\CO _E\right)\gamma ^{-1}\right],$$
where $\gamma \in \wh E^\times$. Let $y_\gamma$ be the point
in $Y'(\BC )$ representing the object $[A_\gamma, C_\gamma]$.
Then $y_\gamma$ depends only on the class of $\gamma$ in 
$E^\times \backslash \wh E^\times /\wh\CO _E^\times.$
Thus we may only consider
$y_\gamma$ with $\gamma$ integral and having components $1$ at places over $N\wp$. Then we have isogenies $\phi _\gamma $ from $[A_1, C_1]$ to $[A_\gamma, C_\gamma]$ given by right multiplication of $\gamma ^{-1}$
on $V_\BZ$.  Let $y_{\gamma, k}$ be the reduction of $y_\gamma$
and let $\phi _{\gamma, k}$ denote the reduction of $\phi _\gamma$.
Then we can choose  isomorphisms in (2.3.1) such that 
for places not dividing $\wp$ they are induced by multiplication
of $\gamma ^{-1}$, and that at place $\wp$, $\phi _{\gamma, k}$
induces the identity on $[\CG^0, \CC^0]$.

Using the Honda-Tate theorem, we  show easily that
$$\phi _{\gamma, k}^{-1}\circ
\End ({y}_k)\circ \phi _{\gamma, k}
=\gamma \wh R(\wp)\gamma^{-1}\cap B(\wp)$$
where $R(\wp)$ (resp.\ $B(\wp)$) denotes $\End (y_{1, k})$
(resp.\ $\End (y_{1, k})\otimes \BQ$), and 
we identify both sides as subrings in 
$$\wh B^{\times, \wp}\simeq \wh B(\wp)^{\times, \wp}.$$
As every order of $ B(\wp)$ of type $( N  \wp,  E)$ is conjugate to one 
of $ \gamma R(\wp) \gamma ^{-1}$, 
 the map in the proposition is surjective.

Let $y_{\gamma _1}$ and $y_{\gamma _2}$ be two Heegner points. 
By  (2.3.1), it is easy to see that the injective map
$$\Isom _{\CF}\left([A_{\gamma _1, k}, C_{\gamma_1, k}], 
[A_{\gamma _2, k}, C_{\gamma_2, k}]\right)\to \End _{\CF_0}(A_{1, k})\otimes \BQ,$$
$$\alpha\to \phi _{\gamma_2, k}^{-1}\alpha \phi _{\gamma _1, k}$$
has the image consisting of elements $b$ such that
$$\left\{ \begin{array}{ll} 
&[\CG^0, \CC^0]\cdot b=[\CG^0, \CC^0],\\
&  [V_{\BZ, p}^{2, \wp}, \left(N_{ E'} ^{-1}/\CO _{ E'}\right
)_p^{2, \wp}]\cdot 
\gamma _1^{-1}b=[V_{\BZ, p}^{2, \wp}, \left(N_{ E'} ^{-1}/\CO _{ E'}\right
)_p^{2, \wp}]\cdot
\gamma _2^{-1},\\
&[\wh V_\BZ   ^p, \left(N_{ E'} ^{-1}/\CO _{ E'}\right)^p]\cdot
\gamma _1^{-1}b
=[\wh V_\BZ   ^p, \left(N_{ E'} ^{-1}/\CO _{ E'}\right)^p]
\cdot\gamma _2^{-1}
.
\end{array}\right.  $$
This is equivalent to  
$$\gamma _1^{-1}b\gamma _2\in \wh {R(\wp)}^\times \CO _{F'}^\times.$$
Thus $y_{\gamma _1, k}$ and $y_{\gamma _2, k}$ are in the same orbit under
$\wh F^\times$ if and only if 
$$\gamma _2\in B(\wp)^\times \cdot \gamma _1\cdot \wh{R(\wp)}^\times
\cdot \wh F^\times.$$
This in turn is equivalent to the fact 
that $\End (y_{1, k})$ and $\End (y_{\gamma _2, k
})$ are conjugate in $B(\wp)$.
\enddemo

2.4.\enspace {\it Liftings of distinguished points}.
\demo{{\rm 2.4.1.} A deformation problem}
Let ${y}_k$ be  a  point of $Y(k)$ which represents
an object $[A_k, C_k]$ of $\CF (k)$.
Let $\CG_k$ denote $A[\wp ^\infty]^2$ and let $\CC_k$ denote the component of
$C$ in $\CG_k$.
Let $\alpha_k : E\to \End _{\CF_0}(A)\otimes \BQ$ be a homomorphism
with order 
$$\CO _c:=\{ x\in  E: \alpha (x)\in \End _\CF([A_k, C_k])\}.$$
Assume:
\begin{itemize}
\item[1.] The order $\CO _c$ has  conductor prime to $N\wp$,
and the restriction of $\alpha_k$ on this order has the positive orientation.
\item[2.] The action of $\CO_c $ on  $\Lie(\CG)$ is 
given by the map 
$$i:\CO_{\alpha_k} \to \CO _{\alpha_k} / q    \to k.$$
\item[3.] The object $[\CG_k, \CC_k]$ is isomorphic to the reduction of 
$[\CG ^0_k, \CC ^0_k]$ with respect to both the actions of
$\iota(\CO_B)$ and $\alpha_k(\CO _c)$. 
\end{itemize}
Let us consider the deformation functor $\CF_\alpha$
\lb{CFalpha} over $\CO _ q^{\ur}$-algebra
with residue field $k$  which sends  an algebra $W$
 to the set of isomorphism classes of objects
$[A, C, \alpha ]$. Here 
$[A, C]$ is an object in $\CF (W)$, and 
 $\alpha :\CO_ c \to \End _{\CF}[A, C]$ is a homomorphism such that 
the following conditions are satisfied:
\begin{itemize}
\item The reduction of $[A, C, \alpha ]$ at $k$ is isomorphic 
to $[A_k, C_k, \alpha _k]$.
\item The Rosati involution induced by $\bar\theta _{A}$ takes 
$\alpha (x)$ to $\alpha (\bar x)$ for any $x\in \CO _ c $.
\item The action of $\alpha (\CO _ c )$ on $\Lie (A)_\wp^2$ is given by the
composition: $$\CO _ c \to \CO _ q^{\ur}\to \CO _S.$$
\end{itemize}

\vglue9pt 

\spn{2.4.2} \proclaim{Proposition}\lb{P:CFalpha}
The functor $\CF _\alpha$ is representable by $\CO_ q^{\ur}${\rm .}
\endproclaim

\vglue9pt 

\demo{Proof}
 The deformation space of $[A_k, C_k, \alpha_k]$ is 
the same as that of\break
$[ \CG  _k, \CC _k, \alpha _k]$. This is isomorphic to
$[\CG _k^0, \CC_k^0, \alpha_k^0]$ by Proposition 2.2.3.
Notice that $\CC_k^0=0$. 
Now the conclusion of Proposition \ref{P:CFalpha} follows from the 
fact that the  formal $E_q$-module $\Sigma $ has universal deformation space
$E_q^\ur$.
\enddemo
 
\vglue9pt 

\spn{2.4.3} \proclaim{{C}orollary} Let ${y}_k$ be  a  point of $Y(k)$ which represents
an object $[A_k, C_k]$ of $\CF (k)${\rm .} Then $y_k$ is a distinguished point 
if and only if there is a homomorphism  
$$\alpha _k:\CO _E\to \End _\CF([A, C])$$
such that the above conditions {\rm 1--3} are satisfied{\rm .}
\endproclaim 

\vglue9pt 

\demo{{\rm 2.4.4.} Canonical liftings}
The universal object over $\CO_ q ^{\ur}$ is called the {\em canonical
lifting}\lb{can.lif} of $[A_k, C_k, \alpha _k]$. 
In this way, if $\wp$ is not split 
in $ E$ then for a fixed  distinguished 
point $y_k\in Y(k)$, the set of positively oriented CM-points  with 
conductor $ c $ prime to 
$   N  \wp$, which reduce to $y_k$ modulo $\wp$,
  is bijective to the set of positively oriented homomorphisms 
$ E\to   B    (\wp)$ with conductor $ c $. 

If $\wp$ is split over $F$, then $\alpha $ is an isomorphism,
and the canonical lifting  of $y_k$ is 
a Heegner point $y$ (of characteristic $0$).\enddemo

\spn{2.4.5} \proclaim{Proposition} Assume $\wp$ is  not split in $ E$ and $\ord _\wp (N)\le 1${\rm .}
Let ${y}  _m=[A _m, C_m]$ be the deformation
of ${y} _k =[A_k, C_k]$ to $\CO _ q     ^{\ur}/q^m$
with respect to $\alpha_k${\rm .}
Then $\End ({y}_m)$ has the same localization as 
$\End ({y}_k )$ at places different from $\wp${\rm ,} and
 $$\End ({y}  _m)_\wp=\CO _{ E, \wp}+q ^{m-1}\End ({y}_k)_\wp.$$
In other words{\rm ,} $\End (y_m)$ is the unique sub\/{\rm -}\/order of $\End (y_k)$
of discriminant $\wp^{b_m}N$ where
$$b_m=\left\{ \begin{array}{ll}  m&\hbox{if $\wp$ is ramified in $E$}\\
2m-1
&\hbox{if $\wp$ is unramified in $E$.}
\end{array}\right. 
$$
Moreover the action on the formal module $\CG _m=A_m[\wp^\infty]^2$
is given by the following composition of canonical homomorphisms\/{\rm :}\/
$$\CO _{ E, \wp}+q ^{m-1}\End ({y}_k)_\wp
\to \CO _{E, \wp}/q^m\to \CO _q^\ur/q^m.$$
\endproclaim

\demo{{P}roof}
By a fundamental theorem of Serre and Tate \cite{Dr}, one can show that
$$\End (y_m)=\End (y_k)\cap \End ([\CG _m, \CC_m])$$
where   $\CC_m$ is the component of
$C_m$ in $\CG _m$. It follows that 
$\End ({y}_m)$ has the same localization as 
$\End ({y}_k )$ at places different from $\wp$, and
 $$\End ({y}  _m)_\wp=\End ([\CG _m, \CC_m])
\simeq 
\End ([\CG^0 _m, \CC^0_m])
$$
where $[\CG ^0_m, \CC^0_m]$ is the restriction of 
$[\CG^0, \CC^0]$ on $\CO _q^\ur/q^m$. We want to use the description
in the proof of Proposition 2.3.2 and Gross' result
\cite{G} to describe $\End ([\CG^0 _m, \CC^0_m])$.

As in the proof of Proposition 2.3.2, let
$D$ denote $\End _{\CO _\wp}(\Sigma _{k})$ and
let $D_m$ denote the suborder 
$\End _{\CO _\wp}(\Sigma_{2, m})$. Then by Gross' result
\cite[Prop.\  3.3]{G}, 
$$D_m=\CO _q+q^{m-1}D.$$
Now in 
\begin{eqnarray*}
\End _{\CO _\wp}(\CG ^0_k)&=&\End _{\CO _\wp}(\Sigma_k\oplus \Sigma_k)
=\RM_2(D),\\
\End _{\CO _\wp}([\CG^0_m, \CC^0_m])&=&\End {\CO _\wp}([\CG ^0_k, \CC ^0_k])
\cap \RM _2(D_m).
\end{eqnarray*}

Using the description in the proof of Proposition 2.3.2,
 we have the following:
\begin{itemize}
\item[1.] If $\wp$ is ramified in $ E$, then $\CC=0$ and
$$\End_{\CO _{B, \wp}}(\CG_m ^0)=D_m \left( \begin{array}{cc}
1&0\\ 0&1 \end{array}\right).$$
\item[2.] If $\wp$ is ramified in $B $, then
$$\End_{\CO_{B, \wp}} (\CG^0_m, \CC^0_m)\simeq
\CO _{ E, \wp}+\CO _{ E, \wp}q^{m}
 \left( \begin{array}{cc}0& \varpi^{-1}\\ \varpi &0 \end{array}\right),$$
where $\varpi$ is an element in $D$ such that
$\varpi ^2=\wp$ and such that $\varpi x=\bar x\varpi$ for $x\in  E_\wp$.
\item[3.] If $\wp$ is unramified  in both $B$ and  $ E$,
 then 
\medbreak
\hfill ${\displaystyle \End _{\CO _{B, \wp}}(\CG_m^0, \CC_m^0)\simeq
\CO _{ E, \wp}+\CO _{ E, \wp}\wp^{m-1}\varpi
 \left( \begin{array}{cc}0&1\\ 1&o \end{array}\right).}$\hfill\qed
\end{itemize}
\enddemo
 
\vglue9pt 

  {\rm 2.4.6.} {\it Quasi-canonical liftings}.
We need also to consider the quasi-canon\-ical lifting. 
Let $y$ be a Heegner point representing $[A, C]$ in $\CF (E_q^\ur)$.
Let $D$ be an admissible submodule of $A$ of order $m=\wp^n$ 
($n\ne 0$) prime to $N$ and let $[A_D, C_D]$ be the quotient constructed in
Proposition 1.4.4. Assume that $D$ is connected
(this is automatically satisfied if $\wp$ is not split in $ E$).
Then $[A_D, C_D]$ is an 
object of $\CF (W)$ where $W$ is a finite extension of $\CO _ q^{\ur}$.
Then $[A, C]$ and $[A_D, C_D]$ have the same reduction modulo $q$. 
Notice that $[A_D, C_D]$ is not a canonical lifting of the reduction 
of $[A_k, C_k]$. We call $[A_D, C_D]$ a {\em quasi-canonical lifting}
\lb{qua.can.lif}
of $[A_k, C_k]$.

\vglue9pt 

\spn{2.4.7} \proclaim{Proposition}
 The objects $[A, C]$ and $[A_D, C_D]$ are not
isomorphic modulo $q^2${\rm .}
\endproclaim

\vglue9pt 

\demo{Proof} By the Honda-Tate theorem we need only check whether or not the associated
divisible groups are isomorphic. It suffices to consider the groups
$A[\wp^\infty]^2$ and $A_D[\wp^\infty]^2$. Then the conclusion follows from 
our precise description for $A[\wp^\infty]^2$ and corresponding results
of Gross (\cite[Prop.\  5.3]{G}) on formal groups of dimension  $1$.
\enddemo

\section{Modular forms and $L$-functions}

In this section we will collect various facts about Hilbert modular forms and
associated $L$-functions.
In Section~3.1, we will recall  definitions of modular forms and
Atkin-Lehner's theory on newforms. In Sections 3.2 and 3.3, we will 
give a newform theory for $X$ using 
Jacquet-Langlands correspondence and some work of Waldspurger. 	In Section 3.4, we will first recall
Hecke's theory of $L$-functions and then 
 prove Theorem B in the introduction. 
In Section 3.5, we will study some standard Eisenstein series and theta series attached
to quadratic characters.

\demo{{\rm 3.1.} Modular forms}
\enddemo

3.1.1. {\it Some definitions}.
Let $k$ be a positive integer,   $N$ an ideal of $\CO_F$,
and  $\omega=\prod \omega _v$  a finite character
of $\BA^\times _F/F^\times$ with   conductor dividing $ N  $  such that
$\omega _v(-1)=(-1)^k$ for $v|\infty$. We want to define the space of 
modular forms of (parallel) weight $k$ and level $ N  $. 
See \cite{B}, \cite{Gel} for general background and references.

 Let
$K_0( N  )$\lb{K_0(N)} denote the  following subgroup of $\GL_2(\wh F)$: 
$$   K  _0( N  )=\left\{ \left( \begin{array}{cc} a&b\\ c&d \end{array}\right)\in \GL _2(\wh F):
c\equiv 0\pmod {\wh N  }\right\}.$$
Let $   K  ^\infty$\lb{K_infty} denote the compact subgroup  $\prod _{v|\infty}
\GL_2(F_v)$ of  matrices of the form
$$r(\theta)=(r(\theta _v),\quad v|\infty)\in \GL_2 (F\otimes \BR)$$
where for $\theta =(\theta _v, v|\infty)\in \BR ^g$, 
$$r(\theta _v)= \left( \begin{array}{cc}
\cos 2\pi\theta _v&\sin 2\pi\theta _v\\ -\sin 2\pi\theta_v
 &\cos 2\pi \theta_v \end{array}\right).$$
Let $Z$ denote the center of $\GL_2$. Extend $\omega $ to a character on 
$Z (\BA_F)   K  _0( N  )   K  ^\infty$ by the formula
$$\omega \left( \left( \begin{array}{cc}z&0\\ 0&z \end{array}\right)
 \left( \begin{array}{cc}a&b\\ c&d \end{array}\right)r(\theta)\right)
=\omega (z)\cdot \prod _{\ord _v( N   )>0}\omega _v(a_v)
\cdot \prod _{v|\infty}e^{2\pi i k\theta_v}.$$

Now by a {\em modular form over $F$ of weight $k$, level $ N  $,  character $\omega$}\lb{mod.for}
we mean a function $\phi$ on $\GL_2(\BA_F)$ satisfying the following conditions:
\begin{itemize}
\item[1.] $\phi (\gamma g)=\phi (g)$ for $\gamma$ in $\GL_2(\BQ)$;
\item[2.] $\phi (gk)=\phi (g)\omega (k)$ for $k$ in $Z(\BA_F)   K  _0( N  )   K   ^\infty$;
\item[3.] $\phi$ is slowly increasing: for every $c>0$, and 
any compact subset $\Omega $ of $\GL_2(\BA_F)$, there are a constant $C$ and $N$ such 
that 
$$\phi \left( \left( \begin{array}{cc}
a&0\\ 0&1 \end{array}\right)g\right)\le C|a|^N$$
for all $g\in \Omega $, and $a\in \BA^\times $ with $|a|\ge c$.
\end{itemize}

Let $\psi$\lb{psi} be a character on $F\backslash \BA _F$ defined by
$$\psi (x)=\exp [2\pi i (\tr _{F/\BQ}(x_\infty)-\tr _{F/\BQ}(x_f)].$$
Then every character on $F\backslash \BA _F$ has the form
$$x\to \psi (\alpha x)$$
with some $\alpha \in F$.

Let $dx$ be a Haar measure on $\BA_F$ which is a product of local 
Haar measures $dx_v$ such that if $v$ is archimedean, $dx_v$ is the usual
Haar measure on $\BR$, and that if $v$ is nonarchimedean,
the volume of $\CO _v$ is $1$. In this way, the volume of 
$\BA _F/F$ is $d_F^{-1/2}$ where $d_F$ denotes $\RN (D_F)$.

For a modular form $\phi$ as above, let
$W_\phi (g)$\lb{W_phi} denote the corresponding
Whittaker function on $\GL _2(\BA)$:
\begin{equation}\lb{E:Whittaker}
W_\phi (g)
=d_F^{-1/2}\int _{F\backslash \BA_F}\phi \left( \left( \begin{array}{cc}
1&x\\ 0&1 \end{array}\right)g\right)
\psi (-x)dx \spq{3.1.1}\end{equation}
where  $W_\phi (g)$ satisfies the same 
condition 2 above as $\phi$, and in addition the following 
property:
\begin{equation}
\lb {E:Wit.pro}W_\phi \left( \left( \begin{array}{cc}
1&x\\ 0&1 \end{array}\right)g\right)
=\psi (x)W_\phi (g). \spq{3.1.2}\end{equation} 
Now $\phi$ has the following Fourier expansion
\begin{equation}\lb{Fou.exp}
\phi(g)=
C_\phi(g)+
\sum _{\alpha\in F^\times}
W_\phi \left( \left( \begin{array}{cc}
a&0\\ 0& 1 \end{array}\right)g\right).\spq{3.1.3}
\end{equation}
where
\begin{equation}
C_\phi(g)\lb{C_phi}=d_F^{-1/2}\int _{F\backslash \BA _F}
\phi \left( \left( \begin{array}{cc}
1&x\\ 0&1 \end{array}\right)g\right)dx. \spq{3.1.4}
\end{equation}

We say a form $\phi$ is {\em cuspidal}\lb{cuspidal}, if
 for almost every $g$,
$$C_\phi(g)=0.$$
Thus cuspidal forms are determined by their
Whittaker functions.

Notice that any double coset in 
$$Z(\BA_F )\GL_2(F)\backslash \GL_2(\BA_F)/K  _0( N  ) K_\infty $$
can be represented by an
element of the form $ \left( \begin{array}{cc} y&x\\ 0&1 \end{array}\right)$
with $y\in \BA _F^\times$, $y_\infty>0$, and $x\in \BA _F$.
We say a form $\phi$ is {\em holomorphic}\lb{holomorphic} if $\phantom{\sum^1}$
$$\omega ^{-1}(y)|y|^{-k/2}\phi 
\left( \left( \begin{array}{cc} y&x\\ 0&1
 \end{array}\right)\right)$$
is holomorphic in
$$x_\infty+iy_\infty\in \CH ^g.$$
\enddemo

\spn{3.1.2} \proclaim{Proposition}
Let $\phi$ be a holomorphic form{\rm .} Then
\begin{itemize}
\ritem{1.} $W_\phi
\left( \left( \begin{array}{cc} y&0\\
0&1 \end{array}\right)\right)\ne 0$
only if $y_\infty>0${\rm .}
\ritem{2.}
There is a function $a$ on the set of fractional ideals which vanishes 
on nonintegral ideals{\rm ,} such that
\begin{eqnarray*}
C_\phi
\left( \left( \begin{array}{cc} y&0\\
0&1 \end{array}\right)\right)
&=&\omega (y)|y|^{k/2}a(0),\\
W_\phi
\left( \left( \begin{array}{cc} y&0\\
0&1 \end{array}\right)\right)
&=&\omega (y)|y|^{k/2}a(y_fD_F)
\psi (iy_\infty),
\end{eqnarray*}
where for  $y\in \BA _F^\times$ 
with $y_\infty >0${\rm ,} $x\in \BA_F$,
and  $ {D_F}   $   the inverse of the different ideal of $F${\rm :}
$$ {D_F}    ^{-1}=\{ x\in F: \tr (x\CO _F)\subset \BZ   \}.$$
\end{itemize}
\endproclaim

\demo{Proof}
From \eqref{E:Wit.pro}, one sees that
$$\phi \left( \left( \begin{array}{cc} y&x\\ 0&1 \end{array}\right)\right)
=\omega(y)|y|^{k/2}
\sum _{\alpha\in F}c(\alpha y)\psi (\alpha x)$$
where $c(y)$ is defined by
\begin{eqnarray*}
C_\phi
\left( \left( \begin{array}{cc} y&x\\
0&1 \end{array}\right)\right)
&=&\omega (y)|y|^{k/2}c(0),\\
W_\phi
\left( \left( \begin{array}{cc} y&x\\
0&1 \end{array}\right)\right)
&=&\omega (y)|y|^{k/2}c(y)\psi (x).
\end{eqnarray*}

As $\phi$  is holomorphic in 
$x_\infty+iy_\infty$, it follows that
$c(y)\ne 0$ only if $y_\infty>0$.
Moreover if $y_\infty>0$, then  
$c(y)$ has the decomposition
$c(y)=c(y_f)\psi (iy_\infty)$.

For any $\alpha \in \wh \CO _{F,+}^\times$, $\beta \in \wh\CO_F$,
since
$$W_\phi\left( \left( \begin{array}{cc} y&x\\ 0&1 \end{array}\right)
 \left( \begin{array}{cc}\alpha &\beta \\ 0&1 \end{array}\right)\right)
=W_\phi\left( \left( \begin{array}{cc}y&x\\ 0&1 \end{array}\right)\right)\omega (\alpha),$$
one has
$$c(\alpha y_f)\psi (\beta y_f)=c(y_f).$$
It follows that $c(y_f)\ne 0$ only if $y_f {D_F}   $ is integral, and
that $c(y_f)$ only depends on the ideal $y_f {D_F}   $. In other words,
there is a function $a(m)$ on the ideals $m$ of $\CO _F$ which vanishes
on nonintegral ideals and such that
\medbreak
\hfill $c(y_f)=a(y_f {D_F}   ).$ \enddemo

\demo{{\rm 3.1.3.} Hecke operators}
Now a holomorphic form is uniquely determined
by $a(m)$. We call $a(m)$ {\em the $  m $-th coefficient}\lb{Fou.coe} of
 $\phi$ and denote it as
$a_\phi(m)$ when $\phi$ is referred.

Now let $  m $ be a nonzero ideal of $\CO_F$. We want to define
{\em the Hecke operator}
 $\RT (  m )$\lb{RT(m)} on the space of cusp forms.
Let $H(  m )$ denote the following subset of $\GL_2(\wh F)$:
$$H(  m )=\left\{ \left( \begin{array}{cc} a&b\\ c&d \end{array}\right)\in 
\RM_2(\wh\CO _F): \quad  (d,  N   )=1, c\in \wh N  , 
(ad-bc)\wh\CO_F=\wh  m \right\}.$$
We define $\RT (  m )$ by the formula:
$$(\RT (  m )\phi)(g) =\RN (  m )^{k/2-1}\int _{H(  m )}\phi (gh)dh$$
where $dh$ is a Haar measure on $\GL_2(\wh F)$ such that $   K  _0( N  )$
 has volume~$1$. 
\enddemo

\spn{3.1.4} \proclaim{Proposition} The Fourier coefficients of 
$\RT (  m )\phi$ are given
by the following formula\/{\rm :}\/
$$a_{\RT (  m )\phi}(\ell)=\sum _{a|  m + \ell  }
\RN (  a   )^{k-1}a_{\phi}(  m \ell   /  a   ^2).$$
\endproclaim

\demo{Proof} We need only prove the corresponding statement for $W_\phi$. The set $H(m)$ is stable under  right multiplication by 
$   K  _0( N  )$ and has   disjoint decomposition:
$$H(  m )=\coprod _{a, b, d}
 \left( \begin{array}{cc} a&b\\ 0&d \end{array}\right)   K   _0( N  )$$
where $(a, d)$ are representatives in the class
$$\wh\CO _F\cap (\wh F)^\times/\wh\CO _F^\times$$
such that $ad$ generates $  m $, and for fixed $(a, d)$,
$b$ are representatives in $\wh\CO _F/a\wh\CO _F$.
Now we have
\begin{eqnarray*}
 \RT (m)W_\phi \left( \left( \begin{array}{cc}
y_f&0\\ 0&1 \end{array}\right)\right) 
&\hskip-3pt =&\RN (m)^{k/2-1}
\sum _{a, b, d} W_\phi\left ( \left( \begin{array}{cc} y_fa/d&y_fb/d\\ 0&1 \end{array}\right)\right)\\
&\hskip-3pt =&\RN (  m  )^{k/2-1}
\sum _{a, b, d} W_\phi\left ( \left( \begin{array}{cc} y_fa/d&0\\ 0&1 \end{array}\right)\right)
\psi (y_fb/d)\\
&\hskip-3pt =&\RN (  m  )^{k/2-1}
\sum _{a} W_\phi\left ( \left( \begin{array}{cc} y_fa/d&0\\ 0&1 \end{array}\right)\right)
\sum_{b, d}\psi (y_fb/d).
\end{eqnarray*}
For fixed $a, d$, if $a_{\phi}(\alpha y_fa/d {D_F}   )\ne 0$ then 
$\alpha y_fa/d {D_F}   $ is an integral ideal. In this case
$b\to \psi (\alpha y_fb/d)$ is a character on $\wh\CO _F/a\wh\CO _F$. 
So the last sum over $b$  is $|a|^{-1}$ if this character is 
trivial; otherwise it is $0$. Notice that this character is trivial
if and only if 
$\alpha y_f d^{-1} {D_F}   $ is an integral ideal.
In terms of Fourier coefficients, we obtain
$$a_{\RT (  m )\phi}(\alpha y_f {D_F}   )
=\RN (  m  )^{k/2-1}\sum_{
a, d\atop
d|\alpha y_f  {D_F}}|a/d|^{k/2}|a|^{-1}
a_{\phi}(\alpha y_fa/d {D_F}   ).$$
For any given nonzero ideal $\ell   $ of $\CO _F$, we always can find 
$\alpha, y$ such that $\alpha y_f {D_F}\break   =\ell   $. 
So the above formula gives
$$a_{\RT (  m )\phi}(\ell   )
=\RN (  m  )^{k/2-1}\sum_{
a, d\atop
d|\ell}|a/d|^{k/2}|a|^{-1}
a_{\phi}(\ell    a/d).$$
 Let $  a   = d\CO _F$; then $\ell    a/d=  m \ell   /  a   ^2$, and
$|a|^{-1}=\RN (  m  /  a   )$ and $|d|^{-1}=\RN (  a   )$. The above
formula, therefore, gives the proposition.
\enddemo 

\vglue9pt 

Setting $\ell=1$ in the formula,   we obtain
$$a_{\RT (  m )\phi}(1)
=a_{\phi} (m).$$

\vglue9pt 

\spn{3.1.5} \proclaim{{C}orollary} \hskip-9pt  If $\phi $ is a nonzero eigenform
for all $\RT (m)$ then $a_\phi(1)\ne 0$
and 
$$a_\phi(1)\RT (m)\phi =a_\phi(m)\phi.$$
\endproclaim

\vglue9pt 

\demo{{\rm 3.1.6.} Newforms and multiplicity one}
The Hecke operators are generated by $\RT (\wp)$
with prime $\wp$ and satisfy the formal identity
$$\sum \frac {\RT (m)}{m^s}=\sum _{\wp|N}
(1-\RT (\wp)\wp ^{-s})^{-1} \prod_{\wp\not |N}
(1-\RT (\wp)\wp ^{-s}+m ^{1-2s})^{-1}.$$ 
It follows that if two eigenforms $\phi_1$ and $\phi_2$
 have the same eigenvalues
under all $\RT(\wp)$, then $\phi_1$ and $\phi_2$ are proportional.
 This will not be
true if we only consider Hecke operators $\RT(  m )$ with $  m $ prime to 
some given ideal $  m '$.
Let $ N  '$ be a factor of $ N  $,  let $d\in \GL_2(\wh F)$ be such that 
$$d^{-1}   K   _0( N  ) d\subset    K   _0( N  '),$$
and let $\phi'$ be a form for $   K   _0( N  ')$. Then the function
$\phi'_d(g)=\phi'(gd)$ is a form for $   K   _0( N  )$. The subspace of
$S_k(K_0( N)  )$ generated by these $\phi'_d$ with $ N  '\ne  N  $ is called
the space of {\em old forms}.\lb{old.for}

We say a form $\phi$  for $   K   _0( N  )$ is {\em new},\lb{new.for} if it is
 perpendicular to the space of old forms.
The space  $S_k^\new (   K  _0 (   N  ) ) $ of new forms is generated by {\em newforms}:\lb{newform}
 eigenforms for $\RT(  m ) \quad ((  m ,    N  )=1)$ whose
first coefficients are $1$. 
Then we have the
 strong multiplicity one theorem: 
\enddemo

\spn{3.1.7} \proclaim{Theorem} 
Let $\phi_i${\rm ,}  $(i=1, 2)${\rm ,} be two newforms of weight $k$
of levels  $   N  _1,    N  _2$ respectively{\rm ,}
such that $a_{\phi_1}(\wp)=a_{\phi _2}(\wp)$ for all but finitely many $\wp${\rm .} Then 
$   N  _1=   N  _2$ and $\phi_1=\phi_2${\rm .}
\endproclaim

\demo{Proof}
See \cite [Th.\ 1.4.4 and 3.3.6]{B}  and \cite{Cas}.
\enddemo

In particular if $\phi$ is a newform of level $   N  $ then 
$w_N(\phi)   =\pm \phi$ since
$w_N(\phi)$ is also a newform 
and shares the same eigenvalues as $\phi$, where 
$$w_N(\phi)(g)=\phi \left(g \left( \begin{array}{cc} 0&1\\ t&0 \end{array}\right)\right)$$
with $t$ a generator of $\wh N$.

One application of this is the rank of the Hecke algebra. Let 
$\BT=\break \BT _k(   K   _0(   N  ))$\lb{BT} denote the subalgebra of
$\End _\BC   (S_k(   K   _0(   N  ))$ generated by $\RT(  m )$ with $(  m ,    N  ) =1$.
Then $\BT$ acts faithfully on 
$$S_   N  :=\oplus _{   N  '|   N  }S_k^\new(   K   _0(   N  ))$$
and there is a nondegenerate bilinear form
$$S_   N  \otimes _\BC   \BT  
 \mathop{\longrightarrow}\limits^{(\cdot, \cdot)} 
\BC  $$
such that 
$$(\phi, \RT(  m ))=a_{\RT (  m )\phi}(1)=a_\phi(m ).$$
Now, we have:
\spn{3.1.8} \proclaim{{C}orollary} 
For any linear map $\alpha :\BT\to \BC  ${\rm ,}
there is a unique form $\phi$ such that 
$$a_\phi(m )=\alpha (\RT (  m ))$$
whenever $(  m ,    N  )=1${\rm .}
\endproclaim

\demo{{\rm 3.2.} Newforms on $X$}
As in the modular curve case, one may  define the notion of modular
form on the curve $X$ defined in the introduction:
$$X=B_+\backslash \CH \times \wh B^\times/\wh F^\times \wh R^\times
\cup\{{\rm cusps}\}.$$
 Here we are only interested in 
forms of weight $2$, which are functions $f$ on $\CH\times \wh B ^\times$ 
such that $f(z)dz$ gives a differential form on $X$. 
For $  m $ prime to $ N  $, 
we  define the action by  the Hecke operator
$\RT (  m )$\lb{RT(m)X} by the following formula:
$$\RT(  m ) \alpha= \sum _{\gamma \in   G_m/ G_1}
\gamma^*\alpha,$$
where $\alpha \in  \Gamma   (X, \Omega ^1)$, and $G_m$ and $G_1$ are 
defined as in Section~1.4.
Let $\BT'$ denote the subalgebra of $\End (\Gamma(X, \Omega _X^1))$ generated
 by images of $\RT(m)$ ($(m, N)=1$).
For every newform $\phi$ 
of level dividing $N$, let $\alpha_\phi$ be 
a character of $\BT$  defined by $\phi$ as in Corollary
3.1.8.

The following theorem translates newforms for $   K  _0( N  )$ into
newforms for~$R ^\times$:
\enddemo

\spn{3.2.1} \proclaim{Theorem} {\rm 1.}
The algebra $\BT'$\lb{BT'} is a quotient algebra of 
the Hecke algebra $\BT$ defined in {\rm 3.1.7.}
\begin{itemize}
\ritem{2.}
If $f$ is a newform of weight $2$ for $K_0(N)$
with trivial character{\rm ,} then the eigen subspace of
$\Gamma(X, \Omega _X^1)$ of $\BT$ with character $\alpha _f$
has dimension $1${\rm .}
\end{itemize}

\endproclaim

\demo{Proof}	Indeed,
as in modular curve case, one can show that $\RT'$ is diagonalizable and every character  $\alpha: \BT'\to \BC$
of $\BT'$ corresponds to an irreducible automorphic representation  of 
$(B\otimes \BA)^\times$. By Jacquet-Langlands theory 
\cite{J-L}, this representation corresponds to a cuspidal representation 
of $\GL _2(\BA _F)$. Thus there is a character $\beta: \BT\to \BC$ such that
$\alpha (\RT (m))=\beta (\RT (m))$. So $\BT'$ is a quotient of 
$\BT$. This proves the first part.

For the second part, let $\pi$ be the 
 cuspidal representation of $\GL_2(\BA_F)$ corresponding to $f$.
Then for each place $\wp$ with $\ord _\wp(N)$ odd, the local 
component $\pi _\wp$ of $\pi$ is special or supercuspidal.
(Otherwise $\pi _\wp$ is principal with   trivial central character.)
So $\pi _\wp=\pi (\mu, \mu ^{-1})$ and the conductor of 
$\pi _\wp$ is the square of the conductor of $\mu$. This 
implies that $\ord _\wp (   N  )$ is even.
(See \cite[p.\ 73]{Gel} for a discussion of conductors.) By  Jacquet-Langlands' 
theory \cite{J-L}, $\pi $ corresponds to a unique admissible representation $\pi'$
of $B^\times(\BA_F)$.
 Let $V'$ be the space of the representation of $\pi'$.
Then the proposition is equivalent to the following:
{\em The space of invariant vectors under $\wh R ^\times$ 
has dimension $1$.}
This is a local problem. In other words, we may check the above
problem for each finite place $\wp$. Thus the proof is reduced to
 the following theorem.
\enddemo

\spn{3.2.2} \proclaim{Theorem} 
Let $F$ be a nonarchimedean local field{\rm ,} $B$  a quaternion algebra 
over $F${\rm ,}  $ E$ an unramified quadratic extension of $F$
embedded in $B${\rm .}
Let $\CO _B$ be a maximal order of $B$ containing $\CO _ E${\rm .}
 Let $(\iota, V)$ 
be an admissible representation of $B^\times$ with  trivial
central characters{\rm .} Assume that the conductor of $\iota$
is $2n$ if $B$ is split{\rm ,} and $2n+1$ if $B$  is not split{\rm .} Then 
the subspace of $V$  of vectors invariant under the action
by $ \Gamma  =(\CO _ E+\wp^n\CO _B)^\times$ is  one dimensional{\rm ,}
where $\wp$ is a uniformizer of $F${\rm .}
\endproclaim

\demo{{P}roof. Case {\rm 1.}   $ E$ is split} The theorem in this case
is a special case of a result  of 
Casselman \cite{Cas}. 
Indeed, in this case we may assume that $\CO _B$ is the 
matrix algebra $\RM _2(\CO _F)$ and $\CO _ E$ is the 
algebra of  diagonal matrices. Let   
 $w= \left( \begin{array}{cc} 0& 1\\ \wp&0 \end{array}\right)$.
Then $w^n \Gamma    w^{-n}= \Gamma   _0(\wp^{2n})$.
In the following we assume that $\CO _ E$ is not split
and $n>0$. 
\enddemo
\vglue9pt

\demo{Case  2} {\it $B$ is split{\rm ,} and $\iota$ is a principal series with  conductor $\wp^{2n}$}.
Then $\iota=\pi (\mu, \mu^{-1})$ with $\mu $ a quasicharacter of 
conductor $\wp^n$, and $\mu ^2(x)\ne |x|^{\pm 1}$.
Recall that $\pi (\mu, \mu ^{-1})$ acts by right translation 
on the space $\CB (\mu , \mu ^{-1})$
of locally constant functions $f$ on $\GL_2(F)$ such that
$$f\left( \left( \begin{array}{cc} a&x\\ 0&b \end{array}\right)g\right)=\mu (a/b)|a/b|^{1/2}
f(g).$$
The restriction on $\GL _2(\CO _F)$ gives an isomorphism from 
$\CB (\mu ,\mu ^{-1})$ to the space of  functions $f$ on $\GL _2(\CO_F)$ 
such that
$$f\left( \left( \begin{array}{cc} a&x\\ 0&b \end{array}\right)g\right)=\mu (a/b)
f(g).$$
The subspace of invariant vectors $f$ for $ \Gamma  $ are functions $f$
on $\GL _2(\CO_F)$ such that
$$f\left( \left( \begin{array}{cc}a&x\\ 0&b \end{array}\right)
g \left( \begin{array}{cc} a'&b'\\ c'&d' \end{array}\right)\right)
=\mu (a/b)f(g)$$
for all $ \left( \begin{array}{cc} a'&b'\\ c'&d' \end{array}\right)$ in $ \Gamma  $.
\smallbreak

Since the embedding of $\CO_ E$ into $\RM_2(\CO_F)$ is unique up to conjugation,
the assertion of the theorem does not depend on the choice of the 
embedding. Now write $\CO_ E=\CO_F+\CO_F\varepsilon $ with 
$\varepsilon ^2\in \CO_F^\times\backslash (\CO_F^\times )^2$.
Define an action of $\RM_2(\CO_F)$ on $\CO_ E$ such that
$$ \left( \begin{array}{cc}a&b\\ c&d \end{array}\right) (x+y\varepsilon )=(dx+cy)+(bx+ay)\varepsilon.$$
Then action of $\CO_ E$ on $\CO_ E$ given by multiplication induces an 
embedding $\alpha$
from $\CO_ E$ into $\RM_2(\CO_F)$. Let $g\in \GL _2(\CO_F)=\Aut _{\CO _F}(\CO_ E)$. 
Let $s=\varepsilon\cdot g^{-1}(\varepsilon)^{-1}$ and $g'=g\cdot \alpha (s)^{-1}$. 
Then $g'$ will fix $\varepsilon$, so it has
the form 
$$g'=\beta (a, x):= \left( \begin{array}{cc}1&x\\ 0&a \end{array}\right).$$ The decomposition
$$g=\beta (a, x)\alpha (s)$$
of this form is obviously unique.
The element $g$ is in $ \Gamma  $ if and only if $$\ord (a-1)\ge n.$$
Now it is easy to see that the space of  functions  invariant under $ \Gamma  $
is the one dimensional space generated by  
\begin{equation}\lb{E:f_0}f_0(\beta (a, x)\alpha (s))=\mu (a)^{-1}. \spq{3.2.1}
\end{equation}
\enddemo
\vglue12pt

\demo{Case  3} {\it $B$ is split{\rm ,} and $\iota$ is a special representation}. Now $\iota$ 
is the quotient representation
 of $\CB (\mu |\cdot |^{-1/2}, \mu |\cdot |^{1/2})$ 
with $\mu ^2=1$, modulo the  one dimensional representation
 $\mu\circ\det (g)$. The restriction of this one dimensional representation
on~$ \Gamma  $ has the form
$$(\mu\cdot\det) (\beta (a, x)\alpha (s))=\mu (\RN _{ E/F}s)).$$
Since $\iota$ has a conductor of even order,
so $\mu $ has a conductor of positive 
order. As $\RN _{ E/F}\CO _ E^\times=\CO_F^\times$,
this one-dimensional representation $\mu\cdot\det$ is nontrivial on~$\Gamma  $. 
It follows that the image of $f_0$ defined in \eqref{E:f_0}
 on the space of $\iota$ 
gives a nonzero generator of the space of invariant vectors for
$ \Gamma  $.  \enddemo
\vglue6pt

\demo{Case {\rm 4.} $B$ {\it is nonsplit or $B$ is split but $\iota$ is supercuspidal}}
The   proof was shown to me  by H. Jacquet and will be given in the next 
subsection.\hfill\qed
\enddemo

\demo{{\rm 3.3.} The supercuspidal case} We prove the theorem in the 
supersingular case in two steps. First we prove that
$V^{\CO _E^\times}$ is one dimensional and then we show that this space is
also invariant under $\Gamma$.

\spn{3.3.1} \proclaim{Proposition} 
The subspace $V^{\CO _ E^\times}$ of $\CO _ E^\times$\/{\rm -}\/invariants in $V$ has dimension $1${\rm .}
\endproclaim

\demo{Proof}
Let $\pi$ be a representation of $\GL _2(\CO_F)$ such that $\pi=\iota$ if
$B$ is split, and $\pi$ is the Jacquet-Langlands correspondence of $\iota$ if
$B$ is nonsplit. Let $m$ be the conductor of $\pi$, so that 
$m=2n$ if $B$ is split and $m=2n+1$ if $B$ is not split.

Since $\iota$ has  trivial central character and $\CO _ E/\CO _F$ is unramified,
$\CO_E^\times$ invariants are simply $E^\times$ invariants. 
 According to Waldspurger, (\cite[Th.\ 2]{W2}),
 $V$ has a nonzero vector invariant under $E^\times$ if and only if
$$\varepsilon \left(\frac 12, \pi\otimes \varepsilon _E\right)
=\varepsilon \left(\frac 12, \pi\right)$$
if $B$ is split, and
$$\varepsilon \left(\frac 12, \pi\otimes \varepsilon _E\right)
=-\varepsilon \left(\frac 12, \pi\right)$$
if $B$ is nonsplit, where $\varepsilon_E$ is the quadratic character of
 $F^\times$ attached to the extension $ E/F$. Moreover by Proposition 1 in 
 \cite{W2},
the space of $E^\times$-invariants 
has dimension $1$ if these conditions are verified.

Now the proposition follows from the   identity:
$$\varepsilon\left(\frac 12, \pi\otimes\varepsilon _E\right)=(-1)^m\varepsilon \left(
\frac 12, \pi\right).$$ 
Let $\psi$ be a nontrivial character of $F$ and let $W$ be a vector in the
Whittaker model $\CW (\pi, \psi)$. As the $L$-function of $\pi$ is 1, one has the 
functional equation:
$$\varepsilon (s, \pi, \psi)
\int W\left[ \left( \begin{array}{cc}a&0\\ 0&1 \end{array}\right)\right]|a|^{s-1/2}d^\times a
=\int \wt W\left[ \left( \begin{array}{cc}a&0\\ 0&1 \end{array}\right)\right]|a|^{1/2-s}
d^\times a$$
where 
$$\wt W(g)=W\left[ \left( \begin{array}{cc} 0&1\\ 1&0 \end{array}\right) {}^tg^{-1} \right].$$
Now assume that $W$ is the essential vector. This means that
$$W\left[ \left( \begin{array}{cc}a&0\\ 0&1 \end{array}\right)\right]
=\left\{ \begin{array}{ll} 
1&\hbox{if $|a|=1$;}\\
0&\hbox{otherwise.}
\end{array}\right. 
$$
Then it follows that
$$\varepsilon (s, \pi, \psi)
=\int \wt W\left[ \left( \begin{array}{cc}a&0\\ 0&1 \end{array}\right)\right]|a|^{1/2-s}
d^\times a.$$

Recall that 
$$\varepsilon (s, \pi, \psi )=q^{(1/2-s)m}\varepsilon \left(\frac 12, \pi\right)
$$ where $q$ is the cardinality of the
residue field of $F$. Thus
$$\wt W\left[ \left( \begin{array}{cc} a&0\\ 0&1 \end{array}\right)\right]\ne 0$$
implies $|a|=q^m$. Consequently,
$$\varepsilon \left(\frac 12, \pi\right)=\int _{|a|=q^m}
\wt W\left[ \left( \begin{array}{cc}a&0\\ 0&1 \end{array}\right)\right]d^\times a.$$
Replacing $\pi$ by $\pi \otimes \varepsilon _E$ and $W(g)$ by $W(g)\varepsilon _E(\det g)$,
we obtain the required equality,
\medbreak
\hfill ${\displaystyle \varepsilon \left(\frac 12, \pi \otimes \varepsilon _E\right)=(-1)^m\varepsilon 
\left(\frac 12, \pi \right).}$
\enddemo
 
\vglue9pt

Let $v\in V$ be a nonzero vector invariant under $E^\times$. 
Let $\varpi$ be a
square root of $\wp$ in $\CO _B$. Then 
$v$ is invariant under $   K  _r=1+\varpi ^r\CO _B$ for some sufficiently large $r$.
 
\spn{3.3.2} \proclaim{Proposition} With the notation and assumption as in Proposition {\rm 3.3.1,}
the smallest $r$ such that $v$ is invariant under $   K  _r$
is $r=m$ if $B$ is split{\rm ,} and $r=m-1$ if $B$ is not split{\rm .}
\endproclaim

\demo{Proof}
We will prove the case that $B$ is not split. The case that $B$ is split 
is similar.
Let $f(g)=(\pi (g)v, v)$ be the coefficient function attached to $v$.
Let $ \Phi $ be the characteristic function of $   K  _r$. Its Fourier
transform is (apart from a positive factor) the function 
$\psi (-\tr (g)) \Psi  $, where $ \Psi  $ is the  the characteristic function
of the set $\varpi ^{1+r}\CO _B$. 
The Godement-Jacquet equation \cite{G-J} reads, apart from a 
 nonzero constant factor,
$$\varepsilon (s, \pi, \psi)=
\int f(g^{-1}) \Psi   (g)\psi (-\tr (g))|\det g |^{1/2-s}d^\times g.$$
Since $\varepsilon (s, \pi, \psi)=q^{m(1/2-s)}\varepsilon (1/2, \pi)$, 
we see that the 
integral does not change if we restrict the domain of the integral to the
set $\CO _B^\times \varpi ^{-m}$. Thus $ \Psi  $ must be nonzero on this set,
 which implies that $r\ge m-1$.
Moreover the nonvanishing of the above integral implies that for at least
one $g\in \CO_B^\times \varpi ^{-m}$ the following integral is nonzero:
\vglue-3pt
$$\int _{   K  _{m-1}}f(k^{-1}g^{-1})\psi (-\tr gk)dk.$$
\vglue3pt\noindent 
Since $\psi (-\tr gk)$ does not depend on $k$, 
\vglue-3pt $$\int _{   K  _{m-1}}f(k^{-1}g^{-1})dk\ne 0.$$
\vglue3pt\noindent This implies that 
\vglue-3pt $$v'=\int_{   K  _{m-1}} \pi (k)v dk\ne 0.$$
\vglue3pt\noindent  As $   K  _m$  is a normal subgroup
of $\CO_B^\times$ and $v$ is invariant under $\CO_ E^\times$, $v'$ is invariant under 
the action of $\CO_ E^\times$. So $v$ is a multiple of $v'$ by the previous 
proposition and $v$ is invariant
under  $   K  _{m-1}$.
\enddemo

\vglue9pt

3.4. $L$-\/{\it functions associated to newforms}.
\demo{{\rm 3.4.1.} Definitions}
Let $\phi$ be a newform for $K_0(N)$ of weight $2$ 
with trivial central character. 
Let $a_\phi(m)$ be the Fourier coefficients of $\phi$. Then the $L$-function
for $\phi$  is defined to be
\vglue-9pt
$$\RL (s, \phi)\lb{L(s,chi,phi)}
=\sum _
{m\in \BN_F}\frac {a_\phi (m)}{\RN (m)^s}
=\prod_{\wp\mid N}\frac 1{1-a_\phi(\wp)\RN (\wp)^{-s}}
\prod_{\wp\nmid N}\frac 1{1-a_\phi(\wp)\RN (\wp)^{1-2s}} $$
\vglue3pt\noindent which is absolutely convergent for $s\in \BC$ with $\Re (s)$ sufficiently
large. 
Recall that \vglue-9pt
$$w_N(\phi)(g):=
\phi \left(g \left( \begin{array}{cc}0&1\\ t&0
 \end{array}\right)\right)=\gamma \phi$$
\vglue3pt\noindent  with $\gamma =\pm 1$, where  
 $t$ is an element of $\BA _F^\times$ such that
\begin{itemize}
\item at archimedean places, $t$ has component $-1$;
\item at finite place, $t$ generates $\wh N$.
\end{itemize}
\enddemo

\spn{3.4.2} \proclaim{Proposition}
The function $L(s, \phi)$ is holomorphic in $s$ and satisfies a functional
equation\/{\rm :}\/
\begin{eqnarray*}
L^*(s, \phi):&=&
d_N^{s/2}d_F^s\left[\frac {\Gamma (s)}
{(2\pi)^s}\right]^gL(s, \phi)\\
&=&\gamma L^*(s, \phi).
\end{eqnarray*}
\endproclaim

\demo{Proof}
Let $d^\times x$ be a Haar measure on $\BA _F^\times$ which is a product of 
local Haar measures $dx_v^\times$ on $F_v^\times$
such that $d^\times x_v =dx/x$ if $v$ is archimedean, and such that 
the volume of $\CO _v^\times$ equals $1$ if $v$ is nonarchimedean.
Let $\Lambda (s,\phi)$ denote the function
$$\Lambda (s, \phi)
=\int _{F^\times \backslash \BA _F^\times}
\phi\left( \left( \begin{array}{cc}y&0\\ 0&1 \end{array}\right)\right)
|y|^{s-1/2}d^\times y$$
where $F_+^*$ (resp.\ $\BA _{F, +}$ ) denotes the subgroup
of $F^\times$ (resp.\ $\BA _F^\times $ ) 
of elements which are totally positive 
at archimedean places.
Then $\Lambda (s,  \phi)$ is absolutely
convergent for all $s\in \BC$ and defines
an entire function   on $\BC$.
Using the Fourier expansion of $\phi$ and Proposition 3.1.2, we have
\begin{eqnarray*}
\Lambda (s,  \phi)
&=&\int _{\wh F^\times}
a_\phi(yD_F)
|y|^{s+1/2} d^\times y\cdot \int _{F_\infty^\times}|y|^{s+1/2}\psi (iy_\infty)
d^\times y\\
&=&d_F^{s+1/2}\RL (s+1/2, \chi, \phi)\cdot 
\left(\frac {\Gamma (s+1/2)}{(2\pi )^{s+1/2}}\right)^g.\end{eqnarray*}

Thus we need only prove the corresponding functional equation for
 $\Lambda (s,\phi)$. By definition of $w_N(\phi)$,    
\begin{eqnarray*}
\phi \left ( \left( \begin{array}{cc}y&0\\ 0&1 \end{array}\right)\right)
&=&\gamma \phi \left ( \left( \begin{array}{cc}y&0\\ 0&1 \end{array}\right)
\cdot  \left( \begin{array}{cc} 0&1\\ t&0 \end{array}\right)\right)\\
&=&\gamma \phi \left (\frac {-1}y\cdot \left( \begin{array}{cc} 0&1\\ -1&0 \end{array}\right)
\cdot  \left( \begin{array}{cc}y&0\\ 0&1 \end{array}\right)
\cdot  \left( \begin{array}{cc} 0&1\\ t&0 \end{array}\right)\right)\\
&=&\gamma \phi\left( \left( \begin{array}{cc} {-t}/y&0\\ 0&1 \end{array}\right)\right).
\end{eqnarray*}
Bringing this to our definition of $\Lambda (s, \chi, \phi)$,  we obtain
\begin{eqnarray*}
\Lambda (s,  \phi)
&=&\gamma \int _{F^\times \backslash \BA _F^\times}
\phi\left( \left( \begin{array}{cc}-t/y&0\\ 0&1 \end{array}\right)\right)
|y|^{s-1/2}d^\times y\\
&=&\gamma 
\cdot \RN (N)^{1/2-s}\cdot
\Lambda (1-s,  \phi).\\
\noalign{\vskip-36pt}
\end{eqnarray*}
 \enddemo

\phantom{warming up}

\demo{{\rm 3.4.3.} Remarks}
Let $\varepsilon$ be the character associated to
the imaginary  quadratic 
extension $E/F$. 
 Let $L(s, \varepsilon, \phi)$ be the twisted 
$L$-series:
$$L(s, \varepsilon, f)
=\sum _m\frac {\varepsilon (m)a_\phi (m)}{\RN (m)^s}.
$$
Then this series is essentially the $L$-series
associated to a new form in the space
of the representation $\pi \otimes \varepsilon$
if $\pi$ is the representation associated
to $\phi$.
Thus it has a  functional equation.

The base change of $L(s, \phi)$ is defined
to be the product:
$$L_E(s, \phi):=L(s, \phi)L(s, \varepsilon, \phi).
$$\lb{L_E(s,phi)}
In Section 6,  using  the Rankin-Selberg  convolution method, 
we 
will prove that $L_E(s, \phi)$
has a functional equation with sign
$\varepsilon (N) (-1)^g$.  
\enddemo

\demo{{\rm 3.4.4.} Proof of Theorem {\rm B}}
Let $J_X$\lb{J_X} denote the Jacobian variety of $X$. 
Let $\CT$\lb{CT} be the $\BZ$-subalgebra  in $\End _\BZ(J_X)$
generated by  $\RT (  m )$ ($(  m ,  N  )=1$).
Then $\CT\otimes \BC  =\BT'$.
For every newform $\phi$ 
of level dividing $N$, let $\alpha_\phi$ be 
a character of $\BT$  defined by $\phi$ as in Corollary
3.1.8, 
let  $\CO _\phi$ be the subalgebra of $\BC$ generated by 
Fourier coefficients $a_\phi (m)$ ($(m, N)=1$),
and 
let $J_\phi$ be the maximal abelian subvariety of $J$
killed by $\ker (\alpha_\phi)$. We say two forms 
$\phi _1$ and $\phi _2$ are {\em conjugate} if 
$\ker (\alpha _{\phi _1})=\ker(\alpha _{\phi _2})$,
or equivalently, there is an automorphism $\sigma$ of 
$\BC$ such that $a_{\phi _1}^\sigma(m)=a_{\phi _2}^2$
for all $m$ prime to $N$. \enddemo

\spn{3.4.5} \proclaim{Lemma} {\rm 1.}  $J_X$ is isogenous to $\oplus_{[\phi]} J_\phi$
where $[\phi]$ runs through the conjugacy classes of newforms 
$\phi$ in $S_N${\rm .} \begin{itemize}
\ritem{2.} If $J_\phi$ is nonzero{\rm ,} then 
$\CO _\phi$ is totally real with finite rank over $\BZ${\rm .}
\ritem{3.} If $\phi$ is a newform of level $N${\rm ,} then 
$\Lie (J_\phi)$ is a free module of rank $1$ over 
$\CO _\phi\otimes \BC${\rm .}
\end{itemize}

\endproclaim

\demo{Proof}
Parts (1) and (3) are reformulations of parts (1) and (2) of
Theorem 3.2.1. 
As $\CT$  acts faithfully on $H^1(J, \BZ)$, the characteristic polynomial
of $\RT (  m )$ is monic and integral. It follows  that the
$a(  m )$ are algebraic integers, and that
the subalgebra  $\CO_\phi$ generated by $a_\phi(  m )$ over $\BZ   $
has finite rank. Also as $\BT'$ is self-adjoint, the  characteristic polynomial
of $\RT (  m )$ has only real roots. So $\CO_\phi$ 
is an order in a totally real number field. 
\enddemo

Now fix a newform $f$ of weight $2$ for $K_0(N)$ with trivial character.
Let $A$ denote $J_f$.
Fix a place $\wp$ not dividing $   N  $. Then $A$ has good reduction
at $\wp$. Let $\ell\ne \wp $ be a prime. 
 Then $A \otimes k(\wp)$ has the same $\ell$-adic 
Tate module as $A$,
that is $H^1(A, \BQ _l)$. The local zeta function of $A$ at $\wp$
is  
$$Z_\wp(t)=\det (1-t\Frob(\wp)|H^1(A, \BQ _\ell)).$$
As $\Frob (\wp)^*$ has the same characteristic polynomial as
$\Frob (\wp)$, 
\begin{eqnarray*}
Z_\wp(t)^2&=&
\det \left[(1-t\Frob (\wp))(1-t\Frob (\wp)^*)
\right]\\
&=&\det (1-t(\Frob (\wp)+\Frob (\wp)^*)
+t^2\Frob (\wp)\Frob (\wp)^*).
\end{eqnarray*}
As $\Frob (\wp)$ has  degree $  \RN(\wp)$, 
we see that $\Frob (\wp)\Frob (\wp)^*=  \RN(\wp)$.
Now the congruence relation
$$a(\wp)=\Frob (\wp)+\Frob (\wp)^*$$
implies
$$Z_\wp(t)^2=\det (1-a(\wp) t+  \RN(\wp)t^2).$$
As $\dim H^1(A, \BQ )=2[\CO _f:\BZ   ]$, 
$$Z_\wp(t)=\RN _{\CO_f /\BZ   }(1-a(\wp) t+\RN(\wp)t^2).$$

Thus Theorem B follows, as the $L$-function of
$A$ is defined as
$$L^{(N)}(s, A)
=\prod_{\wp\nmid N}
Z_\wp(\RN (\wp)^{-s})^{-1}.$$
 \bigbreak

3.5. {\it Eisenstein series
and theta series}.

\demo{{\rm  3.5.1.} Some definitions}
Let $k$ be a positive integer.
Let $\chi$ be a quadratic character on 
$\BA_F^\times /F^\times$ with a square-free conductor 
$D_\chi$ such that $\chi_v(-1)=(-1)^k$,
and that $D_\chi$ is prime to $D_F$. 
We extend $\chi$ to $K_0(D_\chi)$ as
in \S3.1.
For  $s$  a complex number,
we define   a function $H_s$ 
on $\GL_2(\BA _F)$ by
$$H_{s} (g)=\left\{ \begin{array}{ll} \left|\frac ad\right|^{s}
\chi (aur(\theta))&\hbox{if $u\in    K   _0( {D_\chi} )$}\\
0&\hbox{otherwise},
\end{array}\right. $$
where 
every element $g\in \GL _2(\BA_F)$ has the form
$$g= \left( \begin{array}{cc} a&b\\ 0&d \end{array}\right)ur(\theta)$$ 
with $ur(\theta)\in    K  _0(1)   K   ^\infty$,
 the standard maximal subgroup of $\GL_2 (\BA_F)$. 
Let $B$ denote the Borel subgroup (the group of upper triangular 
matrices), then $H_s(g)$ is left invariant under $B(F)$. 

For $\Re (s)>1$, the Eisenstein series 
$$E_s(g)=L(2s, \chi)\sum _{\gamma \in B(F)\backslash \GL_2(F)}H_s(\gamma g)$$
is absolutely convergent and defines a (nonholomorphic and noncuspidal) form 
for $   K   _0( {D_\chi})$ of (parallel) weight $k$, and    character
$\chi$. \enddemo

\spn{3.5.2} \proclaim{Proposition} {\rm 1.} The constant term   of $E_s$ 
at $ \left( \begin{array}{cc}
y&0\\ 0&1 \end{array}\right)$is given by the following formula\/{\rm :}\/
$$C_{E_s}
\left( \left( \begin{array}{cc}
y&0\\ 0&1 \end{array}\right)\right)=\left\{ \begin{array}{ll} 
L(2s, \chi )\chi (y)|y|^{s}&\hbox{if $\chi \ne 1$}\\
\zeta _F(2s)|y|^s+d_F^{-1/2}\zeta _F(2s-1)V_s(0)^g|y|^{1-s} 
&\hbox{if $\chi=1$.}
\end{array}\right. 
$$

{\rm 2.} The Whittaker function at
$ \left( \begin{array}{cc}
y&0\\ 0&1 \end{array}\right)$ of $E_s$
 is $0$
if $y {D_F}$ is not integral\/{\rm ;}  otherwise it is
given by the following formula\/{\rm :}\/
$$W_{E_s}\left( \left( \begin{array}{cc}
y&0\\ 0&1 \end{array}\right)\right)=
\frac{1}{\sqrt{d_Fd_\chi}}\sigma _s(y)|y|^{1-s}
\cdot \prod _{v\mid D_\chi}|y_v\pi _v|^{2s-1}
\varepsilon (y_v)\kappa (v)
$$
with 
$$\sigma _s(y)=\prod_{v\nmid D_\chi\atop v\nmid \infty}
\frac {1-\chi (y_v\delta _v\pi _v)|y_v\delta _v\pi _v|^{2s-1}}
{1-\chi(\pi _v)|\pi _v|^{2s-1}}\cdot \prod _{v\mid \infty}
V_s(y_v),$$
where 
\begin{itemize}
\item $\pi _v$ is a uniformizer of $F_v$ such that $\varepsilon (\pi _v)=1$
if $\pi _v\mid D_\chi${\rm .}
\item   
$\kappa (v)$ is a square root of $(-1)^k$ defined by
$$\kappa (v)=|\pi _v|^{1/2}\sum _{a\in (\CO _v/\pi _v)^\times}
\chi (a/\pi _v)\psi _v(-a/\pi _v).$$
\item $\delta\in \wh F^\times$ is a generator of $D_F${\rm .}
\item \hglue.85in ${\displaystyle V_s(y)=\int _{-\infty}^\infty 
\frac {e^{2\pi i y_vx}}{(x^2+1)^{s-k/2}(x+i)^k}dx.}$ 
\end{itemize}
\endproclaim 

\demo{{P}roof} For $\alpha =0$ or $1$, let
$$c_s(\alpha, y)=d_F^{-1/2}\int _{\BA_F/F}
E_s\left( \left( \begin{array}{cc} y&x\\ 0&1 \end{array}\right)\right)
\psi (-\alpha x)dx.$$
Then 
$$W_{E_s}\left( \left( \begin{array}{cc}y&0\\ 0&1 \end{array}\right)\right)=c_s(1, y),
\qquad
C_{E_s}\left( \left( \begin{array}{cc}y&0\\ 0&1 \end{array}\right)\right)=c_s(0, y).$$

The group $\GL_2(F)$ has the Bruhat decomposition
$$\GL_2(F)=B(F)\coprod\coprod_{u\in F}B(F)w
 \left( \begin{array}{cc} 1&u\\ 0&1 \end{array}\right)$$
where
$$w= \left( \begin{array}{cc}0&-1\\ 1&0 \end{array}\right).$$
Therefore 
\begin{eqnarray*}
c_s(\alpha, y)&=
&L(2s, \chi)d_F^{-1/2}\int _{\BA_F/F}
H_s\left( \left( \begin{array}{cc} y&x\\ 0&1 \end{array}\right)\right)
\psi (-\alpha x)dx\\
&&+\ 
L(2s, \chi)d_F^{-1/2}\int _{\BA_F}
H_s\left(w \left( \begin{array}{cc} y&x\\ 0&1 \end{array}\right)\right)
\psi (-\alpha x)dx.
\end{eqnarray*}

By definition the first term is equal to 
$$L(2s, \chi)d_F^{-1/2}\int _{\BA_F/F}\chi (y)|y|^{s}
\psi (-\alpha x)dx$$
which is 
$$L(2s, \chi)\chi (y)|y|^{s}$$ 
if $\alpha=0$; otherwise it is zero.

To evaluate the second integral, we notice that
$$w \left( \begin{array}{cc}y&x\\ 0&1 \end{array}\right)=
 \left( \begin{array}{cc}1&0\\ 0&y \end{array}\right)
 \left( \begin{array}{cc}0&-1\\ 1&xy^{-1} \end{array}\right).$$
Replacing $x$ by $xy$,  we see that the second integral becomes
$$\int  _{\BA_F}
H_s\left(w \left( \begin{array}{cc} y&x\\ 0&1 \end{array}\right)\right)
\psi (-\alpha x)dx=|y|^{1-s}\prod _vV_s(\alpha_v y_v)$$ 
where for $y\in F_v$,
$$V_s(y)=\int  _{F_v}
H_s\left( \left( \begin{array}{cc} 0&-1\\ 1&x \end{array}\right)\right)
\psi (-xy)dx.$$

\demo{The case where $v$ is archimedean}
If $v$ is archimedean, we have the decomposition 
$$ \left( \begin{array}{cc}0&-1\\ 1&x \end{array}\right)=
 \left( \begin{array}{cc}\frac 1{\sqrt{x^2+1}}& \frac {-x}{\sqrt{x^2+1}}\\
0&{\sqrt{x^2+1}} \end{array}\right)
 \left( \begin{array}{cc}\frac x{\sqrt{x^2+1}}&\frac {-1}{\sqrt{x^2+1}}\\
\frac 1{\sqrt{x^2+1}}&\frac x{\sqrt{x^2+1}} \end{array}\right).
$$
It follows that
\begin{eqnarray}\lb{V_s.arc}
V_s(y)&=&\int _\BR\frac 1{(x^2+1)^{s}}
\left(\frac {x-i}{\sqrt{x^2+1}}\right)^k
e^{-2\pi iyx}dx\spq{3.5.1}\\
\noalign{\vskip4pt}
&=&\int _{-\infty}^\infty \frac {e^{-2\pi i yx}}{(x^2+1)^{s-k/2}(x+i)^k}dx.\nonumber
\end{eqnarray}

\demo{The case where $v$ is nonarchimedean}
If $v$ is nonarchimedean, then
$$ \left( \begin{array}{cc}0&-1\\ 1&x \end{array}\right)\in    K  _v$$
if $x\in \CO _v$; otherwise we have the decomposition
$$ \left( \begin{array}{cc}0&-1\\ 1&x \end{array}\right)
= \left( \begin{array}{cc}x^{-1}&-1\\ 0&x \end{array}\right)
 \left( \begin{array}{cc}1&0\\ x^{-1} &1 \end{array}\right).$$
Thus,
$$H_v\left( \left( \begin{array}{cc}0&-1\\ 1&x \end{array}\right)\right)
=\left\{ \begin{array}{ll} \chi _v(x)|x|^{-2s}&\hbox{if $x\notin \CO _v$;}\\
1&\hbox{if $ x\in \CO _v$, $v\hskip3pt \nmid\hskip2.5pt  D_\chi$}\\
0&\hbox{if $ x\in \CO _v$,
$v\mid D_\chi$.}
\end{array}\right. 
$$
It follows that
\begin{eqnarray*}
V_s(y)&=&\sum _{n\ge 1}
\int _{\CO _v^\times}\chi (x\pi _v^{-n})
|x\pi _v^{-n}|^{-2s}\psi (-xy\pi^{-n})d(\pi ^{-n}x)\\
&&+\ \left\{ \begin{array}{ll} 
\int _{\CO _v}\psi (-yx)dx&\hbox{if $v\!\!\not|\ {D_\chi}   $,}\\
0&\hbox{if $v\mid D_\chi$}
\end{array}\right. \\
&=&\sum _{n\ge 1}\chi (\pi _v)^n|\pi _v|^{2ns-n}
\int _{\CO _v^\times}\chi (x)
\psi (-xy\pi^{-n})dx\\
&&+\ \left\{ \begin{array}{ll} 
1&\hbox{if $v\!\!\not|\ {D_\chi}$, $y\in  {D_F}    ^{-1}$}\\
0&\hbox{otherwise.}
\end{array}\right.  
\end{eqnarray*}\enddemo 

\demo{The case where $v\mid D_\chi$}
If $v$ divides $ {D_\chi}$, then $\int _{\CO _v^\times}\chi (x)
\psi (-xy\pi^{-n})dx$ is nonzero only if $y\ne 0$ and $\ord _v(y)=n-1$.
 In this case
it equals  $\chi (y\pi _v^n)\kappa (v)|\pi _v|^{1/2}$.
  So if $v$ 
divides $ {D_\chi}$, we obtain the following  formula for $V_s(y)$:
\begin{equation}\lb{V_s.ram}
V_s(y)=\left\{ \begin{array}{ll} 
|y|^{2s-1}\chi (y)\kappa (v)|\pi _v|^{2s-1/2}
&\hbox{if $y\ne 0$ and $\ord _v(y)\ge 0$}\\
0&\hbox{otherwise.}
\end{array}\right.  \spq{3.5.2}
\end{equation}
Consequently, if $\chi$ is nontrivial,
$V_s(0)=0$ and  the $0$-th Fourier coefficient of $E_s(g)$ is
$$C_s(y)=L(2s, \chi)\chi (y)|y|^{s}.$$ 
\enddemo

\demo{The case where $v\nmid D_\chi$}
In this case
$$
\int _{\CO _v^\times}\psi (-xy\pi^{-n})dx
=\left\{ \begin{array}{ll} 
1-|\pi _v|&\hbox{if $\ord _v(y {D_F}   )\ge n$;}\\
-|\pi _v|&\hbox{if $\ord _v(y {D_F}   )=n-1$;}\\
0&\hbox{otherwise.}
\end{array}\right. $$
It follows that $V_s(y)$ is nonzero only if $\ord _v(y {D_F}   )\ge 0$ and 
in this case
\begin{eqnarray}\lb{V_s.unr}
V_s(y)&=&\sum _{1\le n\le \ord _v(y {D_F}   )}\chi (\pi _v)^n
|\pi _v|^{2ns-n}(1-|\pi _v|)\spq{3.5.3}\\
\noalign{\vskip4pt}
&&+\ 1-(\chi (\pi _v)|\pi _v|^{2s-1})^{\ord _v(y {D_F}   )+1}|\pi _v|
\nonumber \\
\noalign{\vskip4pt}
&=&(1-\chi _v(\pi_v)|\pi _v|^{2s})\sum _{n=0}^{\ord _v(y_v {D_F}    _v)}
\chi _v(\pi _v)^n|\pi _v|^{2ns-n}.\nonumber
\end{eqnarray}
\vglue-24pt\hfill\qed\enddemo 

\spn{3.5.3} \proclaim{{C}orollary}  If $(F, k, \chi)\ne (\BQ, 2, 1)$ then there is a unique holomorphic form $E_{\chi, k}$ of 
weight $k$ and central character $\chi$ for $K_0(D_\chi)$ such
that the $m^{\rm th}$ Fourier coefficient of $E_{\chi, k}$ is given
by 
$$\sigma _{\chi, k-1}(m)
=\sum _{n\mid m}\chi (n)\RN (n)^{k-1}.
$$
\endproclaim

\demo{Proof}
The function $V_s(y)$ can be analytically extended to a function for all 
$\Re (s)>0$ and has  exponential decay 
with respect to $y$.  
When
$s=k/2$, 
$$V_{k/2}(y)=\left\{ \begin{array}{ll} (-2\pi i)^{k}y^{k-1}e^{-2\pi y}&\hbox{if $y>0$}\\
0&\hbox{if $y<0$.}
\end{array}\right. $$
Now, 
$E_s(g)$ can be analytically continued to a form for $\Re (s)>0$
and $E_{k/2}(g)$ is a holomorphic form whose $  m ^{\rm th}$ Fourier coefficients
are given by
$$
\left\{ \begin{array}{ll} 
L(k, \chi)&\hbox{if $  m  =0$;}\\
A_{\chi, k}
\sigma _{\chi, k-1}(m)&
\hbox{if $  m \ne 0$}
\end{array}\right. 
$$
where \medbreak
\phantom{hi}\hfill ${\displaystyle A_{\chi, k}=\frac {(-2\pi i)^{kg}}{(d_Fd_\chi)^{k-1/2}}\chi (D_F)
\prod _{v\mid D_\chi}\kappa (v).}$
\enddemo

\demo{{\rm 3.5.4.} Remarks} 1. When $k=1$ and $\chi$ is the character attached to an
imaginary quadratic extension
$E/F$,  the form $E_{\chi, k}$ is called the theta series associated to
the extension $E/F$ and is denoted as $\theta _{E/F}$ or 
simply $\theta$; thus
$$E_{1/2}=A_{\varepsilon, 1}\theta .$$
 Notice that in this case
$\sigma _{\chi, k-1}(m)$ is the number of integral ideals in 
$E$ with norm $m'$, where $m'$ is the maximal factor of $m$ prime to
$D_E$. We denote this number
simply by $r(m)$.\lb{r(m)}
\smallbreak
2. When $(F, k, \chi)=(\BQ, 2, 1)$, $E_1(g)$ is holomorphic except for the
constant term.
\enddemo

\vglue-12pt

\section{Global intersections}

In this section we will study the N\'eron-Tate height pairing
$\pair{ z  , \RT (  m )  z  }$ of the Heegner points and the CM-points. 
More precisely, we will first show that 
$\pair{ z  , \RT (  m )   z   }$ is the coefficient of a modular
form $ \Psi  $, and  then express the heights as  the
arithmetic intersections 
 using the arithmetic Hodge index theorem \cite{F}.
Finally  we decompose this number as a sum of
local intersections.
Compared with the case $F=\BQ$, there are
 two major difficulties: one is the absence
of cusps which were used to map the modular curves to their Jacobians;
another is the absence of the Dedekind $\eta$-function which was used 
to compute the self-intersection. 
Therefore  we can only obtain an expression 
of $\pair{ z  , \RT (  m  )  z  }$ as a sum of the local intersections
of CM-points which meet properly at special fibers,
modulo some multiple of the coefficients in the Dirichlet 
series $\zeta _ E(s)$ and $\zeta _F(s)\zeta_F (s-1)$.
At the end of this section, we will use the multiplicity-one theorem to
 show that the modular form $ \Psi  $ actually is 
uniquely determined by our expression.
This section contains most of the new ideas of this paper.

\demo{{\rm 4.1. } Height pairing}
\enddemo

4.1.1. {\it Height pairings as Fourier coefficients}. 
In the introduction, we   defined a Shimura curve  $X$
and a Heegner point $ z  $  in the Jacobian 
$J( E)\otimes \BQ$ of $X$.  In Section~1.4 we   defined the Hecke operator
$\RT (  m  )$ as a correspondence for $  m $ prime to $   N  $. 
As in the modular curve case we want to show that 
$$    
\pair{ z  , \RT_  m   z   },\qquad m \in \BN _F,  \quad (  m  ,    N  )=1 $$
are Fourier coefficients  of a holomorphic cusp form for  $K_0(N)$,
where $\pair{\cdot, \cdot}$ is the N\'eron-Tate 
height pairing on $J(\bar F)\otimes \BQ$.
Actually this is a general fact:

\vglue9pt 

\spn{4.1.2} \proclaim{Lemma}
Let $S_   N  $ denote the  sum of 
$S^\new_2(   K  _0 (   N  '))$ for all $   N  '|   N  ${\rm .}
 For any $x\in\Jac (X)(\bar F)${\rm ,} there is a unique 
element $f_x$ in $S_   N  $ such that 
$\pair{x, \RT (  m  ) x}$ is the $  m^{\rm th}$ coefficient  in the 
Fourier expansion of $f$ at $\infty$ for all $  m \in \BN _F$ prime to $   N  ${\rm .}
\endproclaim

\vglue9pt 

\demo{Proof}
Now $\BT' $ also acts on
$J(\bar F)\otimes \BC  $. So $\RT\to \pair {x, \RT x}$ gives 
a linear function  on $\BT'$ and, therefore, on $\BT$.   Now the conclusion follows from
Corollary 3.1.8 and Lemma 3.4.5.
\enddemo

\vglue9pt 

\demo{{\rm 4.1.3.} Height pairings as intersection pairings}
Let $ \Psi  $\lb{Psi} denote the form $f_{z}$ defined in the  
lemma. The purpose of this section is to show that
$ \Psi  $ is determined by the local arithmetic intersections of 
some CM-divisors.

We have constructed an integral  model $\CX$ for
$X$ over $\CO _F$. However this model is not fine enough for
the computation of intersection numbers. Instead of $X$ we will
consider $\wt X$\lb{wtX} which is the Shimura curve corresponding 
to a smaller group $\wt K$ such that the corresponding curve
has a regular model. For example, we may take
$\wt K:=(1+N_E\wh \CO _{B, \wp})^\times \cap\   U  $ where
$  U  $ is an open compact subgroup  of $G(\BA _f)$ which is maximal
at
places dividing $   N     {D_E}   $. When $  U  $ is sufficiently small,
 $\wt X$
has a regular model over $\wt\CX$ over $\CO _F$. As $  U  $ is maximal at
places dividing $   {D_E}   $, $\wt\CX\times \Spec\CO _ E$ is also regular.
Let $\pi: \wt X_E\to  X_E$
 be the  projection induced by the 
inclusion $\wt   K  \to    K  $, and  let $\wt  z  $ be the pullback of
$ z   $ on $\wt X_E$. Then $\wt  z  $ has degree $0$ on each
irreducible component  of $\wt X_E$. The projection formula
for heights gives
$$\pair { z  , \RT (  m  ) z   }=\pair{\wt z  , \RT (  m  )\wt z  }/\deg
\pi.$$
Here the pairing on the right-hand side is the N\'eron-Tate pairing
on the Jacobian of $\wt X\otimes  E$, which by definition is the product
of Jacobians of irreducible components.

We may write $\pair{\wt z  , \RT (  m  )\wt z  }$
as an intersection of arithmetic divisors on $\wt\CX\otimes \CO _ E$
(\cite{F2}, \cite{G-S1}, \cite{G-S2}). More precisely, let $\wh z  $ be the arithmetic
divisor on $\wt \CX\otimes \CO _ E$ which has curvature $0$ on the Riemann
surface $\wt X(\BC  )$ and has zero degree on each irreducible
component $C$ of the special fibers of $\wt\CX\otimes \CO _ E$;
then the Hodge index theorem gives
$$\pair{\wt z  , \RT (  m  )\wt z  }=-(\wh  z  ,\RT (  m  )\wh z  ).$$
Here the right-hand side is the arithmetic intersection. \enddemo

\demo{{\rm 4.1.4.} A formula for $\wh z$}
We write a formula for $\wh z  $ and
let $\eta$\lb{eta} be the divisor 
$$\eta =u^{-1}\sum _x [x]$$
where $u=[\CO _E^\times: \CO _F^\times]$ and $x$ runs through the 
set of positively oriented Heegner points on $X$.
Let $\wt\eta$ be the pull-back of $\eta$ on $\wt X\otimes  E$. 
 Let $\bar \eta$ denote the 
Zariski closure of  $\wt\eta$ on $\wt\CX\otimes \CO _ E$. For each infinite
 place
$\tau$ of $F$, $X_\tau(\BC  )$ is a Riemann surface compactified from 
a quotient of $\CH$. Let $d\mu$ be a volume form on $\wt X_ E(\BC  ) $
 such that on each irreducible component $X_i$ of $\wt X_ E(\BC  )$,
$d\mu $ has volume $1$, and  the pull-back of $d\mu$ on 
$\CH$ is proportional to the Poincar\'e metric $dxdy/y^2$
for $x+yi\in \CH$. Let $g$ denote  
Green's function on $X(\BC   )$ with respect to the Poincar\'e 
volume form $d\mu$:
$$\pp g=\delta _{\eta_i}-\deg(\eta _i)d\mu, \quad \hbox{where}\quad \eta _i=
\wt\eta|_{X_i}.$$
Let ${\wh \eta}$ denote the arithmetic divisor $(\bar \eta, g)$.

Let $\xi$\lb{xi} be the class in $\Pic (X)\otimes\BQ$ which has  component
$\xi_i$ on each geometrically connected component $X_i$ defined as
in the introduction. Then $z$\lb{z} is the class of $\eta-h\xi$
where $h$ is a number such that $z$ has degree $0$ on each irreducible
component of $X$. Let $\wt\xi$ be the pull-back of $\xi$ on $\wt X_ E$. Then 
$\wt \xi$ is the class of 
the bundle $\Omega ^1_{\wt X}[{\rm cusps}]$
divided by its degree.
We will
find an extension of $\wt \xi$ to an arithmetic class
 $\wh \xi$ 
whose curvature is a  multiple of $d\mu$ on each 
component $X_i$. 
We need only do this locally at each place $v$ of $\CO _F$. 

Choose $F'$ as before.  Let $\wt X'$
be  the Shimura curve defined over $F'$ associated to the open compact 
subgroup $K'=\wt   K  \cdot  J $ of 
$\wh  B    ^{'\times}$, where $ J $ is an 
 open compact subgroup  of $\wh \CO _{F'}^\times$
which is maximal at places dividing $   N  $.
Choose $  U  $ and $ J $ sufficiently small so that $\wt\CF:=\CF_{K'}$
\lb{wtCF} is representable.
Let $\CA$ be the universal abelian variety
over $\wt X'$ and  let $\CL_{F'}$ denote $\det (\Lie \CA)^\vee$.
 Then 
by  the Kodaira-Spencer map, 
$\CL_{F'}$ equals   the canonical
bundle $\Omega ^1[{\rm cusps}]$ on $\wt X'$.

If $v$ is an infinite place $\tau$ of $F$, then $\wt X _{\tau}$
can  be embedded into $\wt X'_{\tau}$. The bundle $\CL_v $ has 
a Peterson-Weil metric $\norm\cdot$: for a point $x\in X_\tau (\BC  )$
 representing
an abelian variety $A$, and for an element $\alpha\in \CL_v = \Gamma  
(A, \Omega _A^{4g})$, 
$$\norm {\alpha}^2=(-i)^{g^2}\int _{A (\BC)} \alpha \wedge \bar \alpha .$$
 So we obtain a 
metric on $\xi$; this is nothing else but the 
standard hyperbolic metric up to a constant multiple.

If $v$ is a finite place $\wp$,  we assume that $F'$ is split at 
$\wp$ and that $ J $ is maximal at places dividing $\wp$.
We assume that $\CF _{K', \wp}$ is representable by a 
regular scheme $\wt\CX'$ over $\CO _\wp$. Then we can define a bundle
$\CL_v $ on 
 $\wt\CX'$ the same way.  Let $\CO _\wp^\ur$ be the 
completion of the maximal unramified extension of $\CO _\wp$; then
 $\wt\CX _{\CO _\wp^\ur}$ can be 
embedded into $\wt \CX'_{\CO _\wp^\ur}$. Now  the restriction of $\CL_v$ on $\wt \CX_{\CO _\wp^\ur}$
defines an extension of $\Omega ^1[{\rm cusps}]$.
If $\wp$ does not divide $N$, then this integral structure is the same as that induced
by $\Omega ^1$ on $\wt \CX$ at $v$.

Let $\CL$\lb{CL} be the extension of $\Omega ^1[{\rm cusps}]$ on $\wt \CX$ such that 
$\CL _\wp=\CL\otimes \CO _\wp^\ur$ for every $\wp$. 
Let $\wh\xi$ be the arithmetic divisor class 
of the  hermitian line bundle $(\CL, \|\cdot \|)$ dividing by its degree.
Then $\wh \eta-h\wh \xi$ has curvature zero.

 For simplicity of notation and
computation, we will assume that $ E/F$ is not unramified. In this
case $\eta$ will have the same degree on each geometrically 
connected component of $X$, and 
so will $\wt \eta$. 
  Now we can write $$\wh  z  :=\wh \eta-h\wh \xi+Z$$
where $h$ is a number such that $\wh z  $ has degree $0$ on each 
geometrically connected component of the generic fiber,
and  $Z$\lb{Z} is  a vertical divisor of $\wt\CX\otimes \CO _ E$ 
such that $\wh  z  $  
 has the degree $0$ on any irreducible component of the 
special fibers of $\wt\CX\otimes \CO_ E$. 
In the following subsections we will compute 
$\RT (  m  )\eta$, $\RT (  m  )\wh \xi$, and  $\RT (  m  )Z$ respectively.
\enddemo

\demo{{\rm 4.2.} Computing $\RT (  m  )\eta$}

\spn{4.2.1} \proclaim{Proposition} For $ c $ prime to $   N  ${\rm ,} let 
$$\eta_ c =u_ c ^{-1}\sum _x x,$$\lb{eta_c}
where $u_ c $ is the cardinality of $\CO_ c ^\times/\CO _F^\times${\rm ,}
and where the sum runs through the set of positively oriented 
{\rm CM-}\/points of conductor $c${\rm .} Then for $m$ prime to $N${\rm ,}
$$\RT (  m  ) \eta_1=\sum_{
 c \in   \BN _F\atop
 c  |  m}
r(  m  / c )\eta_ c  $$
where  $r(  m  )$ denotes the number of integral ideals
in $\CO _ E$ with norm $  m ${\rm .}
\endproclaim

\demo{Proof} The map $(\sqrt{-1}, g)\to g$ identifies the
set of CM-points with the set
$$ E^\times \backslash \wh   B    ^\times/\wh    R    ^\times.$$
For any $\CO _F$-module $M$, write
$M^\flat$ for $M\otimes \CO ^\flat$, 
where $\CO ^\flat$ is the product
$\prod _{\wp\nmid   N  }\CO _\wp$.
Also 
write $E^\sharp$  for the group of elements in  $E^\times$
which is a unit  at $\wp$ for any place 
$\wp$  dividing $   N  $. Then the set of positively oriented CM-points 
is identified with
$$E^\times \backslash \prod _{\wp\mid N}
E_\wp^\times R_\wp^\times\cdot B^{\flat, \times} /\wh R^\times 
=  E^\sharp \backslash B    ^{\flat,\times} /
R    ^{\flat, \times}.$$

As $  B    $ is unramified off $   N  $, there is an
isomorphism $\CO _  B    ^\flat\simeq \End _{ \CO _F}(\CO _ E^\flat)$
of the left $\CO _ E^\flat$-algebras.  Now 
 the correspondence $g\to  g\CO _ E^\flat$ 
gives a bijection between the set of positively oriented  CM-points
and the set of
classes of $ \CO ^\flat$-lattices in $ E^\flat$:
$$  E^\sharp \backslash \{  \CO^\flat -\hbox{ lattices in }
 E^\flat\}$$
where $ E^\sharp$ acts on the lattices by  left multiplication.
It is not difficult to show that if a CM-point has   order $\CO_ c $ then
the corresponding lattice class has the form $g\CO_ c ^\flat$ with 
$g$ an element in $ E^\flat$. This shows that the set of CM-points
of conductor $c $ is bijective to 
$$ E^\sharp\backslash E^{\flat, \times} /\CO_c^{\flat,\times}.$$
More precisely, let $S_c$ denote a subset of $E^\flat$ representing
the above set; then $\eta _c$ has the expression
$$\eta_c =u_c^{-1}\sum _{\gamma \in S_c}[\sqrt{-1}, \gamma]$$
where $S_c$ is considered as a subset of $\wh B^\times$
by setting components $1$ at places dividing $N$.

The action of $\RT (  m   )$ on CM-points can be described as follows.
 If $x$
is represented by a lattice  $L$ in $ E^\flat$ then
$\RT (  m  )x $ is the sum of classes of all sublattices $M$
of norm $  m $ (this means that the product of elementary factors of\break $\CO _F$-module $L/M$ is $  m $).

 Let $[g\CO _ c ^\flat]$ be a lattice class
with $g\in  E^{\flat, \times}$.  Then the multiplicity of
$[g\CO _{c}^\flat]$ in $\RT (  m  ) \eta_1 $ is equal to
$u_1^{-1}$ times the number of pairs
$$(\gamma, k)\in S_1
\times E^\sharp/\CO _c^\times$$
such that $kg\CO _c^b$ is a sublattice of $\gamma \CO _E^b$ of norm
$m$, or equivalently
$$\gamma^{-1} gk\in \wh \CO _E^b, \qquad 
\RN (\gamma^{-1} gk)=m/c.$$
Now the surjective map 
$$S_1\times E^\sharp/\CO _c^\times \to 
(E^\flat)^\times /\CO _E^{\flat,\times}, 
\quad (g, k)\to \gamma ^{-1}gk \pmod {\CO _E^{\flat, \times}}$$ 
is $[\CO _E^\times :\CO _c^\times ]$ to $1$. Thus the multiplicity is equal to
\medbreak
\phantom{hi} \hfill ${\displaystyle u_1^{-1}[\CO _E^\times :\CO _c^\times ]\#\{\gamma \in \wh \CO _E^\flat/\CO
_E^{\flat,
\times} :\quad
\RN (\gamma )=m/c\}=u_c^{-1}r(m/c).}$
\enddemo

\vglue9pt

Let $\eta_ c ^0$\lb{eta_c^0} denote the sum of $\eta_  a   $ for all 
$  a    | c $ and $  a   \ne \CO _F$,
and define
\begin{equation}
\RT ^0(  m  )\lb{RT^0} \eta=
\sum _{ c  |  m }\varepsilon( c ) \eta_{  m  / c }^0. \spq{4.2.1}
\end{equation}
Then   $\RT ^0(  m  )\eta$ is disjoint to $\eta$. As
$r(  m  )=\sum _{  n |  m }\varepsilon (  n),$
we obtain:

\spn{4.2.2} \proclaim{{C}orollary}  If $  m $ is prime to $   N     {D_E}   ${\rm ,} then
$$\RT (  m  )\eta=\RT^0 (  m  )\eta+r(  m  )\eta.$$
\endproclaim

\vglue9pt

\demo{{\rm 4.3.} Computing  $\RT (  m  )\wh\xi$}
\enddemo

4.3.1. {\it Some definitions}.
Let  $\pi :U\to V$ be a finite flat morphism of integral schemes.
Let ${\bf Pic}(U)$, ${\bf Pic}(V)$ be categories of line bundles
on $U$ and $V$ respectively.
Then we can define the  pull-back functor $\pi^*: {\bf Pic} (V)\to 
{\bf Pic} (U)$ as
usual,  and the  norm functor $\RN _\pi: {\bf Pic}(U)\to {\bf Pic} (V)$ 
as follows. If $L$ is a line bundle on $U$ then $ \RN _\pi(L)$ is a 
line bundle on $V$ which is locally generated by $\RN _\pi (\ell)$
with $\ell$ a  section of $\CL$  such that 
$$\RN_\pi (f\ell)={\rm Norm} (f)\RN _\pi (\ell)$$
where ${\rm Norm} $ is the norm map $f_*\CO _U\to \CO _V$  for the algebra extension $\CO _V\to f_*\CO _U$. It
follows from the  definition that
if $L=\CO_U (D)$ for a divisor $D$ on $U$,
 then $N_{\pi}(L)$ is canonically isomorphic to
$\CO _V(\pi _*D).$

If $W$ is an integral subscheme of  $U\times V$ such that the projection
from $W$ to $U$ is finite and flat then we can define a functor 
$W: {\bf Pic} (V)\to {\bf Pic} (U)$ as $W(L)=\RN _{\pi _U}\pi^* _V(L)$
where $\pi _U$, $\pi _V$ are projections from $W$ to $U$ and $V$
respectively. We may extend this definition linearly to any 
correspondence $W$ of $U\times V$ such that $W$ has 
all irreducible components finite and flat over $U$.

 It is easy to see that at the generic fiber
$$\RT (  m  ) \xi=\sigma _1(  m  ) \xi.$$
The following proposition gives the corresponding formula for
$\RT (  m ) \wh \xi$.

\spn{4.3.2} \proclaim{Proposition} There is a morphism 
$$\psi_m: \RT (  m  )\CL \to \CL ^{\sigma_1(  m  )}$$
 such that the following conditions are verified\/{\rm :}
\begin{itemize} 
\ritem{1.} Let $ c \in \BN _F$ be such that
$$\psi_m\left(\RT (  m  )\CL\right)= c  \CL ^{\sigma _1(  m  )}.$$
Then  for each finite place $\wp${\rm ,} 
$$\ord _\wp ( c )=2\sigma _1(  m\wp^{-\ord _\wp (m)})
\sum _{i=0}^n i\RN(\wp^{n-i}).$$
\ritem{2.}
Let $\psi$  be the following function on $\wt X(\BC)$
$$\norm \psi_m (x):=\frac {\norm {\psi_m \beta }}{\norm \beta},$$
where  $\beta$ is a nonzero element  in $\RT (  m  ) (\CL)(x)${\rm .} Then 
$$\norm \psi_m (x) =\RN (  m  )^{2\sigma_1 (  m  )}.$$
\end{itemize}
\endproclaim

\demo{Proof} We need only prove the corresponding  statement on $\wt\CX'$.
For this we extend $\RT (  m  )$ to $\wt\CX '$ by the formula (1.4.1).
By Proposition 1.4.2, we have 
the following modular interpretation 
for $\RT (  m )$:
For any object $[A, C]$ of $\wt\CF (S)$,  
$$\RT (  m  )[A, C]=\sum _D [A_D, C_D]$$
where $D$ runs through the set of  admissible submodules of $A$ of
order $m$, $A=A/D$, $C_D=C+D/D$.
Let $\CX _  m $ be the subscheme of $\CX'\times \CX'$ which represents the 
isogenies $A_1\to A_2$
with admissible kernel  of order $  m $; then 
$\RT (  m  )$ is induced by $\CX _  m $.

Let $\pi: \CA _1\to \CA _2$ be the universal isogeny over 
$\CX _  m $, and let $p_1, p_2$\break be the projection of 
$\CX _  m $ to $\CX'$.  Then $p_i^*\CL =\det\Lie (\CA_i)^\vee$.
The morphism\break $\pi ^*: \Lie (\CA_2)\to \Lie (\CA_1)$
therefore induces a morphism of line bundles on 
$\CX'$:
$$ \RN _{p_1}(p_2^*\CL) \to \RN _{p_1}(p_1^*\CL).$$
Notice that 
by definition $\RT (  m  ) \CL =\RN_{p_1}p_2^*(\CL )$,
and $\RN _{p_1}p_1^*\CL=\CL ^{\sigma_1 (  m  )}$.
Therefore, we obtain a morphism of line bundles:
$$\psi_  m  :\RT (  m  )\CL \to \CL ^{\sigma _1(  m  )}.$$

To prove (1),  we need only check the proposition locally at 
each finite place $\wp$ prime to $ N  $. Write $  m =  m '\wp^n$ with $(  m ',\wp)=1$. Then 
$\psi_  m $ is factored as a composition of $\psi _{  m '}$ and
$\psi _{\wp^n}$:
\[\RT (  m  )(\CL )=\RT (\wp^n)\RT (  m ')\CL \lrar{
 {\RT (\wp^n)\psi _{  m '}} }
\RT (\wp^n)\CL ^{\sigma _1(m ')} \lrar{{\psi _{\wp^n}^{\otimes \sigma _1(  m ')}}}
\CL ^{\sigma _1(  m  )}.\]
As $\RT (m')$ is \'etale at $\wp$, it follows that if  $\psi _{\wp^n}$ has order $t$ at $\wp$, then  $\psi _  m $ has order
$\sigma _1(m')t$.

Let $x: \Spec W\to \CX _\wp$ be a strictly henselian point represented by an 
abelian variety $A$ with  ordinary reduction. 
  Then 
$$\RT (\wp^n)(\CL _x)=\otimes _  D\det\Lie (A/  D)^\vee$$
where $  D$ runs through the set    of admissible submodules of 
$A$ of order $  m $.
Fix an isomorphism $\CO _{B, \wp}\simeq \RM_2(\CO _\wp)$.
Let $G$ denote the $\CO _\wp$-module
$ \left( \begin{array}{cc} 1&0\\ 0&0 \end{array}\right)A[\wp^\infty]^2$
of dimension $1$; then $\det\Lie (A)=\Lie (G)^{\otimes 2} $.
It follows that
$$\RT (\wp^n)\CL =\prod _H(\Lie G/H)^{\otimes{-2}}$$
where $H$ runs through the set of submodules  of $G$ of order $\wp^n$,
and that the morphism $\psi: \RT (\wp^n)\CL \to \CL ^{\sigma _1(\wp^n)}$
is induced by morphisms $\pi ^*: \Lie (G/H)^\vee\to \Lie (G)^\vee$.
Let 
$$0\to H'\to H\to H''\to 0$$
be the formal-\'etale decomposition.
Then 
$$\Lie (G)^\vee\bigg /\pi ^*\Lie (G/H)^\vee\simeq 0^*(\Omega _{H})
\simeq 0^*(\Omega _{H'})$$
where $0$ is the $0$-section of $G$.

Now $G$ has a decomposition $G=\Sigma _1\oplus F_\wp/\CO _\wp$ 
where $\Sigma _1$ is a formal\break $\CO _\wp$-module of height $1$.
It follows that $H$ has the form
$H=\Sigma _1[\wp^i]\otimes G_{n-i, \lambda} $
where $0\le i\le n$, $\lambda \in \wp^{i-n}\CO _{\wp}/\CO _{\wp}$,
and $G_{n-i, \lambda}$ is the subgroup with the generic fiber
$\{(\lambda x, x): x\in   \wp^{i-n}\CO _{\wp}/\CO _{ \wp}\}$.
Thus  
$$0^*(\Omega _H)\simeq 0^*(\Omega _{\Sigma _1[\wp ^i]})
\simeq \Lie (\Sigma _1)^\vee/\wp ^i\Lie (\Sigma _1)^\vee
\simeq \CO _\wp/\wp ^i.$$
It follows that  the quotient of $\psi$ has the order
$$\sum _{i=0}^n 2i\RN(\wp^{n-i}).$$

It remains to prove (1).  Recall that $\RT (  m  )\CL (x)$ is equal to
$$\otimes _D\det\Lie (A/D)^\vee$$
where $D$ runs through the set of  admissible submodules 
of order $  m $. As
$\psi$ is induced by the maps
$$\pi _D^*: \Omega ^1_{A/D}\to \Omega ^1_A,$$
 the norm of $\psi$ is the product of the norms  of
$$\det \pi _D^*: \det  \Gamma   (\Omega ^1_{A/D})
\to \det  \Gamma   (\Omega ^1_A)$$
which is $(\deg \pi _D)^{1/2}=\RN (  m  )^2$.
Then for any infinite place $\tau$,

\phantom{hi}\hfill ${\displaystyle\norm {\psi}_\tau(x)=\RN (  m  )^{2\sigma_1 (  m  )}.}$
\enddemo

\vglue9pt

\spn{4.3.3} \proclaim{{C}orollary} 
Let $\phi$ be a function on the set of elements of $\BN_F$ prime to 
$   {\rm N}  $ with values in the group of 
arithmetic divisors on $\CO _F$ defined by the formula
$$\RT(  m  )  \wh \xi=\sigma _1(  m  )(\wh \xi+\phi(  m  )).$$
Then $\phi $ is quasi\/{\rm -}\/additive\/{\rm :} for any $  m '$ and
$  m ''$ such that $(  m ',   m '')=1$
then  
$$\phi (  m '  m '')=\phi (  m ')+\phi (  m '').$$
\endproclaim

\demo{Proof}
We decompose $\phi(m)=\sum \phi (m)_v [v]$ where $v$ runs through 
all places of $F$. Then by the proposition,
$$\phi (m)_v=
\left\{ \begin{array}{ll}  
c\sigma (\wp^{\ord _\wp(m)})^{-1}\sum _{i=0}^n i\RN(\wp^{n-i})&\hbox{if $v=\wp$ is finite},\\
c\log \RN (m)&\hbox{if $v$ is infinite} 
\end{array}\right. 
$$
where $c$ is some fixed constant.
Thus $\phi $ is additive for coprime 
$m$'s.
\enddemo

\demo{{\rm 4.4.} Computing $\RT (  m  )Z$}
\enddemo

4.4.1. {\it Decompositions}.
For each finite place $\wp$ of $F$, let $V_\wp$ denote the group
of $\BQ$-divisors of $\wt \CX$ supported in the fiber over $\wp$
modulo the subgroup of $\BQ$-divisors of connected components.
Then we have the decomposition
$$Z=\sum _\wp Z_\wp$$
where $Z_\wp$ are elements in $V_\wp$.
We want to study $\RT (  m )Z_\wp$ for $  m $ prime to $ N   {D_E}   $.
If we choose different models $\wt\CX$, then the decomposition is
preserved by the pull-back maps.  So we assume that $\wt \CX$ has the 
same level structure as $\CX$ at the place $\wp$.

\vglue9pt

\spn{4.4.2} \proclaim{Proposition} Assume that $\wp$ is split in $B${\rm .}
Then 
$$\RT (  m )Z_\wp=\sigma _1(  m  )Z_\wp.$$
\endproclaim

\vglue9pt

\demo{{P}roof} By definition
$\RT (  m  ) Z$ is a unique solution to the equations
$$(\RT (  m  )\wh\eta -h\RT (  m  )\wh \xi+\RT (  m  )Z, P)=0$$
for any irreducible vertical divisor $P$ on $\wt\CX\otimes \CO_ E$. 
 As $\wt\CX\otimes  E$ is smooth at the places 
not dividing $   N  $, we need only check that the differences 
$$Z_1=\RT (  m  )\hat\eta-\sigma _1(  m  )\hat\eta\quad \hbox{and}\quad
Z_2=\RT (  m  )\wh\xi-\sigma _1(  m  )\wh\xi$$ both have degree $0$ on each 
irreducible component   of $\wt\CX$ over $\wp$ dividing $   N  $.
For $Z_2$ this follows from Proposition 4.3.2. It remains to study 
$Z_1$.
\enddemo
\vglue9pt

\demo{Case 1}  {\it $\wp$ does not divide $ N  $}.
 In this case, each geometrically 
connected 
component of $\wt\CX_\wp$ has only one irreducible component. Thus
$Z_\wp=0$.
\enddemo

\vglue9pt

\demo{Case 2} {\it $\wp$ split in $ E$}. Let $\wt   K  _0$ 
denote the level structure obtained by 
replacing the level structure $   K   _p$ by the maximal one $\CO _{B, \wp}^\times$.
Let $\wt\CX _0$ denote the corresponding Shimura curve. Then the 
natural map $\wt \CX\to \wt \CX _0^{\phantom{|}}$ induces a bijection 
on the set of connected components. Over $\wt\CX _0$ we have the 
divisible $\CO _\wp$-module $ \CG   ^1$ of height $2$ and $\wt \CX $ classifies
the ``cyclic" submodules $C$ of $ \CG  ^1$ of order $\wp ^{\ord _\wp( N  )}$.
For a fixed irreducible component $D$ of the special fiber of $\wt\CX _0$
over $\wp$, by Proposition 1.3.2,
the set of irreducible components of $\wt\CX $ over $D$
is indexed by the types of the subgroups over the ordinary points over $D$.
By Proposition 2.2.3, all divisors $\eta _ c $ will have  ordinary 
reduction at $\wp$ and the corresponding subgroups are of the same type:
either all \'etale or all formal. It follows that all CM-divisors $\eta _ c $
with  positive orientation will reduce to the same irreducible component
of $\wt\CX _\wp$ over $D$. This implies that $Z_1$ has  degree $0$ on
each irreducible component of $\wt\CX _\wp$. \enddemo

\demo{Case 3} $\wp$ {\it is split in $B$ and inert in $ E$}. 
We claim that each 
connected component of $\wt \CX $ over $\wp$ has only one irreducible component.
With the notation as Section 1.3.1 and Proposition 1.3.2, 
the set of  irreducible components of $\wt \CX _\wp$
over $D$ is indexed by $\BP ^1(F_\wp)/   K   ^\wp$ where 
$   K   ^\wp=F_\wp ^\times  R  ^\times _\wp$. As $ R  $ contains $\CO _{ E, \wp}$,
and $F_\wp ^\times \CO _{ E, \wp}^\times = E_\wp ^\times$,
it suffices to show that $\BP ^1(F_\wp)$ has only one orbit under the 
action of $ E_\wp ^\times$ for any embedding $ E_\wp\to \RM_2(F_\wp)$.
Up to a conjugation, we may identify $\BP ^1(F_\wp)$ as the set of 
surjective $F_\wp$-homomorphisms, from $ E _\wp $ to $F_\wp$ and the action
of $ E_\wp ^\times $ is given by   multiplication on $ E_\wp$. It follows that
$\BP ^1(F_\wp)$ has one element $\tr :  E _\wp\to F_\wp$. As
the pairing 
$$ E _\wp\times  E_\wp\to F_\wp,\qquad
(x, y)\to \tr (xy)$$
is nondegenerate, any other surjective $F_\wp$-homomorphism 
$\phi:  E _\wp\to F_\wp$ will have the form
$$\phi (x)=\tr (ax)$$
where $a$ is a nonzero element of $ E_\wp$. In other words, 
$\phi =a(\tr)$, or the action of $ E _\wp$ on $\BP ^1(F_\wp)$ 
is transitive. Consequently, each connected component of $\wt\CX _\wp$
has only one irreducible component. As in case 1, we have $Z_\wp=0$.\hfill\qed
\enddemo

It remains to consider the case where $\wp$ is not split in $B$. 
The conclusion of the previous proposition will definitely not be true.
But we have the following:

\spn{4.4.3} \proclaim{Proposition} Assume that $\wp$ is not split in $B${\rm .} Then
for any element $D\in V_\wp${\rm ,}  
there is a holomorphic $V_\wp$\/{\rm -}\/valued cusp form $f$ 
of weight $2$ and 
 level $   K   _0( N^\wp)$
such that for all but finitely many  $  m $, the $  m $\/{\rm -}\/coefficient of 
$f$ is given by $\RT (  m  )D${\rm ,}
where $ N^\wp$ denotes $N   \wp ^{-\ord _\wp ( N   )}${\rm .}
\endproclaim

\demo{Proof} Actually by Proposition 1.3.4, 
the set of  irreducible components of
 $\wt \CX $ is identified with 
$$S_{\wt K_0}=B(\wp )^\times \backslash \wh {B (\wp )}^\times /\wh F^\times 
\GL _2(\CO _\wp) \wt K   ^\wp.$$
The group $V_\wp$ is therefore identified with a subgroup of the space $\wt V$
of  complex  functions on $S_{\wt   K   _0}$. 
By Jacquet-Langlands  theory \cite{J-L},
the action of the Hecke correspondences is  factored through the
action of the Hecke algebra of   holomorphic cusp forms of 
 weight $2$ and  level 
$\wh F^\times \GL _2(\CO _\wp)\wt   K   ^\wp$, so the proposition is true 
with the level structure $   K  _0( N  )$ replaced by 
$K_0(N^\wp)_N\wt   K   ^N$.

Using the  pull-back of divisors,  we notice that the minimal level of the 
forms which have  $\RT (m)D$ as Fourier coefficients does exist and
does not depend  on the choice of $\wt K$. Thus this  minimal level must be $K_0(N^\wp)$.
\enddemo

\demo{{\rm 4.4.4.} Some definitions}
Let $\CS $\lb{CS}
 denote the vector space of complex-valued functions on $\BN_F$
modulo  an equivalence relation  so that 
two functions $a$ and $b$ are equivalent if and only if
there is some element $M$    such that $a(\ell )=b(\ell)$ for any 
 $\ell$ prime to $M$.
The strong multiplicity theorem 3.1.7 implies that the map
$$f \longrightarrow \wt f: n  \to a_f( n  )$$ is an embedding
from $S_N $ into $\CS$.
We say a function $h$ in $\CS$ is {\em quasi-multiplicative}\lb{qua.mul}
if there is an $M\in \BN _F$ such that
 $$f(  m  n  )=f(  m  )f( n  )$$
 for all  $  m ,  n  \in \BN_F$ such that
$$(  m ,  n  )=(  m  n  , M)=1.$$
 For a quasi-multiplicative function $f$, a function $h$ is called
{\it an $f$\/{\rm -}\/derivative}  if
$$h(  m  n  )= f(  m  )h( n  )+f( n  )h( m  )$$
for all $(  m ,  n  )$ as above.

Let $\sigma _1$ and 
$r$ denote the elements in $\CS$ defined by: $   m \to \sigma _1(  m  )$
and $  m \to r(  m  )$ respectively, and let $\CD_N$\lb{CD} be the subspace of $\CS$ generated by
$\sigma _1$, $r$, $\sigma _1$-derivatives, and $r$-derivatives,
and the Fourier coefficients
 corresponding to the old cusp forms of weight $2$.
Then we have:
 
\spn{4.4.5} \proclaim{Proposition} Let $\wh \Psi  $ denote the image of $ \Psi  $ in $\CS${\rm .} 
 Then in $\CS${\rm ,}  
$$\wh \Psi (m)  = -\intn{\wh\eta, \RT ^0(m)\wh\eta}/\deg \pi\pmod {\CD_N}.$$ 
\endproclaim

\demo{Proof}
By discussions in   4.1.3 and 4.1.4, 
for $  m $ prime to $   N     {D_E}   $, 
$$\wh \Psi   (  m  )=-\intn {\wh\eta-h\wh\xi+Z,  \RT(  m  )
 (\wh \eta -h\wh \xi+Z)}/\deg \pi.$$
Now we have shown:
\begin{itemize}
\item  $\RT (  m  )Z_\wp=\sigma _1(  m  )Z_\wp$ if $\wp$ is split in $B$,
and $  m \to \RT (  m )Z_\wp$ is given by an old cusp form of weight $2$
if $\wp$ is not split in $B$;
\item  $\RT (  m  )\wh\xi=\sigma _1(  m  )(\xi +\psi 
(  m  ))$ with $\psi$ quasi-additive;
\item   $\RT (  m  )\wh \eta=r(  m  )\wh\eta
+\RT ^0(  m  )\wh\eta$. 
\end{itemize}
It follows that 
 
\vglue-6pt 
\phantom{hi}\hfill ${\displaystyle \wh  \Psi   (  m  )=-\intn{\wh\eta, \RT ^0(  m  )\wh\eta}/
\deg \pi \pmod {\CD_N}.}$
\enddemo
  \vglue6pt
\demo{{\rm 4.5.} A uniqueness theorem}
Now we are going to prove that the relation in Proposition 4.4.5
 determines 
a newform projection of $ \Psi  $ uniquely: 

\bigbreak

{\elevensc Proposition 4.5.1}.
{\it Let $f$ be an element in the space $S_   N  $
such that in} $\CS${\rm ,}
$$\wh f\equiv 0\pmod {\CD_N}.
$$
{\it Then $f$ is an old cusp form of weight $2${\rm .}}

\demo{{P}roof}
We start from the following:

\spn{4.5.2} \proclaim{Lemma}
Let $\alpha _1, \cdots , \alpha _\ell$ 
be distinct nonzero quasi\/{\rm -}\/multiplicative 
elements in $\CS${\rm .} Then the equation 
$$(c_1\alpha _1+h_1 )+\cdots +(c_\ell \alpha _\ell + h_\ell)=0$$
in $\CS$ does not have a nonzero solution
 $$x=(c_1,h_1,\cdots, c_\ell , h_\ell),$$ 
where for each $i${\rm ,}  $c_i$ is a  constant and 
$h_i$ is an   $\alpha _i$\/{\rm -}\/derivative{\rm .}
\endproclaim

\demo{Proof} Assume that the lemma is not true, then we will have one solution
$x_0=(c_1, h_1, \cdots, c_l, h_l)$. Let $M$ be an element in $\BN _F$  such that
$$(c_1\alpha _1( n  )+h_1( n  ))+\cdots+ (c_\ell \alpha _\ell ( n  )+ 
h_\ell ( n  ))=0$$
for any $ n  $ prime to $M$. Let $  m $ be any ideal prime to $M$;
then for any $ n  $ prime to $  m  M$, we have 
$$(c_1\alpha _1(  m  n  )+h_1(  m  n  ))+ (c_2\alpha _2(  m  n  )+ h_2(  m n  ))+
\cdots=0.$$
So we have a new solution
$$x_1=(c_1\alpha _1(  m  )+h_1(  m  ), \alpha _1(  m  )h_1, \cdots,
c_\ell \alpha _\ell (  m  )+h_\ell (  m  ), \alpha _\ell (  m  )h_\ell ).$$
If $h_1(  m  )\ne 0$ then we obtain a solution
$$x'=x_1-\alpha _1(  m  )x_0=(h_1(  m  ), 0, \cdots) ,$$
in which $h_1=0$ and $c_1\ne 0$. Doing this for each $i$, 
  we obtain a solution 
in which every $h_i=0$ but some $c_i$ will not be $0$.
We need only to show that $\alpha _1, \cdots, \alpha _\ell$ are linearly independent. 
This  is similar to the proof of the linear independence 
of the characters of a group. 
\enddemo

Now go back to the proof of our proposition. Decompose $f$ into a sum of newforms
of levels dividing $   N  $ and  forms of  type 
$\phi\left(g \left( \begin{array}{cc}1&0\\ 0&d \end{array}\right)\right)$, where
$d\ne 1$ is a divisor of $ N  $ in $\wh F^\times$
and $\phi$ is a newform of level
$ d^{-1}  N  $.
Then the above lemma implies the proposition if we 
can show that $\sigma_1$ and $r$ are distinct and not in the image 
of newforms. For any quasi-multiplicative $a$ in $\CS $, we define
the Dirichlet series
$$L(s, a)=\sum _ n   a( n  )\RN( n  )^{-s}$$ 
which is well-defined modulo finitely many factors.
Then it is easy to see up to finitely many factors,
$$L(s, r)=\zeta_ E(s),\qquad L(s, \sigma _1) =\zeta _F(s-1)\zeta _F(s).$$
So $L(s, r)$ has a pole at $s=1$ and 
$L(s, \sigma _1)$ has a pole
at $s=2$.
If $r$ is multiple of $\sigma_1$ in $\CS$, then
$L(s, r)$ should be equal to a multiple of $L(s, \sigma _1)$ 
up to finitely many Euler factors. This is impossible as they have different 
poles. The same argument shows that
$\sigma_1$ and $r$ should not be  equal to $\hat f$ for any cusp form $f$,
as $L(s, \wh f)=L(s, f)$ is holomorphic at $s=2$ and $s=1$.\hfill\qed
\enddemo 
 
\demo{{\rm 4.5.3.} Remarks}
The number $(\wh \eta, \RT ^0(m)\wh\eta)/\deg \pi$ does not depend on
the choice of $\pi:\wt X\to X$. Let us denote it by
$(\eta, \RT ^0(m)\eta)$. As two divisors $\wh \eta$ and 
$\RT ^0(m)\wh\eta$ are disjoint at the generic fibers,  it has 
decomposition
$$(\eta, \RT ^0(m)\eta)
=\sum _v(\eta, \RT ^0(m)\eta)_v$$
where $v$ runs through the set of all places of $F$, 
 and
$$(\eta, \RT ^0(m)\eta)_v=
\sum _{w|v}(\wh \eta, \RT ^0(m)\wh \eta )_w/\deg \pi$$
where $w$ runs through all places of $E$ over $v$.

Assume that $m$ is prime to $ND_E$, then by
(4.2.1),  the computation of $ \Psi  $ modulo oldforms is reduced to 
the computation of 
$$(\eta, \eta_ c ^0)_v:=(\wh \eta, \wh \eta ^0)/\deg \pi.$$
In the following section we will first compute
the local intersections then  add them together.
\enddemo

\section{Local intersections}

In this section we are going to compute $(\eta, \RT (  m  )^0\eta)_v$
where $v$ is a place of $E$.
 We  follow the method of Gross-  Kohnen-Zagier \cite{G-Z2}.
 However  we are working on CM-points with the discriminants not necessarily
coprime. Again, we need to assume that every factor of $2$ is split
in $E$.

\demo{{\rm 5.1.} Archimedean intersections}  In this subsection
we want to compute the infinite local intersections 
$(\eta, \eta ^0_ c )_v$
where $ c \in \BN_F$ is prime to $   N     {D_E}   $
and $\pi: \wt X\to X$ is  some covering of $X$ constructed as in Section 4.1.
First of all let us assume that $v$ is over $\tau$, the embedding
 chosen  in the introduction.

\demo{{\rm 5.1.1.} Intersections as  Green\/{\rm '}\/s functions}
Let $ R  _i$'s be nonconjugate orders of $B$ of type $( N  ,  E)$. 
Then $X(\BC  )$ is a union of Riemann surfaces
$$X_i= R  _{i +}^\times \backslash \CH =
 R  _i^\times \backslash \CH^\pm.$$
We want to compute the intersections 
$(\eta_i, \eta ^0_{ c , i})_v$ separately,
where $\eta _i$ and $\eta ^0_{ c , i}$ are restrictions of $\eta$ and
$\eta ^0_{ c }$ on $X_i$.
Let $g_i(x, y)$  be  Green's function on the compactification of
 $X_i$ with respect to the Poincar\'e
metric $d\mu$:
$$\frac{\partial _x\bar\partial _x}{\pi i}g_i(x, y)
=\delta _y-d\mu (x).$$
We can linearly extend $g_{i}$ to a function on the set of disjoint pairs
of divisors. 

\spn{5.1.2} \proclaim{Lemma}
  ${\displaystyle (\eta_i, \eta _{ c , i}^0)_v=g_{i}(\eta _{i}, 
\eta _{c, i}^0).}$  
\endproclaim

\demo{Proof}
Let $\wt X_i$ be the part of $\wt X$ which projects into $X_i$. Then
the Riemann surface  $\wt X_i$ is a union of Riemann
surfaces of the form: $\wt X_{i, j}=\bar \Gamma   _{i, j}\backslash \CH$,
where $\bar \Gamma   _{i, j}\subset \PGL _2(\BR)^+$ acts freely on $\CH$.
Let $\eta _{i, j}$ and $\eta ^0_{c, i,j}$ be the restrictions of $\wt\eta_i$
and $\wt\eta _{ c , i}^0$ on $\wt X_{i, j}$; then by construction of
$\wh \eta_i$, as  in 4.1.4,
$$(\wh \eta_i, \wh\eta _{ c , i}^0)_v=\sum _{j}g_{i, j}(\eta _{i, j}, 
\eta _{c, i, j}^0)/\deg \pi,$$ 
 where $g_{i, j}(x, y)$  are  Green's functions 
on the compactifications of $\wt X_{i, j}$'s with respect to the Poincar\'e
metric. Now the lemma follows easily from the projection formula
$$g_i(x, y)=\sum _jg_{i, j}(\pi ^{-1}(x), \pi ^{-1}(y))/\deg \pi$$
for any distinct $x$ and $y$ in $X_i(\BC)$.
 \enddemo

\demo{{\rm 5.1.3.} Construction of Green\/{\rm '}\/s functions}
Now  $g_{i}(x, y)$ can be constructed  as follows \cite{G2}: 
for $s\in \BC  $ with $\Re(s) >1$, define the function 
$$G_s(z, w)=Q_{s-1}\left(1+\frac {|z-w|^2}{2\Im z\Im w}\right)$$
where  $Q_{s-1}(u)$ is  the Legendre function 
of the second kind:
$$Q_{s-1}(u)=\int _0^\infty (u+\sqrt{u^2-1}\cosh t)^sdt.$$
\lb{Q_{s-1}(u)}
Then the function on $X_i$,
$$g_{s,i}(z, w)=\sum _{\gamma \in \Gamma  _{i}}
G_s(z, \gamma w)$$
is convergent and has a simple pole at $s=1$ with residue $1/\chi_{i}$
 where $\Gamma _i$ is the image of $R_i^\times$ in 
${\rm PSL} _2(\BR)$ and $\chi_{i}$ is the Euler characteristic of $X_{i}$. 
Then we have an identity
$$g_{i}(z, w)=\lim _{s\to 1}\left(g_{s, i}(z,w)-\frac 1{s(s-1)\chi_i}\right).$$

 Let $x, y$ be two points on $X_\tau (\BC  )$ represented by $z,w$ on $\CH$.
Let $u_x$ and $u_y$ be the orders of stabilizers of $x$ and $y$ in $\Gamma _i$
respectively.
Let $P_x$ and $P_y$ be the sets of points on $\CH$ mapping to $x$ and $y$, respectively.
Then,
$$g_{i}(x/u_x, y/u_y)=\lim _{s\to 1}
\left[\sum _
{(z, w)\in \Gamma_i \backslash P_x\times P_y}
g_s(z, w)-\frac {u_x^{-1}u_y^{-1}}{s(s-1)\chi_i}\right].$$
Applying this to components in $\eta _i$ and $\eta _{ c , i}^0$, we see that 
Lemma 5.1.2 gives
\begin{equation}(\eta_i, \eta _{ c , i} ^0)_v
=\lim _{s\to 1}
\left[\sum _
{(z, w)\in \Gamma_i  \backslash P_i\times P_{ c , i}^0}
g_s(z, w)-\frac {\deg \eta \deg \eta _ c ^0}{s(s-1)\chi_i}\right],\hskip.5in \spq{5.1.1}
\end{equation}
where $P_i$ and $P_{ c , i}^0$ are the sets of points in $\CH$ mapping to 
components of $\eta_i$ and~$\eta _{ c , i}^0$.

\demo{{\rm 5.1.4.} Descriptions of {\rm CM-}\/points}
We may identify $\CH^\pm$ with $$\Hom _\BR(\BC  , \RM_2(\BR))$$ such that
if $z=g(\sqrt{-1})\in \CH^\pm$ with $g\in \GL _2(\BR)$, then the corresponding
element $\phi _z: \BC  \to \RM_2(\BR)$ takes $a+bi$ to
$g \left( \begin{array}{cc}a&b\\ -b&a \end{array}\right)g^{-1}$. In this way,
the CM-points on $X_i$ are those points induced by a homomorphism
$\phi: K\to B$ with order given by $\phi^{-1}( R _i)$.
 
For two points $z$ and $w$ in $\CH^\pm$ corresponding to two 
homomorphisms $\phi _z$ and $\phi _w$ in $\Hom (\BC  , \RM_2(\BR ))$
it is easy to check that
$$1+\frac {|z-w|^2}{2\Im z\Im w}=-\frac 12 \tr (i_zi_w),$$
where $i_z=\phi _z(i)$ and $i_w=\phi _w(i)$.
It follows that $z$ and $w$ are in the  same connected component and $z\ne w$
if and only 
if   $-\frac 12\tr (i_zi_w)>1$.
Let $\CP_i $ (resp.\ $\CP _{ c , i} ^0$) 
denote the inverse image of $\eta_i$ and $\eta_{ c , i} ^0$ on $\CH^\pm$,
and let $\CP _{ c , i}$ denote the union of $\CP_i$ and $\CP _{ c , i} ^0$.
 Then we have
$$
(\eta_i, \eta_{ c , i}^0)_v
=
\lim _{s\to 1}
\left\{
\sum_{
(z, w)\in \CP_i \times \CP _{ c , i}^0/   R _i   ^\times
\atop
-\frac 12\tr (i_zi_w)>1}
Q_{s-1}\left(-\frac 12\tr (i_zi_w)\right)+
\frac {\deg \eta_i\deg \eta_{ c , i}^0}{s(s-1)\chi}
\right\}.
$$

For an element $ c \in \BN _F$, define
\begin{equation}
u_{v, s}( c ,i   )=\sum_{
(z, w)\in \CP \times \CP _{ c ,i}/   R    ^\times\atop
-\frac 12 \tr (i_zi_w)> 1}
Q_{s-1}\left(-\frac 12\tr (i_zi_w)\right). \spq{5.1.2}
\end{equation}
Then formula (5.1.1) gives
\begin{equation}
(\eta_i, \eta_{ c ,i} ^0)_{\tau _1}=\lim_{s\to 1}
\left(u_{v, s}( c ,i)-u_{v, s}(1  ,i)+
\frac {\deg \eta_i\deg \eta_{ c , i} ^0}{s(s-1)\chi}\right).\hskip.5in \spq{5.1.3}
\end{equation}
\enddemo

\demo{{\rm 5.1.5.} Linking numbers} 
Each  pair  $(z, w)\in \CP_i \times \CP _ {c, i} $ of $X_i$ 
determines two homomorphisms $\phi _z$ and $\phi _w$ from 
$K$ to $B$ such that 
$\phi _z(\CO _K)\subset  R _i$ and $\phi _w(\CO _ c  )\subset  R  _i$ and such
that $\phi _z$ and $\phi _w$ have the same orientation.
Let $a, b$ be two totally positive elements of $ c $ and $ {D_E}   $ respectively
such that both $a$ and $b$ are prime to $ N  $ and that $\sqrt{-b}$ is in
$ E$. Let $e_1=\phi_z(\sqrt {-  b})$ and
$e_2=\phi _w(c\sqrt {-  b})$ in $   R _i  $.
 As $\phi _z$ and $\phi _w$ have  positive 
orientation, we have
$$(ae_1-e_2)^2\equiv 0\pmod {4   N  }.$$
In other words there is an $n\in  N   a^{-1}b^{-1}$
 such that $$\tr (e_1e_2)=-2ab+4n ab.$$ 
It is easy to verify that $n$ is independent 
of the choice of $c$ and $d$.
So $n\in  N    c ^{-1} {D_E}    ^{-1}$.
We call $n$ the {\em linking number} of $z$ and $w$ 
(or $\phi_z$ and $\phi _w$) and denote it by
$n(z, w)$ ( or $n(\phi _z, \phi _w)$).

As
$$-\frac 12\tr (i_zi_w)=1-2\tau _1(n),$$
 formula (5.1.2) becomes
\begin{equation}
u_{v, s}( c , i   )
=
\sum_{
n\in  N   c ^{-1} {D_E}   ^{-1}\atop
 \tau (n)<0}
\varrho _{v}( c , n,i)
Q_{s-1}(1-2\tau _1(n)), \spq{5.1.4}
\end{equation}
where $\varrho _v( c ,  n,i)$ is the number of  conjugacy classes
 of pairs $(z, w)\in \CP_i\times \CP _{ c , i}$ such that
$n(z, w)=n.$ \enddemo

\demo{{\rm 5.1.6.} Summing up} We need to sum up formula (5.1.3) for all
$i$, but only for  $\varrho _v( c ,  n,i)$'s and residues.
Let $P(n)_i$ denote the set of conjugacy classes of 
pairs $(\phi_1, \phi_2)\in  \Hom ( E,  B)^2   $ such that
$$\phi_1(\CO _ E)\subset  R  _i,\quad
\phi _2(\CO _ c )\subset  R  _i,
\qquad
n(\phi _1, \phi _2)=n.$$ 
Then any  pair $(\phi _1, \phi _2)$ 
defines two CM-points $(z, w)\in \CH\times \CH$
with conductor $1  $ and $ c $ respectively.
These two points are  in $\CP_i\times \CP _{ c , i}$ if and only if the 
morphism $\phi _z$ defined by $z$ has positive  orientation.
The orientation group $\CW$ acts freely on $\cup _iP(n)_i$ which, therefore, has  cardinality 
$\varrho _\tau( c , n):=2^s\sum_i\varrho _v( c , n, i)$.
\lb{varrho_v(c,n)}

Now we want to treat the residue term in formula (5.1.3). \enddemo

\spn{5.1.7} \proclaim{Lemma} The numbers $\deg \eta _i${\rm ,} $\deg \eta _{c, i}${\rm ,} $\chi _i$
do not depend on $i${\rm ,} if they are nonzero.
\endproclaim

\demo{Proof} The natural projection from $X_\tau (\BC)$ onto the set of 
its connected components is given by the determinant map
$$X(\BC)\to \pi _0(X(\BC):=F_+^\times\backslash \wh F^\times /\wh F^{\times, 2}
\wh \CO _F^\times,$$
$$(z, g)\in \CH \otimes \wh B^\times\to \det g.$$
Now $\eta_c$ is $u_c^{-1}$ times the sum of CM-points   represented by
$(\sqrt{-1}, g)$ with $g$ in $E^\times \backslash \wh E^\times
/\wh F^\times \wh \CO _c^\times$. The determinant map restricted on these 
CM-points is given by the norm homomorphism
$$\RN _{E/F}: E^\times \backslash \wh E^\times
/\wh F^\times \wh \CO _c^\times\to F_+^\times\backslash\wh F^\times /\wh F^{\times, 2}
\wh \CO _F^\times.$$
Thus the pre-image of every point in $\pi _0(X(\BC)$ has the same cardinality
if it is not empty. This implies that $\deg \eta _{c, i}$, 
therefore,
 $\deg\eta _i$ and  $\deg \eta ^0_{c, i}$ do not depend on $i$ if they 
are nonzero. 

It remains to show that $\chi_i$ does not depend on $i$. Recall that $\chi_i$
is the volume of $X_i$ with respect to the measure $dxdy/y^2$ on $\CH$ 
times
an absolute constant. We need only show that the volume of $X_i$ does not
depend on $i$. For this we use Hecke's correspondence $\RT (m)$. By 
the definition
of $\RT (m)$ in  Section 1.4, the induced action of $\RT (m)$ on $\pi _0(X(\BC))$ is
given by $[x]\to \sigma _1(m)[mx]$, where $[x]$ denote a point represented
by $x\in \wh F^\times$. On the other hand, $\RT (m)$ changes volume form
$dxdy/y^2$ to $\sigma _1(m)dxdy/y^2$. Thus all connected components of
$X_\tau (\BC)$ must have the same volume.
\enddemo

\demo{{\rm 5.1.8.} Intersection on other archimedean places}
Now we want to compute the archimedean intersection for
places of $ E$ over $\tau _2, \cdots, \tau _g$. For this we need to describe
the conjugation $X _{\tau _k}(\BC   )$ of $X(\BC   )$ over $F$.
Let $  B   (\tau _k)$ denote a quaternion algebra obtained 
from $  B   $ by switching invariants at $\tau _1$ and $\tau _k$.
Fix an order $   R    (\tau _k)$ of $  B   (\tau _k)$ of type $(   N  ,    E)$;
then 
$$X_{\tau _k}(\BC   )\simeq   B(\tau _k)^\times\backslash
\CH^\pm\times \wh B (\tau _k)^\times/\wh R    (\tau _k)^\times.$$
So the above  formulas (5.1.2)--(5.1.4) and Lemma 5.1.7
 for $(\eta, \eta_ c  ^0)$ work for each $\tau _k$.

More precisely for each infinite place $\tau _k$, let
$\varrho _{\tau _k}( c , n)$ be defined as above for $B(\tau _k)$;
then we have the following:

\spn{5.1.9} \proclaim{Proposition} For each infinite place $v$ of $ E$ over an
infinite place $\tau _k$ of $F${\rm ,} the local intersection
$(\eta, \eta_ c  ^0)_v$ is given by the formula
$$\lim_{s\to 1}
\left(u_{\tau _k, s}(c)-u_{\tau _k, s}(\CO _F)+
\frac {\deg \eta\deg \eta_ {c}  ^0}{s(s-1)\chi}\right),$$
where $\chi$ is a constant independent of $c$ and $\tau _k${\rm ,} and 
$$
u_{\tau _k, s}(c   )
=
\sum_{
n\in  N   c ^{-1} {D_E}    ^{-1}\atop
 \tau _k(n)<0}2^{-s}
\varrho _{\tau _k}(c, n)
Q_{s-1}\left(1-2\tau _k(n)\right).
$$
\endproclaim

\demo{{\rm 5.2.} Nonarchimedean intersections} 
In this section we want to
compute the intersection of $\eta$ and $\eta_ c  ^0$ 
at  a  place $v$ of $ E$ over a prime $\wp$ of $F$. \enddemo

5.2.1. {\it Some intersection settings}.
Let $ q    $ denote the prime of $\CO _ E$ corresponding to $v$,
and let $\CO _ q    ^{\ur}$ be the completion of the maximal unramified 
extension of $\CO_ q    $ with a uniformizer $\pi$. Let $E_q^\ur$
denote its field of fractions. Then 
 $\bar X:=\wt\CX\otimes \CO _ q     ^{\ur}$ 
can be embedded into $\bar X'=\wt\CX'\otimes \CO _ q     ^{\ur}$. For any 
$\CO _ q     ^{\ur}$-scheme $S$, the set $\Hom _{\CO _ q    ^{\ur}}(S, \bar \CX')$ 
 parametrizes
isomorphism classes of objects $[A, C, \kappa  _0]$ where $[A, C]$ is an object 
of $\CF (S)$, and $\kappa _0$ is a level structure defined by
the compact subgroup $  U  \times  J $.

Let  $x$ and $y$ be  two  integral components 
of $\eta$ and $\eta_ c  ^0$, respectively, over $E_q     ^{\ur}$. 
Then $x$ is the image of a morphism from $E_q     ^{\ur}$ to $X$
and $y$ is the image of a morphism from $F(W)$ to $X$
where $F(W)$ is the fraction field of a finite extension $W$
of $E_q    ^{\ur}$. Let us denote 
$$(x, y)_ q    =(\overline {\pi ^*(x)}/u_x, \overline{\pi ^*(y)}/u_y)/\deg \pi,$$
where $\overline {\pi ^*(x)}$ and $\overline {\pi ^*(y)}$ are the Zariski
closures of $\pi ^*(x)$ and $\pi ^*(y)$.

Assume that $  U  $ and $ J $ are maximal at places dividing
$ c $; then 
$$\pi ^*(x)=u_x\sum x_i,\qquad \pi ^*(y)=u_y\sum y_j$$
where the $x_i$ are points of $\wt X$ defined over $E_q     ^{\ur}$ and
the $y_j$ are points defined over $F(W)$. Now, we have
\begin{equation}
(x, y)_ q    =\frac 1{\deg \pi} \sum _{(i, j)}(\bar x_i, \bar y_j). \spq{5.2.1}
\end{equation}

\demo{{\rm 5.2.2.} Moduli interpretation}
The  schemes 
$\overline{\pi ^*(x)}=u_x\sum \overline{x_i}$,
and $\overline {\pi ^*(y)}=u_y\sum \overline {y_j}$ 
 represent  objects 
$[A, C, \kappa _i]$ and  $[A', C', \kappa  _j']$,
where $[A, C]$ and $[A', C']$ are objects represented by the Zariski closures
$\bar x$ and $\bar y$
of $x$ and $y$ in $\CX$, respectively, 
and $\kappa _i$ and $\kappa  _j$ are  level structures on them
for the group $  U  \cdot J $.

Now let us study the local intersection $\intn {\eta, \eta _ c ^0}_ q    $ in two cases: $\wp\, \nmid\, c$ and $\wp\!\mid\!
c$.
\enddemo

\vglue9pt

\demo{Case 1} $\wp\, \nmid\,  c $. Let $x$ and $y$ be integral components of 
$\eta $ and $\eta ^0_ c $ over $E_q    ^{\ur}$. Then all
   $x_i$ and $y_j$ are sections of $\CX$ over $\CO _ q    ^{\ur}$. 
Let $z_1, z_2, \cdots,$ be the inverse images of the reduction $z$ of $\bar x$ 
on $\bar X$. Let $[A^0, C^0, \kappa  ^0_k]$ be the corresponding 
objects.

If $\wp$ is split in $ E$ and  $\bar x_i$ and $\bar y_j$ intersect at 
some $z_k$, then
  both $\bar x_i$ and $\bar y_j$ are canonical liftings of $z_k$
with the same multiplication by $\End (x)$, so that  $x_i=y_j$.
This is impossible and now $(x, y)_ q    =0$.

If  $\wp$ is not split in $ E$ 
and  $\bar x_i$ and $\bar y_j$ intersect at some $z_k$ in the special fiber, 
then there are two
embeddings $\alpha _x:\End (x)\to \End (z)$
and $\alpha _y:\End (y)\to \End (z)$. 
With respect to $\alpha _x$ and $\alpha _y$,
 $x$ and $y$ are canonical liftings. 
 
Fix  isomorphisms
\begin{equation}
\left\{ \begin{array}{ll} 
&\End (x)\simeq \CO _ E, \\
& \End (y)\simeq \CO _ c, \\  
&\End (z)\simeq    R    (\wp),
\end{array}\right.  \spq{5.2.2}
\end{equation}
 where $   R    (\wp)$ is an order
 of type $(N(\wp),    E)$ in
the quaternion algebra  $  B   (\wp )$. 
We require that the first two isomorphisms
 satisfy the conditions of   Proposition 2.1.3.
Let $$n=n(\alpha _x, \alpha _y)\in  N(\wp)   c ^{-1} {D_E}   ^{-1}$$ 
be the link number defined as in   5.1.5. 
\enddemo

\spn{5.2.3} \proclaim{Lemma} Assume that $\ord _\wp (N)\le 1${\rm .}
 Then the intersection of $x$ and $y$ is given by 
 $(x, y)_ q=m(n)$  where
$$m(n)=\left\{ \begin{array}{ll} 
\ord _\wp (n\wp)&\hbox{if $\wp\mid D_E$}\\
\left[\ord _\wp (n\wp/N)/2\right]&\hbox{if $\wp\nmid D_E$.}
\end{array}\right. 
$$\endproclaim
 
\demo{Proof}  In this case the component of $C$ at $\wp$ is $0$.
Thus the formal deformation of the formal group
gives a formal neighborhood of  $z_i's$ in $\wh X$.
By (5.2.1), it is not 
difficult to show that 
$(\bar x, \bar y)$ equals  the maximal integer $m$ such that
\begin{itemize}
\item[1.]
$\End (x_m)$ contains the images of $\alpha _x$ and $\alpha _y$
where $x_m$ is the 
restriction of $x$ on $\CO _ q     ^{\ur}/q^m\CO _ q     ^{\ur}$;
\item[2.] $\alpha _x=\alpha _y \pmod {q^{m-1}\pi}$
in $\CO _{B(\wp)}$ as $x$ and $y$ have the same orientation at $\wp$, 
where 
$$
\pi=\left\{ \begin{array}{ll} q&\hbox{if $\wp \mid ND_F$}\\
\varpi &\hbox{otherwise},
\end{array}\right. 
$$ 
and where $\varpi$ is a uniformizer of $B(\wp)$
\end{itemize}

By Proposition 2.4.5, $\End (x_m)$ is the unique  suborder of
$R(\wp)$ of type 
$(E, \wp ^{b_m}N)$ where
$$b_m=\left\{ \begin{array}{ll} 
2m-1&\hbox{if $\wp\, \nmid\,   {D_E}   $}\\
m  &\hbox{if $\wp\hskip-1pt \mid\hskip-1pt   {D_E}   $.}
\end{array}\right. 
$$ 
On the other hand,
the algebra $\CO _{x, y}$ generated by the images of
$\alpha _x, \alpha_y$ has  discriminant $D_n:= c ^2  {D_E}   ^2 n (1-n)$.
Thus $m(n)$ is the largest number such that $\ord _\wp (D_n)\ge b_m$.
Now, the first condition is
equivalent to
\begin{equation}
m\le 
\left\{ \begin{array}{ll} \ord_\wp(n(1-n)\wp^2)&
\hbox{if $\wp\mid D_E$}\\
\left [\frac 12\ord _\wp(n(1-n)\wp/N)\right]
&\hbox{otherwise.}
\end{array}\right.  \spq{5.2.3}
\end{equation}

For the second condition we let $t$ be an element in $\CO _q$
such that 
$$\CO _q=\CO _\wp +\CO _\wp t, \qquad
 t^2\in \CO _\wp.$$
Then the second condition is equivalent to
$$\alpha _x(t)-\alpha _y(t)=0\pmod {q^{m-1}\pi}.$$
Let $\mu$ be an element in $\CO _{B(\wp)}$ such that the following
conditions hold:
$$\CO _{B(\wp)}=\CO _q+\CO _q \mu,\qquad \mu ^2\in \CO_\wp$$
$$\mu x=\bar x\mu\qquad \hbox{for all } x\in \CO _q.$$
Consider $t$ as an element in $R(\wp)$ via $\alpha _x$ then
$\alpha _y(t)$ will have the form
$$\alpha _y(t)=t(\alpha+\beta \mu), \qquad \alpha ^2-\beta \bar\beta
\mu ^2=1,$$
where $\alpha \in \CO _\wp$, $\beta \in \CO _q$.  
Now the second condition
is equivalent to the following:
\begin{equation}\alpha -1=0\pmod{q^{m-1}\pi t^{-1} },\qquad
\beta \mu  =0\pmod{q^{m-1}\pi t^{-1}}. \hskip.4in \spq{5.2.4} 
\end{equation}  
By the definition of $n$,
$$\tr (\alpha _x(t)\alpha _y(t))=2t^2-4t^2n.$$
Thus 
$$\alpha -1=-2n,\qquad \beta\bar\beta \mu^2=4n(n-1)$$
 and (5.2.4) is equivalent to
$$
n=0\pmod{q^{m-1}\pi t^{-1}},
\qquad 
n(1-n)=0\pmod{(q^{m-1}\pi t^{-1})^2},
$$ 
or equivalently,
\begin{eqnarray*}
m&\le &\ord _q (tq/\pi)
+\ord _q (\wp)\min \left\{\ord _\wp (n), 
\frac 12\ord _\wp (n(n-1))\right\}\\
&\le &\ord _q (tq/\pi)
+\ord _q (\wp) 
\frac 12\ord _\wp (n).
\end{eqnarray*}
Thus the second condition is equivalent to
\begin{equation}
m\le \left\{ \begin{array}{ll} \ord _\wp (n\wp)&\hbox{if $\wp\mid D_F$}\\
\frac 12\ord _\wp(n)&\hbox{if $\wp\mid N$}\\
\noalign{\vskip2pt}
\frac 12\ord _\wp(n\wp)&\hbox{otherwise.}
\end{array}\right.  \spn{5.2.5}
\end{equation}
The lemma follows from (5.2.3), (5.2.5),  and the fact that  
$\ord _\wp (n)>0$ if $\wp$ is unramified in $E$, as 
$n\in N(\wp) c^{-1}D_E^{-1}$.
\enddemo

\vglue9pt

Conversely if $\alpha _1$ and $\alpha _2$ are two  homomorphisms
from $\CO _ E$ and $\CO _ c $  to $\End (z)$, respectively,
which have   positive orientation, 
then by Proposition 1.5.1, 
we can find  objects $[A, C]$ and $[A', C']$ which are canonical liftings
of $[A^0, C^0]$ with respect to $\alpha _1$ and $\alpha _2$. This 
defines a component $x$ for $\eta$ and a component $y$ for $\eta _ c ^0$. Now for each $z_k$,
the level structure $\kappa  ^0_k$ can be uniquely extended to
level structure on $[A, C]$ and $[A', C']$ so that  we obtain some
sections $x_i$ and $y_j$ which intersect at $z_k$. 
It is easy to see that the number of $z_k$ is 
$\deg \pi /c_z$
where $c_z=\#[   R    (\wp)^\times/\CO _F^\times].$
Now the total intersection of $(\eta, \eta_ c ^0)_ q    $ at $z$
is given by
$$\sum _n\varrho (z, c, n)m(n),$$
where $\varrho (z, c, n)$ is the number  of $   R    (\wp)$-conjugacy classes
 of
pairs $(\phi_1, \phi_2)$  as above 
with link number $n$.

Write $\varrho _\wp( c , n)$ as the sum of $\varrho (   R    (\wp), c ,  n)$
over all nonconjugate orders $   R   (\wp)$ of $  B    (\wp)$ of type 
$(N(\wp),    E)$, where
$\varrho (   R   (\wp), c, n)$ is the  number  of $   R    (\wp)$-conjugacy
 classes of
pairs $(\alpha_1, \alpha_2)$  of  homomorphisms from 
$\CO _ E$ and $\CO _ c $ to $   R    (\wp)$ with the same orientation and
 link number $n$  in  $N(\wp)c  ^{-1} {D_E}    ^{-1}$.
As  in the archimedean case, 
$$\sum _z\varrho (z,  c , n)=2^{-s(\wp)}\varrho _\wp(  c , n)$$
where $s(\wp)$ is the number of prime factors of 
$ N(\wp)$ not dividing $   {D_E}   $.
Then we obtain 
\begin{equation}
(\eta, \eta _ c ^0)_v=u_\wp(c )-u_\wp (1), \spq{5.2.6}
\end{equation}
where 
$$
u_\wp( c  )=2^{-s(\wp)}\sum _{n\in  N   c ^{-1} {D_E}   ^{-1}}
\varrho _\wp ( c , n)
 m(n),$$ 
with $m(n)$ given by formula (5.2.6). Here $2^{-1}$ appears in the formula
because of  the symmetry between $n$ and $1-n$.

\phantom{rain}

\demo{Case 2}  $\wp| c $. Write $ c = c '\wp^s$ with $s=\ord _\wp( c )$. Then 
$\eta_ c  ^0$ can be written as
$$\eta_ c  ^0=\sum _{x'\in \eta}x'(s)+\sum _{y'\in \eta_{ c '}^0}(y'+y'(s))$$
 where $x'(s)$ and $y'(s)$ are sums  of  quasi-canonical liftings 
of the reductions of  $x'$ and $y'$ of levels up to $s$.
 If $x$ is a  section of $\eta$ then a component $\bar x_i$
of  $\overline{\pi ^*x}$ has an  intersection with a component
$\overline {x'(s)_j}$ of  
$\overline{\pi ^*(x'(s))}$ if and only if $x=x'$, and then
 $(\bar x_i, \overline{x'(s)}_j)=s$.
Similarly, a component $\bar x_i$ of $\overline{\pi ^*x}$ has an 
intersection with a component $\overline{y(s)}_i$ of $\overline {\pi^*y(s)}$ 
if and only if $\bar x_i$ has an
intersection with $\bar y_j$, and then $(\bar x_i, \overline{y(s)}_j)_ q    =s$.
 It follows that 
$$(\eta, \eta_ c  ^0)_ q    = sh_1,$$
if $\wp$ is split in $ E$,  and that
$$
(\eta, \eta_ c  ^0)_ q    =sh_2+ u_\wp( c )-u_\wp (1),$$
if $\wp$ is inert in $ E$, 
where $h_1, h_2$ are constants independent of $ c $, and where
$$u_\wp( c )=2^{-s(\wp)}\sum _{n'\in  N   c ^{-1} {D_E}   ^{-1}}
\varrho _\wp ( c ', n')
(s+m(n')).$$

In summary  we have  proved the following: \enddemo
 
\spn{5.2.4} \proclaim{Proposition} Assume that $ c $ is prime to $   N     {D_E}   ${\rm .}
\begin{itemize}
\ritem{1.} If $\varepsilon(\wp) =1${\rm ,} then  
$$(\eta, \eta_ c  ^0)_v=\ord _\wp( c )h_1$$
where $h_1$ is a constant independent of $ c ${\rm .}
\ritem{2.} If $\varepsilon(\wp)=0${\rm ,} then 
$$(\eta, \eta_ c ^0)_v=u_\wp ( c  )-u_\wp (1  ),$$ 
where $u_\wp$ is given by the formula
$$
u_\wp=2^{-s(\wp)}\sum _{n\in  N   c ^{-1} {D_E}   ^{-1}}
\varrho _\wp ( c , n)
m(n)  .$$
\ritem{3.} If $\varepsilon(\wp)=-1${\rm ,} then
$$(\eta, \eta_ c  ^0)_v=\ord _\wp( c )h_2+u_\wp( c  )-u_\wp (\CO _F)$$ 
where $h_2$ is a constant independent of $ c ${\rm ,}
 and where $u_\wp$ is given by the formula
$$
u_\wp ( c )=2^{-s(\wp)}\sum _{n'\in  N   c ^{'-1} {D_E}   ^{-1}}
\varrho _\wp ( c ', n')
(\ord _\wp (c)+m(n')),$$
where  $ c '= c \wp^{-\ord _\wp ( c )}${\rm .}
\end{itemize}

\endproclaim

\demo{{\rm 5.3.} Clifford algebras}
\enddemo

5.3.1. {\it Determining the ramification type}.
Let $\Delta$ be a quaternion algebra over $F$ with embeddings
$\phi _1$ and $\phi _2$ from $ E$ into $\Delta$.
Let  $d$ be a nonzero element in $ F^\times$ such that $\sqrt{-d}\in  E$,
and denote 
$$e_1=\phi _1(\sqrt{-d}),\qquad e_2=\phi_2(\sqrt{-d}).$$
 Let $m\in F^\times$ be defined by
$$e_1e_2+e_2e_1=2dm.$$
Then $m$ does not depend on the choice of $d$.
We want to  describe  the places at which $\Delta$ is ramified
in terms of $m$. 
 
\spn{5.3.2} \proclaim{Proposition}
Let $v$ be a place of $F${\rm .} The algebra $\Delta$ is ramified at $v$ if and
only if $\varepsilon _v(m^2-1)=-1${\rm .}
\endproclaim

\demo{Proof}
Let $\Delta^0$
denote the vector space of trace $0$ elements in $\Delta$. Then $\Delta^0$ is
ramified at a place $v$ of $F$ if and only if $\Delta^0\otimes F_v$ has
no nonzero element with square $0$. Now $\Delta^0$ is a vector space
over $F$ generated by $e_1$, $e_2$ and $e_1e_2-md$
and it is easy to check that for any $x, y, z$ in $F_v$,
$$[xe_1+ye_2+z(e_1e_2-md)]^2=-  d x^2+2mdxy-d y^2+(m^2-1)d^2z^2.$$
This form is linearly equivalent to
$$- d x^2+(m^2-1)y^2-z^2.$$
Thus  $\Delta$ is ramified at $v$ if and only if 
$(m^2-1)$ is  not a norm from $ E_v$, or equivalently,
$\varepsilon _v(m^2-1)=-1$. 
\enddemo

\demo{{\rm 5.3.3.} Counting orders}
Let $ c $ be a nonzero ideal of $\CO _ E$ prime to $ N   {D_E}   $ .
 Let $S$ denote the $\CO _F$-subalgebra in $\Delta$
generated by $\phi _1(\CO _ E)$ and $\phi _2(\CO _ c )$.
 Then $S$ is finite  over $\CO _F$
 if and only if $m+1\in 2 c ^{-1} {D_E}   ^{-1}$.
The discriminant of $S$ is
$D_S =(m^2-1) c ^2 {D_E}   ^2$.

Let $\ell   $ be an ideal of $\CO _F$ such that the following 
conditions are satisfied:
\begin{itemize}
\item[1.]  $\ord_v(\ell   )$ is even if $v$ is split in $\Delta$
 and inert in $ E$;
\item[2.] $\ord _v(\ell   )$ is odd if $v$ is ramified in $\Delta$ 
and inert in $ E$;
\item[3.] $\ord _v(\ell   )$ is $0$ if $v$ is split in $\Delta$ 
and ramified in $ E$;
\item[4.] $\ord _v(\ell   )$ is $1$ if $v$ is ramified in both $\Delta$ 
and $ E$.
\end{itemize}
In the following we want to compute the number of orders in $\Delta$ of type
$(\ell,    E)$ containing $S$.
The above conditions imply the existence of the  orders in $\Delta$ of type 
$(\ell,    E)$. Indeed,  conditions 1 and 2 imply that there is an 
ideal $\ell _E$ in $\CO _E$ with norm $\ell/D_\Delta$ where $D_\Delta$
is the product of primes in $F$ over which $\Delta$ is ramified. 
Let $\CO _\Delta$ be any maximal order of $\Delta$ containing $\phi_1(\CO _E)$. Then 
$$\phi _1(\CO _E)+\phi _1(\ell _E)\CO _\Delta $$
is  an order in $\Delta$ of type $(\ell, E).$
Let $\varrho (S)$\lb{varrho(S)} denote the number 
of  orders in $\Delta$ of discriminant 
$\ell   $ containing 
$S$. \enddemo

\spn{5.3.4} \proclaim{Proposition}
Assume $m+1\in 2 c ^{-1} {D_E}   ^{-1}${\rm .}
There is an order of discriminant $\ell   $ containing $S$  only if 
$D_S$ is divisible by $\ell   ${\rm .} If $\ell   |D_S   ${\rm ,} then 
$$\varrho (S)=r( D_S   /\ell   )
\cdot
\prod_{v|( D_S   ,    {D_E}   )\atop \varepsilon _v(m^2-1)=1}
2.$$
\endproclaim

\demo{Proof}
Since the correspondence $\CO\to \wh{\CO}$ gives a bijection between 
the set of orders of $\Delta$ and the orders of $\wh{\Delta}$, it follows that 
$\varrho (S)$ equals  the product of the numbers $\varrho _v(S)$ of 
orders on $\Delta_v$ of type $(\ell   ,    E)$ containing $S_v$ for all finite places $v$
of $F$. Fix a finite place $v=\wp$. We want to compute $\varrho _v(S)$
case by case. Let $W$ denote the ring 
$\phi _1(\CO _{ E, v})$
contained in $S$.
Recall that the discriminant of $S$ is $ D_S  $. 

If $\varepsilon (v)=-1$,
or  $v$ is ramified in $\Delta$, 
then there is a unique order  in $\Delta_v$ of 
discriminant $\wp^{\ord _v(\ell   )}$ containing $W$.
 This order contains
$S_v$ if and only if $\ord _v(\ell   )\le \ord _v( D_S  )$. 
 In other words, $\varrho _v(S)=1$ if $\ord _v ( D_S  )\le \ord _v(\ell   )$ and 
$\varrho _v(S)=0$ if $\ord _v( D_S  )<\ord _v(\ell   )$.

If $\varepsilon (\wp)=1$ then $W\simeq \CO _\wp^2$ as a ring.
 It follows that
$S_\wp$ is an Eichler order  conjugate
to an order of the form 
$$\left\{ \left( \begin{array}{cc} a&b\\ \wp ^{\ord _\wp( D_S  )}c& d \end{array}\right)
\biggm| a, b, c, d\in \CO _{\wp}\right\}.$$
This order is contained in an order of the discriminant $\wp^{\ord _\wp(\ell   )}$
if and only if $\ord _\wp(\ell   )\le \ord _\wp( D_S  )$.
If this is the case, then each  order of $\Delta$ of 
the type $(\ell   ,    E)$ 
containing this order has the  form
$$\left\{ \left( \begin{array}{cc}a&\wp ^{-k}\ell    b\\ \wp^{k}c& d \end{array}\right)\biggm| a,b, c, d, 
\in \CO _{\wp}\right\}$$
with $0\le k\le \ord_\wp ( D_S  )-\ord _\wp (\ell   )$.
 Hence $\varrho _\wp(S)=1+\ord _\wp( D_S  /\ell   )$.

If $\varepsilon (\wp)=0$ and $\wp$ is not ramified on $\Delta$,
then $S$ is generated by $e_1$ and $e_2$ such that
$$e_1e_2+e_2e_1=2md$$
and $W$ is generated by $e_1$, where $e_i^2=d$.
The correspondence $I\to \End _{\CO _\wp}(I)$ gives
a bijection  between 
the set of   maximal orders of $\Delta_\wp$ containing $W$ and
 the set of ideals in $W$  modulo  an equivalence relation:
$I_1\sim I_2$ if and only if $I_1=I_2\alpha$ for an
$\alpha \in F^\times$.  We have two maximal orders
$\End _{\CO _\wp}(W)$ and $\End _{\CO _\wp}(We_1)$
corresponding to ideals $W$ and $We_1$. 
One of them must contain $S$, say the first one. We want to prove 
that the second one also contains $S$. It suffices to show that the second
one contains $e_2$, or in other words $e_2(We_1)\subset We_1$.
First of all, as $e_2W\subset W$ and 
$$e_1e_2+e_2e_1=2md,$$ we have 
$$e_2(We_1)\subset 2mdW+We_1.$$
Secondly, as $\wp|( D_S  ,    {D_E}   )$ with $ D_S  =(m^2-1)d^2\CO _F$,
 we must have $\ord \wp(m)\ge 0$. So $e_1|md$ at the place $\wp$. It follows that 
$e_2(We_1)\subset We_1$. So we have proved that $\varrho _\wp(S)=2$.
This completes the proof of the proposition.
\enddemo

We will use Proposition 5.3.4 to compute the embeddings from $S$ into
orders of type $(\ell, E)$:

\spn{5.3.5} \proclaim{Proposition}
Let $\CO _1, \ldots, \CO _h$ be a  representing set of all conjugacy classes
of  orders in $\Delta$ of type $(\ell,    E)${\rm .} 
Then
$$\sum _{i=1}^h\#\{\phi: S\to \CO _i \pmod {\CO _i^*}\} =2^{t(\ell   )}\varrho (S),$$
where $\CO_i^\times $ acts on the set of embeddings from 
$S$ into $\CO_i$ by conjugations{\rm ,} and $t(\ell   )$ is the number
of finite places dividing $\ell   ${\rm .}
\endproclaim

\demo{Proof}
The proof is almost exactly like the modular curve case treated by
Gross,  Kohnen, and Zagier \cite{G-Z2}. We omit the details.
\enddemo

\demo{{\rm 5.4.} The final formula}
For each place $w$ of $F$, let $(\eta, \eta _ c ^0)_w$ denote the total
intersection over the places over $w$:
\begin{equation}
(\eta, \eta _ c ^0)_w =\sum _{v|w}(\eta, \eta _ c ^0)_v\log \RN (v)
\spq{5.4.1}
\end{equation}
where the $v$ are places of $ E$, and $\log \RN (v)$ is set to be $2$ if 
$v$ is a complex place.
In this section we want to compute the local intersection
\begin{equation}
(\eta, \RT (  m  )^0\eta)_w=\sum_{ c  |  m }
\varepsilon ( c )(\eta, \eta_{  m / c }^0)_w.\spq{5.4.2}
\end{equation}
\enddemo

\demo{{\rm 5.4.1.} The archimedean case\/{\rm :} Computation}
By Proposition 5.1.9, for a place $\tau _i$,
$(\eta, \RT (  m  )^0\eta )_{\tau _i}$ is equal to
\begin{equation}
2\lim _{s\to 1}
\left(U_{\tau _i, s}(  m  )-r(  m  )U_{\tau _i, s}(\CO _F)
+\frac {(\deg \eta)^2\deg \RT ^0(  m  )}{s(s-1)\chi}\right)\hskip.5in \spq{5.4.3}
\end{equation}
where $U_{\tau _i, s}(  m  )$ is
\begin{eqnarray*}&&2^{-s}\sum _{ c  |  m }
\varepsilon ( c )
\sum_{
n\in  N   c   m ^{-1} {D_E}   ^{-1}\atop
\tau _i(n)<0
}
\varrho _{\tau _i}(  m / c ,n)Q_{s-1}\left(1-
2\tau _i(n)\right)\\
&&\qquad\quad =\ 2^{-s}\sum_{
n\in  N    m ^{-1} {D_E}   ^{-1}\atop
\tau _i(n)<0
}
\sum_{
 c |  m \atop
{ c |n  m  {D_E}    N  ^{-1}}}
\varepsilon ( c )\varrho _{\tau _i}(  m / c , n)
Q_{s-1}\left(1-2\tau _i(n)\right).
\end{eqnarray*}

 \spn{5.4.2} \proclaim{Lemma} Let $n\in  N   c^{-1} {D_E}   ^{-1}$ be such that
$\tau _i(n)<0$.
Then  the following assertions hold\/{\rm :}
\begin{itemize} 
\ritem{1.}  $\varrho _{\tau _i}( c , n)\ne 0$ if and only if the following 
are satisfied\/{\rm :}\/
\begin{itemize}
\ritem{(a)}  $0<\tau _j(n)<1$ for $j\ne i${\rm ;}
\ritem{(b)} $\varepsilon _\wp(n(n-1))=1$ for any $\wp | {D_E}   ${\rm ;}
\ritem{(c)} $r(n(n-1) c ^2  N   ^{-1})\ne 0${\rm .}
\end{itemize}
\ritem{2.} Assume  conditions {\rm (a)} and {\rm (b).} Then
$$\varrho _{\tau _i}( c , n)=2^s
r(n(n-1) c ^2  N   ^{-1})\delta (n)$$
where $\delta (n)=\prod _{\wp\mid ( {D_E}, n\wp)}2${\rm .}
\end{itemize}

 \endproclaim
 
\demo{Proof}
Let $\Delta$ and $S$ be a Clifford algebra and 
an order defined as in 5.3.1 and 5.3.3 with $m=2n-1$. 
Then $S$ has discriminant 
$$D_S  =n(n-1) c ^2 {D_E}   ^2.$$

By   5.1.6,  Propositions 5.3.5 and 5.3.4, 
$\varrho_{\tau _i}(c,n)\ne 0$
is equivalent to   the
following:
\begin{itemize}
\item   $\Delta$ is isomorphic to $B(\tau _i)$, or 
equivalently by Proposition 5.3.2,
$$\varepsilon _v(n(n-1))=-1$$ if and only if
$v$ is ramified in $B(\tau _i)$.
\item  $D_S$ is divisible by $\ell=N$, or equivalently 
$n(n-1) c ^2 {D_E}   ^2 N   ^{-1}$ is an integer.
\end{itemize}

Recall that  $B(\tau _i)$ is ramified exactly at 
archimedean place $\tau _j$ ($j\ne i$) and finite places $\wp$ 
such that $\varepsilon _\wp (N)=-1$.
Thus these two conditions are equivalent to the conditions (a), (b), and (c)
because of the following:
\begin{itemize}
\item For an infinite place $\tau _j$, 
$$\varepsilon _{\tau _j}(n(n-1))<0
\Longleftrightarrow\tau_j(n(n-1))<0;$$
\item $(D_S, N)=1$ so that $B(\tau _i)$ is unramified at all places dividing $D_E$;
\item $r(n(n-1) c ^2 N   ^{-1})\ne 0$ if and only if 
$n(n-1) c ^2 {D_E}   ^2 N   ^{-1}$ is an integer and 
$$\varepsilon _\wp(n(n-1) c ^2 {D_E}   ^2 N   ^{-1})=1$$ for all finite 
place $\wp\nmid D_E$.
\end{itemize}
This proves the first assertion in the lemma.

By assertion 1, the equality in assertion 2  follows if 
$\varrho _{\tau _i}(c, n)=0$. Otherwise, 
by Propositions 5.3.5 and 5.3.4, 
 $\varrho _{\tau _i}(c, n)=2^s\varrho(S)$ and $\varrho(S)$ is given by
$$r( D_S   /\ell   )\cdot
\prod_{ v|( D_S   ,    {D_E}   )\atop \varepsilon _v(m^2-1)=1}
2.$$
Now for any $\wp\mid D_E$, condition (a) implies $\varepsilon _\wp (m^2-1)=1$,
and $\wp\mid D_S$ is equivalent to $\ord _\wp(n)\ge 0$. Thus we have 
assertion 2.
\enddemo

\spn{5.4.3} \proclaim{Lemma}
Let $  a   $ and $ b $ be two nonzero ideals{\rm .} Then
$$\sum _{ c  |(  a   ,  b )} \varepsilon ( c )r(\frac {  a    b }{ c ^2})
=r(  a   )r( b ).$$
\endproclaim

\vglue9pt

\demo{Proof} 
It is easy to reduce to the case
where $  a   =\wp^m$ and $ b =\wp^n$  both are powers of a prime
ideal in $\CO _F$. 
In this case the lemma is obvious.
\enddemo

 Applying  Lemmas 5.4.2 and 5.4.3 to formula (5.4.3), 
we, therefore, obtain  the following:

\spn{5.4.4} \proclaim{Proposition} Let $\tau _i$ be an infinite place of $F${\rm .} Then in $\CS/{\CD_N}${\rm ,}
$(\eta, \RT (  m  )^0\eta)_{\tau _i}$ is given by the limit as $s\to 1$ of
$$2
\sum_{\begin{array}{c}
{\scriptstyle n\in  N    m ^{-1} {D_E}   ^{-1},\tau _i(n)<0,}\\
{\scriptstyle  \quad 0<\tau _j(n) <1, \forall j\ne i}\\
{\scriptstyle \varepsilon _\wp (n(n-1))=1, \forall \wp| {D_E}  }
\end{array}}
\delta (n)r(n c     N  ^{-1})r((n-1)  m    ) 
Q_{s-1}\left(1-2\tau _i(n)\right).$$
\endproclaim

Here  in
$\CS/{\CD_N}${\rm ,} the limit makes sense{\rm ,}
as the term
$$\frac {(\deg \eta)^2\deg \RT ^0(  m  )}{s(s-1)\chi}$$
 is an element in ${\CD_N}$ as a function of $  m ${\rm .}

\demo{{\rm 5.4.5.} The nonarchimedean case}
Now let us treat the nonarchimedean case. Fix a prime $\wp$.
To compute $(\eta, \RT (  m  ) ^0\eta)_\wp$, there are three cases:
\enddemo

\demo{Case 1} $\varepsilon (\wp)=1$. 
  By Proposition 5.2.3,
\begin{equation}\lb{epsilon1}
(\eta, \RT (  m  )^0\eta)_\wp=2h_1j_\wp (  m  ), \spq{5.4.4} \end{equation}
where $h_1$ is a constant independent of $m$, and where
\begin{equation}\lb{epsilon1'}
j_\wp (  m  )=
\sum _{ c |  m }
\varepsilon ( c )\ord _\wp(  m / c )\log \RN (\wp)
=\frac 12r(  m  )\ord _\wp(  m  )\log \RN (\wp). \hskip.25in \spq{5.4.5}
\end{equation}
\enddemo

\demo{Case 2} $\varepsilon (\wp)=0$.  As $  m $ is prime to $ {D_E}   $, 
by 5.2.3, one has
\begin{equation}\lb{epsilon0}
(\eta, \RT (  m  ) ^0\eta)_\wp
=(U _\wp(  m  )-R (  m  )U_\wp(1  ))\log \RN (\wp). \spq{5.4.6}
\end{equation}
where 
$$U_\wp(  m  )
=2^{-s(\wp)}\sum _{n\in  N(\wp)    m ^{-1} {D_E}   ^{-1}}
\sum_{c |  m \atop  c |n  m  {D_E}    N  ^{-1}}
\varepsilon ( c )\varrho _\wp(m/c , n)
m(n)$$
where $m(n)$ is as defined by formula (5.2.3).
By the   proof of Lemma 5.4.2 we have:

\spn{5.4.6} \proclaim{Lemma} Let $n\in  N(\wp)    c ^{-1} {D_E}   ^{-1}${\rm .} Then   the following
assertions hold\/{\rm :}
\begin{itemize}
\ritem{1.} 
$\varrho _\wp(  c ,  n)\ne 0$ if and only if the following are satisfied
\begin{itemize}
\ritem{(a)}   For all infinite places $\tau _i$,
$0\le \tau _i(n )\le 1${\rm ;}
\ritem{(b)} $\varepsilon _\wp(n(n-1))=-1$, and
$\varepsilon _ q    (n(n-1))=1$ for all $ q    | {D_E}   $, $ q    \ne \wp${\rm ;}
\ritem{(c)} $r(n(n-1) c ^2  N  ^{-1})\ne 0${\rm ;}
\end{itemize}
\ritem{2.} If conditions {\rm (a)} and {\rm (b)} are satisfied then 
 $$\varrho _\wp( c , n)=2^{s(\wp)} 
\delta (n)r(n(n-1) c ^2  N  ^{-1}).$$
\end{itemize}
\endproclaim

Applying Lemmas 5.4.6 and 5.4.3, we see that $U_\wp(  m  )$ is equal to
 \begin{equation}\lb{epsilon0'}
\sum_{\begin{array}{c}
{\scriptstyle n\in  N    m ^{-1} {D_E}   ^{-1},0<n<1}\\
{\scriptstyle \varepsilon _\wp(n(n-1))=-1}\\
{\scriptstyle \varepsilon _ q    (n(n-1))=1, \forall \wp\ne  q    | {D_E}   }\end{array}}
\delta(n)r(n  m      N  ^{-1}/\wp)r((n-1) m     )
m(n) \hskip.75in \spn{5.4.7}
\end{equation} 
where the inequality $0<n<1$ means $0<\tau _i(n)<1$ for all $\tau _i$. \enddemo

\demo{Case 3} $\varepsilon (\wp)=-1$. 
 Again by Proposition 5.2.3, 
\begin{equation}\lb{epsilon-1}
  (\eta, \RT (  m  ) ^0\eta)_\wp
=h_2j_\wp (  m  )
+(U _\wp(  m  )-R (  m  )U_\wp(1))\log \RN (\wp).\hskip.4in \spq{5.4.8}
\end{equation}
Here $h_2$ is a constant independent of $  m $, and 
\begin{eqnarray}
\quad j_\wp (  m  )&=&
2\sum _{ c |  m }\varepsilon ( c )\ord _\wp(  m / c )\log \RN (\wp)
\spq{5.4.9} \\
&=&r(  m  )\ord _\wp(  m  )\log \RN (\wp)
+\ord _\wp(  m \wp)r(  m  /\wp)\log \RN (\wp),\lb{epsilon-1'} \nonumber
\end{eqnarray}
and $U_\wp(  m  )$ is equal to
$$2^{1-s(\wp)}\sum_{ c |  m }
\varepsilon ( c )\sum_
{n\in    m ^{'-1} c ' {D_E}   ^{-1} N  \wp}
\varrho _\wp(\frac {  m '}{ c '}, n)
\left[\ord _\wp(\frac   m  c )+m(n)\right],$$
where $  m '=  m '\wp^{-\ord _\wp (  m  )}$ and  $ c '= c \wp^{-\ord _\wp ( c )}$.
Changing the order of  sums and writing $ c = c '\wp^t$, we see that
$U_\wp(  m  )$ is equal to
\begin{eqnarray*}
&&\hskip-.5in 2^{1-s(\wp)}\sum _
{n\in   N    m ^{'-1} {D_E}   ^{-1}\wp}
\sum_{ c |  m '\atop  c |n   m ' {D_E}    N  ^{-1}}\varepsilon ( c ')
\varrho _\wp(\frac {  m '}{ c }, n)\\
& &\hskip.5in \cdot\sum _{t=0}^{\ord _\wp (  m  )}
(-1)^t\left[m(n)+\ord_\wp(m) -t\right].
\end{eqnarray*}
The last two sums are independent. Let us evaluate them separately. 

As in the other two case, we have the following lemma: \enddemo

\spn{5.4.7} \proclaim{Lemma} Let $c$ be an integer prime to $\wp$  and let 
$n\in   N    c ^{-1} {D_E}   ^{-1}\wp${\rm .}
Then  the 
following two assertions hold\/{\rm :}
\begin{itemize}
\ritem{1.}
 $\varrho_\wp ( c , n)\ne 0${\rm ,} if and only if the following conditions are satisfied\/{\rm :}\/
\begin{itemize}
\ritem{(a)} $0<n<1$\/{\rm ;}\/
\ritem{(b)} $\varepsilon _\ell(n(n-1))=1$, for all $\ell | {D_E}   ${\rm ;}
\ritem{(c)} $r(n(n-1)  c ^2   N  ^{-1}\wp^{-1})\ne 0${\rm .}
\end{itemize}
\ritem{2.} Moreover if conditions {\rm (a)} and {\rm (b)} are satisfied then
 $$\varrho _\wp( c , n)=
2^{s(\wp)}\delta (n)r(n(n-1)  c ^2   N  ^{-1}\wp^{-1}).$$
\end{itemize}
\endproclaim

Applying  Lemmas 5.4.7 and 5.4.3, we obtain
$$
2^{-s(\wp)}\sum_{  c |  m '\atop
 c |n  m ' {D_E}    N  ^{-1}} \varepsilon ( c )
\varrho _\wp(\frac {  m '}{ c }, n) 
= r(n  m '    N  ^{-1}/\wp)
r((n-1)  m '    ).
$$

The second sum (in Case 3) can be evaluated directly:
\begin{eqnarray*}
&&\hskip-.5in \sum _{t=0}^{\ord _\wp (  m  )}
(-1)^t\left[m(n)+\ord _\wp(m)-t\right]\\
&&\qquad\quad = \left\{ \begin{array}{ll} 
\left[m(n)+\frac 12 \ord _\wp(m)\right]&\hbox{if $\ord _\wp(  m  )$ is even,}\\
\noalign{\vskip2pt}
\frac 12\ord _\wp (  m \wp)&\hbox{if $\ord _\wp (  m  )$ is odd.}
\end{array}\right. \end{eqnarray*} 
Thus $U_\wp(m )$ is equal to 
\begin{eqnarray}\lb{epsilon-1''}
&&\hskip-.25in\sum_{\begin{array}{c}
{\scriptstyle n\in    m ^{'-1} {D_E}   ^{-1} N(\wp)  }\\
{\scriptstyle 0< n<1}\\
{\scriptstyle \varepsilon_\ell(n(n-1))=1, \forall \ell| {D_E}   }\end{array}
} 
r(n m '/N  ^{-1}\wp)r((n-1)  m '    )\delta (n)\cdot
\spq{5.4.10} \\
&&\hskip.5in\cdot
\left\{ \begin{array}{ll} 
2\left[m(n)+\frac 12 \ord _\wp(m)\right]&
\hbox{if $\ord _\wp(  m  )$ is even,}
 \\
\ord _\wp (  m \wp)&\hbox{if $\ord _\wp (  m  )$ is odd.}
\end{array}\right. \nonumber 
\end{eqnarray}

In summary we obtain the following:

\spn{5.4.8} \proclaim{Proposition} Assume that $\varepsilon (\wp)=1$ if either 
$\ord _\wp (N)>1$ or $\wp\mid 2${\rm .} Then 
 the local intersection $(\eta, \RT (  m  ) ^0 \eta)_\wp$
is given by the following formulas\/{\rm :}\/
\begin{itemize}

\ritem{1.} If $\varepsilon (\wp)=1$ then
$(\eta, \RT (  m  ) ^0\eta)_\wp \pmod {\CD_N}$
is equal to 
$$h_1r(m)\ord _\wp (m)\log \RN (\wp)$$
where $h_1$ is a constant independent of $m$ and $\wp${\rm .}

\ritem{2.} If $\varepsilon (\wp)=0$ then
$(\eta, \RT (  m  ) ^0 \eta)_\wp$ is equal to
$$(U_\wp(m)-R(m)U_\wp(1))\log \RN (\wp)$$
where $U_\wp (m)$
is equal to
$$
\sum_{\begin{array}{c}
{\scriptstyle n\in   m ^{-1} {D_E}   ^{-1} N(\wp)  , 0<n <1}\\
{\scriptstyle \varepsilon _\wp(n(n-1))=-1}\\
{\scriptstyle \varepsilon_ q    (n(n-1))=1, \forall \wp\ne  q    | {D_E}}\end{array}   
}
\delta(n)r(n  m /N \wp )r((n-1)  m       )\ord _\wp(n \wp). 
$$

\ritem{3.} If $\varepsilon (\wp)=-1$ then  $(\eta, \RT (  m  ) ^0\eta)_\wp $
is equal to
\begin{eqnarray*}
&&h_2r(m)\ord _\wp(m)\log \RN (\wp)+h_2\ord _\wp(m\wp)r(m/\wp)\log \RN (\wp)\\
&&\phantom{h_2r(m)\ord _\wp(m)\log \RN (\wp)} +(U_\wp(m)-U_\wp(1))\log \RN (\wp)
\end{eqnarray*}
where $h_2$ is a constant independent of $m, \wp$, and $U_\wp(m)$ is equal to
\begin{eqnarray*}
&&\hskip-.75in \sum_{\begin{array}{c}
{\scriptstyle n'\in    m ^{'-1} {D_E}   ^{-1} N(\wp)  } \\
{\scriptstyle 0< n'<1}\\
{\scriptstyle \varepsilon_\ell(n'(n'-1))=1, \forall \ell| {D_E}   }\end{array}
} 
r(n'  m '    N  ^{-1}/\wp)r((n'-1)  m '    )\delta (n') 
 \\
&&\qquad\quad \cdot
\left\{ \begin{array}{ll} 
2\left[\frac 12 \ord _\wp(n'\wp m)\right]&\hbox{if $\ord _\wp(  m  )$ is even,}
 \\
\ord _\wp (  m \wp)&\hbox{if $\ord _\wp (  m  )$ is odd.}
\end{array}\right. \nonumber\\
&
\end{eqnarray*}
\end{itemize}
\endproclaim

\demo{Proof}
All these  statements follow  from formulas (5.4.4)--(5.4.10) and  Lemma 5.2.3,
with the fact that in the case
$\varepsilon (\wp)=-1$, the term for an $n'$ in (5.4.10) has nonzero 
contribution only if $\ord _\wp(n'N(\wp))$ is even. Thus 
$$m(n')=\left[\frac 12 \ord _\wp(n'\wp)\right]$$ 
even when $\ord _\wp (N)$ is odd.
\enddemo

\section{Derivatives of $L$-series}

In this section, we will compute $L_ E'(f,s)$ 
using the method of Gross and Zagier shown in \cite {G-Z1}.
 We will start with a formula  which expresses
$L_ E(f, 1/2\break +s)$ as an inner product of $f$ with a 
nonholomorphic form $ \Phi   _s(z)$. Then we compute the 
Fourier expansion for $ \Phi  _s(z)$ and get a formula
for some multiple 
$\wt \Phi $ of  $\frac \partial {\partial s}\big|_{s=1/2}
\wt \Phi  _s(z).$
Finally  the holomorphic projection
of $\wt \Phi $ gives a holomorphic form $ \Phi $ with
Fourier coefficients given explicitly.

\demo{{\rm 6.1.} The Rankin-Selberg method}  Let $f$ be a new form for $   K   _0( N  )$ and let $ E$ be an imaginary
quadratic  extension of $F$ as before. Then the base change $L$-function of $f$ to $ E$ is defined to
be $L_ E(s, f)=L(s, f)L(s, \varepsilon, f)$ where $\varepsilon$ is the character
attached to $ E/F$. See   Section 3.4 for definitions.
 For any nonzero ideal $  m $ let $r(  m )$ denote the 
number of integral ideals in $\CO_ E$ with the norm $  m $. Using Proposition 3.1.4, 
one  shows that
\begin{equation}
L_ E(f, s)=L^ N  (2s-1,\varepsilon )\sum _{  m \in \BN _F}
a(  m )r(  m  )\RN (  m  )^{-s} \spq{6.1.1}
\end{equation}
where $L^ N  (s, \varepsilon)$ denotes the series
$$\sum_{m \in \BN_F\atop (  m ,  N   {D_E}   )=1}
\varepsilon (  m )\RN(  m  )^{-s},$$
where $ {D_E}   $ is the conductor of $\varepsilon$.
In other words, $L_ E(f, s)$ is essentially the Rankin-Selberg convolution 
of $L(s,f)$ with $\zeta _ E(s)$. We want to express this convolution
 as an inner product of $f$ with a modular form.  We will construct such a
form using the Eisenstein series defined in 
Section 3.5.
\enddemo

\demo{{\rm 6.1.1.}  $L_ E(f, s)$ as an inner product of $f$ with some 
other form} We need to define a Haar measure on $Z(\BA _F)
\backslash G(\BA _F)$. Let $dk=\otimes dk_v$ be the Haar measure on $   K  _0(1)$ with volume 
$1$ on each component. Recall  that  $dx=\otimes dx_v$ is defined in 3.1.1
to be a measure on $\BA _F$ 
such that $dx_v$ is the usual Euclidean measure if $v$ is infinite,
and that $\CO_v$ has volume $1$ if $v$ is finite. Also recall  that
$d^\times x=\otimes d^\times x_v$ is defined in the proof of 3.4.2 to
be a Haar measure on $\BA _F^\times$ such that 
$d^\times x_v=|dx_v/x_v|$ if $v$ is infinite, and that 
$\CO _v^\times$ has volume $1$ if $v$ is finite.
Now $dg$ on $G (\BA_F)$ is defined by the formula
$$\int _{Z(\BA _F)\backslash G(\BA _F)}f(g)dg=
\int _{\BA _F^\times}\int _{\BA _F}\int _   K  
f\left( \left( \begin{array}{cc}y&x\\ 0&1 \end{array}\right)k\right)dkdx\frac {d^\times y}
{|y|}.$$
For any two function $f$ and $g$ on $Z(\BA _F)G(F)\backslash G (\BA _F)$
the integral $f\bar g$ (if it is absolutely convergent) is denoted as 
$(f, g)$.

Let $E_s$ \lb{E_s}
be the Eisenstein series defined in 3.5.1 with
$\chi=\varepsilon$ associated to the extension $E/F$. 
Let  $E_{s,  N  }$\lb{E_{s,N}}
 be an Eisenstein series defined by the formula
\begin{equation}
E_{s,  N  }(g)=E_s\left(g  \left( \begin{array}{cc}1&0\\ 0&\pi _ N 
 \end{array}\right)\right) \spq{6.1.2}
\end{equation}
 where $\pi _ N  $ is an idele with components $1$ at places not
dividing $ N  $ and such that $\pi _ N  $ generates $\wh N  $. Then  
$E_{s, N  }$ is a form of level
$K_0( {D_E}    N ) $. 
\enddemo

\spn{6.1.2} \proclaim{Proposition} 
Let $f$ be a new cusp form of weight $2${\rm ,}
 for $K_0(N)$ a
 trivial central character{\rm .}
 Then 
$$(f, E_{1/2}E_{s, N})=A(s)
L_E(s+1/2, f)$$
where 
$$A(s)=\left[\frac {\Gamma (s+1/2)}{2^{2s}\pi ^{s-1/2}}\right]^g
d_F^{s+1/2}d_N^sd_E^{-1/2}\mu (ND_E)$$
where $\mu (ND_E)$ is the volume of $K_0(ND_E)${\rm .}
\endproclaim

\demo{{P}roof}
For each factor $e$ of $N$, let 
$E_s^e$ be the Eisenstein series defined in the same way as $E_s$ 
in 3.5.1 with
factor $L(2s, \varepsilon)$ replaced by $L^e(2s, \varepsilon)$ and with 
$H_s$ replaced by the following $H_s^e$:
$$H_{s}^ e (g)=\left\{ \begin{array}{ll} \left|\frac ad\right|^{s}
\varepsilon (akr(\theta))&\hbox{if $k\in    K   _0( {D_E}    e )$}\\
0&\hbox{otherwise.}
\end{array}\right. $$
For $\Re (s)>1$, 
$E_s^{ e }(g)$
is absolutely convergent and defines a (nonholomorphic) form 
for $   K   _0( {D_E}    e )$ of (parallel) weight $1$ with   character
$\varepsilon$. If $e=\CO _F$, $E^e_s=E_s$.
\enddemo

\spn{6.1.3} \proclaim{Lemma} Let $f$ be a cusp form of weight $2$ for $   K   _0( N   {D_E}   )$ with 
trivial central character{\rm .} Let $\theta$ be the theta series with Fourier coefficients
$r(m)$ defined as in {\rm 3.4.5.} Then
$$(f, \theta E^{ N  }_{\bar s})=\mu (ND_E)d_F^{s+1}
\left[\frac { \Gamma   (s+1/2)}{(4\pi
)^{s+1/2}}\right] ^g
L_ E (s+1/2, f).$$
\endproclaim

\demo{Proof}
By definition of $E ^ N  _s$, up to a factor $L(2s, \varepsilon)$,
$(f, \theta E_{\bar s}^ N  )$
is given by
 \begin{eqnarray*}
 &&\hskip-.5in 
 \int _{Z(\BA _F)B(F)\backslash G(\BA _F)}f \overline {\theta H^ N 
_{\bar s}} dg  \\
\noalign{\vskip4pt} &&\hskip.25in
 = \int _{\BA _{F, +}^\times/F_+}\int _{\BA _F/F}\int _   K   
(f\overline {\theta H^ N  _{\bar s}})\left( \left( \begin{array}{cc}
y&x\\ 0&1 \end{array}\right)k\right)dkdx\frac {d^\times y}{|y|},
 \end{eqnarray*}  
where $\BA ^\times _{F, +}$ denote 
 ideles with  positive components at
the infinite places.
By definition of $H_s^ N  $, the inner integral over $K$ is 
$$\mu (ND_E)
(f\bar\theta)\left( \left( \begin{array}{cc}y&x\\ 0&1 \end{array}\right)\right)
|y|^{s}\varepsilon (y).$$
Using Fourier expansions of $f$ and $\theta$ in (3.1.3) and Proposition 3.1.2, 
\begin{eqnarray*}
&&\hskip-.5in \int _{\BA _F/F} 
(f\bar\theta)\left( \left( \begin{array}{cc}y&x\\ 0&1 \end{array}\right)\right)dx\\
\noalign{\vskip4pt}
&&\hskip.25in =\ d_F^{1/2}|y|^{3/2}\varepsilon (y)\sum _{\alpha>0}
a(\alpha y_f  {D_F}   )r(\alpha y_f {D_F}   )\psi (2\alpha y_\infty i),
\end{eqnarray*} 
where $a(  m )$ are the Fourier coefficients of $f$ as defined in Proposition 3.1.2. 
Combining these,  $(f, \theta E^ N  _{\bar s})$ up to a factor $L(2s, \varepsilon)$, 
is equal to 
$$\mu (ND_E)d_F^{1/2}
\int _{\BA _{F, +}^\times }|y|^{s+1/2}\sum _{\alpha>0}
a( y_f {D_F}    )r( y_f {D_F}    )\psi (2y_\infty i)d^\times y.
$$
This integral is the product of the integrals over infinite ideles
$\prod _{v|\infty}F _{v, +}$
and over finite ideles $\wh F^\times$.
The integral over infinite ideles gives
$$\left[\frac { \Gamma   (s+1/2)}{(4\pi )^{s+1/2}}\right]^g$$
while the integral over the finite ideles gives
$$d_F^{s+1/2}\sum _  m  \frac {a(  m )r(  m )}{\RN (  m )^{s+1/2}}.$$
The lemma follows from (6.1.1).
\enddemo

\vglue9pt 
The next lemma gives a comparison between $E^e_s$ and $E_{s, n}$.

\vglue9pt 

\spn{6.1.4} \proclaim{Lemma}
  ${\displaystyle E_{s,  N  }=d_N^{s}\sum _{  a   | N  }\frac {\varepsilon (  a   )}
{\RN (  a   )^{2s}}E_s^{ N  /  a   }.}$ 
\endproclaim

\vglue9pt 

\demo{Proof} Let $H_{s,  N  }$ be the function on $G(\BA _F)$ defined in the same
way as $E_{s,  N  }$:
$$H_{s,  N  }(g)=H_s\left(g  \left( \begin{array}{cc}1&0\\ \noalign{\vskip3pt} 0&\pi _ N   \end{array}\right)\right).$$
It suffices to prove the corresponding statement for $H_{s,  N  }$
on $   K  _0(1)$:
$$H_{s,  N  }=d_N^{s}\sum _{  a   | N  }\frac {\varepsilon (  a   )}
{\RN (  a   )^{2s}}\frac {L^{ N  /  a   }(2s, \varepsilon)}
{L(2s, \varepsilon)}H_s^{ N  /  a   }.$$ 
We do   this by testing their values on elements $k=\left( \begin{array}{cc}
a&b\\ \noalign{\vskip3pt} c&d \end{array}\right)$ of $   K  _v$ for finite places 
$v$ not dividing $ {D_E}   $.
Now we have the decomposition:
$$k \left( \begin{array}{cc}1&0\\ \noalign{\vskip3pt} 0&\pi _ N   \end{array}\right)
= \left( \begin{array}{cc} \frac 1d (ad-bc)&b\pi _N\\ \noalign{\vskip3pt} 0&d\pi _N \end{array}\right)
 \left( \begin{array}{cc}1&0\\ \noalign{\vskip3pt} \frac c{d\pi _ N  }&1 \end{array}\right)$$
if $ N  |c$, and
$$k \left( \begin{array}{cc}1&0\\ \noalign{\vskip3pt} 0&\pi _ N   \end{array}\right)
= \left( \begin{array}{cc}\frac {\pi _ N  } c(ad-bc)&a\\ \noalign{\vskip3pt} 0&c \end{array}\right)
 \left( \begin{array}{cc}0&-1\\ \noalign{\vskip3pt} 1&\frac {d\pi _ N  }c \end{array}\right)$$
if $c|\frac {\pi _ N  }{\pi _v}$. It follows that
$$H_{s,  N  }\left(k\right)
=\left\{ \begin{array}{ll}  |\pi _ N  |^{-s}&\hbox{if $ N  |c$}\\ \noalign{\vskip3pt}
\varepsilon_v \left(\frac {\pi _ N  } c\right)\left|\frac {\pi _ N  } {c^2}
\right|^{s}
&\hbox{if $c|\frac {\pi _ N  } {\pi_v}$.}
\end{array}\right. 
$$
Let $\wp_v$ be the prime ideal of $\CO _F$ corresponding to $v$.
By definition of $H_s^ e $, 
$$H_s^{ e }(k)-H_s^{\wp e }(k)=\left\{ \begin{array}{ll} 1&
\hbox{if $\ord _v( N  )=\ord _v ( e )$}\\ \noalign{\vskip3pt}
0&\hbox{otherwise.}\end{array}\right. $$
It follows that
\begin{eqnarray*}
&&\hskip-.65in \RN ( N  _v)^{s}H_{s,  N  }(k) \\ \noalign{\vskip4pt}
&&= H_s^ N   (k) +\sum _{
1\le i\le \ord _v( N  )}
\varepsilon (\wp _v^i)\RN(\wp_v^i)^{2s}
(H_s^{ N  /\wp^i}(k)-H_s^{ N  / \wp^{i-1}}(k))\\
\noalign{\vskip4pt}
&&=
\sum _{0\le i\le \ord _v( N  )}
\frac {L^{ N  / \wp^i}(2s, \varepsilon)}{L(2s, \varepsilon)}\varepsilon 
( \wp ^i)\RN ( \wp^i )^{2s}
H_s^{ N  / \wp ^i}(k).\\
\end{eqnarray*}
\vglue-40pt
\enddemo

\phantom{dream}
\vglue9pt

Now go back to the proof of Proposition 6.1.2. 
By Proposition 3.5.4,
$$E_{1/2}=\frac{(2\pi )^g}{\sqrt{d_Fd_E}}
\theta.
$$ 
By the     two lemmas above, 
the proof of the proposition is reduced 
to showing that 
$$(f, E_{1/2}E_s^e)=0$$
for any factor $e\ne N$ of $N$. 

Let $\tr _{D_E}$ be the trace 
operator from the space of cusp forms of level $K_0(ND_E)$ to
$K_0(N)$: 
for any form
$\phi$ of level $ {D_E}    N  $, 
\begin{equation}\lb{tr_{D_E}}
(\tr _ {D_E}   \phi)(g)
=\sum _{\gamma \in    K  _0( N  ) /   K   _0( N   {D_E}   )}
\phi (g\gamma ). \spq{6.1.3}
\end{equation}
Then 
$$(f, E_{1/2}E_s^e)=[K_0(N):K_0(ND_E)]^{-1}
(f, \tr _{D_E}(E_{1/2}E_s^e).$$
As representatives of $   K   _0( N   )/   K   _0( N   {D_E}   )$ will also serve as representatives
for $   K   _0( e )/   K   _0( e  {D_E}   )$ for any $ e | N  $,
$\tr _{D_E}(E_{1/2}E_s^e)$ is a form of level $K_0(e)$.
Thus it is orthogonal to $f$ as $f$ is a newform.\hfill\qed

\demo{{\rm 6.1.5} Definition of $\Phi _s$}
We define 
\begin{equation}\lb{Phi_s}
 \Phi  _s(g)=\tr _ {D_E}   \left( \frac1{2^{\#\{v: v| {D_E}   \}}}
\sum _{ e | {D_E}   }\RN (e)^{s-1/2} \Phi ^e_s
\right). \spq{6.1.4}
\end{equation}
Here, for  $ e  $  a divisor of
 $ {D_E} $, 
\begin{equation}\lb{Phi^e_s}
\Phi _s^e(g)=(E_{1/2}E_{s,  N  })
(g\gamma _ e ) 
\end{equation}
 where $\gamma _ e  $ is an element of $\GL_2
 (\BA_F)$ which has components $1$
at places not dividing $ e $ and at a place $v$ dividing $ e $ it has 
the component $ \left( \begin{array}{cc} 0&-1\\ \pi _v&0 \end{array}\right)$ where
$\pi_v$ is normalized such that $\varepsilon (\pi _v)=1$, 
and where $\tr _{D_E}$ is as defined in (6.1.3).
It is easy to check that
 $ \Phi  _s$ is a form of weight $1$
for $K_0( N ) $ with trivial character.
\enddemo

\spn{6.1.6} \proclaim{{C}orollary}  Let $f$ be a newform of weight $2$ for $ K_0(N)  $ with
 trivial central character{\rm .} Then 
$$(f,  \Phi  _{\bar s})=B(s)L_ E (f, 1/2+s)$$
where
$$B(s)=\left[\frac { \Gamma   (s+1/2)}{2(4\pi )^{s-1/2}}\right]
^gd_F^{s+1/2}d_N^{s}d_E^{-1/2}\mu (N).
$$
\endproclaim 

\demo{Proof} Fix any factor $ e  $ of $ {D_E}   $.
Write $f^ e $ for the form $g\to f(g\gamma_ e ^{-1})$. Then again
$f^e$ is a form of level $K_0(N)$.
By definition,
$(f,  \Phi^e  _{\bar s})$ is equal to
$$[K_0(N):K_0(ND_E)](f^ e , E_{1/2}E_{\bar s, N}).
$$
By Proposition 6.1.2, this is 
$$A(s)[K_0(N):K_0(ND_E)]L_E(s+1/2, f^e).$$

As 
$$f^ e  \left( \left( \begin{array}{cc}y&x\\ 0&1 \end{array}\right)\right)
=f \left( \left( \begin{array}{cc}y/\pi _ e  &x\\ 0&1 \end{array}\right)\right),$$
if $f$ has the  Fourier coefficients $a(  m  )$ then $f^ e  $ will have the 
Fourier coefficients 
$$a(f^ e ,   m  )=\RN ( e  )a(  m / e ).$$
It follows that
$$L_ E(f^ e , s)=\RN ( e )^{1-s}L_ E(f, s).$$
The proposition follows.
\enddemo

\demo{{\rm 6.2.} Fourier coefficients}
\enddemo

6.2.1. {\it Strategy}.
In this section we want to compute the Fourier coefficients 
$c_s(\alpha, y)$  ($\alpha\in F$)  for $ \Phi  _s$ defined by
(6.1.4), where 
\begin{equation}
c_s (\alpha, y)=d_F^{-1/2}\int _{\BA _F/F}
\Phi_s \left( \left( \begin{array}{cc} y&x\\ 0&1 \end{array}\right)\right)
\psi (-\alpha x)dx. \spq{6.2.1}
\end{equation}
It suffices to compute $c_s(\alpha, y)$ for $\alpha =0$ or $1$, as
for $\alpha \in F^\times$, 
$$c_s(\alpha, y)=c_s(1, \alpha y).$$
We proceed with the following steps:
\begin{itemize}
\item[1.] Compute the Fourier coefficients $c^e_s(\alpha, y)$
 for $E_s(g\gamma _ e )$.
This will give the Fourier coefficients for $\Phi^e_s(g)$ defined by
(6.1.5).
\item[2.]
For a factor $ g $ of $ {D_E}   $ and 
 an integral adele $a$ which is $0$ at places not dividing $ g $, let 
 $\gamma _{ g , a}$ denote the  element  in $\GL_2 (\BA)$ which has the component $1$ at
places $v$  not dividing $ g $, otherwise it is given by
$ \left( \begin{array}{cc} a_v&-1\\ 1&0 \end{array}\right).$
Then  
$$\left\{\gamma _{ g , a}\bigg|\quad  g | {D_E},   \quad a\pmod  g 
\right\}$$ 
forms
 a set of representatives for $   K   _0( N  )/   K   _0( {D_E}    N  )$.
Compute the Fourier coefficients $c_s^{e, g}(\alpha, y)$ for
\begin{equation}
 \Phi _s^{e, g}:=\sum _{a\pmod  g }  \Phi  ^ e  (g\gamma _{ g , a}). \spq{6.2.2}
\end{equation}
\item[3.] Compute the Fourier coefficients of $\Phi _s$ using the following 
expression: 
\begin{equation}
\Phi _s= 2^{-\#\{v: v| {D_E}   \}}
\sum _{ e, g | {D_E}   }\RN (e)^{s-1/2} \Phi _s^{e,g}. \spq{6.2.3}
\end{equation}
\end{itemize}

\spn{6.2.2} \proclaim{Lemma} 
The Fourier coefficient $c_s^ e  (\alpha, y)$ of $E_s(g\gamma_ e)$
 is zero if $\alpha y {D_F}   $ is nonintegral{\rm .}
Otherwise{\rm ,} it is given by the following expressions\/{\rm :}
$$
c_s^e(0, y)=\left\{ \begin{array}{ll} \varepsilon (y)L(2s, \varepsilon) |y|^s
&\hbox{if $e=1$}\\
\frac {(-1)^g}{
d_F^{1/2}d_E^{s}}
V_s(0)^gL(2s-1, \varepsilon) |y|^{1-s}&\hbox{if $e=D_F$}\\
0&\hbox{otherwise},
\end{array}\right. 
$$
$$c_s^e(1, y)=\frac {(-1)^g}{\sqrt{d_Fd_E}}
\RN (e)^{1/2-s}\sigma _s( y)|y|^{1-s}
\prod _{v\mid {D_E}   / e }|y_v\pi _v|^{2s-1}\varepsilon (-y_v)
\kappa (v),$$
where 
$$\sigma _s(y)= \prod_{v\nmid D_E\atop v\nmid\infty}
\frac {1-\varepsilon (y_v\delta _v \pi _v)
|y_v\delta _v \pi _v|^{2s-1}}{1-\varepsilon (\pi _v )|\pi _v|^{2s-1}}
\cdot\prod _{v\mid \infty}V_s(y_v),$$
and for $y\in \BR${\rm ,} 
$$V_s(y)=\int _{-\infty}^\infty\frac {e^{-2\pi iyx}}{(x^2+1)^{s-1/2}
(x+i)}dx.$$
\endproclaim

\demo{Proof}  The case $e=1$ was covered in Proposition 3.5.2.
So we assume $e\ne 1$ in the following.
Again by Bruhat decomposition, $c^e(\alpha, y)$ is equal to
\begin{eqnarray*}
&&\hskip-.25in L(2s, \varepsilon)d_F^{-1/2}\int _{\BA _F/F}
H_s \left( \left( \begin{array}{cc} y&x\\ 0&1 \end{array}\right)\gamma_ e  \right)
\psi (-\alpha x)dx\\
&&\hskip.25in +\ L(2s, \varepsilon)d_F^{-1/2}\int _{\BA _F}
H_s\left(w  \left( \begin{array}{cc} y&x\\ 0&1 \end{array}\right)\gamma _ e \right)
\psi (-\alpha x)dx.
\end{eqnarray*}
Let $v$ be a place dividing $ e $, then $ \left( \begin{array}{cc} y&x\\ 0&1 \end{array}\right)
\gamma _ e $ has the component
$$ \left( \begin{array}{cc}y_v&x_v\\ 0&1 \end{array}\right) \left( \begin{array}{cc}0&-1
\\ \pi _v&0 \end{array}\right)
= \left( \begin{array}{cc}y_v&x_v\pi _v\\ 0&\pi _v \end{array}\right) \left( \begin{array}{cc}0&-1
\\ 1&0 \end{array}\right).$$
It follows that
$$H_s \left( \left( \begin{array}{cc} y&x\\ 0&1 \end{array}\right)\gamma_ e  \right)=0$$
and that $c^ e (\alpha, y)$ is given by
\begin{eqnarray*}
&&\hskip-.25in L(2s, \varepsilon)d_F^{-1/2}\int _{\BA _F}
H_s\left(w  \left( \begin{array}{cc} y&x\\ 0&1 \end{array}\right)\gamma_ e  \right)
\psi (-\alpha x)dx\\
&& \hskip.25in =\  L(2s, \varepsilon)d_F^{-1/2}
|y|^{1-s}\prod _{v\nmid  e }V_s(\alpha_vy_v)\cdot
\prod _{v\mid e }V_s'(\alpha _vy_v),
\end{eqnarray*}
where $V_s$ is as defined   in the
proof of Proposition 3.5.2, and $V_s'$ is given by
the formula
$$V_s' (y)=\int _{F_v}H_s\left(w \left( \begin{array}{cc}
1&x\\ 0&1 \end{array}\right) \left( \begin{array}{cc}0&-1\\ \pi _v&0 \end{array}\right)\right)\psi (-xy)dx.$$
Now
$$w \left( \begin{array}{cc}1&x\\ 0&1 \end{array}\right) \left( \begin{array}{cc}
0&-1\\ \pi _v&0 \end{array}\right)= \left( \begin{array}{cc}
-\pi _v&0\\ x\pi _v&-1 \end{array}\right)$$
has the decomposition
$$
 \left( \begin{array}{cc}-\pi _v&0\\ 0&1 \end{array}\right)
 \left( \begin{array}{cc}1&0\\ x\pi _v &-1 \end{array}\right),$$
if $\ord _v(x)\ge 0$, and
$$ \left( \begin{array}{cc}-x^{-1}&\pi _v\\ 0&-x\pi _v \end{array}\right)
 \left( \begin{array}{cc}0&1\\ -1&x^{-1}\pi _v ^{-1} \end{array}\right),
$$
if $\ord _v(x)<0$. It follows that
$$H_s\left( \left( \begin{array}{cc}
-\pi _v&0\\ x\pi _v&-1 \end{array}\right)\right)
=\left\{ \begin{array}{ll} |\pi _v|^{s}\varepsilon _v(-1)&\hbox{if $\ord _v(x)\ge 0$,}\\
0&\hbox{otherwise,}
\end{array}\right. 
$$
and
$$V_s'(y)=\left\{ \begin{array}{ll} |\pi _v|^{s}\varepsilon _v(-1)&
\hbox{if $\ord _v(y)\ge 0$,}\\
0&\hbox{otherwise.}
\end{array}\right. 
$$

Now the lemma follows easily from this formula
and the formulas for $V_s$ derived in the proof 
of Proposition 3.5.2.
\enddemo

\spn{6.2.3} \proclaim{Lemma}The Fourier coefficient $c_s^{ e ,  g }(\alpha, y)$ 
($\alpha =0, 1$)
of $\Phi  ^{ e, g }_s$ is zero if $\alpha y {D_F}   $ is nonintegral{\rm .} Otherwise{\rm ,} 
it is given by
$$
c_s^{ e ,  g }(\alpha, y)=\RN ( g )
\varepsilon ( N  )
\sum _{n\in F}c_{1/2}^{ e * g } (\alpha-n, \pi _ g  y)
c_s^{ e * g } (n, \pi _ g  y/\pi _ N  )
$$
where $ e * g $ denotes $ e  g /( e ,  g )^2${\rm .}
\endproclaim

\demo{Proof}
By definition,
$$
 \Phi  ^{ e ,  g }_s\left( \left( \begin{array}{cc} y&x\\ 0&1 \end{array}\right)\right)
= \sum _{a\pmod  g }
 \Phi  ^{ e }_s\left( \left( \begin{array}{cc} y&x\\ 0&1 \end{array}\right)
\gamma _a\right).$$
From the following decomposition at any place $v$ dividing $ g $,
$$ \left( \begin{array}{cc} a_v&-1\\ 1&0 \end{array}\right)=\frac 1{\pi _v}
 \left( \begin{array}{cc} \pi _v&a_v\\ 0&1 \end{array}\right)
 \left( \begin{array}{cc} 0&-1\\ \pi_v &0 \end{array}\right),$$
we see that
$$ \left( \begin{array}{cc}y&x\\ 0&1 \end{array}\right)\gamma _a=
\frac 1{\pi _ g } \left( \begin{array}{cc} \pi _ g  y&ay+x\\ 0&1 \end{array}\right)
\gamma _ g .$$
It follows that
$$
 \Phi  ^{ e ,  g }_s\left( \left( \begin{array}{cc} y&x\\ 0&1 \end{array}\right)\right)
= \sum _{a\pmod  g }
 \Phi  ^{ e * g }_s\left( \left( \begin{array}{cc} \pi _ g  y&ay+x\\ 0&1 \end{array}\right)
\right).$$
Thus 
the Fourier coefficients of $\Phi _s ^{ e, g}$ are given by
$$a_s^{ e * g } (\alpha, \pi _ g  y)\sum _{a\pmod  g }
\psi (\alpha ya),$$
or in other words, $c_s^{ e ,  g }(\alpha, y)$ is nonzero only 
if $\alpha y$ is
integral at places dividing $g$. In this case it is given by 
$$c_s^{ e ,  g }(\alpha, y)=\RN ( g )a_s^{ e * g } (\alpha, \pi _ g  y)
$$
where $a_s^e(\alpha, y)$ is the Fourier coefficient of $\Phi _s^e$
which can be expressed as 
\medbreak
\phantom{bizarre} \hfill ${\displaystyle a_s^e(\alpha, y)
=\varepsilon (N)
\sum _{n\in F}c_{1/2}^e(\alpha -n, y)c_{s, N}^e(n, y/\pi_N).}$
\enddemo

\spn{6.2.4} \proclaim{Proposition}
The Fourier coefficient $ c_s(\alpha, y)$ ($\alpha=0, 1$)
of $ \Phi _s$ is nonzero only if
$\alpha y {D_F}   $ is integral{\rm .} In this case it is given by
$$
c_s(\alpha, y)=\frac {\varepsilon ( N  )d_F^{1-s}}{d_Ed_F}
\sum_{n\in F
} a_s^n(\alpha, y)
$$
where $a_s^n(\alpha, y)$ is as given by the following formulas if it is nonzero{\rm .}
\begin{itemize}
\ritem{1.}  If $n\ne 0$ and $n\ne \alpha${\rm ,} then $a_s^n(\alpha, y)\ne 0$ only
if $ny {D_E}     {D_F}    N  ^{-1}$ is integral{\rm .} In this case
$a_s^n(\alpha, y)$ is equal to 
\begin{eqnarray*}
&&\hskip-.5in |y|^{3/2-s}\delta (ny)
 \prod _{v\mid {D_E}    }\frac 
{1+|n_vy_v\pi _v|^{2s-1}\varepsilon _v((n-\alpha)n)}2\\
&&\qquad\quad\cdot
\sigma _{1/2}((\alpha -n)y)\sigma _s(ny/\pi _ N  ),
\end{eqnarray*}
where for an idele $y$,
$\delta (y)=2^{\#\{ v\mid {D_E}   ,\ord _v(y)\ge 0\}}$.
\ritem{2.}  If $n=0${\rm ,} $\alpha =1${\rm ,} then $a_s^n(\alpha, y)$ is equal to
\begin{eqnarray*}
&&\hskip-.25in \sigma _{1/2}(y)
d_F^{1/2}d_E^{1/2}d_N^{2s-1}\varepsilon ( N  )
i^g L(2s, \varepsilon)|y|^{1/2+s}\\ \noalign{\vskip4pt}
&&\qquad\quad +\ \sigma _{1/2}(y)V_s(0)^gL(2s-1, \varepsilon )|y|^{3/2-s}.
\end{eqnarray*}
\ritem{3.} If $n=\alpha=0${\rm ,} then $a_s^n(\alpha, y)$ is equal to
\begin{eqnarray*}
&&\hskip-.25in \varepsilon (N)d_Fd_Ed_N^{2s-1}L(1, \varepsilon)L(2s,\varepsilon)|y|^{1/2+s}\\
\noalign{\vskip4pt}
&&\qquad\quad +\ V_{1/2}(0)V_s(0)L(0, \varepsilon)L(2s-1, \varepsilon)|y|^{3/2-s}.
\end{eqnarray*}
\ritem{4.} If $n=1${\rm ,} $\alpha =1${\rm ,} then $a_s^n(\alpha, y)$  is equal to
\begin{eqnarray*}
&&\hskip-.25in \left[
d_F^{1/2}d_E^{1/2} L(1, \varepsilon)i^g|\pi _ {D_E}    y_ {D_E}   | ^{2s-1}
\varepsilon ^ {D_E}    (y)+L(0, \varepsilon )V_{1/2}(0)^g\right]\\
\noalign{\vskip4pt}
&&\qquad \cdot \sigma _s(y/\pi _ N  )|y|^{3/2-s},
\end{eqnarray*}
where $\varepsilon ^ {D_E}   (y)$ denotes $\varepsilon (y)\prod _{v\mid
 {D_E}   }\varepsilon _v(y)${\rm .}
\end{itemize}
\endproclaim
 
\demo{{P}roof} By formula (6.2.3), $c_s(\alpha, y)$ is equal to 
$$\frac {1}{\delta (1)}
\sum _{ e ,  g }\RN ( e )^{s-1/2}c_s^{ e ,  g }(\alpha, y)
=\frac {\varepsilon ( N  )d_N^{1-s}}{d_Ed_F}
\sum _{n\in F}a_s^n(\alpha, y)
$$
where $a_s^n(\alpha, y)$ is equal to
\begin{equation}
\frac {d_Fd_Ed_N^{s-1}}{\delta (1)}
\sum _{ e ,  g }\RN ( g )\RN ( e  )^{s-1/2}
c_{1/2}^{ e * g }(\alpha-n, \pi _ g  y)
c_s^{ e * g }(n, \pi _ g  y/\pi _ N  ). \spq{6.2.4}
\end{equation}
\enddemo

\vglue9pt 

\demo{Case 1} $n\ne 0, \alpha$. 
If $a_s^n(\alpha, y)\ne 0$, one must have 
$$\ord _v(ny\pi _v)\ge 0\quad\hbox{for each}\quad  v\mid {D_E}   .$$
Assume this is the case and let $ g  _0$ be the factor of $ {D_E}   $ consisting of places
$v$ such that
$|n_vy_v\pi _v|=1$.
Then
$$c_{1/2}^{ e * g }(\alpha-n, \pi _g y)
c_s^{ e * g }(n, \pi _ g  y/\pi _ N  )\ne 0$$
only if $ g  _0| g $ and in this case by Lemma 6.2.2,  it equals 
\begin{eqnarray*}
&&\hskip-.25in 
\frac {|\pi _ N   |^{s-1}}
{d_Ed_F}|\pi _ g  y|^{3/2-s}\sigma _{1/2}((\alpha-n)y)
\sigma _s(ny/\pi _ N  )\RN ( g * e )^{1/2-s}\\
&&\hskip.25in \cdot \prod _{v| {D_E}   /( g * e )}
|n_vy_v\pi _{ g , v}\pi _v|^{2s-1}\varepsilon _v((\alpha-n)n).
\end{eqnarray*}
It follows that $a_s^n(\alpha, y)$ is equal to
\begin{eqnarray}
&&\quad  \frac {|y|^{3/2-s}}{\delta (1)}
\sigma _{1/2}((\alpha-n)y)
\sigma _s(ny/\pi _ N  )\spq{6.2.5} \\
&&\hskip.55in \cdot \sum_{
 g _0| g | {D_E}   \atop
 e | {D_E}   }
\frac {\RN ( g  e )^{s-1/2}}{\RN ( g * e )^{s-1/2}}
\prod _{v| {D_E}   /( g * e )}
|n_vy_v\pi _{ g , v}\pi _v|^{2s-1}\varepsilon _v((\alpha-n)n).\nn
\end{eqnarray}
Notice that
$$
\frac {\RN ( g  e )}{\RN ( e * g )}
\prod _{v| {D_E}   /( g * e )}
|\pi _{ g , v}|^{2}=1.
$$
Now the last sum is
$$
\sum_{
 g _0| g | {D_E}   \atop
 e | {D_E}   }
\prod _{v| {D_E}   /( g * e )}
|n_vy_v\pi _v|^{2s-1}\varepsilon _v((\alpha-n)n).$$
When $ e $ is exchanged for  $( {D_E}   / e )* g $, this sum equals
$$
\delta (1/ g _0)
\prod _{v| {D_E}   }
\left[1+|n_vy_v\pi _v|^{2s-1}\varepsilon _v((\alpha-n)n)\right].$$
Bringing this to (6.2.4), we obtain the formula for $a_s^n(\alpha, y)$ 
in the proposition.
\enddemo

\vglue9pt 

\demo{Case 2} $n=0$, $\alpha=1$.   In this case,
$a_s^0(1, y)$ is equal to
\begin{eqnarray*}
&&\hskip-.25in \frac {d_Fd_Ed_N^{s-1}}{\delta (1)}
\sum _{ g }\RN ( g )^{1/2+s}c_{1/2}^1(1, \pi _ g  y)
c_s^1(0, \pi _ g  y/\pi _ N  )\\  
&&\hskip.25in +\ \frac {d_Ed_Fd_N^{s-1}}{\delta (1)}
\sum _{ g }d_E^{s-1/2}\RN ( g )^{3/2-s}
c_{1/2}^ {D_E}    (1, \pi _ g  y)c_s^ {D_E}    (0, \pi _ g  y/\pi _ N  ).
\end{eqnarray*}
The formula in the proposition 
follows, as $c_{1/2}^1(1, \pi _ g  y)c_s^1(0, \pi _ g  y/\pi _ N  )$
is equal to
$$\frac {i^g\varepsilon ( N  ) d_N^{s}}
{d_E^{1/2}d_F^{1/2}}\sigma _{1/2}(y)L(2s, \varepsilon) |y\pi _ g |^{1/2+s}
\varepsilon^ {D_E}    (y)$$
and $c_{1/2}^ {D_E}   (1, \pi _ g  y)c_s^ {D_E}    (0, \pi _ g  y/\pi _ N  )$
is equal to
$$
\frac {d_N^{1-s}}{d_Ed_F}d_E^{1/2-s}V_s(0)^g\sigma _{1/2}(y) 
L(2s-1, \varepsilon )|y\pi _ g |^{3/2-s}
$$
where $\varepsilon ^ {D_E}   (y)$ denotes $\varepsilon (y)\prod _{v\mid
 {D_E}   }\varepsilon _v(y)$ which equals $1$ if $\sigma _{1/2}(y)\ne 0$.
\enddemo

\demo{Case 3} $n=\alpha=0$.  In this case, 
$a_s^0(1, y)$ is equal to
\begin{eqnarray*}
&&\hskip-.25in \frac {d_Fd_Ed_N^{s-1}}{\delta (1)}
\sum _{ g }\RN ( g )^{1/2+s}c_{1/2}^1(0, \pi _ g  y)
c_s^1(0, \pi _ g  y/\pi _ N  )\\
\noalign{\vskip4pt}
&&\hskip.25in +\ \frac {d_Ed_Fd_N^{s-1}}{\delta (1)}
\sum _{ g }d_E^{s-1/2}\RN ( g )^{3/2-s}
c_{1/2}^ {D_E}    (0, \pi _ g  y)c_s^ {D_E}    (0, \pi _ g  y/\pi _ N  ).
\end{eqnarray*} 
The formula in the proposition follows, as 
$c_{1/2}^1(0, \pi _ g  y)c_s^1(0, \pi _ g  y/\pi _ N  )$
is equal to 
$$\varepsilon (N)d_N^sL(1, \varepsilon)L(2s, \varepsilon)|\pi _gy|^{1/2+s}$$
while 
$c_{1/2}^ {D_E}   (0, \pi _ g  y)c_s^ {D_E}    (0, \pi _ g  y/\pi _ N  )$
is equal to 
$$\frac {d_N^{1-s}}{d_Fd_E}d_E^{-s}
V_{1/2}(0)V_s(0)L(0, \varepsilon)L(2s-1, \varepsilon)|y|^{3/2-s}.$$
\enddemo

\vglue9pt

\demo{Case 4} $n=\alpha=1$.
This case can be treated  similarly.
We  have  $a_s^1(1, y)$ equal to
\begin{eqnarray*}
&&\hskip-.25in \frac{d_Fd_Ed_N^{s-1}}{\delta (1)}
\sum _{ g } \RN ( g )^{1/2+s}
c_{1/2}^1(0, \pi _ g  y)c_s^1(1, \pi _ g  y/\pi _ N  )\\
\noalign{\vskip4pt}
&&\hskip.25in +\  \frac{d_Fd_Ed_N^{s-1}}{\delta (1)}
\sum _ g  d_E^{s-1/2}\RN ( g )^{3/2-s}c_{1/2}^ {D_E}    (0, \pi _ g  g)
c_s^ {D_E}    (1, \pi _ g  y/\pi _ N  ),
\end{eqnarray*}
where $c_{1/2}^1(0, \pi _ g  y)c_s^1(1, \pi _ g  y/\pi _ N   )$ is equal to
$$
\frac {d_N^{1/2-s}i^gL(1, \varepsilon )}
{d_F^{1/2}d_E^{1/2}}
\sigma _s (y/ \pi _ N  )|y|^{1-s}|\pi_ {D_E}    y_ {D_E}   | ^{2s-1}\varepsilon ^ {D_E}    (y) 
|\pi _ g |^{1/2+s},
$$
and $c_{1/2}^ {D_E}    (0, \pi _ g  y)c_s^ {D_E}    (1, \pi _ g  y/\pi _ N  )$
is equal to 
\medbreak
\phantom{long} \hfill ${\displaystyle \frac {L(0, \varepsilon )d_N^{1-s}}{d_Ed_F}
d_E^{1/2-s} V_{1/2}(0)^g\sigma _s(y/\pi _ N  )|\pi _ g  y|^{3/2-s}.
}$\hfill\qed
\enddemo

\demo{{\rm 6.3.} Functional equations and derivatives}

\spn{6.3.1} \proclaim{Proposition}
The Fourier coefficient $c_s^e(\alpha, y)$ of $E(g\gamma _e)$
 has the following
functional equation\/{\rm :}
\begin{eqnarray*}
\wt c{\hskip1pt}_s^{e} (\alpha, y)&:=&(d_Fd_Ed_e)^{s-1/2}
\left[\Gamma   (s+1/2)\pi^{1/2-s}\right]^{g} c_s^e(\alpha, y)\\
&=&i^g\varepsilon (y)\prod _{v\mid D_E/e}\varepsilon _v(-1)
\wt c_{1-s}^{D_E/e}(\alpha, y).
\end{eqnarray*}
\endproclaim

\demo{Proof} If $\alpha=1$, then by Lemma 6.2.2, up to a factor independent
of $s$ and $e$, $\wt c{\hskip1pt}_s^e(1, y)$ is given by
\begin{eqnarray*}
&&\hskip-.5in \prod_{
v \mid\hskip-2.85pt/\hskip 0.5pt D_E\atop
v \mid\hskip-2.85pt/\hskip 0.5pt \infty
}
|y_v\delta _v|^{1/2-s}\frac {1-\varepsilon (y_v\delta _v \pi _v)|y_v\delta _v\pi _v|^{2s-1}}
{1-\varepsilon (\pi _v )|\pi _v|^{2s-1}}\\
\noalign{\vskip6pt} 
&&\hskip.25in \cdot
 \prod _{v\mid \infty}\left[\Gamma   (s+1/2)\pi^{1/2-s}\right]
|y_v|^{1/2-s}V_s(y_v)\\
\noalign{\vskip6pt} 
&&\hskip.25in \cdot \prod _{v\mid e}|y_v\pi _v|^{1/2-s}\\
\noalign{\vskip6pt} 
&&\hskip.25in \cdot \prod _{v\mid D_E/e}|y_v\pi _v|^{s-1/2}
\varepsilon (-y_v)\kappa (v).
\end{eqnarray*}
By Proposition 3.3 in \cite[p.\ 278]{G-Z}, with $k=1$,
 $V_s(t)$ ($t\ne 0$) has a functional equation
$$V_s^*(t):=(\pi|t|)^{1/2-s}\Gamma (s+1/2)V_s(t)
=\sgn (t)V^*_{1-s}(t).$$
(Notice that  $V_s$ as defined in Lemma 6.2.2 
is $V_{s+1/2}$ in \cite{G-Z}.)
Thus the functional equation in the lemma follows from the local equations
and the equality
$$\prod _{v\mid D_E}\kappa (v)=i^g\varepsilon (D_F).$$

Now we want to treat the case where $\alpha=0$. By Lemma 6.2.2, we 
need only consider the case where $e=1$ or $e=D_E$. 
Recall that $L(s, \varepsilon)$ has a  functional equation:
\begin{eqnarray*}
L^*(s, \varepsilon )&:=&\left(d_Ed_F \right)^{s/2}
 \left[\Gamma  (s/2+1/2)\pi ^{1/2-s/2}\right]^gL(s, \varepsilon)\\
\noalign{\vskip6pt}
&=&L^*(1-s, \varepsilon).
\end{eqnarray*}
(This can be proved by using  functional equations for both $\zeta _E$
and $\zeta _F$ and the identity
 $\zeta _E(s)=L(s, \varepsilon)\zeta _F(s).$)
Thus $\wt c_s^1(0, y)$ is equal to
\begin{eqnarray}
&&\quad   (d_Fd_E)^{s-1/2}
\left[\Gamma   (s+1/2)\pi^{1/2-s}\right]^{g}L(2s, \varepsilon )
\varepsilon (y)|y|^s\spq{6.3.1} \\
\noalign{\vskip6pt}
&&\hskip.25in =\
 (d_Fd_E)^{-1/2}L^*(2s, \varepsilon)\varepsilon(y)|y|^s
=
(d_Fd_E)^{-1/2}L^*(1-2s, \varepsilon)\varepsilon (y)|y|^s. \nonumber
\end{eqnarray}

On the other hand, by Proposition (3.3) in \cite[p.\ 277]{G-Z},
\begin{eqnarray*}
V_s(0)&=&-\pi i 2^{2-2s}\Gamma (2s-1)/\Gamma (s-1/2)\Gamma (s+1/2)\\
\noalign{\vskip4pt}
&=&-i\pi ^{1/2}\Gamma (s)/\Gamma (s+1/2).
\end{eqnarray*} 
Thus $\wt c_s^{D_E}(0, y)$ is equal to
\begin{eqnarray*}
&&\hskip-.75in (d_Fd_E^2)^{s-1/2}
\left[\Gamma   (s+1/2)\pi^{1/2-s}\right]^{g}
\frac {(-1)^g}{d_F^{1/2}d_E^s}V_s(0)^gL(2s-1, \varepsilon)|y|^{1-s}\\ \noalign{\vskip3pt}
&&\hskip.5in =(d_Fd_E)^{s-1}\left[i\pi ^{1-s}\Gamma (s)\right]^gL(2s-1, \varepsilon)
|y|^{1-s}\\ \noalign{\vskip3pt}
&&\hskip.5in = i^g(d_Fd_E)^{-1/2}L^*(2s-1, \varepsilon).
\end{eqnarray*}
Combining this with (6.3.1), we have shown 
$$\wt c_{1-s}^{D_E}(0, y)=i^g\varepsilon (y)\wt c_{s}^1(0, y).$$
Thus, the lemma is proved in this case.
\enddemo

\vglue9pt 

\spn{6.3.2} \proclaim{{C}orollary}  The function $\Phi _s$ satisfies the following functional 
equation\/{\rm :}
\begin{eqnarray*}
\Phi _s^*&:=&(d_Fd_E)^{s-1/2}\left[\Gamma   (s+1/2)\pi^{1/2-s}\right]^{g}
\Phi _s\\ \noalign{\vskip3pt}
&=&(-1)^g\varepsilon (N)\Phi _{1-s}^*.
\end{eqnarray*}
\endproclaim 
\vglue9pt 

\demo{Proof}
We need only prove the following functional equation for $\Phi_s^e$
defined in (6.1.4):
\begin{eqnarray*}
\Phi _s^{e,*}&:=&(d_Fd_E\RN (e))^{s-1/2}\left[\Gamma   (s+1/2)\pi^{1/2-s}\right]^{g}\Phi _s^e\\ \noalign{\vskip3pt}
&=&(-1)^g\varepsilon (N)\Phi _{1-s}^{e, *}.
\end{eqnarray*}
As both sides are modular forms for $K_0(ND_E)$ with  trivial
 character, it suffices to check the functional equation for 
its Fourier coefficients. But this follows from Lemma 6.3.1, as
the Fourier coefficients of $\Phi _s^e$ are  expressed in the form
\medbreak
\phantom{safe}\hfill
${\displaystyle a_s^e(\alpha, y)
=\varepsilon (N)
\sum _{n\in F}c_{1/2}^e(\alpha -n, y)c_{s, N}^e(n, y/\pi_N).}$
\enddemo 
 
\vglue9pt

\spn{6.3.3} \proclaim{Theorem} 
The function $L_E(s, f)$ satisfies the following functional equation\/{\rm :}
\begin{eqnarray*}
L^*(s, f)&:=&(d_F^2d_Ed_N)^{s-1}\left[\Gamma (s)(2\pi )^{1-s}\right]^{2g}
L_E(s, f)\\ \noalign{\vskip3pt}
&=&(-1)^g\varepsilon (N)L_E^*(2-s, f).
\end{eqnarray*}
\endproclaim

\demo{Proof} This follows from Corollaries 6.3.2 and 6.1.6.
\enddemo

\spn{6.3.4} \proclaim{Proposition}
Assume that $\varepsilon (N)=(-1)^{g-1}${\rm .} Then $\Phi _{1/2}=0$
and the Fourier coefficient $ c'(\alpha, y)$ $(\alpha=0, 1)$
of 
$$\Phi':=\frac{\partial}{\partial s} \Phi _s\big |_{s=1/2}$$
\lb{Phi'}
 is nonzero only if
$\alpha y {D_F}   $ is integral. In this case it is given by
$$
c'(\alpha, y)=\frac {\varepsilon ( N  )d_N^{1/2}}{d_Ed_F}
\sum_{ n\in F
} b^n(\alpha, y)
$$
where $b^n(\alpha, y)$ is give by the following formulas 
if it is nonzero\/{\rm :}
\begin{itemize}
\ritem{1.}  If $n\ne 0$ and $n\ne \alpha$, then $b^n(\alpha, y)\ne 0$ only
if $ny {D_E}     {D_F}    N  ^{-1}$ is integral and 
$(\alpha -n)y$ is totally positive{\rm .} In this case
$b^n(\alpha, y)$ is equal to 
$$(-4\pi ^2)^g|y|\psi (i\alpha y_\infty)
\delta (ny)r((\alpha -n)yD_F)
\sum _v b^n_v(\alpha, y)$$
where $v$ runs through all places of $F${\rm ,}
$\delta (y)=2^{\#\{ v\mid {D_E}   ,\ord _v(y)\ge 0\}}${\rm ,}
 and $b^n_v$ is given by the following
formulas\/{\rm :}
\begin{itemize}
\ritem{(a)} If $v$ is an infinite place, then $b_v^n(\alpha, y)$ is nonzero
only if 
\begin{itemize}
\item[$\bullet$] $ny$ is negative at place $v$ and positive at other infinite places{\rm ,}
\item[$\bullet$]  $\varepsilon _\wp ((n-\alpha)n)=1$
for every place $\wp$ of $D_E${\rm .}
\end{itemize}
In this case{\rm ,} $b_v^n(\alpha, y)$ is equal to 
$$
r(nyD_F/N )q(4\pi |ny_v|)$$
where 
$$q(t)=\int _1^\infty e^{-xt}\frac {dx}x,
\qquad (t>0).$$
\ritem{(b)} If $v$ is a finite place ramified in $E${\rm ,} then 
$b_v^n(\alpha, y)$ is nonzero only if 
\begin{itemize}
\item[$\bullet$] $ny$ is totally
positive{\rm ,}
\item[$\bullet$] $\varepsilon_v((n-\alpha)n)=-1$ but $\varepsilon _\wp ((n-\alpha)n)=1$
for every place $\wp$ of $D_E${\rm .}
\end{itemize}
In this case{\rm ,} $b_v^n(\alpha , y)$ is equal to 
$$
-r(nyD_F/N )\log |n_vy _v\pi _v/\pi _{N, v}| .
$$
\ritem{(c)} If $v$ is a finite place inert in $E${\rm ,} then
$b_v^n(\alpha, y)\ne 0$ only if 
\begin{itemize}
\item[$\bullet$] $ny$ is totally
positive{\rm ,}
\item[$\bullet$] $\varepsilon_\wp((n-\alpha)n)=1$ 
for every place $\wp$ of $D_E${\rm ,}
\item[$\bullet$] $\ord _v(nyD_F/N)$ is odd.
\end{itemize}
In this case, $b_v^n(\alpha , y)$ is equal to 
$$
-r(nyD_F/(N\wp_v) )\log |n_vy_vD_{F, v}\wp _v|
$$
where $\wp_v$ is the prime corresponding to $v${\rm .}
\ritem{(d)} If $v$ is a finite place split in $E${\rm ,} then
$b_v^n(\alpha, y)=0${\rm .}
\end{itemize}

\ritem{2.}  If $n=0${\rm ,} $\alpha =1${\rm ,} then $b^n(\alpha, y)$ is nonzero only
if $y$ is totally positive{\rm .} In this case{\rm ,} it is equal to
$$r(yD_F)\psi (iy_\infty)|y|(c_1+c_2\log |y|)
$$
where $c_1$ and $c_2$ are constants{\rm .}
\ritem{3.} If $n=\alpha=0${\rm ,} then $b^n(\alpha, y)$ is equal to
$$|y|(c_3+c_4\log |y|)$$
where $c_3$ and $c_4$ are constants{\rm .}
\ritem{4.} If $n=1${\rm ,} $\alpha =1${\rm ,} then $b^n(\alpha, y)$  is equal to
$$|y|\psi (iy_\infty)
\left[c_5r(yD_F/N)\log |\pi _{D_E}y_{D_E}|+c_6r'(yD_F/N)
\right]$$
where $c_5, c_6$ are constants{\rm ,} and for a nonzero integral ideal $m${\rm ,}
$$r'(m)=\sum _{n\mid m}\varepsilon (n)\log \RN (n).$$
\end{itemize}
\endproclaim 
 
\demo{{P}roof} The vanishing of $\Phi _s$ at $1/2$ follows from
Theorem 6.3.3. To compute the Fourier coefficients of 
$\Phi '$ we use formulas in Proposition 6.2.4 with 
$$b^n(\alpha, y)=\frac {\partial}{\partial s}a_s^n(\alpha, y)\big|_{s=1/2}.$$
Notice that $a_s^n$ vanishes at $s=1/2$. 
This can be checked
from its expression, or from formula (6.2.4) and Proposition
6.3.1.
\enddemo

\demo{The case where $n\ne 0$, $n\ne \alpha$}
In this case,  $a_s^n$ is a product of 
$$y^{3/2-s}\delta (ny)\sigma _{1/2}((\alpha -n)y)
\cdot \prod _v\sigma _{s, v}^n(\alpha, y/\pi _N)$$
where $v$ runs through   the set of all places of $F$, and
$$\sigma _{s, v}^n(\alpha, y)
=\left\{ \begin{array}{ll}  \frac {1+\varepsilon ((n-\alpha)n)|n_vy_v\pi _v|^{2s-1}}2
&\hbox{if $v|D_E$}\\
V_s(n_vy_v)&\hbox{if $v\mid\infty$}\\
\frac {1-\varepsilon (n_vy_v \delta _v\pi _v)|n_vy_v\delta _v\pi _v|^{2s-1}}
{1-\varepsilon (\pi _v)|\pi _v|^{2s-1}}
&\hbox{otherwise.}
\end{array}\right. 
$$
If $b^n(\alpha, y)\ne 0$, then $\sigma _{1/2}((\alpha -n)y)\ne 0$
and one and only one factor of $\sigma _{s, v}$ vanishes at $s=\frac 12$.
If this is the case, then 
$$\sigma_{1/2}((\alpha -n)y)=r((\alpha -n)yD_F)\prod _{v\mid\infty}
V_{1/2}((\alpha _v-n_v)y_v).$$
 By Proposition 3.3 in \cite[p.\ 278]{G-Z},
we know that $V_{1/2}(t)$ for $t\in \BR$, is given by
$$V_{1/2}(t)=\left\{ \begin{array}{ll} 0&\hbox{if $t<0$}\\
-2\pi i e^{-2\pi t}&\hbox{if $t>0$}.
\end{array}\right. 
$$
Thus, if $b^n(\alpha, y)\ne 0$, then $(\alpha-n)y$ must be totally positive
and $r((\alpha -n)y)\ne 0$. In this case  $b^n$ has the expression in the 
proposition
with $b_v^n(\alpha, y)$ equal to 
$$\psi (-iny_\infty)\frac{\partial}{\partial s} \sigma _{v, s}^n(\alpha,
 y)\big|_{s=1/2}
\cdot \prod _{w\ne v}\sigma _{w, 1/2}^n(\alpha, y).$$ 
The proposition in this case can be checked case by case.
Notice that when $v$ is archimedean, we have used the identity 
$$\frac {\partial}{\partial s}V_s(t)\big|_{s=1/2}
=-2\pi iq(t)e^{-2\pi it} \qquad (t<0)$$
in Proposition 3.3 in \cite {G-Z}. \enddemo 

\demo{The cases where  $n=0$ or 
$n=\alpha$} These cases can be verified from the expressions in
Proposition 6.2.4.\hfill\qed
\enddemo

The same proof will also give the following:
\spn{6.3.5} \proclaim{Proposition} Assume that $\varepsilon (N)=(-1)^g$. Then 
the Fourier coefficient $ c_{1/2}(\alpha, y)$ ($\alpha=0, 1$)
of $ \Phi _{1/2}$ is nonzero only if
$\alpha y {D_F}   $ is integral{\rm .} In this case it is given by
$$
c_{1/2}(\alpha, y)=\frac {\varepsilon ( N  )d_F^{1/2}}{d_Ed_F}
\sum_{n\in F
} a_{1/2}^n(\alpha, y)
$$
where $a_{1/2}^n(\alpha, y)$ 
is given by the following formulas if it is nonzero\/{\rm :}\/
\begin{itemize}
\ritem{1.}  If $n\ne 0$ and $n\ne \alpha${\rm ,} then $a_{1/2}^n(\alpha, y)\ne 0$ only
if the following holds\/{\rm :}
\begin{itemize}
\ritem{(a)} $ny {D_E}     {D_F}    N  ^{-1}$ is integral{\rm ,}
\ritem{(b)} both $(\alpha-n)y$ and $ny$ are totally positive{\rm ,}
\ritem{(c)} $\varepsilon _v(n(n-1))=1$ for all $v\mid D_E${\rm .}
\end{itemize} 
In this case
$a_{1/2}^n(\alpha, y)$ is equal to 
$$(-4\pi ^2)^g|y|\delta (ny)r((\alpha-n)yD_F)r(nyD_F/N)
\psi (i\alpha y_\infty)
$$
where for an idele $y${\rm ,}
$\delta (y)=2^{\#\{ v\mid {D_E}   ,\ord _v(y)\ge 0\}}${\rm .}
\ritem{2.}  If $n=0$, $\alpha =1$, then $a_{1/2}^n(\alpha, y)$ is equal to
$$c_1r(yD_F)|y|\psi (iy_\infty)$$
where $c_1$ is a constant independent of $y${\rm .}
\ritem{3.} If $n=\alpha=0${\rm ,} then $a_{1/2}^n(\alpha, y)$ is equal to
$$c_2|y|$$
where $c_2$ is a constant independent of $y${\rm .}
\ritem{4.} If $n=1${\rm ,} $\alpha =1${\rm ,} then $a_{1/2}^n(\alpha, y)$  is equal to
$$c_3r(yD_F/N)|y|\psi (iy_\infty)$$
where $c_3$ is a constant independent of $y${\rm .}
\end{itemize}

\endproclaim

\demo{{\rm 6.4.} Holomorphic projections}
\enddemo

6.4.1. {\it Asymptotic formula for $\Phi'$
near cusps}.
Assume that $\varepsilon ( N   )=(-1)^{g-1}$. 
The form $\Phi' $ defined in Proposition 6.3.4 is not holomorphic. 
We want to find a {\em holomorphic
projection $ \Phi $}\lb{Phi} which is a holomorphic cusp form for
$   K   _0(   N    )$ and has the property that for any newform $f$, 
$f$ has the same scalar product with both  $\wt \Phi $ and  $ \Phi $.

As in the case $F=\BQ$ treated by Gross and Zagier \cite[Chap.\  IV,
 \S6]{G-Z}, $\Phi' $ 
satisfies the growth condition 
\begin{equation}
\Phi'  \left(  \left( \begin{array}{cc} y&x\\ 0&1 \end{array}\right) g\right)
=a_g|y|\log |y|+b_g|y|+
O_g(|y|^{1-\varepsilon}) \spq{6.4.1}
\end{equation}
 for each $g\in \GL(\BA _F)$.  By Proposition 6.3.4, 
the asymptotic formula (6.4.1) is true for $g=1$, and  
$$
\Phi'  \left(  \left( \begin{array}{cc} y&x\\ 0&1 \end{array}\right) g\right)
=c_3|y|\log |y|+c_4|y|+
O(|y|^{1-\varepsilon})
$$
where $c_3$ and $c_4$ are constants independent of $g$
as defined in Proposition 6.3.4.

For any $e|N$, let $g_e$ denote an element 
in $\GL _2(\BA _F)$ which has  components $1$ at places 
not dividing $N$ and has components  
 $ \left( \begin{array}{cc}1&0\\
\pi _v^{\ord _v(e)}&1 \end{array}\right)$
at each place $v$ dividing $N$. 
Using the same method, we may compute the Fourier coefficients 
of $\Phi'(gg_e)$ and will obtain the  formula
similar to  (6.4.1). As $\GL _2(\BA _F)$ is a union of the form $B(\BA)
g_eK_0(N)$, formula (6.4.1) is true for every $g\in \GL _2(\BA _F)$.

We have to subtract some Eisenstein 
series to make $a_g=b_g=0$ for every $g\in \GL _2(\BA _F)$.  
Let  $E_{2, s}(g)$ be the Eisenstein series constructed in
 3.5.1 with $k=2$, $\chi=1$, and $N=1$.
Then $E_{2, s}(g)$ is perpendicular to all holomorphic cusp forms.
 Proposition 3.5.2 implies the following asymptotic formula:
\begin{equation}
E_{2, s}\left( \left( \begin{array}{cc}y&x\\ 0&1 \end{array}\right)
\right)
=\zeta _F(2s)|y|^s+c(s)|y|^{1-s}+O(1) \spq{6.4.2}
\end{equation}
as $y\to \infty$, where $c(s)$ and $O(1)$ are holomorphic in $s$ near
$s=1$. Define by continuation
$$E(g)=E_{2, s}(g)|_{s=1}
\quad\hbox{and}\quad 
F(g)=\frac {\partial}{\partial s}\bigg|_{s=1}
E_{2, s}(g).$$

For each $e\mid N$, let $h_e$ be an element of $\GL _2(\BA _F)$
which has components $1$ at places not dividing $N$ and has 
components $ \left( \begin{array}{cc}1&0\\ 0&\pi _v^{\ord _v(e)}
 \end{array}\right)$ at places $v$ dividing $N$.

\spn{6.4.2} \proclaim{Lemma}
There are some pairs of numbers $(\alpha _e, \beta _e)$ $(e\mid N)$
such that the form
$$\wt \Phi(g):=\Phi'  (g)-\sum _ e 
\left[\alpha _eF\left(gh_e\right)
+\beta_e E\left(gh_e\right)
\right]
$$
has the same holomorphic projection as $\Phi' ${\rm ,}  and $ \wt\Phi$
satisfies the bound
$$
\wt \Phi \left( \left( \begin{array}{cc}y&x\\ 0&1 \end{array}\right)
g\right)=O(|y|^{1-\varepsilon})$$
 as $y\to\infty${\rm ,} for every $g\in \GL _2(\BA _F)${\rm .}
\endproclaim

\demo{Proof} We need only find $(\alpha _e, \beta _e)$'s so that the
equation in the lemma holds for $g=g_f$'s,
as $\GL _2(\BA _F)$ is a union of $B(\BA _F)\gamma _eK_0(N)$.
Now $ \left( \begin{array}{cc}y&x\\ 0&1 \end{array}\right)g_fh_e$ has the 
decomposition at a place $v$ of $N$:
$$
\left\{ \begin{array}{ll} 
 \left( \begin{array}{cc}y_v&x_v\pi _v^{m_v}\\
0&\pi _v^{m_v} \end{array}\right)\cdot
 \left( \begin{array}{cc}1&0\\ \pi _v^{n_v-m_v}&1 \end{array}\right)
&\hbox{if $n_v\ge m_v$}\\
\noalign{\vskip4pt}
 \left( \begin{array}{cc}y_v\pi _v^{m_v-n_v}& y_v+x_v\pi _v^{n_v}\\
0&\pi _v^{n_v} \end{array}\right)
\cdot  \left( \begin{array}{cc}0& -1\\ 1&\pi _v^{m_v-n_v} \end{array}\right)
&\hbox{if $m_v> n_v$}
\end{array}\right. 
$$
where $m_v=\ord _v(e)$, $n_v=\ord _v(f)$.
Thus  (6.4.2)  implies
$$
E_{2, s}\left( \left( \begin{array}{cc}y&x\\ 0&1 \end{array}\right)
g_fh_e\right)
=C_N(e,f)^s\zeta _F(2s)|y|^s+c(s)
C_N(e,f)^{1-s}|y|^{1-s}+O(1)
$$
where $C_N(e, f)=\RN (e,f)^2/\RN (e)$.
It follows that,
\begin{eqnarray*}
E\left( \left( \begin{array}{cc} y&x\\ 0&1 \end{array}\right)g_fh_e\right)
&=&\zeta _F(2)C_N(e,f)|y|+O(1),\\
\noalign{\vskip4pt}
F\left( \left( \begin{array}{cc} y&x\\ 0&1 \end{array}\right)g_fh_e\right)
&=&\zeta _F(2)C_N(e, f)|y|\log |y|+O(\log |y|).
\end{eqnarray*}

Now the asymptotic formula  in the lemma is equivalent to
\begin{eqnarray*}
\noalign{\vskip4pt}
\sum _{e\mid N}\alpha _e C_N(e, f)&=&\zeta _F(2)^{-1}a_{g_f},\\ \noalign{\vskip4pt}
\sum _{e\mid N}\beta _e C_N(e, f)&=&\zeta
_F(2)^{-1}b_{g_f},\\
\noalign{\vskip-4pt}
\end{eqnarray*}
 for all $f\mid N$, where $a_{g_f}$ and $b_{g_f}$ are constants 
in (6.4.1). Thus it suffices to show that the matrix
$C_N=(C_N(e, f))_{e, f\mid N}$ is invertible. It is easy to see
that $C_N$ is multiplicative for coprime $N$'s in the sense
of tensor products. Thus it suffices to prove that $C_N$ is
invertible for $N=\wp ^n$ to be a power of a prime. 
In this case, $C_{\wp ^n}$ has determinant $\pm (\RN (\wp)^2-1)^n$.
This completes the proof of the lemma.
\enddemo

\vglue9pt

\spn{6.4.3} \proclaim{Lemma} Let $\wt\phi$ be a form which has   growth $O(y^{1-\varepsilon})$
near each cusp{\rm .}  Let $\wt c(y)$ denote the Whittaker function
at $ \left( \begin{array}{cc}y&0\\ \noalign{\vskip4pt}  0&1 \end{array}\right)$
of $\wt\phi$\/{\rm :}\/
$$\wt c(y)
=d_F^{-1/2}
\int _{\BA _F/F}
\wt\phi \left( \left( \begin{array}{cc} y&x\\ \noalign{\vskip4pt}  0&1 \end{array}\right)\right)
\psi (-x)dx.$$
 Then the Fourier coefficient  of the 
holomorphic projection $\phi$ of $\phi^*$
is given by
$$a(  m )=(4\pi )^g\lim _{s\to 1}\int _{\BR _+^g}
|t|^{-1}\wt c(ty_\infty)
\psi (i y_\infty)|y_\infty|^{s-2}dy_\infty $$
where $t$ is a generator of $mD_F^{-1}$ in $\wh F^\times${\rm .}
\endproclaim

\vglue9pt

\demo{Proof} 
For $  m $  a nonzero ideal of $\CO _F$, let $P_{  m , s} (g)$
be the $  m^{\rm th}$ Poincar\'e series defined by
$$P_{  m , s} (g)=\sum _{\gamma \in Z(F)U(F)\backslash \GL _2(F)}
H_  m  (\gamma g)$$
where $U$ denotes the algebraic group of matrices
$ \left( \begin{array}{cc} 1&x\\ \noalign{\vskip3pt}0&1 \end{array}\right)$ and
$H_  {m, s} $ denotes a function on $Z(\BA_F)\backslash G(\BA_F)$ such that
for $y\in \BA _F^\times, x\in \BA_F, r(\theta)k\in    K  $, 
$$H_{  m , s} \left( \left( \begin{array}{cc}
y&x\\ \noalign{\vskip3pt}0&1 \end{array}\right) r(\theta) k\right)
=|y|^{s}\psi (2\theta+x+iy_\infty)$$
if $y_\infty >0$, $k\in    K  _0( N  )$, and $y_f {D_F}    =  m $;
 otherwise,
it is zero.  Then the Petersson product $(\wt\phi, P_{  m , \bar s})$
is equal to

\begin{eqnarray*}
\noalign{\vskip-12pt}
&&\hskip-.25in \int _{Z(\BA_F)G(F)\backslash G(\BA_F)}
\wt \phi \overline {P_{  m , \bar s}}dg
=\int _{Z(\BA_F) U(F)\backslash G(\BA_F)}
\wt\phi \overline {H_{  m , \bar s}}dg\\ \noalign{\vskip3pt}
&&\hskip.25in =\ \int _{\BA_F^\times}\int _{\BA_F/F}\int _   K   
(\wt\phi \overline {H_{  m , s}})
\left( \left( \begin{array}{cc} y&x\\ \noalign{\vskip3pt} 0&1 \end{array}\right)k \right) dk dx \frac {d^\times y}
{|y|}\\ \noalign{\vskip3pt}
&&\hskip.25in =\ \mu (N)
\int _
{y_\infty\in \BR _+^g}\int_ 
{t\wh \CO _F^\times}
\int _{\BA _F/F}\wt\phi
\left( \left( \begin{array}{cc} y&x\\ \noalign{\vskip3pt} 0&1 \end{array}\right) \right) |y|^{s-1}
\psi (-x+iy_\infty) dx d^\times y\\ \noalign{\vskip3pt}
&&\hskip.25in =\ \mu (N)d_F^{1/2}
\int _{\BR _+^g}\int _{t\wh\CO _F^\times}
\wt c(y)\psi (i y_\infty)|y|^{s-1}d^\times y.
\end{eqnarray*}
Thus we have 
\begin{equation}(\wt\phi, P_{  m , \bar s})
=|t|^{s-1}\mu (N)d_F^{1/2}
\int _{\BR _+^g}\wt c(ty_\infty)
\psi (i y_\infty)|y_\infty|^{s-2}dy. \hskip.4in \spq{6.4.3}
\end{equation}
If we replace $\wt\phi$ by $\phi$ with the  Whittaker function
$$c(y)=|y|a(y_fD_F)\psi (iy_\infty),$$
then  
\begin{equation}
(\phi, P_{  m , s})=|t|^{s}\mu (N)d_F^{1/2}
\left[\frac { \Gamma   (s)}{(4\pi)^{s}}\right]^ga(  m  ) .\spq{6.4.4}
\end{equation} 
As  $P_  m =\lim _{s\to 0}P_{  m , s}$  is a holomorphic
form,  
\begin{equation}
(\phi, P_  m )=(\wt\phi, P_  m  ). \spq{6.4.5}
\end{equation} 
The lemma follows from (6.4.3)--(6.4.5).
\enddemo

We want to apply this lemma to $\wt \Phi $.

\spn{6.4.4} \proclaim{Lemma} Let $a(m)$ be the Fourier coefficient of the holomorphic 
projection of $\Phi'${\rm .} Then for $m$ prime to $ND_E${\rm ,}
$$a(  m  ) \pmod{\CD_N}
=(4\pi)^g\lim _{s\to 1}
\int _{\BR _+^g} |t|^{-1}c'(1,ty_\infty)\psi (iy_\infty)
|y_\infty|^{s-2} dy_\infty
$$
where $t$ is a generator of $\wh m\wh D_F^{-1}$ in $\wh F^\times${\rm .}
\endproclaim

\demo{Proof}
Let $a(y)$, $b(y)$, and $\wt c(y)$ 
be Fourier coefficients of  $E(g)$, $F(g)$, and $ \wt \Phi  (g)$;
 then 
$$\wt c(y)=c'(1,y)-\alpha _1a(y)-\beta_1b(y)$$
where $c'(1, y)$ is the Whittaker function of 
$\Phi'$, and $\alpha _1, \beta _1$ are constants as in Lemma 6.4.2.
By Lemma 6.4.3, $a(  m  )$ is equal to

\begin{equation}
(4\pi)^g\lim _{s\to 0}
\int _{(\BR ^+)^g}|t|^{-1} c'(1,ty)\psi (iy_\infty)|y|^{s-2}dy
-\alpha_1c_s(m)-\beta_1b_s(m) \qquad\spq{6.4.6}
\end{equation}
where
 $$b_s(  m  )=(4\pi)^g
\int _{(\BR ^+)^g}|t|^{-1} b( ty)\psi (i y_\infty)|y_\infty|^{s-2}
dy_\infty$$
and 
$$c_s(  m  )=(4\pi)^g
\int _{(\BR ^+)^g} |t|^{-1}c(ty_\infty)
\psi (iy_\infty)|y_\infty|^{s-2}dy_\infty.$$

Write $\sigma _s(  m )=\sum _{  a   |  m }\RN (  a   )^{s}$ and
$\sigma '(  m )=\frac {\partial}{\partial s}\big|_{s=1}\sigma_s (  m )$.
One can show from the Fourier expansion of $E_{2, s}$ that 
$$
b_s(  m  )=k_1\sigma _1(  m  )+o(s-1),
$$
and
$$
c_s(  m  )
=k_2 \sigma _1(  m  )+k_3\sigma '(  m  )+\frac {k_4}{s-1}+k_5
+o(s-1).
$$
Here the  $k_i$'s are constants independent of $  m $. 
Thus $c_s(m)$, $b_s(m)$ only contribute elements in $\CD_N$ in
(6.4.6). The lemma follows.
\enddemo

Applying the formula for Fourier coefficients of $\wt  \Phi $,
 we obtain the following:

\spn{6.4.5} \proclaim{Proposition} Let $a(  m  ) $ be the Fourier coefficients for the holomorphic 
projection $ \Phi $ of $\Phi'  ${\rm .} Then for $  m $ prime to $N{D_E}   ${\rm ,}
$$a(m)\pmod{\CD_N}=-\frac {(2\pi)^{2g}d_N^{1/2}}{d_Ed_F}
\sum _va_v(  m )$$
 where $\RN (v)=1$ if $v$ is archimedean and $a_v(m)$ is  given by
the following formulas\/{\rm : }\/
\begin{itemize}
\ritem{1.} If $v|\infty${\rm ,} then $a_v(m)$ is equal to the constant term in
the Taylor expansion in $s-1$ of
$$
\sum_{\begin{array}{c} 
{\scriptstyle n\in  N    m  ^{-1} {D_E}    ^{-1},
n_v<0}\\ \noalign{\vskip3pt}
{\scriptstyle 0<n_w<1 \forall v\ne w|\infty}\\
{\scriptstyle \varepsilon _w(n(n-1))=1 \forall w| {D_E} } \end{array} 
}\delta (n)r((1-n)m)
r(n  m  {D_E}   / N  )
p_s(|n_v|)
$$
where 
  the sum is over the set of places of $F${\rm ,}
$$\delta (n)=2^{\#\{ v\mid {D_E}   ,\ord _v(n)\ge 0\}},$$
and
$$p(s,t)=\int _1^\infty (1+tx)^{-s}\frac {dx}x,
\qquad (t>0).$$

\ritem{2.} If $v=\wp\!\!\not|\ \infty, \varepsilon (v)=0${\rm ,} $a_v(  m )$ is equal to
$$\sum_{\begin{array}{c} {\scriptstyle n\in  N    m  ^{-1} {D_E}   ^{-1} }\\
{\scriptstyle \varepsilon_v ((n-1)n)=1}\\
{\scriptstyle \varepsilon_w ((n-1)n)=1
\forall v\ne w| {D_E}   }\\
{\scriptstyle 0< n<1}\end{array}}\delta (n)r((1-n)m)
r(n  m  / N)
\ord _v(n m\wp)\log \RN (v).
$$ 

\ritem{3.} If $v=\wp \!\!\not|\ \infty, \varepsilon (v)=-1${\rm ,} $a_v(  m )$ is equal to
$$\sum_{\begin{array}{c} {\scriptstyle n\in  N    m  ^{-1} {D_E}   ^{-1}}\\
{\scriptstyle \varepsilon_w ((n-1)n)=1
\forall v| {D_E} }\\
{\scriptstyle 0< n<1}\end{array} }\delta (n)r((1-n)m)
r(n  m  / N\wp)
\ord _v(n m\wp/N)\log \RN (v).
$$ 

\ritem{4.} If $v\!\!\not|\ \infty, \varepsilon (v)=1${\rm ,}
$$a_v(  m )
=0.$$
\end{itemize}
\endproclaim

\demo{{P}roof}
By Proposition 6.3.4 and 6.4.4, the Fourier coefficients
$a(m)\pmod {\CD_N}$ are given by 
\begin{equation}
 \frac {\varepsilon ( N  )d_N^{1/2}}{d_Ed_F}
(4\pi )^g\lim _{s\to 1}\sum_{n\in F
} b^n_s(m) \spq{6.4.7}
\end{equation}
where $$b^n_s(m)=\int _{\BR _+^g}|t|^{-1}b^n(1, ty_\infty)
\psi (iy_\infty)|y_\infty|^{s-2}dy_\infty.$$
From formulas of $b^n(1, ty_\infty)$ one can show that
if $n=0$ or $1$,
$b^n_s(m)$ is a linear combination of a 
multiple of $r(m)$ and its derivatives.
Thus, modulo $\CD_N$, we may assume that $n\ne 0, 1$ in (6.4.7).
Moreover  $b^n_s(m)\ne 0$ only
if $n {D_E} mN  ^{-1}$ is integral and 
$1 -n$ is totally positive. In this case
$b^n_s(m)$ is equal to 
\begin{equation}
(-4\pi ^2)^g\delta (n)r((1-n)m)
\sum _v b_{v, s}^n(m) \spq{6.4.8}
\end{equation}
where $v$ runs through all places of $F$,
 and  
$$b^n_{v, s}(m)=\int _{\BR _+^g}
b^n_v(1, ty_\infty)\psi (2iy_\infty)|y_\infty|^{s-1}dy_\infty.$$
Now we compute $b^n_{v, s}(m)$ case by case using Proposition 6.3.4.
\enddemo

\demo{Case 1}
 $v\mid \infty$.
In this case, $b_{v, s}^n(m)$ is nonzero
only if 
\begin{itemize}
\item $n$ is negative at place $v$, but positive at other infinite places,
\item  $\varepsilon _\wp ((n-1)n)=1$
for every place $\wp$ of $D_E$.
\end{itemize}
In this case, $b_{v, s}^n(m)$ is equal to 
\begin{eqnarray}
&&\hskip-.25in r(nm/N )\int _{\BR _+^g}q(4\pi |ny_v|)\psi (2iy_\infty)|y_\infty|^{s-1}dy_\infty\spq{6.4.9} \\
&&\hskip.25in =r(nm/N)\left[\frac {\Gamma (s)}{(4\pi )^s}\right]^g \nonumber
p(s, |n_v|).
\end{eqnarray}
\enddemo

\vglue9pt
\demo{Case 2}  $v\ \nmid\ \infty, \varepsilon (v)=0$.
In this case,  
$b_{v, s}^n(m)$ is nonzero only if 
\begin{itemize}
\item $n$ is totally
positive,
\item $\varepsilon_v((n-1)n)=-1$ but $\varepsilon _\wp ((n-1)n)=1$
for every place $\wp$ of $D_E$.
\end{itemize}
In this case, $b_{v, s}^n(m)$ is equal to 
\begin{equation}
r(nm/N )\ord _\wp(nm\wp)
\log \RN (\wp)\left[\frac {\Gamma (s)}{(4\pi )^s}\right]^g. \spq{6.4.10}
\end{equation}
\enddemo

\vglue9pt
\demo{Case 3} $v\ \nmid\ \infty, \varepsilon (v)=-1$.
In this case,  
$b_{v, s}^n(m)\ne 0$ only if 
\begin{itemize}
\item $n$ is totally
positive,
\item $\varepsilon_\wp((n-1)n)=1$ 
for every place $\wp$ of $D_E$.
\item $\ord _v(nm/N)$ is odd.
\end{itemize}
In this case, $b_{v, s}^n(m)$ is equal to 
\begin{equation}
r(nm/(N\wp_v) )\ord _\wp(nm\wp/N)\log \RN(\wp)
\left[\frac {\Gamma (s)}{(4\pi )^s}\right]^g
. \spq{6.4.11}
\end{equation}

\vglue9pt
\demo{Case 4} $v$ is a finite place split in $E$.
In this case, 
\begin{equation}
b_{v, s}^n(m)=0. \spq{6.4.12}
\end{equation}

The proposition follows from (6.4.7)--(6.4.12) with 
\medbreak 
\phantom{tomorrow}\hfill ${\displaystyle a_v(m)=-(4\pi )^{g}\lim _{s\to 1}\sum_{ n\in F\atop
n\ne 0, 1} b_{v, s}^n(m).}$\hfill\qed
\enddemo

The same proof will also give the following:

\spn{6.4.6} \proclaim{Proposition} Assume that $\varepsilon (N)=(-1)^g${\rm .} Let $b(m)$ denote  
the Fourier coefficient  of the holomorphic projection of 
 $\Phi _{1/2}${\rm .} Then for $m$ prime to $ND_E${\rm ,}
$$
b(m)\pmod{\CD_N}=\frac {(2\pi )^{2g}d_F^{1/2}}{d_Ed_F}
\sum_{\begin{array}{c}
{\scriptstyle n\in ND_E^{-1}m^{-1}}\\
{\scriptstyle 0<n<1}\\
{\scriptstyle \varepsilon _v(n(n-1))=1, \forall v\mid D_E}\end{array}
} \delta (n)r((1-n)m)r(nm/N).
$$
\endproclaim

\section{Proof of the main theorems}

In this section we will finish the proofs of the theorems stated in 
the introduction. We need only prove
Theorem C and   A. 

\demo{{\rm 7.1.} Proof of Theorem {\rm C}}
Recall that  $\Phi$ is the holomorphic form of weight
$2$ for $K_0(N)$ with trivial character as constructed in 6.4.1,
which is the holomorphic projection of $\frac {\partial}{\partial s}
\Phi _s\big|_{s=1/2}$ where $\Phi _s$ is a form constructed as in 
 6.1.5. By Corollary 6.1.6, we thus have
\begin{equation}
(f, \Phi)=B(1/2)L_E'(f, 1) \spq{7.1.1}
\end{equation}
where
$$B(1/2)=2^{-g}d_Fd_N^{1/2}d_E^{-1/2}\mu (N).$$

Recall also that we have constructed a form $\Psi$ in
4.1.3 whose Fourier coefficients 
are height pairings of CM-points $\pair{z, \RT (m)z}$, 
where $z$ is the class
of $\eta$ in the Jacobian of $X$ and the pairing here is the Neron-Tate
height pairing.

The proof of Theorem C will be easily reduced to the following:
\enddemo

\spn{7.1.1} \proclaim{Proposition} With the notation of {\rm  4.4.4,}
$$\wt\Phi =\frac {(2\pi)^{2g}d_N^{1/2}}{d_Ed_F}\wt \Psi \pmod{\CD_N}.$$
\endproclaim 

\demo{{P}roof}
We need to show that both sides have the same value for all
$m\in \BN _F$
prime to $ND_E$, modulo ${\CD_N}$. By Proposition 4.4.5  and  4.5.3,
modulo $\CD_N$,
$\wt \Psi (m)$ is equal to the sum of $-(\eta, \RT ^0(m)\eta)_v$.
On the other hand,   
we have studied the  Fourier coefficients $a(m) \pmod {\CD_N}$
in Proposition 6.4.5 by decomposing it into a sum of local
terms $a_v(m)$. Now we  only need to show that
$$
\sum _{v}(\eta, \RT ^0(m)\eta)_v=\sum _va_v(m)\pmod {\CD_N}.
$$
We will prove this by splitting the sum according to types of $v$. 
More precisely we want to show
\begin{equation}
\sum _{v\in S}(\eta, \RT ^0(m)\eta)_v=\sum _{v\in S}a_v(m)\pmod {\CD_N},
\spq{7.1.2}
\end{equation}
where $S$ is one of the following:
\begin{itemize}
\item[1.] $S_\infty$: the set of infinite places;
\item[2.] $S_0$: the set of finite  places ramified in $E$;
\item[3.] $S_1$: the set of finite places split in $E$;
\item[4.] $S_{-1}$: the set of finite places inert in $E$.
\end{itemize}
\enddemo

\demo{The case of archimedean places}
Since $S_\infty$ is finite, we need only show the individual 
identity
$$(\eta, \RT ^0(m)\eta)_v=a_v(m) \pmod {\CD_N}.$$
In view of Proposition 5.4.4 and 6.5.4, we need only show that the quantity
$$E(s)=
\sum_{\begin{array}{c} {\scriptstyle n\in  N    m  ^{-1} {D_E}    ^{-1},
n_v<0}\\
 {\scriptstyle 0<n_w<1 \forall v\ne w|\infty}\\
 {\scriptstyle \varepsilon _w(n(n-1))=1 \forall w| {D_E}   }\end{array}
}\delta (n)r((1-n)m)
r(n  m  {D_E}   / N  )\varepsilon _s(|n_v|)$$
has limit $0$ as $s\to 1$,
where $$\varepsilon _s(t)=
p_s(t)-2Q_{s-1}(1+2t) \qquad (t>0).
$$
One can show that
$$\varepsilon _1(t)=0, \qquad \varepsilon _s(t)=O(t^{-1-s})$$
as $t\to \infty$. Thus $E(s)$ is absolutely convergent
for $\Re (s)>0$, 
and has limit $0$ as $s\to 1$. \enddemo

\demo{The case of ramified places}
Again $S_0$ is finite, we need only show the individual identity
$$(\eta, \RT ^0(m)\eta)_v=a_v(m) \pmod {\CD_N}.$$
This follows directly from Proposition 5.4.8 and (6.4.5).\enddemo

\demo{The case of split places}
In this case $S_1$ is not finite. But by Proposition 5.4.8, 
the sum
$$\sum _{v\in S_1}(\eta, \RT ^0(m)\eta)_v
=r(m)\sum _{v\in S}\ord _v(m)\log \RN (v)$$
has only finitely many nonzero terms and defines an element in
${\CD_N}$. On the other hand $\sum _{v\in S_1}a_v(m)=0$.
\enddemo

\demo{The case of inert places} 
Again by Proposition 5.4.8, the left-hand side of (7.1.2)
is equal to
\begin{equation}
\sum _{\wp\in S_{-1}}(U_\wp (m)-U_\wp (1)R(m))\log \RN (\wp)
\pmod {\CD_N}. \spq{7.1.3}
\end{equation}
We need to compare $(U_\wp (m)-U_\wp (1)R(m))\log \RN (\wp)$
with $a_\wp (m)$.
For $\ell \in \BN _F$ prime to $\wp$ define
\begin{eqnarray*}
k_\wp(\ell)&=&
\sum_{\begin{array}{c}
{\scriptstyle n\in    \ell ^{-1} {D_E}   ^{-1} N  }\\
{\scriptstyle 0< n <1}\\
{\scriptstyle \varepsilon _ q     (n(n-1))=1, \forall  q    | {D_E} }\end{array}  
}
r(n  \ell  N  ^{-1}/\wp)r((n-1)  \ell)\delta (n),
\\
k'(\ell)&=&
\sum_{\begin{array}{c}
{\scriptstyle n\in    \ell ^{-1}D_E^{-1} N  }\\
{\scriptstyle 0< n <1}\\
{\scriptstyle \varepsilon _ q     (n(n-1))=1, \forall  q    | {D_E}  }\end{array} 
}
r(n \ell  N  ^{-1})r((n-1)  \ell /\wp)\delta (n).
\end{eqnarray*}
\enddemo

\spn{7.1.2} \proclaim{Lemma} Let $m'=m\wp^{-\ord _\wp(m)}${\rm .}
\begin{itemize}
\ritem{1.} If $\ord _\wp(m)$ is even{\rm ,} then
$$U_\wp(m)\log \RN (\wp)-a_\wp (m)=0.$$
\ritem{2.} If $\ord _\wp (m)$ is odd{\rm ,} then
\begin{eqnarray*}
U_\wp(m)&=&\ord _\wp (m\wp)k_\wp(m'),\\
a_\wp(m)&=&\ord _\wp (m\wp)\log \RN (\wp)k'_\wp(m').
\end{eqnarray*}
\end{itemize}
\endproclaim

\demo{Proof}
If $\ord _\wp (  m  )$ is even, then the only nonzero terms in 
$a_\wp(  m  )$ are for those
$n$ which lie in 
 $  m ^{'-1} {D_E}   ^{-1} N  $, where 
$m'=m\wp^{-\ord _\wp (m)}$. 
This is clear if $\wp\ \nmid\ m$. Otherwise, $\wp\nmid N$
and then $r (nm/N\wp)\ne 0$ will imply that $\ord _\wp(n)$ is 
odd. But then $r((n-1)m)\ne 0$ will imply that $\ord _\wp(n)$
is nonnegative. Thus $U_\wp(m)\log \RN (\wp)=a_\wp (m)$.

If $\ord _\wp(  m  )$ is odd, 
then the only nonzero terms in $a_\wp(  m  )$ are for those
$n$ which have zero order at $\wp$. Indeed, $r(nm/N\wp)\ne 0$ implies 
$\ord _\wp(n)$ is even, but $r((1-n)m)\ne 0$ implies $\ord _\wp(n)=0$.
Actually, $\ord _\wp(1-n)$ is positive and even. Thus
\medbreak
\phantom{higher} \hfill ${\displaystyle a_\wp(  m  )=\sum_{\begin{array}{c}
{\scriptstyle n\in    m ^{'-1} {D_E}   ^{-1} N  }\\
{\scriptstyle 0< n' <1}\\
{\scriptstyle \varepsilon _ q    (n(n-1))=1, \forall  q    | {D_E}}\end{array}   
}
r(nm ' N  ^{-1})r((n-1)  m '  /\wp)\delta (n)\ord _\wp (  m \wp).
}$
\enddemo

\vglue12pt

\spn{7.1.3} \proclaim{Lemma} Let $\ell\in \BN_F$ be prime to $\wp$. Then 
$$k(\ell)-k(1)=k'(\ell)-k'(1).$$
\endproclaim

\demo{Proof}
From the proof of Proposition 5.4.8, it is not difficult to see that 
$k_\wp(\ell)-k_\wp(1)$
is the local intersection of $\eta$ and $\RT^0 (\ell)\eta$ 
over $\wp$  without
counting multiplicities. See formula (5.4.10) with
$m=\ell$ and $m(n')=1$.
But this does not give a description for $k'_\wp(\ell)-k'_\wp(\ell)$.
So we need to give a description in a different setting. 

 Let $R(\wp)$ be an order of $B(\wp)$ of type 
$(E, N(\wp))$ and consider the projection map
$$\pi :C=E^\times \backslash \wh B(\wp)^\times/\wh R(\wp)^\times
\to S=B^\times \backslash \wh B(\wp)^\times /\wh R(\wp)^\times.$$

The set $C$ can be considered as the set of CM-points, and the 
set $S$  as supersingular points,
the reduction of CM-points. We may define conductors for elements 
in $C$, and 
orientations for elements in $C$ with conductor prime to $N\wp$
 for each place dividing $N\wp$.
 The group 
$$\CW=\{b\in \wh B(\wp )^\times: b^{-1}\wh R(\wp) b\}
/\wh R (\wp )^\times$$ 
acting on $C$ does not change reductions and conductors
but is free and transitive on orientations. For each place $v$ dividing 
$N\wp$, we call the orientation defined by  $1$ the positive orientation.

By a $\BQ$-divisor in 
$C$ we just mean an element in the free abelian group
$\BQ[C]$. For $\ell$ prime to 
$N (\wp)$ we can also define the 
Hecke operator $\RT (\ell)$ on $\BZ   [C]$.
Let $\eta (\wp)$ (resp.\ $\eta (\wp)'$) be the set of elements in the 
first set with trivial conductor and positive orientations
at all places of $N\wp$ (resp.\ positive orientations at places 
dividing $N$ but negative orientation at place $\wp$).
Then $\eta(\wp)$ and $\eta (\wp)'$ have the exact same reduction
because of the action of $\CW$. 
Now, $k(\ell)-k(1)$ (resp.\ $k'(\ell)-k'(1)$)
is the intersection number of $\eta(\wp)$ (resp.\ $\eta (\wp)'$)
 and $\RT (\ell  )^0(\eta(\wp))$
under the specialization map. Thus they are same
since $\eta (\wp)$ and $\eta (\wp)'$ have the same reductions.
\enddemo

We return to the proof of (7.1.2) for $S_{-1}$.
By (7.1.3), the difference of two
sides of (7.1.2) for $S_{-1}$ is equal to
\begin{eqnarray*}
\sum _{\varepsilon (\wp)=-1}(U_\wp(m)\log \RN (\wp)-a_\wp(m))&-&\sum _{\varepsilon (\wp)=-1}(U_\wp (1)\log \RN
(\wp)-a_\wp (1))r(m)\\
&-&\sum _{\varepsilon (\wp)=-1}a_\wp (1)r(m).
\end{eqnarray*}
The first two terms vanish by the above two lemmas.
The last sum is absolutely convergent, and thus defines an element 
in $\CD_N$. 
\enddemo

\spn{7.1.4} \proclaim{{C}orollary} 
For any 
newform $f$ for  $K  _0(   N    )${\rm ,} 
$$L_ E'(f, 1)=\frac {(8\pi ^2)^g}{d_F^2\sqrt{d_ E}}
[   K   _0(1):   K  _0( N  )](f,  \Psi   ).$$
\endproclaim 

\demo{Proof}
By Propositions 4.5.1 and 7.1.1, 
$$\Phi =
\frac {(2\pi)^{2g}d_N^{1/2}}{d_Ed_F}\Psi 
+\hbox{an old form}.$$
Now the conclusion follows from formula (7.1.1).
\enddemo

\demo{{\rm 7.1.5.} Proof of Theorem {\rm C}}
The ideal is exactly as in \cite[p.\ 308]{G-Z}.
By Lemma 3.4.5, we may decompose $z$ in $\Jac (X)\otimes \BC$
into eigenvectors $z_\phi$ with the same eigenvalues as $\phi$
under Hecke operators $\RT (m)$ with $m$ prime to $N$:
$$z=\sum _{\phi\in S_N}z_\phi,\qquad \RT(m)z_\phi=a_\phi (m)z_\phi.$$
As Hecke operators on $\Jac (\BC)\otimes \BC$ are self-adjoint 
with respect to the Neron-Tate height pairing,
one has the decomposition
$$\Psi =\sum _{\phi\in S_N}\pair{z_\phi, z_\phi}\phi.$$
Now Theorem C follows from this equality and 
Corollary 7.1.4. \enddemo

\demo{{\rm 7.2.} Proof of Theorem {\rm A}}
By Theorems B and C, it suffices to prove the following generalization
of a theorem of   Kolyvagin:
\spn{7.2.1} \proclaim{Proposition} Assume that the Heegner point $ y   _f$ in $A$ is
not a torsion point{\rm .}  Then 
\begin{itemize}
\item $A(F)$ has rank given by
$$\rank A(F)=[\CO_f:\BZ   ]\ord _{s=1}L(s, f),$$
\item $\Sha (A)$ is finite.
\end{itemize}
\endproclaim

In view of Kolyvagin's method for other cases 
 (\cite{G3}, \cite{K}, \cite{K-L}, \cite{K-L2}),
 we need only  to construct certain Euler system  of CM-points.
We consider square-free elements $n\in \BN_F$ which are square-free and
prime to $   N       {D_E}   $ and such that every prime factor $\ell$ is inert in $K$.
For every such $n$, we choose a CM-point 
$x_n$ of the conductor $n$  such that 
$$ \hbox{$x_n$ is included in $\RT (\ell )x_m$}$$
if $n=m\ell$ with $\ell$ a prime. Then 
$x_n$ is defined over $ E_n$, the ring class field of the conductor $n$
over $ E$.

\spn{7.2.2} \proclaim{Lemma}
If $n=m\ell$ as above{\rm ,} then
$$u_n^{-1}\sum _{\sigma \in \Gal ( E_n/ E_m)}x_n^\sigma =u_m^{-1}\RT (\ell )x_m,$$
where $u_n$ is the cardinality of the group
$\CO _ c ^\times/\CO _F ^\times${\rm .}
\endproclaim

\demo{Proof}
By an argument similar to the proof of Proposition 4.2.1, one 
can show that $P:=\frac {u_m}{u_n}\RT (\ell )x_m$ is a divisor with
integral coefficients.  
It follows that $P=Q$ because of the following facts:
\medbreak
 \noindent\hglue.25in $\bullet$     $P$ includes
the divisor $Q=\sum _{\sigma \in \Gal ( E_n/ E_m)}x_n^\sigma$;
\smallbreak  \noindent\hglue.25in $\bullet$    $\deg P=\deg Q$;
\smallbreak \noindent
 \hglue.25in $\bullet$   $Q$ is irreducible
over $ E_m$.
 \enddemo

As in the case $F=\BQ$, this lemma   implies that
the collection of $x_n$ forms
an Euler system \cite{K}.
\bigbreak
\def\ensp{\enspace}
\def\ni{\\ \noalign{\vskip2pt}}
\centerline{\bf Index}
\bigbreak
{\ninepoint 

\begin{tabular}{ll}
\phantom{endomorphism of a distinguished point, and more }&\phantom{doesn't matter}\\
\noalign{\vskip-6pt}
$\cdot^1$, $\cdot^2$, \ensp 13, 20 &$H_c$, \ensp 29\ni 
$\cdot^{2,1}$, $\cdot^{2,2}$, \ensp 23&Heegner point, \ensp 30\ni
$\cdot^{\rho}$, $\cdot_\rho$, \ensp 13& holomorphic form, \ensp 44\ni 
$\alpha$, \ensp 30 &holomorphic projection, \ensp 106\ni
admissible submodules, \ensp 19\ni  %
\phantom{hi there}&$J$,\ensp 7\ni
$B, B'$, \ensp 6, 7  &$J_X$,\ensp 55\ni %
\phantom{hi there} &$j,j_z$, \ensp 8\ni
$C_\phi(g)$, \ensp 44 &\ni  %
CM-points, \ensp 28 &$k$, \ensp 33\ni%
canonical lifting, \ensp 40& $K, K'$,\ensp 7\ni%
conductor of a CM-point, \ensp 28&$K_0,K'_0$, \ensp 25\ni%
cuspidal form, \ensp 44 &  $K_0(N), K^\infty$, \ensp 43\ni%
cyclic submodule structure, \ensp 25\ni%
\phantom{hi}&$\lambda$, \ensp 7\ni%
$\delta$, \ensp 8&  $L(s,\varepsilon,\phi), L_E(s,\phi)$,\ensp 54\ni%
$D_F$, $D_E$, \ensp 6& $L(s,A)$,\ensp 56\ni%
$d_F,d_E,d_N$,\ensp 6& ${\cal L}$, \ensp 63\ni
$\cal D$, \ensp 70& linking numbers,\ensp 75, 78\ni %
derivative, \ensp 70\ni %
det, \ensp 6&$m(n)$,\ensp 78\ni%
distinguished point, \ensp 38 &$M_K, M_{K'}$,\ensp 7\ni%
Drinfeld base, \ensp 15&${\cal M}_K, {\cal M}_{K',\wp}$, \ensp 16\ni%
\phantom{hi}& $\mu_\wp$,\ensp 13\ni%
$\varepsilon$, \ensp 6, 24&modular form,\ensp 43\ni%
$\eta,\tilde\eta, \hat\eta,\eta_c,\eta^0_c$, \ensp 62, 63 \ni%
$E, E'$, $E^{\rm ur}_q$, \ensp 24, 26, 33&$\nu(K')$,\ensp 11\ni%
$E_s, E_{s,N}$, \ensp 90 &$N,N_B,\tilde N,N_E$,\ensp 24, 25\ni 
End$_{\cal F}(A,C)$, End$_{\cal F}(A)$, \ensp 30 &$N_{F'}$,\ensp 26\ni%
endomorphism ring of a CM-point, \ensp 28&$N(\wp)$,\ensp 37\ni%
endomorphism of a distinguished point,\ensp 38&$\Bbb N_F$, \ensp 6\ni%
\phantom{hi}&newform, 47\ni %
$F,F'$, \ensp 6, 7\ni%
${\cal F}, {\cal F}_0$, \ensp 30&${\cal O}_B, {\cal O}_{B'}$,\ensp10, 11\ni %
${\cal F}^0_{K'}$, ${\cal F}_{K'}$, ${\cal F}_{K',\rho}$, \ensp 9, 11, 14&${\cal O}^{\rm ur}_q$, \ensp 33\ni
$\tilde{\cal F}$, \ensp 62& old form,\ensp 47\ni%
\phantom{hi}& ordinary reduction,\ensp 17\ni
$G_m, G_1$, \ensp 18&orientation,\ensp 28, 38\ni%
$[{\cal G},{\cal C}]$, $[{\cal G}^0, {\cal C}^0]$, \ensp 33, 34\ni %
${\cal G}_K$, \ensp 16&$\wp,\wp_E$, \ensp 13, 25\ni%
\phantom{hi} &$\wp_{F'},\wp_{E'}$, \ensp 26\ni %
$\Bbb H$, \ensp 6 &$\psi_F$, \ensp 9\ni 
$\cal H^{\pm}$, \ensp 7&$\Phi_s, \Phi^e_s$,\ensp 93\ni \end{tabular}
 
{\ninepoint 
\begin{tabular}{ll}
\phantom{endomorphism of a distinguished point, and more }&\phantom{doesn't matter}\\
\noalign{\vskip-6pt}
$\Phi',\Phi$,\ensp 103, 106&$Z$, \ensp 63\ni
\phantom{dumb}&$z,\tilde z,\hat z$, \ensp 61, 62\ni
 quasi-canonical liftings,\ensp 42\ni
 quasi-multiplicative,\ensp 70\ni
 $Q_{s-1}(u)$,\ensp 73\ni
\phantom{hi}\ni
\noindent $\rho$, \ensp 24\ni
$\varrho_v(c,n)$,\ensp 75, 80\ni
$\varrho(S)$,\ensp 82\ni
$R$,\ensp 24\ni
$R(\wp)$,\ensp 37\ni
$r(m)$,\ensp 60\ni
\phantom{Hi}\ni
$\sigma_1(m)$,\ensp 18\ni
$\Sigma_h$,\ensp 34\ni
$S_N$,\ensp 48\ni
${\cal S}$,\ensp 69\ni
special module,\ensp 15\ni
supersingular reduction,\ensp 17\ni
\phantom{s}\ni
$\tau$,\ensp 6, 7, 28\ni
$t(\ell),t'(\ell)$,\ensp 8, 30\ni
tr$_{D_E}$,\ensp 93\ni
$\Bbb T$,\ensp 48\ni
$\Bbb T'$,\ensp 48\ni
${\cal T}$,\ensp 55\ni
${\rm T}(m)$,\ensp 46, 48\ni
${\rm T}^0(m)$, \ensp 65\ni
type $(N,E)$,\ensp 64\ni
\phantom{still}\ni
${\cal U}$, \ensp 11\ni
\phantom{going}\ni
$V,V_{\Bbb E}, V_{\Bbb Z}$, \ensp 8, 11\ni
\phantom{strong}\ni
${\cal W}$, \ensp 29\ni
$W_\phi(g)$, \ensp 44\ni
\phantom{after}\ni
$\xi,\tilde\xi, \hat\xi$, \ensp 62, 63\ni
$X$, \ensp 25\ni
$\tilde X,\tilde{\cal X}$, \ensp 61\ni
\phantom{all}\ni
$Y,Y',Y_c$, \ensp 29,  30\ni
${\cal Y}, {\cal Y}_k, y, \bar{y},y_k$, \ensp 33\ni
\phantom{these}\ni
\end{tabular}
}

\bye
\begin{theindex}
\item $\cdot ^1$, $\cdot ^2$, \pr{^1}, \pr{cdot^1ell}
\item $\cdot ^{2,1}$, 
$\cdot ^{2,2}$, \pr{cdot^{2,1}}
\item $\cdot ^\wp$, $\cdot_\wp$,
\pr {cdot^wp}

\indexspace
\item $\alpha$, \pr{alpha}
\item admissible submodules
\pr{D:admissible}

\indexspace
\item $B$, $B'$, \pr{B}, \pr{B'}

\indexspace
\item $C_\phi (g)$, \pr{C_phi}
\item CM-points, \pr{CMpoint}
\item canonical lifting, \pr{can.lif}
\item conductor of a CM-point,
\pr{conductor}
\item cuspidal form, \pr{cuspidal}
\item cyclic submodule structure,
\pr{cyc.sub.mod.str}

\indexspace

\item $\delta$, \pr{delta}
\item $D_F$, $D_E$, \pr{D_F}
\item $d_F$, $d_E$, $d_N$, \pr{d_F}
\item $\CD$, \pr{CD}
\item derivative, \pr{derivative}
\item $\det$, \pr{det}
\item distingushed point, \pr{dis.poi}
\item Drinfeld base, \pr{Drinfeldbase}

\indexspace

\item $\varepsilon$, \pr{epsilon}, \pr{epsilon(m)}
\item $\eta$, $\wt\eta$, $\wh \eta$, $\eta_c$, $\eta _c^0$, \pr{eta},
  \pr{eta_c}
\item $E$, $E'$, $E_q^\ur$, \pr{E}, \pr{E'}, \pr{E_q^ur}
\item $E_s$, $E_{s, N}$, \pr{E_s}
\item $\End_\CF (A, C)$,  $\End _{\CF _0}(A)$, \pr{End_0}
\item endomorphism ring of a CM-point,
\pr{end.rin}
\item endomorphism of a distinguished point,
\pr{Endwp}

\indexspace

\item $F$, $F'$, \pr{F}, \pr{F'}
\item $\CF$, $\CF _0$,
\pr{CFR}
\item $\CF^0_{K'}$, $\CF_{K'}$, $\CF _{K', \wp}$, \pr{CF^0K'}, \pr{CFK'}, \pr{CFK'wp}
\item $\wt\CF$, \pr{wtCF}

\indexspace
\item $G _m$, $G_1$,  \pr{G_m}
\item $[\CG, \CC]$, $[\CG ^0, \CC^0]$, \pr{CG}, \pr{CG^0CC^0}
\item $\CG _{K'}$, 
 $\CG _K$, \pr{CGK}

\indexspace

\item $\BH$, \pr{BH}
\item $\CH^\pm$, \pr{CHpm}
\item $H_c$, \pr{H_c}
\item Heegner point, \pr{Hee.poi}
\item holomorphic form, \pr{holomorphic}
\item holomorphic projection, \pr{Phi}

\indexspace

\indexspace

\item $J$, \pr{J}
\item $J_X$, \pr{J_X}
\item $j$,  $j_z$, \pr{j_z}

\indexspace

\item $k$, \pr{k}
\item $K$,  $K'$, \pr{K},
\pr{K'}
\item $K_0$, 
 $K_0'$, \pr{K0'}
\item $K_0(N)$, $K_\infty$, \pr{K_infty}

\indexspace

\item $\lambda$, \pr{lambda}
\item $L(s, \chi, \phi)$, $L_E(s, \phi)$, 
\pr{L(s,chi,phi)}, \pr{L_E(s,phi)}, 
\item $L(s, A)$,\pr{L(s,A)}
\item $\CL$, \pr{CL}
\item linking numbers, \pr{lin.num}, \pr{lin.num.2}

\indexspace

\item $m(n)$, \pr{m(n)}
\item $M_K$, $M_{K'}$, \pr{MK}, \pr{MK'}
\item $\CM _K$, $\CM _{K',\wp}$, \pr{CMK}, \pr{CMK'wp}
\item $\mu _p$, \pr{mup}
\item modular form, \pr{mod.for}

\indexspace

\item $\nu(K')$, \pr{nuK'}
\item $N$, $N_B$,
 $\wt N$, 
 $N_E$, \pr{NE}
\item $N_{F'}$,  \pr{NE'}
\item $N(\wp)$, \pr{N(wp)}
\item $\BN _F$, \pr{BN_F}
\item new form,  newform, \pr{newform}

\indexspace

\item $\CO _B$, 
 $\CO _{B'}$, \pr{COB'}
\item $\CO _q^\ur$, \pr{CO_q^ur}
\item old form, \pr{old.for}
\item ordinary reduction, \pr{ocase}
\item orientation,
\pr{orientation}, \pr{orientationwp}

\indexspace

\item $\wp$, $\wp_E$, \pr{wp}, \pr{wpE}
\item $\wp _{F'}$, 
 $\wp_{E'}$, \pr{wpE'}
\item $\psi _F$, \pr{psiF}
\item $\Phi _s$, $\Phi_s^e$, \pr{Phi_s}
\item $\Phi'$, $\Phi$, \pr{Phi'}, \pr{Phi}

\indexspace
\item quasi-canonical liftings,
\pr{qua.can.lif}
\item quasi-multiplicative,
\pr{qua.mul}
\item $Q_{s-1}(u)$, \pr{Q_{s-1}(u)}

\indexspace
\item $\rho$, \pr{rho}
\item $\varrho _v(c, n)$, \pr{varrho_v(c,n)},\pr{varrho_wp(c,n)}
\item $\varrho (S)$, \pr{varrho(S)}
\item $R$, \pr{R}
\item $R(\wp)$, \pr{R(wp)}
\item $r(m)$, \pr{r(m)}

\indexspace
\item $\sigma _1(m)$, \pr{sigma_1(m)}
\item $\Sigma _h$, \pr{Sigma_h} 
\item $S_N$, \pr{S_N}
\item $\CS$, \pr{CS}
\item special module, \pr{specialmodule}
\item supersinglular reduction, \pr{sscase}

\indexspace

\item $\tau$, \pr{tau}, \pr{tau'}, \pr{tauE}
\item $t(\ell)$, $t'(\ell)$, \pr{tell},\pr{t'}
\item $\tr_{D_E}$, \pr{tr_{D_E}}
\item $\BT$, \pr{BT}
\item $\BT'$, \pr{BT'}
\item $\CT$, \pr{CT}
\item $\RT (m)$, \pr{OT(m)}, \pr{RT(m)}, \pr{RT(m)X}
\item $\RT ^0(m)$, \pr{RT^0}
\item type $(N, E)$, \pr{typeNE}

\indexspace

\item $\CU$, \pr{CU}

\indexspace

\item $V$,  $V_\BR$, $V_\BZ$, \pr{VBR}, \pr{VBZ}

\indexspace
\item $\CW$, \pr{CW}
\item $W_\phi(g)$, \pr{W_phi}

\indexspace

\item $\xi$, $\wt\xi$, $\wh \xi$, \pr{xi}
\item $X$, \pr{X}
\item $\wt X$, $\wt\CX$, \pr{wtX}

\indexspace
\item $Y$, $Y'$, $Y_c$,\pr{Y}, \pr{Y'}, \pr{Y_c}
\item $\CY$,  $\CY _k$, $y$,  $\bar y$,  $y_k$, \pr{y_k}

\indexspace
\item $Z$, \pr{Z}
\item $z$, $\wt z$, $\wh z$, \pr{z}

\indexspace
\end{theindex}

\end{thebibliography}
\bigskip


\begin{references}


\bibitem{B-C}
\name{J.-F. Boutot} and \name{H. Carayol}, Uniformisation $p$-adique des courbes
de Shimura: les th\'eor\`emes de \v Cerednik et de   Drinfeld,
{\it Ast\'erisque\/} {\bf 196--197} (1991),  45--158.

\bibitem{B}
\name{D. Bump},
{\it Automorphic Forms and Representations},
{\it Cambridge Stud.\  in Adv.\ Math\/}.
{\bf 55}, Cambridge Univ.\  Press, Cambridge (1997).

\bibitem{B-F-H}
\name{D.~Bump, S.~Friedberg}, and \name{J.~Hoffstein},
 Nonvanishing theorems for $L$-functions of modular forms and
their   derivatives, {\it Invent.\ Math\/}.\ {\bf 102} (1990), 543--618.

\bibitem{C}
\name{H. Carayol},  Sur la mauvaise r\'eduction des courbes de
Shimura,
{\it Compositio Math\/}.\ {\bf 59} (1986),   151--230.

\bibitem{Cas}
\name{W. Casselman},  On some results of Atkin and Lehner,
{\it Math.\ Ann\/}.\ {\bf 201} (1973), 301--314.

\bibitem{D}
\name{P. Deligne},  Travaux de Shimura, {\it S{\rm \'{\it e}}minaire Bourbaki\/}, Exp.\
No.\ 389,
in {\it Lecture Notes  in Math\/}.\ {\bf 244}, 123--165, 
Springer-Verlag, New York, 1971.

\bibitem{Dr}
\name{V. G. Drinfeld},  Coverings of $p$-adic symmetric regions,
{\it Funct.\ Anal.\ Appl\/}.\ 
{\bf 10} (1976),  29--40.

\bibitem{F}
\name{G. Faltings},  Endlichkeitss\"atze f\"ur abelsche Variet\"aten
\"uber Zahlk\"opern, 
{\it Invent.\ Math\/}.\ {\bf 73} (1983), 349--366.

\bibitem{F2}
\bibline,  Calculus on arithmetic surfaces,
{\it Ann.\ of  Math\/}.\ {\bf 119} (1984), 387--424.

\bibitem{Gel}
\name{S. Gelbart},
{\it Automorphic Forms on Ad\`ele Groups},
{\it Ann.\  of Math.\ Studies\/} {\bf 83},
Princeton Univ.\  Press, Princeton, NJ  (1975).

\bibitem{G-S1} 
\name{H. Gillet} and \name{C. Soul\'e},
Arithmetic intersection theory,
{\it I.H.E.S.\ Publ.\ Math\/}.\ {\bf 72} (1990), 94--174.

\bibitem{G-S2}
\bibline,
 Characteristic class for algebraic vector bundles with Hermitian 
metrics, I, II, {\it Ann.\ of Math\/}.\ {\bf 131} (1990),
163--203, 205--238.


\bibitem{G-J}
\name{R. Godement} and \name{H. Jacquet},
{\it Zeta Functions of Simple Algebras},
{\it Lecture Notes in Math\/}.\ {\bf 260}
Springer-Verlag, New York  (1972).

\bibitem{Go-Zh}
\name{D. Goldfeld} and \name{S. Zhang},
Holomorphic kernels of Rankin-Selberg convolutions,
manu\-script, 1988. 

\bibitem{G}
\name{B. H. Gross},  On canonical and quasi-canonical liftings, {\it Invent.\
Math\/}.\
{\bf 84} (1986), 321--326.



\bibitem{G2}
\bibline,  Local heights on curves, in {\it Arithmetic Geometry\/}
(Storrs, Conn., 1984), 327--339,
Springer-Verlag, New York (1986).

\bibitem{G3}
\bibline,  Kolyvagin's work on modular elliptic curves,
in {\it $L$-functions and Arithmetic\/} (Durham, 1989), 235--256,
Cambridge Univ.\  Press, Cambridge  (1991). 

\bibitem{G4}
\bibline,  Heegner points on $X_0(N)$, in 
{\it Modular Forms\/} (Durham, 1983),  
87--105, Ellis Horwood,  Chichester (1984).

\bibitem{Gr5}
\bibline, Local orders, root numbers, and modular curves,
{\it Amer.\ J.\ Math\/}.\ {\bf 110} (1988),  1153--1182.

\bibitem{G-Z}
\name{B. H. Gross} and \name{D. B. Zagier},
On singular moduli, {\it J. Reine.\ Angew.\ Math\/}.\ 
{\bf 355} (1985),  191--220.

\bibitem{G-Z1}
\bibline,  Heegner points and derivatives of
$L$-series,
{\it Invent.\ Math\/}.\ {\bf 84} (1986),  225--320.



\bibitem{G-Z2}
\name{B. H. Gross, W.~Kohnen},  and \name{D.~B.~Zagier},
 Heegner points and derivatives of
$L$-series II,
{\it Math.\ Ann\/}.\ {\bf 278} (1987),  497--562.

\bibitem{H}
\name{E. Hecke}, {\it Vorlesung {\rm \"{\it u}}ber die Theorie der Algebraischen Zahlen},
Akademische Verlagsgesellchaft, Leipzig (1923).

\bibitem{J-L}
\name{H. Jacquet} and \name{R. Langlands},
{\it Automorphic Forms on $\GL _2$},
{\it Lecture  Notes in Math\/}.\ {\bf 114}, 
Springer-Verlag, New York  (1971).

\bibitem{K-M}
\name{N. Katz} and \name{B. Mazur},
{\it Arithmetic Moduli of Elliptic Curves},
{\it Ann.\ of  Math.\ Studies} {\bf 108}, Princeton Univ.\ Press,
Princeton, NJ (1985).


\bibitem{Ke}
\name{K. Keating},
Intersection numbers of Heegner divisors on Shimura curves,
preprint.

\bibitem{Kn}
\name{M. Kneser},  Hasse principle for $H^1$ of simply connected groups,
in {\it Algebraic Groups and Discontinuous Subgroups}, 
{\it Proc.\ Sympos.\  in Pure Math\/}.\ {\bf IX} (Colorado, 1965), 
159--163, A.M.S., 
Providence, RI  (1966).

\bibitem{K}
\name{V. A. Kolyvagin},
{\it Euler Systems},
The Grothendieck Festschrift,  {\it Progr.\ in Math\/}.\ {\bf 87}, Birkh\"auser
Boston, Boston, MA (1990).

\bibitem{K-L}
\name{V. A. Kolyvagin} and \name{D. Yu. Logachev},
 Finiteness of the Shafarevich-Tate group and the group of rational
points for some modular abelian varieties,
{\it Leningrad Math.\ J\/}.\ {\bf 1} (1990),  1229--1253.

\bibitem{K-L2}
\bibline,
Finiteness of $\sSha$ over totally real fields,
{\it Math.\ USSR Izvestiya\/} {\bf 39} (1992), 829--853.

\bibitem{Ku}
\name{S. Kudla},
 Central derivatives of Eisenstein series and height pairings,
{\it Ann.\ of  Math\/}.\ {\bf 146} (1997), 545--646.


\bibitem{Mu}
\name{D. Mumford},
{\it Geometric Invariant Theory}, Springer-Verlag, New York  (1965).


\bibitem{R}
\name{D. P. Roberts},
Shimura curves analogous to $X_0(N)$,
Ph.D. Thesis, Harvard University (1989).

\bibitem{Se}
\name{J-P. Serre},
{\it Abelian $\ell$-adic Representation and Elliptic Curves},
W. A. Benjamin, Inc.,
New York (1968).

\bibitem{S}
\name{G.  Shimura},  Construction of class fields and zeta function of 
algebraic curves,
{\it Ann.\ of  Math\/}.\ {\bf 85} (1967),  58--159.

\bibitem{S2}
\bibline, {\it Introduction to the Arithmetic Theory of Automorphic 
Forms}, {\it Publications of the Math.\  Society of Japan\/} {\bf 11},
Princeton Univ.\  Press, Princeton, NJ  (1971).

\bibitem{T}
\name{J. Tate}, {\it Classes d\/{\rm '}\/Isog{\rm \'{\it e}}nie des Vari{\rm \'{\it e}}t\'es Ab{\rm \'{\it e}}liennes sur
un Corps Fini}, {\it S{\rm \'{\it e}}minaire Bourbaki\/}, no.\ 352, 1968--1969.

\bibitem{Tu}
\name{J. Tunnell}, 
 Local $\epsilon$-factors and characters of $\GL(2)$,
{\it Amer.\  J.\ Math\/}. {\bf 105} (1983),  1277--1307.


\bibitem{V}
\name{M.-F.  Vign\'eras}, {\it Arithm{\rm \'{\it e}}tique des Alg{\rm \'{\it e}}bres de Quaternions},
{\it Lecture Notes in Math\/}.\ {\bf 800}, Springer-Verlag, New York,
1980.

\bibitem{W}
\name{J.-L.~Waldspurger},  Sur les values de certaines fonctions $L$ automorphes en
leur centre de symm\'etre, {\it Compositio Math\/}.\ {\bf 54} (1985),
173--242.

\bibitem{W2}
\bibline, Correspondances de Shimura et quaternions,
{\it Forum Math\/}.\ {\bf 3} (1991),  219--307.

\bibitem{W-M}
\name{W. Waterhouse} and \name{J. Milne},  Abelian varieties over finite fields,
{\it Proc.\ Sympos.\ Pure Math\/}.\ (SUNY, Stony Brook, 1969), 
{\bf 20}, 53--64, A.M.S., Providence, RI
(1971).


\end{references}
\end{document}